\documentclass[10pt]{amsart}

\usepackage{xypic,amscd,amssymb,latexsym}
\usepackage{amsmath}	
\usepackage{graphicx}

\usepackage{color}
\usepackage[pdftex]{hyperref}
\usepackage{cancel}

\usepackage{soul}

\usepackage{mathrsfs}

\usepackage{url}
\usepackage{pifont}

\DeclareMathAlphabet{\mathpzc}{OT1}{pzc}{m}{it}

\usepackage{setspace}

\usepackage{multirow}

\usepackage{changepage}
\usepackage{xcolor}

\usepackage{geometry}

\usepackage{tikz}

\definecolor{pinky}{rgb}{1.0, 0, 1.0}

\newcommand{\red}[1]{{\color{red}#1}}

\begin{document}
\pagenumbering{arabic}

\renewcommand{\circled}[1]{\scalebox{0.8}{\textcircled{\fontsize{8}{9}\selectfont \raisebox{-0.1pt}{#1}}}}

\newcommand{\circledd}[1]{\scalebox{0.8}{\textcircled{\fontsize{6}{7}\selectfont \raisebox{0.5pt}{#1}}}}

\newcommand*\circleddd[1]{\scalebox{0.4}{\tikz[baseline=(char.base)]{
            \node[shape=circle,draw,inner sep=1pt,very thick] (char) {#1};}}}

\newcommand{\vbar}[0]{\textcolor{gray}{$\, |$}}

\newtheorem{theorem}{Theorem}[section]
\newtheorem{proposition}[theorem]{Proposition}
\newtheorem{lemma}[theorem]{Lemma}
\newtheorem{corollary}[theorem]{Corollary}
\newtheorem{remark}[theorem]{Remark}
\newtheorem{definition}[theorem]{Definition}
\newtheorem{question}[theorem]{Question}
\newtheorem{claim}[theorem]{Claim}
\newtheorem{conjecture}[theorem]{Conjecture}
\newtheorem{defprop}[theorem]{Definition and Proposition}
\newtheorem{example}[theorem]{Example}
\newtheorem{deflem}[theorem]{Definition and Lemma}

\newcommand{\II}[0]{{\rm \it II}}

\def\qed{{\quad \vrule height 8pt width 8pt depth 0pt}}

\newcommand{\cplx}[0]{\mathbb{C}}

\newcommand{\fr}[1]{\mathfrak{#1}}

\newcommand{\vs}[0]{\vspace{2mm}}

\newcommand{\til}[1]{\widetilde{#1}}

\newcommand{\mcal}[1]{\mathcal{#1}}

\newcommand{\myul}[1]{\underline{#1}}

\newcommand{\ol}[1]{\overline{#1}}

\newcommand{\wh}[1]{\widehat{#1}}

\newcommand{\smallmatthree}[9]{
\left(
\begin{smallmatrix}
#1 & #2 & #3 \\
#4 & #5 & #6 \\
#7 & #8 & #9
\end{smallmatrix}
\right)
}

\author[H. Kim]{Hyun Kyu Kim}
\email{hkim@kias.re.kr}

\address{School of Mathematics, Korea Institute for Advanced Study (KIAS), 85 Hoegi-ro, Dongdaemun-gu, Seoul 02455, Republic of Korea}

\numberwithin{equation}{section}

\title[Naturality of ${\rm SL}_3$ quantum trace maps for surfaces]{Naturality of ${\rm SL}_3$ quantum trace maps for surfaces}

\thanks{This paper will appear in Quantum Topol. \quad DOI: 10.4171/QT/214}

\begin{abstract}
Fock-Goncharov's moduli spaces $\mathscr{X}_{{\rm PGL}_3,\frak{S}}$ of framed ${\rm PGL}_3$-local systems on punctured surfaces $\frak{S}$ provide prominent examples of cluster $\mathscr{X}$-varieties and higher Teichm\"uller spaces. In a previous paper of the author \cite{Kim}, building on the works of others, the so-called ${\rm SL}_3$ quantum trace map is constructed for each triangulable punctured surface $\frak{S}$ and an ideal triangulation $\Delta$ of $\frak{S}$, as a homomorphism from the stated ${\rm SL}_3$-skein algebra of the surface to a quantum torus algebra that deforms the ring of Laurent polynomials in the cube-roots of the cluster coordinate variables for the cluster $\mathscr{X}$-chart for $\mathscr{X}_{{\rm PGL}_3,\frak{S}}$ associated to $\Delta$. We develop quantum mutation maps between special subalgebras of the cube-root quantum torus algebras for different triangulations and show that the ${\rm SL}_3$ quantum trace maps are natural, in the sense that they are compatible under these quantum mutation maps. As an application, the quantum ${\rm SL}_3$-${\rm PGL}_3$ duality map constructed in the previous paper is shown to be independent of the choice of an ideal triangulation.
\end{abstract}

\maketitle

\vspace{-7mm}

\tableofcontents

\section{Introduction}
\label{sec:introduction}

\subsection{Naturality of quantum ${\rm SL}_3$-${\rm PGL}_3$ duality maps for Fock-Goncharov cluster varieties}

Let $\frak{S}$ be a {\em generalized marked surface} (or a {\em decorated surface}), obtained from a compact oriented smooth real surface $\ol{\frak{S}}$ with possibly empty boundary by removing a non-empty finite set of points called {\em marked points}, where we choose at least one marked point from each boundary component of $\ol{\frak{S}}$. So each component of the boundary of $\frak{S}$ is diffeomorphic to an open interval; we call it a {\em boundary arc} of $\frak{S}$. A marked point in the interior of $\ol{\frak{S}}$ is called a {\em puncture} of $\frak{S}$. If $\partial \ol{\frak{S}} = {\O}$, $\frak{S}$ is called a {\em punctured surface}. Let ${\rm G}$ be a split reductive algebraic group over $\mathbb{Q}$, such as ${\rm SL}_n$ or ${\rm PGL}_n$, where $n\ge 2$. The moduli space $\mathscr{L}_{{\rm G},\frak{S}}$ of ${\rm G}$-local systems on $\frak{S}$ has been a central object of study in many areas of mathematics and physics. Some enhanced versions $\mathscr{A}_{{\rm G},\frak{S}}$, $\mathscr{X}_{{\rm G},\frak{S}}$ and $\mathscr{P}_{{\rm G},\frak{S}}$ of ${\rm G}$-local systems with certain kinds of boundary data are defined and studied by Fock and Goncharov \cite{FG06} and by Goncharov and Shen \cite{GS19}; the spaces $\mathscr{X}_{{\rm G},\frak{S}}$ and $\mathscr{P}_{{\rm G},\frak{S}}$ are equipped with Poisson structures, and they coincide with each other in the case when $\frak{S}$ is a punctured surface. One of the crucial properties of these enhanced moduli stacks is that they have structures of {\em cluster varieties} \cite{FG06} \cite{GS19}, which first appeared in the early 2000s and are gaining more interest especially recently, where these moduli spaces associated to surfaces and  algebraic groups form a very important class of examples. 

\begin{figure}[htbp!]
\vspace{-4mm}
\begin{center}
\hspace*{-0mm}
\raisebox{-0.3\height}{\scalebox{1.0}{
\begingroup%
  \makeatletter%
  \providecommand\color[2][]{%
    \errmessage{(Inkscape) Color is used for the text in Inkscape, but the package 'color.sty' is not loaded}%
    \renewcommand\color[2][]{}%
  }%
  \providecommand\transparent[1]{%
    \errmessage{(Inkscape) Transparency is used (non-zero) for the text in Inkscape, but the package 'transparent.sty' is not loaded}%
    \renewcommand\transparent[1]{}%
  }%
  \providecommand\rotatebox[2]{#2}%
  \newcommand*\fsize{\dimexpr\f@size pt\relax}%
  \newcommand*\lineheight[1]{\fontsize{\fsize}{#1\fsize}\selectfont}%
  \ifx\svgwidth\undefined%
    \setlength{\unitlength}{368.50393701bp}%
    \ifx\svgscale\undefined%
      \relax%
    \else%
      \setlength{\unitlength}{\unitlength * \real{\svgscale}}%
    \fi%
  \else%
    \setlength{\unitlength}{\svgwidth}%
  \fi%
  \global\let\svgwidth\undefined%
  \global\let\svgscale\undefined%
  \makeatother%
  \begin{picture}(1,0.23076923)%
    \lineheight{1}%
    \setlength\tabcolsep{0pt}%
    \put(0,0){\includegraphics[width=\unitlength,page=1]{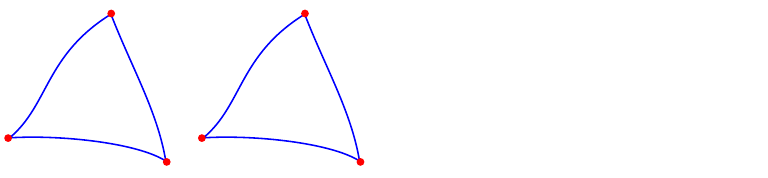}}%
    \put(0.02611683,0.17113519){\makebox(0,0)[lt]{\lineheight{1.25}\smash{\begin{tabular}[t]{l}$\Delta$\end{tabular}}}}%
    \put(0.2811286,0.17011185){\makebox(0,0)[lt]{\lineheight{1.25}\smash{\begin{tabular}[t]{l}$Q_\Delta^{[2]}$\end{tabular}}}}%
    \put(0,0){\includegraphics[width=\unitlength,page=2]{m-triangulation.pdf}}%
    \put(0.5538532,0.1863939){\makebox(0,0)[lt]{\lineheight{1.25}\smash{\begin{tabular}[t]{l}$Q_\Delta^{[3]}$\end{tabular}}}}%
    \put(0,0){\includegraphics[width=\unitlength,page=3]{m-triangulation.pdf}}%
    \put(0.80215431,0.1863939){\makebox(0,0)[lt]{\lineheight{1.25}\smash{\begin{tabular}[t]{l}$Q_\Delta^{[4]}$\end{tabular}}}}%
    \put(0,0){\includegraphics[width=\unitlength,page=4]{m-triangulation.pdf}}%
  \end{picture}%
\endgroup%
}} 
\end{center}
\vspace{-5mm}
\caption{$n$-triangulation quiver, for one triangle}
\vspace{-2mm}
\label{fig:n-triangulation}
\end{figure}

\vspace{2mm}

Here we focus on $\mathscr{A}_{{\rm SL}_n,\frak{S}}$, $\mathscr{X}_{{\rm PGL}_n,\frak{S}}$ and $\mathscr{P}_{{\rm PGL}_n,\frak{S}}$. Let us recall the quivers relevant to the cluster variety structures on them. Choose an {\em ideal triangulation} $\Delta$ of $\frak{S}$, i.e. a mutually disjoint collection of simple paths in $\frak{S}$ running between marked points, called {\em ideal arcs}, dividing $\frak{S}$ into {\em ideal triangles}, which are regions bounded by three ideal arcs. We assume that the valence of $\Delta$ at each puncture of $\frak{S}$ is at least two, which means that we do not allow `self-folded' triangles; see \S\ref{subsec:surfaces_and_triangulations} for a discussion on this condition. For each ideal triangle of $\Delta$, consider the quiver as in Fig.\ref{fig:n-triangulation} depending on $n$, and glue them throughout the surface to obtain a single quiver, called the {\em $n$-triangulation quiver} $Q_\Delta^{[n]}$ for $\Delta$ \cite{FG06} \cite{FG06b} (one must cancel the length 2 cycles formed by dashed arrows). In the present paper we mainly deal with the case of $n=3$ only, so $Q_\Delta^{[3]}$ will be denoted by $Q_\Delta$ in the main text. For any quiver $Q$, denote the set of all nodes of $Q$ by $\mathcal{V}(Q)$, and its {\em signed adjacency matrix} by $\varepsilon = \varepsilon_Q$, which is a $\mathcal{V}(Q) \times \mathcal{V}(Q)$ matrix whose entries $\varepsilon_{vw}$ are defined as
$$
\varepsilon_{vw} = (\mbox{the number of arrows from $v$ to $w$}) - (\mbox{the number of arrows from $w$ to $v$}), \qquad v,w\in \mathcal{V}(Q).
$$
As in Fig.\ref{fig:n-triangulation} we allow dashed arrows, which can be viewed as `half' arrows; each of them contributes by $\frac{1}{2}$ when counting the number of arrows.

\vs

\begin{figure}[htbp!]
\hspace*{-5mm}
\raisebox{-0.5\height}{\scalebox{0.7}{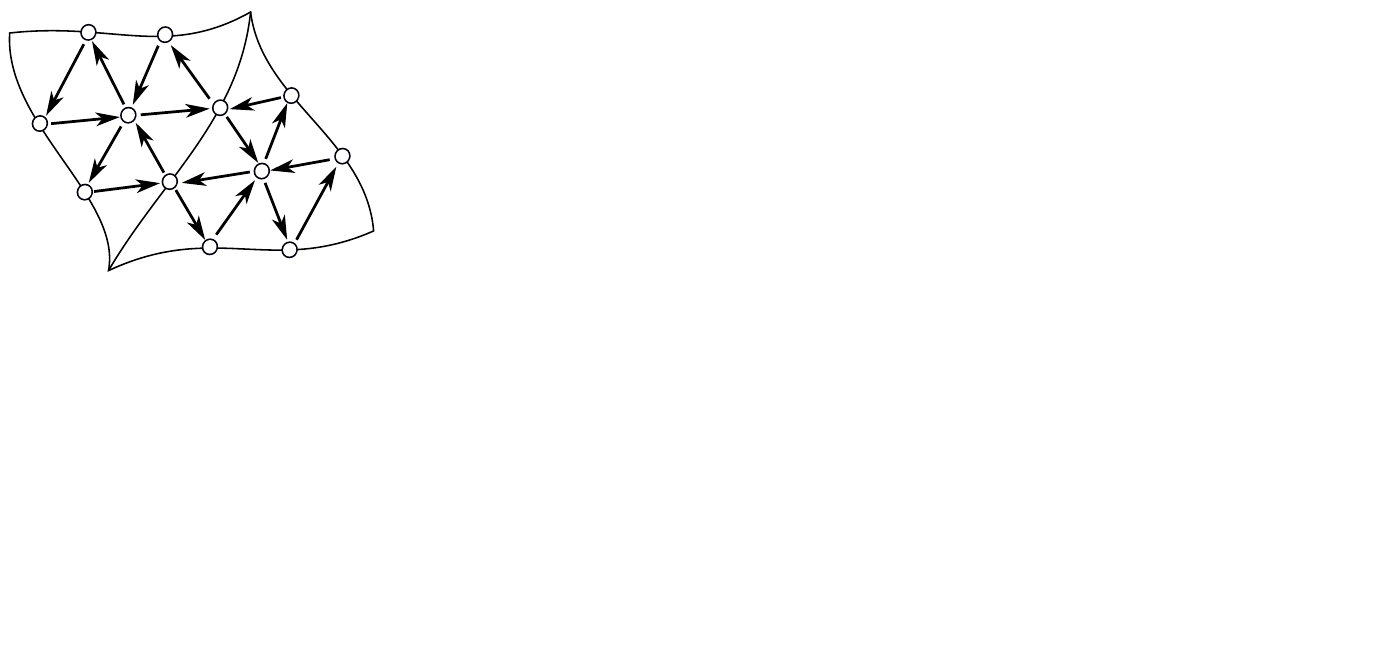}}
\vspace{-2mm}
\caption{\resizebox{140mm}{!}{The sequence of four mutations for a flip at an arc, transforming $Q_\Delta=Q_\Delta^{[3]}$ to $Q_{\Delta'}=Q_{\Delta'}^{[3]}$}}
\vspace{-2mm}
\label{fig:mutations_for_a_flip}
\end{figure}

Per each choice of an ideal triangulation $\Delta$ of $\frak{S}$, it is known that there exist birational maps \cite{FG06} \cite{GS19}
$$
\mathscr{A}_{{\rm SL}_n,\frak{S}} \dashrightarrow (\mathbb{G}_m)^{\mathcal{V}(Q)} \quad\mbox{and}\quad
\mathscr{P}_{{\rm PGL}_n,\frak{S}} \dashrightarrow (\mathbb{G}_m)^{\mathcal{V}(Q)},
$$
called {\em cluster $\mathscr{A}$- and $\mathscr{X}$-charts} for $\Delta$, respectively\footnote{The charts for the space $\mathscr{P}_{{\rm PGL}_n,\frak{S}}$ are called ``$\mathscr{X}$-charts", instead of ``$\mathscr{P}$-charts", where the terminology ``$\mathscr{X}$-charts" comes from the theory of cluster $\mathscr{X}$-varieties.}, with $Q$ being the $n$-triangulation quiver $Q_\Delta^{[n]}$ for $\Delta$. Here, $\mathbb{G}_m={\rm Spec}(\mathbb{Q}[x^{\pm 1}])$ is the multiplicative group scheme, whose set of ${\bf k}$-points is ${\bf k}^*$, for a field ${\bf k}$. We denote the above cluster $\mathscr{X}$-chart of $\mathscr{P}_{{\rm PGL}_n,\frak{S}}$ associated to $\Delta$ by the symbol $\Gamma_\Delta$. The transition maps between two such charts for different ideal triangulations are given by compositions of certain sequences of {\em cluster $\mathscr{A}$- and $\mathscr{X}$-mutation} formulas. Let us elaborate a little more. Given a cluster $\mathscr{A}$-chart with the underlying quiver $Q$, with the cluster $\mathscr{A}$-coordinate variables $A_v$ for the nodes $v$ of $Q$, through the {\em mutation} $\mu_k$ at the node $k$ one obtains another cluster $\mathscr{A}$-chart with the quiver $\mu_k(Q) = Q'$ such that $\mathcal{V}(Q')=\mathcal{V}(Q)$ whose signed adjacency matrix $\varepsilon'=\varepsilon_{Q'}$ is given in terms of the original matrix $\varepsilon = \varepsilon_Q$ through the {\em quiver mutation formula}
$$
\varepsilon'_{vw} = \left\{
\begin{array}{ll}
- \varepsilon_{vw}  & \mbox{if $k \in \{v,w\}$}, \\
\varepsilon_{vw} + \frac{1}{2}(\varepsilon_{vk} | \varepsilon_{kw}| + |\varepsilon_{vk}| \varepsilon_{kw}) & \mbox{if $k \not\in \{v,w\}$},
\end{array}
\right.
$$
and with the cluster $\mathscr{A}$-variables $A'_v$ for $v\in \mathcal{V}(Q')=\mathcal{V}(Q)$ given by the {\em cluster $\mathscr{A}$-mutation} formulas
$$
A'_v = \left\{
\begin{array}{ll}
A_v & \mbox{if $v\neq k$} \\
A_k^{-1}( \prod_{w\in \mathcal{V}(Q)} A_w^{[\varepsilon_{wk}]_+}  + \prod_{w\in \mathcal{V}(Q)} A_w^{[-\varepsilon_{wk}]_+}) & \mbox{if $v =k$},
\end{array}
\right.
$$
where $[\sim]_+$ is the {\em positive part}, i.e. $[a]_+ = a$ if $a\ge 0$ and $[a]_+ = 0$ if $a<0$. Similarly, a cluster $\mathscr{X}$-chart with the quiver $Q$ and the cluster $\mathscr{X}$-variables $X_v$, $v\in \mathcal{V}(Q)$, transforms via the mutation $\mu_k$ at the node $k$ to a cluster $\mathscr{X}$-chart with the quiver $\mu_k(Q)=Q'$ and the cluster $\mathscr{X}$-variables $X'_v$ given by
\begin{align}
\label{eq:intro_X-mutation}
X'_v = \left\{
\begin{array}{ll}
X_k^{-1} & \mbox{if $v=k$}, \\
X_v(1+X_k^{-{\rm sgn}(\varepsilon_{vk})})^{-\varepsilon_{vk}} & \mbox{if $v\neq k$},
\end{array}
\right.
\end{align}
where ${\rm sgn}(\sim)$ is the sign, i.e. ${\rm sgn}(a)=1$ if $a>0$ and ${\rm sgn}(a)=-1$ if $a<0$. For the current situation for the $n$-triangulation quivers, note that when the ideal triangulations $\Delta$ and $\Delta'$ of $\frak{S}$ are related by a {\em flip} at an arc, i.e. differ exactly by one arc, it is known that the $n$-triangulation quivers $Q_\Delta^{[n]}$ and $Q_{\Delta'}^{[n]}$ are related by a certain sequence of $\frac{1}{6}(n-1)n(n+1)$ mutations (see \cite{FG06}). When $n=3$ for example, first mutate $Q_\Delta^{[3]}$ at the two nodes lying in the arc of $\Delta$ that is to be flipped, then mutate at the two nodes lying in the interiors of the two triangles of $\Delta$ having the to-be-flipped arc as a side, to land in $Q_{\Delta'}^{[3]}$; see Fig.\ref{fig:mutations_for_a_flip}. 
The cluster $\mathscr{A}$-charts of $\mathscr{A}_{{\rm SL}_n,\frak{S}}$ for $\Delta$ and $\Delta'$ are related by the composition of the cluster $\mathscr{A}$-mutations for this same sequence of $\frac{1}{6}(n-1)n(n+1)$ mutations, and the cluster $\mathscr{X}$-charts of $\mathscr{P}_{{\rm PGL}_n,\frak{S}}$ for $\Delta$ and $\Delta'$ by the composition of the cluster $\mathscr{X}$-mutations for the same mutation sequence. Moreover, $\mathscr{P}_{{\rm PGL}_n,\frak{S}}$ is equipped with a canonical Poisson structure \cite{FG06} \cite{GS19}, given on each cluster $\mathscr{X}$-chart $\Gamma_\Delta$ for an ideal triangulation $\Delta$ by
$$
\{X_v,X_w\} = \varepsilon_{vw} X_v X_w, \qquad \forall v,w\in \mathcal{V}(Q_\Delta^{[n]}),
$$
where $\varepsilon = (\varepsilon_{vw})_{v,w\in \mathcal{V}(Q_\Delta^{[n]})}$ is the signed adjacency matrix of $Q_\Delta^{[n]}$.

\vs

One major line of research is on the quantization of the Poisson moduli space $\mathscr{P}_{{\rm PGL}_n,\frak{S}}$, or more precisely, of its cluster $\mathscr{X}$-variety structure. One first needs to construct a corresponding {\em quantum cluster $\mathscr{X}$-variety}, as a `non-commutative scheme'. There is such a formulation by Fock and Goncharov for a general cluster $\mathscr{X}$-variety \cite{FG09a} \cite{FG09b}. For each (classical) cluster $\mathscr{X}$-chart $\Gamma$ with the underlying quiver $Q$ (not necessarily the above discussed special cluster $\mathscr{X}$-chart $\Gamma_\Delta$ which has $Q_\Delta^{[n]}$ as its underlying quiver), consider the {\em Fock-Goncharov algebra} $\mathcal{X}^q_\Gamma$ defined as the associative algebra over $\mathbb{Z}[q^{\pm 1/18}]$ defined by
$$
\mbox{generators} ~:~ {\bf X}_v^{\pm 1}, ~ v\in \mathcal{V}(Q), \quad
\mbox{relations} ~:~ {\bf X}_v {\bf X}_w = q^{2\varepsilon_{vw}} {\bf X}_w {\bf X}_v, \quad {\bf X}_v{\bf X}_v^{-1} = {\bf X}_v^{-1}{\bf X}_v=1, \quad \forall v,w\in \mathcal{V}(Q).
$$
For the purpose of the present introduction section, a reader may think of this algebra as being defined over $\mathbb{Z}[q^{\pm 1/2}]$ instead of  $\mathbb{Z}[q^{\pm 1/18}]$; see Rem.\ref{rem:coefficient_ring}. This non-commutative algebra, which is an example of a quantum torus algebra, is what deforms the classical ring of functions on the chart $\Gamma$, namely the Laurent polynomial ring $\mathbb{Z}[\{X_v^{\pm 1} \, | \, v\in \mathcal{V}(Q)\}]$, in the direction of the above Poisson structure. For a mutation $\mu_k : \Gamma \leadsto \Gamma' = \mu_k(\Gamma)$, one would associate a {\em quantum mutation map} between the skew fields of fractions of the Fock-Goncharov algebras
$$
\mu^q_{\Gamma\Gamma'} = \mu^q_k ~:~ {\rm Frac}(\mathcal{X}^q_{\Gamma'}) \to {\rm Frac}(\mathcal{X}^q_\Gamma)
$$
so that it recovers the classical mutation formula as $q\to 1$, and that satisfies the consistency relations satisfied by their classical counterparts; namely, $\mu^q_k \mu^q_k = {\rm id}$ should hold for each initial cluster $\mathscr{X}$-chart $\Gamma$, $\mu^q_j \mu^q_k \mu^q_j \mu^q_k  = {\rm id}$ when $\varepsilon_{jk}=0$, and $\mu^q_j \mu^q_k \mu^q_j  \mu^q_k \mu^q_j = P_{(jk)}$ when $\varepsilon_{jk}=\pm 1$, where $P_{(jk)}$ stands for the label exchange $j \leftrightarrow k$ (Prop.\ref{prop:quantum_mutation_relations}). Such quantum mutation maps are found in \cite{FG09b}, based on earlier works such as \cite{BZ} \cite{CF99}, constituting a version of a quantum cluster $\mathscr{X}$-variety. In particular, the quantum isomorphism $\mu^q_{\Gamma\Gamma'}  : {\rm Frac}(\mathcal{X}^q_{\Gamma'}) \to {\rm Frac}(\mathcal{X}^q_\Gamma)$ can be constructed for each pair of cluster $\mathscr{X}$-charts $\Gamma$ and $\Gamma'$ in a consistent manner, by composing those for the mutations connecting $\Gamma$ and $\Gamma'$. For the case of $\mathscr{P}_{{\rm PGL}_n,\frak{S}}$, denote by
\begin{align}
\label{eq:intro_Phi_q}
\Phi^q_{\Delta\Delta'}  := \mu^q_{\Gamma_\Delta \Gamma_{\Delta'}} ~:~ {\rm Frac}(\mathcal{X}^q_{\Gamma_{\Delta'}}) \to {\rm Frac}(\mathcal{X}^q_{\Gamma_\Delta})
\end{align}
the quantum isomorphism for the cluster $\mathscr{X}$-charts $\Gamma_\Delta$ and $\Gamma_{\Delta'}$ for $\mathscr{P}_{{\rm PGL}_m,\frak{S}}$ associated to two ideal triangulations $\Delta$ and $\Delta'$ of $\frak{S}$. In particular, when $\Delta$ and $\Delta'$ are related by the flip at an arc, $\Phi^q_{\Delta\Delta'}$ is a composition of $\frac{1}{6}(n-1)n(n+1)$ number of quantum mutation maps $\mu^q_k$; see Def.\ref{def:Phi_q_i} for $n=3$.

\vs

We are interested in the problem of constructing a {\em deformation quantization map}, which is a map connecting the classical cluster $\mathscr{X}$-variety and the corresponding quantum cluster $\mathscr{X}$-variety. More precisely, it is an assignment to each `quantizable' classical observable function a corresponding quantum observable. One first needs to decide which classical functions to quantize, and the natural candidates would be the {\em universally Laurent} functions, i.e. the functions that are Laurent in all cluster $\mathscr{X}$-charts. In our case of $\mathscr{P}_{{\rm PGL}_n,\frak{S}}$, these form the ring denoted by $\mathcal{O}_{{\rm cl}}(\mathscr{P}_{{\rm PGL}_n,\frak{S}})$, which is proved in \cite{Shen} to equal the ring $\mathcal{O}(\mathscr{P}_{{\rm PGL}_n,\frak{S}})$ of regular functions on $\mathscr{P}_{{\rm PGL}_n,\frak{S}}$. Then a deformation quantization map would be a map
$$
\mathcal{O}(\mathscr{P}_{{\rm PGL}_n,\frak{S}}) \to \mathcal{O}^q(\mathscr{P}_{{\rm PGL}_n,\frak{S}})
$$
satisfying some conditions, where $\mathcal{O}^q(\mathscr{P}_{{\rm PGL}_n,\frak{S}})$ stands for the ring of all quantum universally Laurent elements, i.e. the intersection of all quantum Laurent polynomial rings $\mathcal{X}^q_\Gamma \subset {\rm Frac}(\mathcal{X}^q_\Gamma)$, where ${\rm Frac}(\mathcal{X}^q_\Gamma)$ for different $\Gamma$'s are identified via the quantum mutation maps $\mu^q_{\Gamma\Gamma'}$ in a consistent manner. One standard approach would be to first establish a {\it duality map}
$$
\mathbb{I} ~:~ \mathscr{A}_{{\rm SL}_n,\frak{S}}(\mathbb{Z}^T) \to \mathcal{O}(\mathscr{P}_{{\rm PGL}_n,\frak{S}}),
$$
whose existence was originally conjectured by Fock and Goncharov in \cite{FG06}, and whose image forms a basis of $\mathcal{O}(\mathscr{P}_{{\rm PGL}_n,\frak{S}})$, enumerated by the set $\mathscr{A}_{{\rm SL}_n,\frak{S}}(\mathbb{Z}^T)$ of $\mathbb{Z}^T$-points of $\mathscr{A}_{{\rm SL}_n,\frak{S}}$, where $\mathbb{Z}^T$ is the {\em semi-field of tropical integers} (see \S\ref{subsec:compatibility_of_quantum_duality_maps}), and then to establish a quantum duality map
$$
\mathbb{I}^q ~:~ \mathscr{A}_{{\rm SL}_n,\frak{S}}(\mathbb{Z}^T) \to \mathcal{O}^q(\mathscr{P}_{{\rm PGL}_n,\frak{S}}),
$$
which deforms $\mathbb{I}$ in a suitable sense. For a discussion on the domain set $\mathscr{A}_{{\rm SL}_n,\frak{S}}(\mathbb{Z}^T)$, we refer the readers to \cite{Kim} and to \S\ref{subsec:compatibility_of_quantum_duality_maps} of the present paper. Then one would construct a deformation quantization map by sending each basis element $\mathbb{I}(\ell)$ for $\ell \in \mathscr{A}_{{\rm SL}_n,\frak{S}}(\mathbb{Z}^T)$ to the corresponding element $\mathbb{I}^q(\ell)$. 

\vs

The setting of $n=2$ is referred to as the quantum Teichm\"uller theory; for punctured surfaces $\frak{S}$, a classical duality map $\mathbb{I}$ is constructed by Fock and Goncharov \cite{FG06}, and a quantum duality map $\mathbb{I}^q$ by Allegretti and the author \cite{AK}, based on Bonahon and Wong's ${\rm SL}_2$ quantum trace map \cite{BW}. These constructions heavily use geometry and topology of the surface $\frak{S}$. For other $n\ge 2$, and in fact for a much more general class of cluster $\mathscr{X}$-varieties, a duality map $\mathbb{I}$ is constructed by Gross, Hacking, Keel and Kontsevich \cite{GHKK} \cite{GS18}, and a quantum duality map $\mathbb{I}^q$ by Davison and Mandel \cite{DM}. These general constructions are very powerful when proving properties, but lack geometric intuition on surface geometry, and are quite difficult to compute. Even for the simplest possible punctured surfaces like the once-punctured torus, a direct computation has not been established yet, for a crucial ingredient called a `consistent scattering diagram' has not been described in a manner that can be used in a direct computation. In the meantime, a geometric and  straightforward-to-compute duality map for $n=3$ in the case of punctured surfaces $\frak{S}$ is constructed by the author in \cite{Kim}. Moreover, in \cite{Kim}, an {\em ${\rm SL}_3$ quantum trace map} is developed, and is used to construct a quantum duality map too. More precisely, as for the quantum duality maps, for each ideal triangulation $\Delta$ of a triangulable punctured surface $\frak{S}$, a map
\begin{align}
\label{eq:intro_I_q}
\mathbb{I}^q_\Delta ~:~ \mathscr{A}_{{\rm SL}_3,\frak{S}}(\mathbb{Z}^T) \to \mathcal{X}^q_{\Gamma_\Delta}
\end{align}
is constructed, and several nice properties are proved. One of the most important and fundamental properties for these $\mathbb{I}^q_\Delta$ is the {\it naturality}, or the compatibility under the change of ideal triangulations. This naturality, which was not proved and merely left as a conjecture in \cite{Kim}, is the major motivation of, as well as the major consequence of, the main theorem of the present paper. 
\begin{theorem}[main application: naturality of ${\rm SL}_3$-${\rm PGL}_3$ quantum duality maps]
\label{thm:intro_main_application}
Let $\frak{S}$ be a triangulable punctured surface. For any two ideal triangulations $\Delta$ and $\Delta'$ of $\frak{S}$ (without self-folded triangles), the ${\rm SL}_3$-${\rm PGL}_3$ quantum duality maps in eq.\eqref{eq:intro_I_q} for $\Delta$ and $\Delta'$, constructed in \cite{Kim}, are related by the quantum coordinate change map $\Phi^q_{\Delta\Delta'}$, i.e.
$$
\mathbb{I}^q_\Delta = \Phi^q_{\Delta\Delta'} \circ \mathbb{I}^q_{\Delta'}.
$$
\end{theorem}
This theorem, which implies that the above proposed deformation quantization map for the space $\mathscr{P}_{{\rm PGL}_3,\frak{S}}=\mathscr{X}_{{\rm PGL}_3,\frak{S}}$ for a triangulable punctured surface $\frak{S}$ is independent of the choice of an ideal triangulation $\Delta$ of $\frak{S}$, can be regarded as the principal result of the present paper, for a reader whose primary area is the theory of cluster varieties.

\subsection{Naturality of ${\rm SL}_3$ quantum trace maps}

\begin{figure}[htbp!]
\vspace{-4mm}
\begin{center}
\hspace*{-5mm}\begin{tabular}{ c | c | c }
\raisebox{-0.4\height}{
\begingroup%
  \makeatletter%
  \providecommand\color[2][]{%
    \errmessage{(Inkscape) Color is used for the text in Inkscape, but the package 'color.sty' is not loaded}%
    \renewcommand\color[2][]{}%
  }%
  \providecommand\transparent[1]{%
    \errmessage{(Inkscape) Transparency is used (non-zero) for the text in Inkscape, but the package 'transparent.sty' is not loaded}%
    \renewcommand\transparent[1]{}%
  }%
  \providecommand\rotatebox[2]{#2}%
  \newcommand*\fsize{\dimexpr\f@size pt\relax}%
  \newcommand*\lineheight[1]{\fontsize{\fsize}{#1\fsize}\selectfont}%
  \ifx\svgwidth\undefined%
    \setlength{\unitlength}{90.70866142bp}%
    \ifx\svgscale\undefined%
      \relax%
    \else%
      \setlength{\unitlength}{\unitlength * \real{\svgscale}}%
    \fi%
  \else%
    \setlength{\unitlength}{\svgwidth}%
  \fi%
  \global\let\svgwidth\undefined%
  \global\let\svgscale\undefined%
  \makeatother%
  \begin{picture}(1,0.28125)%
    \lineheight{1}%
    \setlength\tabcolsep{0pt}%
    \put(0.23728242,0.11890283){\color[rgb]{0,0,0}\makebox(0,0)[lt]{\lineheight{1.25}\smash{\begin{tabular}[t]{l}$=[3]_q {\O}=$\end{tabular}}}}%
    \put(0,0){\includegraphics[width=\unitlength,page=1]{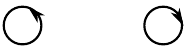}}%
  \end{picture}%
\endgroup%
} & \raisebox{-0.4\height}{
\begingroup%
  \makeatletter%
  \providecommand\color[2][]{%
    \errmessage{(Inkscape) Color is used for the text in Inkscape, but the package 'color.sty' is not loaded}%
    \renewcommand\color[2][]{}%
  }%
  \providecommand\transparent[1]{%
    \errmessage{(Inkscape) Transparency is used (non-zero) for the text in Inkscape, but the package 'transparent.sty' is not loaded}%
    \renewcommand\transparent[1]{}%
  }%
  \providecommand\rotatebox[2]{#2}%
  \newcommand*\fsize{\dimexpr\f@size pt\relax}%
  \newcommand*\lineheight[1]{\fontsize{\fsize}{#1\fsize}\selectfont}%
  \ifx\svgwidth\undefined%
    \setlength{\unitlength}{116.22047244bp}%
    \ifx\svgscale\undefined%
      \relax%
    \else%
      \setlength{\unitlength}{\unitlength * \real{\svgscale}}%
    \fi%
  \else%
    \setlength{\unitlength}{\svgwidth}%
  \fi%
  \global\let\svgwidth\undefined%
  \global\let\svgscale\undefined%
  \makeatother%
  \begin{picture}(1,0.2195122)%
    \lineheight{1}%
    \setlength\tabcolsep{0pt}%
    \put(0.40460658,0.0870799){\color[rgb]{0,0,0}\makebox(0,0)[lt]{\lineheight{1.25}\smash{\begin{tabular}[t]{l}$=-[2]_q$\end{tabular}}}}%
    \put(0,0){\includegraphics[width=\unitlength,page=1]{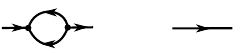}}%
  \end{picture}%
\endgroup%
} & \raisebox{-0.4\height}{
\begingroup%
  \makeatletter%
  \providecommand\color[2][]{%
    \errmessage{(Inkscape) Color is used for the text in Inkscape, but the package 'color.sty' is not loaded}%
    \renewcommand\color[2][]{}%
  }%
  \providecommand\transparent[1]{%
    \errmessage{(Inkscape) Transparency is used (non-zero) for the text in Inkscape, but the package 'transparent.sty' is not loaded}%
    \renewcommand\transparent[1]{}%
  }%
  \providecommand\rotatebox[2]{#2}%
  \newcommand*\fsize{\dimexpr\f@size pt\relax}%
  \newcommand*\lineheight[1]{\fontsize{\fsize}{#1\fsize}\selectfont}%
  \ifx\svgwidth\undefined%
    \setlength{\unitlength}{155.90551181bp}%
    \ifx\svgscale\undefined%
      \relax%
    \else%
      \setlength{\unitlength}{\unitlength * \real{\svgscale}}%
    \fi%
  \else%
    \setlength{\unitlength}{\svgwidth}%
  \fi%
  \global\let\svgwidth\undefined%
  \global\let\svgscale\undefined%
  \makeatother%
  \begin{picture}(1,0.30909091)%
    \lineheight{1}%
    \setlength\tabcolsep{0pt}%
    \put(0.3074048,0.1579046){\color[rgb]{0,0,0}\makebox(0,0)[lt]{\lineheight{1.25}\smash{\begin{tabular}[t]{l}$= $\end{tabular}}}}%
    \put(0,0){\includegraphics[width=\unitlength,page=1]{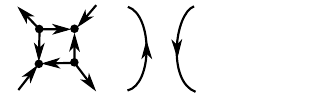}}%
    \put(0.62490477,0.1579046){\color[rgb]{0,0,0}\makebox(0,0)[lt]{\lineheight{1.25}\smash{\begin{tabular}[t]{l}$+$\end{tabular}}}}%
    \put(0,0){\includegraphics[width=\unitlength,page=2]{skein_rel7_corrected.pdf}}%
  \end{picture}%
\endgroup%
} \\
{\rm (S1)} & {\rm (S2)} & {\rm (S3)} \\  \hline
& \raisebox{-0.4\height}{
\begingroup%
  \makeatletter%
  \providecommand\color[2][]{%
    \errmessage{(Inkscape) Color is used for the text in Inkscape, but the package 'color.sty' is not loaded}%
    \renewcommand\color[2][]{}%
  }%
  \providecommand\transparent[1]{%
    \errmessage{(Inkscape) Transparency is used (non-zero) for the text in Inkscape, but the package 'transparent.sty' is not loaded}%
    \renewcommand\transparent[1]{}%
  }%
  \providecommand\rotatebox[2]{#2}%
  \newcommand*\fsize{\dimexpr\f@size pt\relax}%
  \newcommand*\lineheight[1]{\fontsize{\fsize}{#1\fsize}\selectfont}%
  \ifx\svgwidth\undefined%
    \setlength{\unitlength}{150.23622047bp}%
    \ifx\svgscale\undefined%
      \relax%
    \else%
      \setlength{\unitlength}{\unitlength * \real{\svgscale}}%
    \fi%
  \else%
    \setlength{\unitlength}{\svgwidth}%
  \fi%
  \global\let\svgwidth\undefined%
  \global\let\svgscale\undefined%
  \makeatother%
  \begin{picture}(1,0.28301887)%
    \lineheight{1}%
    \setlength\tabcolsep{0pt}%
    \put(0.23913077,0.12612741){\color[rgb]{0,0,0}\makebox(0,0)[lt]{\lineheight{1.25}\smash{\begin{tabular}[t]{l}$=q^{-2/3}$\end{tabular}}}}%
    \put(0,0){\includegraphics[width=\unitlength,page=1]{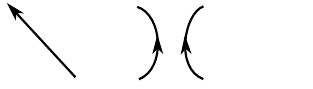}}%
    \put(0.65847036,0.12612741){\color[rgb]{0,0,0}\makebox(0,0)[lt]{\lineheight{1.25}\smash{\begin{tabular}[t]{l}$+q^{1/3}$\end{tabular}}}}%
    \put(0,0){\includegraphics[width=\unitlength,page=2]{skein_rel8_corrected.pdf}}%
  \end{picture}%
\endgroup%
} & \raisebox{-0.4\height}{
\begingroup%
  \makeatletter%
  \providecommand\color[2][]{%
    \errmessage{(Inkscape) Color is used for the text in Inkscape, but the package 'color.sty' is not loaded}%
    \renewcommand\color[2][]{}%
  }%
  \providecommand\transparent[1]{%
    \errmessage{(Inkscape) Transparency is used (non-zero) for the text in Inkscape, but the package 'transparent.sty' is not loaded}%
    \renewcommand\transparent[1]{}%
  }%
  \providecommand\rotatebox[2]{#2}%
  \newcommand*\fsize{\dimexpr\f@size pt\relax}%
  \newcommand*\lineheight[1]{\fontsize{\fsize}{#1\fsize}\selectfont}%
  \ifx\svgwidth\undefined%
    \setlength{\unitlength}{150.23622047bp}%
    \ifx\svgscale\undefined%
      \relax%
    \else%
      \setlength{\unitlength}{\unitlength * \real{\svgscale}}%
    \fi%
  \else%
    \setlength{\unitlength}{\svgwidth}%
  \fi%
  \global\let\svgwidth\undefined%
  \global\let\svgscale\undefined%
  \makeatother%
  \begin{picture}(1,0.28301887)%
    \lineheight{1}%
    \setlength\tabcolsep{0pt}%
    \put(0.23913077,0.12612741){\color[rgb]{0,0,0}\makebox(0,0)[lt]{\lineheight{1.25}\smash{\begin{tabular}[t]{l}$=q^{2/3}$\end{tabular}}}}%
    \put(0,0){\includegraphics[width=\unitlength,page=1]{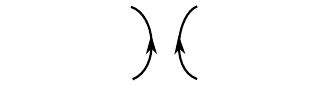}}%
    \put(0.65847036,0.12612741){\color[rgb]{0,0,0}\makebox(0,0)[lt]{\lineheight{1.25}\smash{\begin{tabular}[t]{l}$+q^{-1/3}$\end{tabular}}}}%
    \put(0,0){\includegraphics[width=\unitlength,page=2]{skein_rel9_corrected.pdf}}%
  \end{picture}%
\endgroup%
} \\ 
& {\rm (S4)} & {\rm (S5)} \\ 
\end{tabular}
\end{center}
\vspace{-3mm}
\caption{${\rm SL}_3$-skein relations, drawn locally (${\O}$ means empty) in $\frak{S}$, with the framing pointing toward the eyes of the reader; the regions bounded by a loop, a $2$-gon, or a $4$-gon in (S1), (S2), (S3) are contractible, and $[m]_q = \frac{q^m - q^{-m}}{q-q^{-1}} ~\in~\mathbb{Z}[q^{\pm 1}]$}
\vspace{-2mm}
\label{fig:A2-skein_relations_quantum}
\end{figure}

We now describe a more general statement, which we shall formulate as the actual main theorem. We first need to introduce the {\em ${\rm SL}_3$-skein algebra} \cite{S01} \cite{S05} \cite{FS}. For a generalized marked surface $\frak{S}$, consider the 3-dimensional manifold $\frak{S} \times {\bf I}$ called the {\em thickened surface} of $\frak{S}$, where
$$
{\bf I} = (-1,1)
$$
is the open interval in $\mathbb{R}$ whose elements are called {\em elevations}. Each boundary arc of $\frak{S}$ corresponds to a {\em boundary wall} $b\times {\bf I}$. 
An {\em ${\rm SL}_3$-web} $W$ in $\frak{S} \times {\bf I}$ (Def.\ref{def:SL3-web}), which goes back to \cite{Kuperberg} in its simplest case, is a disjoint union of oriented simple loops in $\frak{S} \times {\bf I}$, oriented edges in $\frak{S}\times {\bf I}$ ending at boundary walls, and oriented 3-valent graphs in $\frak{S}\times {\bf I}$ which may have endpoints at boundary walls, such that $W$ meets boundary walls transversally at 1-valent endpoints, the endpoints of $W$ lying in each boundary wall have mutually distinct elevations, and each 3-valent vertex is either a source or a sink. Also, $W$ is equipped with a framing. A {\em state} of $W$ is a map $s:\partial W \to \{1,2,3\}$, and $(W,s)$ is called a {\em stated ${\rm SL}_3$-web}. A {\em (reduced) stated ${\rm SL}_3$-skein algebra} $\mathcal{S}^\omega_{\rm s}(\frak{S};\mathbb{Z})_{\rm red}$ (Def.\ref{def:stated_SL3-skein_algebra}) \cite{Higgins} \cite{Kim} is defined as the free $\mathbb{Z}[\omega^{\pm 1/2}]$-module freely spanned by all isotopy classes of stated ${\rm SL}_3$-webs in $\frak{S} \times {\bf I}$, mod out by the {\em ${\rm SL}_3$-skein relations} in Fig.\ref{fig:A2-skein_relations_quantum} and the {\em boundary relations} in Fig.\ref{fig:stated_boundary_relations}, where $\omega$ is related to $q$ as
\begin{align}
\label{eq:intro_q_and_omega}
q = \omega^9,
\qquad \mbox{or more precisely, \qquad $q^{\pm 1/18} = \omega^{\pm 1/2}$}.
\end{align}
The product of $\mathcal{S}^\omega_{\rm s}(\frak{S};\mathbb{Z})_{\rm red}$ is defined by superposition, i.e. $[W_1,s_1] \cdot[W_2,s_2] = [W_1\cup W_2, s_1\cup s_2]$ when $W_1\subset\frak{S}\times (0,1)$ and $W_2\subset \frak{S}\times (-1,0)$, where $[W,s]$ denotes the element of $\mathcal{S}^\omega_{\rm s}(\frak{S};\mathbb{Z})_{\rm red}$ represented by the stated ${\rm SL}_3$-web $(W,s)$; we stack the former on top of the latter.

\vs

When $\frak{S}$ is a punctured surface, $\mathcal{S}^\omega_{\rm s}(\frak{S};\mathbb{Z})_{\rm red}$ can be understood as the {\em ${\rm SL}_3$-skein algebra} $\mathcal{S}^\omega(\frak{S};\mathbb{Z})$ which is defined just by the ${\rm SL}_3$-skein relations in Fig.\ref{fig:A2-skein_relations_quantum} for isotopy classes of ${\rm SL}_3$-webs in $\frak{S}\times {\bf I}$, without the boundary relations in Fig.\ref{fig:stated_boundary_relations} or states (Def.\ref{def:stated_SL3-skein_algebra}). 
It is known from \cite{S01} \cite{S05} that $\mathcal{S}^\omega(\frak{S};\mathbb{Z})$ is a quantum algebra deforming $\mathcal{O}(\mathscr{L}_{{\rm SL}_3,\frak{S}})$, the coordinate ring of the ${\rm SL}_3$-character stack. Similarly for the ${\rm SL}_2$ case, the ${\rm SL}_3$-skein algebras play a crucial role in the cluster-variety-theoretic study of the moduli spaces $\mathscr{A}_{{\rm SL}_3,\frak{S}}$ and $\mathscr{P}_{{\rm PGL}_3,\frak{S}}$, where the bridge to the world of cluster varieties is the family of maps
$$
{\rm Tr}^\omega_\Delta = {\rm Tr}^\omega_{\Delta;\frak{S}} ~:~ \mathcal{S}^\omega_{\rm s}(\frak{S};\mathbb{Z})_{\rm red} \to \mathcal{Z}^\omega_\Delta,
$$
called the {\em ${\rm SL}_3$ quantum trace maps} (Thm.\ref{thm:SL3_quantum_trace}) associated to each triangulable generalized marked surface $\frak{S}$ and its ideal triangulation $\Delta$, where $\mathcal{Z}^\omega_\Delta$ is a cube-root version of the Fock-Goncharov algebra, namely a $\mathbb{Z}[\omega^{\pm 1/2}]$-algebra defined by
$$
\mbox{generators} ~:~ {\bf Z}_v^{\pm 1}, ~ v\in \mathcal{V}(Q_\Delta^{[3]}), \quad
\mbox{relations} ~:~ {\bf Z}_v {\bf Z}_w = \omega^{2\varepsilon_{vw}} {\bf Z}_w {\bf Z}_v, \quad {\bf Z}_v{\bf Z}_v^{-1}={\bf Z}_v^{-1}{\bf Z}_v=1,\quad \forall v,w\in \mathcal{V}(Q_\Delta^{[3]}),
$$
where $\varepsilon=(\varepsilon_{vw})_{v,w\in \mathcal{V}(Q_\Delta^{[3]})}$ is the signed adjacency matrix of the quiver $Q_\Delta^{[3]}$. The usual Fock-Goncharov algebra $\mathcal{X}^q_\Delta:=\mathcal{X}^q_{\Gamma_\Delta}$ embeds into $\mathcal{Z}^\omega_\Delta$ as
$$
{\bf X}_v^{\pm 1} ~\mapsto~ {\bf Z}_v^{\pm 3}, \qquad \forall v\in \mathcal{V}(Q_\Delta^{[3]}).
$$
The ${\rm SL}_3$ quantum trace maps are constructed in \cite{Kim}, building on the work by Douglas \cite{Douglas1} on loops, as an ${\rm SL}_3$ analog of Bonahon and Wong's ${\rm SL}_2$ quantum trace \cite{BW}; see \cite{LY23} for a generalization to ${\rm SL}_n$ quantum trace maps. The motivating property is that it deforms the classical map from $\mathcal{S}^1(\frak{S};\mathbb{Z})$ to $\mathcal{O}(\mathscr{L}_{{\rm SL}_3,\frak{S}})$ (\cite{S01}) or to $\mathcal{O}(\mathscr{X}_{{\rm SL}_3,\frak{S}})$ (\cite{Kim}) in a certain sense; in particular, ${\rm Tr}^1_\Delta$ should yield the {\em trace-of-monodromy} functions along oriented loops, which are studied essentially by Fock and Goncharov \cite{FG06}. Another characterizing property of ${\rm Tr}^\omega_\Delta$ is the axiom about cutting and gluing (Thm.\ref{thm:SL3_quantum_trace}(QT1)), which says that the ${\rm SL}_3$ quantum trace maps are compatible with the process of cutting the surface $\frak{S}$ along an ideal arc of a triangulation $\Delta$. Although several favorable properties of the ${\rm SL}_3$ quantum trace maps are shown and used crucially in \cite{Kim}, one fundamental property was just conjectured but not proved in \cite{Kim}: namely, the naturality, or the compatibility under the change of ideal triangulations, which is the main theorem of the present paper.

\vs

\begin{figure}[htbp!]
\vspace{-4mm}
\begin{center}
\hspace*{-0mm}
\raisebox{-0.3\height}{\scalebox{0.9}{
\begingroup%
  \makeatletter%
  \providecommand\color[2][]{%
    \errmessage{(Inkscape) Color is used for the text in Inkscape, but the package 'color.sty' is not loaded}%
    \renewcommand\color[2][]{}%
  }%
  \providecommand\transparent[1]{%
    \errmessage{(Inkscape) Transparency is used (non-zero) for the text in Inkscape, but the package 'transparent.sty' is not loaded}%
    \renewcommand\transparent[1]{}%
  }%
  \providecommand\rotatebox[2]{#2}%
  \newcommand*\fsize{\dimexpr\f@size pt\relax}%
  \newcommand*\lineheight[1]{\fontsize{\fsize}{#1\fsize}\selectfont}%
  \ifx\svgwidth\undefined%
    \setlength{\unitlength}{113.38582677bp}%
    \ifx\svgscale\undefined%
      \relax%
    \else%
      \setlength{\unitlength}{\unitlength * \real{\svgscale}}%
    \fi%
  \else%
    \setlength{\unitlength}{\svgwidth}%
  \fi%
  \global\let\svgwidth\undefined%
  \global\let\svgscale\undefined%
  \makeatother%
  \begin{picture}(1,0.875)%
    \lineheight{1}%
    \setlength\tabcolsep{0pt}%
    \put(0,0){\includegraphics[width=\unitlength,page=1]{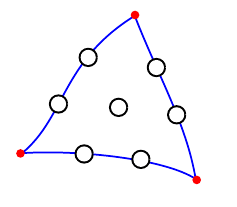}}%
    \put(0.25936061,0.11386532){\makebox(0,0)[lt]{\lineheight{1.25}\smash{\begin{tabular}[t]{l}$v_{e_1,2}$\end{tabular}}}}%
    \put(0.54607736,0.09319353){\makebox(0,0)[lt]{\lineheight{1.25}\smash{\begin{tabular}[t]{l}$v_{e_1,1}$\end{tabular}}}}%
    \put(0.03425753,0.42402765){\makebox(0,0)[lt]{\lineheight{1.25}\smash{\begin{tabular}[t]{l}$v_{e_2,1}$\end{tabular}}}}%
    \put(0.18603104,0.69997803){\makebox(0,0)[lt]{\lineheight{1.25}\smash{\begin{tabular}[t]{l}$v_{e_2,2}$\end{tabular}}}}%
    \put(0.81074466,0.37537559){\makebox(0,0)[lt]{\lineheight{1.25}\smash{\begin{tabular}[t]{l}$v_{e_3,2}$\end{tabular}}}}%
    \put(0.69476653,0.62447421){\makebox(0,0)[lt]{\lineheight{1.25}\smash{\begin{tabular}[t]{l}$v_{e_3,1}$\end{tabular}}}}%
    \put(0.46239583,0.48630232){\makebox(0,0)[lt]{\lineheight{1.25}\smash{\begin{tabular}[t]{l}$v_t$\end{tabular}}}}%
  \end{picture}%
\endgroup%
}} 
\end{center}
\vspace{-5mm}
\caption{Labels of the nodes of a $3$-triangulation quiver in a triangle}
\vspace{-2mm}
\label{fig:3-triangulation_node_names}
\end{figure}
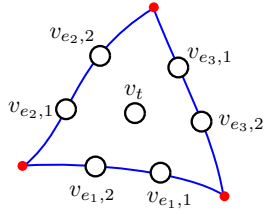

The first major step toward this naturality statement is to find a sensible formulation of it, which is already non-trivial because the values of the ${\rm SL}_3$ quantum trace maps are Laurent polynomials in the cube-root variables ${\bf Z}_v$'s, instead of the usual quantum cluster $\mathscr{X}$-variables ${\bf X}_v$'s. The transformation formulas for the latter variables ${\bf X}_v$ under the quantum mutation maps $\mu^q_k$ are certain non-commutative rational formulas deforming eq.\eqref{eq:intro_X-mutation} (see Def.\ref{def:FG_quantum_mutation}); in general, one would not expect that each ${\bf Z}_v$ would transform by rational formulas. It is only the elements of some subalgebra of ${\rm Frac}(\mathcal{Z}^\omega_\Delta)$  
that do transform via rational formulas. 
This subalgebra which we find in the present paper, as well as the characterizing condition for its elements, is called {\em balanced}, as they are the ${\rm SL}_3$ analog of Hiatt's balancedness condition for ${\rm SL}_2$ \cite{Hiatt} \cite{BW}, used for Bonahon and Wong's ${\rm SL}_2$ quantum trace \cite{BW}. The description of the ${\rm SL}_3$ balancedness condition is 
more complicated than that of ${\rm SL}_2$, and is inspired by the properties of the values of the {\em tropical coordinates} of {\em ${\rm SL}_3$-laminations} in $\frak{S}$ (Def.\ref{def:SL3-lamination}) \cite{DS1} \cite{Kim}. 
\begin{definition}[Def.\ref{def:Delta-balanced_elements}; \cite{Kim}]
\label{def:intro_balanced}
Let $\Delta$ be an ideal triangulation of a triangulable generalized marked surface $\frak{S}$. Let $\mathcal{V} = \mathcal{V}(Q_\Delta^{[3]})$. An element $(a_v)_{v\in \mathcal{V}} \in (\frac{1}{3}\mathbb{Z})^\mathcal{V}$ is said to be \ul{\em $\Delta$-balanced} if for each ideal triangle $t$ of $\Delta$, the following holds: denoting the sides of $t$ by $e_1,e_2,e_3$ (with $e_4:=e_1$), and the nodes of $Q_\Delta^{[3]}$ lying in $t$ by $v_{e_\alpha,1}$, $v_{e_\alpha,2}$ (for $\alpha=1,2,3$), and $v_t$ as in Fig.\ref{fig:3-triangulation_node_names},  one has
\begin{enumerate}\itemsep0em
\item[\rm (1)] the numbers $\sum_{\alpha=1}^3 a_{v_{e_\alpha,1}}$ and $\sum_{\alpha=1}^3 a_{v_{e_\alpha,2}}$ belong to $\mathbb{Z}$;

\item[\rm (2)] for each $\alpha = 1,2,3$, the number $a_{v_{e_\alpha,1}} + a_{v_{e_\alpha,2}}$ belongs to $\mathbb{Z}$;

\item[\rm (3)] for each $\alpha=1,2,3$, the number $-a_{v_t}+a_{v_{e_\alpha,2}} + a_{v_{e_{\alpha+1},1}}$ belongs to $\mathbb{Z}$.
\end{enumerate}
\end{definition}

We refer the readers to \cite{DS2} (see also \cite{DS1}) for a combinatorial and representation-theoretic formulation of the balancedness condition in Def.\ref{def:intro_balanced} in terms of the so-called Knutson-Tao rhombi \cite{KT99} \cite{GS15}.

\vspace{0mm}

\vspace{-1,5mm}

\begin{definition}[Def.\ref{def:balanced_subalgebras}--\ref{def:balanced_fraction_algebra}]
Let $\Delta$, $\frak{S}$ and $\mathcal{V}$ be as in Def.\ref{def:intro_balanced}. The \ul{\em $\Delta$-balanced cube-root Fock-Goncharov algebra} $\wh{\mathcal{Z}}^\omega_\Delta$ is the subalgebra of $\mathcal{Z}^\omega_\Delta$ spanned by the Laurent monomials $\prod_v {\bf X}_v^{a_v} := \prod_v {\bf Z}_v^{3a_v}$ with the powers forming a $\Delta$-balanced element $(a_v)_{v\in\mathcal{V}} \in (\frac{1}{3}\mathbb{Z})^\mathcal{V}$.  The \ul{\em $\Delta$-balanced fraction algebra} for $\Delta$ is the subalgebra $\wh{\rm Frac}(\mathcal{Z}^\omega_\Delta)$ of the skew field of fractions ${\rm Frac}(\mathcal{Z}^\omega_\Delta)$ consisting of all elements that can be written as ${\bf P}{\bf Q}^{-1}$ with ${\bf P}\in \wh{\mathcal{Z}}^\omega_\Delta\subset \mathcal{Z}^\omega_\Delta$ and $0\neq {\bf Q}\in \mathcal{X}^q_\Delta\subset\mathcal{Z}^\omega_\Delta$.
\end{definition}
In fact, one can identify $\wh{\rm Frac}(\mathcal{Z}^\omega_\Delta)$ with the skew field of fractions ${\rm Frac}(\wh{\mathcal{Z}}^\omega_\Delta)$ (Lem.\ref{lem:balanced_fraction_skew-field}). We note that $\mathcal{X}^q_\Delta \subset \wh{\mathcal{Z}}^\omega_\Delta$, as well as ${\rm Frac}(\mathcal{X}^q_\Delta)\subset \wh{\rm Frac}(\mathcal{Z}^\omega_\Delta)$.

\vspace{1mm}

We show that the quantum mutation maps in eq.\eqref{eq:intro_Phi_q} can be extended to these balanced fraction algebras.
\begin{proposition}[the balanced cube-root version of quantum coordinate change maps; \S\ref{subsec:the_balanced_algebras_and_quantum_coordinate_change_maps_for_them}]
Let $\frak{S}$ be a triangulable generalized marked surface. There is a family of algebra isomorphisms between the balanced fraction algebras
$$
\Theta^\omega_{\Delta\Delta'} ~:~ \wh{\rm Frac}(\mathcal{Z}^\omega_{\Delta'}) \to \wh{\rm Frac}(\mathcal{Z}^\omega_\Delta)
$$
defined for each pair of ideal triangulations $\Delta$ and $\Delta'$, that extend the maps $\Phi^q_{\Delta\Delta'}$, recover the classical coordinate change maps as $\omega^{1/2}\to 1$, and satisfy the consistency $\Theta^\omega_{\Delta\Delta''} = \Theta^\omega_{\Delta\Delta'} \Theta^\omega_{\Delta'\Delta''}$.
\end{proposition}
The formula for $\Theta^\omega_{\Delta\Delta'}$ is directly inspired by $\Phi^q_{\Delta\Delta'}$. When $\Delta$ and $\Delta'$ are related by the flip at an arc, while $\Phi^q_{\Delta\Delta'}$ is given as the composition of four quantum mutation maps $\mu^q_k$ (Def.\ref{def:Phi_q_i}) among the skew fields of fractions ${\rm Frac}(\mathcal{X}^q_\Gamma)$ of the Fock-Goncharov algebras $\mathcal{X}^q_\Gamma$, we construct $\Theta^\omega_{\Delta\Delta'}$ (Def.\ref{def:Theta}) as the composition of four balanced cube-root quantum mutation maps $\nu^\omega_k$ (Def.\ref{def:nu_omega_k}) defined among certain subalgebras $\wh{\rm Frac}_k(\mathcal{Z}^\omega_\Gamma)$ of ${\rm Frac}(\mathcal{Z}^\omega_\Gamma)$ (Def.\ref{def:Gamma_k_balanced}; see eq.\eqref{eq:FG_algebra_for_Gamma} and Def.\ref{def:FG_algebra} for $\mathcal{Z}^\omega_\Gamma$). One key thing to check is, in case one starts with an element of the $\Delta'$-balanced fraction algebra $\wh{\rm Frac}(\mathcal{Z}^\omega_{\Delta'})$ and tries to apply the four balanced cube-root quantum mutation maps $\nu^\omega_k$, whether the element at each step belongs to the domain $\wh{\rm Frac}_k(\mathcal{Z}^\omega_\Gamma)$ of $\nu^\omega_k$ (Lem.\ref{lem:Theta_well-definedness2}). We also check that the resulting element lies in the $\Delta$-balanced fraction algebra $\wh{\rm Frac}(\mathcal{Z}^\omega_\Delta)$ (Lem.\ref{lem:Theta_well-definedness2}). Then, to show the consistency relations (Prop.\ref{prop:consistency_for_Theta}) we resort to the known results on identities of classical and quantum cluster mutations (Lem.\ref{lem:classical_consistency_relations_for_flips}--\ref{lem:classical_consistency_relations_for_mutations}, Prop.\ref{prop:quantum_mutation_relations}, Prop.\ref{prop:quantum_consistency_of_flips}; \cite{FG06} \cite{KN}).

\vs

We can now state the main result of the present paper, which was conjectured in \cite{Kim} and also partially in an earlier work of Douglas \cite{Douglas1}.
\begin{theorem}[main theorem, Thm.\ref{thm:main}: naturality
of the ${\rm SL}_3$ quantum trace maps]
\label{thm:intro_main}
Let $\frak{S}$ be a triangulable generalized marked surface. For any two ideal triangulations $\Delta$ and $\Delta'$ of $\frak{S}$ (without self-folded triangles), the ${\rm SL}_3$ quantum trace maps for $\Delta$ and $\Delta'$ of \cite{Kim} are related by the balanced cube-root quantum coordinate change maps $\Theta^\omega_{\Delta\Delta'}$, i.e.
$$
{\rm Tr}^\omega_\Delta  = \Theta^\omega_{\Delta\Delta'} \circ {\rm Tr}^\omega_{\Delta'}.
$$
\end{theorem}
To prove this, we first establish the compatibility (Prop.\ref{prop:compatibility_of_Theta_under_cutting}) of the balanced coordinate change maps $\Theta^\omega_{\Delta\Delta'}$ with the cutting of the surfaces along internal arcs of ideal triangulations (Def.\ref{def:cutting_process}), and then using the cutting/gluing property of ${\rm Tr}^\omega_\Delta$ we reduce the situation to the case when $\frak{S}$ is a quadrilateral. Then, in fact we also use the cutting/gluing property of ${\rm Tr}^\omega_\Delta$ with respect to an ideal arc isotopic to a boundary arc $e$, so that $e$ cuts out a {\em biangle} (Def.\ref{def:ideal_triangulation}), to observe that it suffices to just check the above theorem for simple oriented edges living over a quadrilateral surface (\S\ref{subsec:the_base_case}). Still, a direct computational check would be quite involved, and we use several tricks to reduce the amount of the computations (\S\ref{subsec:classical_compatibility_and_Weyl-ordering}--\ref{subsec:checking_quantum_Laurent}). Namely, we use the equivariance of ${\rm Tr}^\omega_\Delta$ under the elevation-reversing map on $\mathcal{S}^\omega_{\rm s}(\frak{S};\mathbb{Z})_{\rm red}$ and the $*$-structure on $\mathcal{Z}^\omega_\Delta$ (Prop.\ref{prop:elevation_reversing_and_star}; \cite{Kim}), together with basic observations on $*$-invariant Laurent monomials of $\mathcal{Z}^\omega_\Delta$ (Lem.\ref{lem:star-invariance_and_Weyl-ordering_new}, Lem.\ref{lem:classicalization_determines_quantum}), so that what remains to check is whether ${\rm Tr}^\omega_\Delta([W,s])$ stays being Laurent and multiplicity-free as in Lem.\ref{lem:star-invariance_and_Weyl-ordering_new} in the cube-root quantum variables after applying the balanced cube-root quantum mutation {maps $\nu^\omega_k$ at  special nodes, which we verify carefully in \S\ref{subsec:reducing_the_amount_of_computation}--\S\ref{subsec:checking_quantum_Laurent}.

\vs

As mentioned above, perhaps the most interesting and important consequence of the main theorem, Thm.\ref{thm:intro_main}, is Thm.\ref{thm:intro_main_application} which we prove in \S\ref{subsec:compatibility_of_quantum_duality_maps} and which is about the quantum duality map for the space $\mathscr{X}_{{\rm PGL}_3,\frak{S}}$ for a triangulable punctured surface $\frak{S}$. However, we expect that Thm.\ref{thm:intro_main} itself would also serve as the first step toward a much wider range of future research topics, such as the representation theory for various versions of the ${\rm SL}_3$-skein algebras. One now has a consistent way of relating these algebras with various versions of Fock-Goncharov algebras, which are quantum torus algebras and hence admit a straightforward representation theory. Thus one might seek for the ${\rm SL}_3$ analogs of Bonahon and Wong's series of works on the similar topic for ${\rm SL}_2$ (see e.g. \cite{BW16}), which might also find applications in 3-dimensional topological quantum field theories or 2-dimensional conformal field theories. As suggested by a referee, we remark here that the main theorems of the present paper still hold when $q$ is a root of unity, which is the setting of several potential applications of these theorems, including the results obtained in \cite{BW16} \cite{BL20} (also \cite{KL}).

\vs

\noindent{\bf Acknowledgments.} This work was supported by the National Research Foundation of Korea(NRF) grant
funded by the Korea government(MSIT) (No. 2020R1C1C1A01011151). H.K. has been supported by a KIAS Individual Grant (MG047203, MG047204) at Korea Institute for Advanced Study. H.K. thanks Linhui Shen for helpful discussions. H.K. thanks the anonymous referees for their valuable remarks and suggestions to greatly improve the paper.

\section{${\rm SL}_3$ quantum trace maps}

In the present section, we recall the ${\rm SL}_3$ quantum trace maps from \cite{Kim}, as well as basic necessary notions from references therein.

\subsection{Surfaces and triangulations}
\label{subsec:surfaces_and_triangulations}

\begin{definition}[\cite{Le17} \cite{Le18}]
A \ul{\em generalized marked surface} $(\Sigma,\mathcal{P})$ is a pair of a compact oriented smooth surface $\Sigma$ with possibly empty boundary $\partial \Sigma$ and a non-empty finite subset $\mathcal{P}$ of $\Sigma$, such that each component of $\partial \Sigma$ contains at least one point of $\mathcal{P}$. Elements of $\mathcal{P}$ are called the \ul{\em marked points}, and the elements of $\mathcal{P}$ not lying in $\partial\Sigma$ are called the \ul{\em punctures}. When $\partial \Sigma = {\O}$, we say that $(\Sigma,\mathcal{P})$ is a \ul{\em punctured surface}.
\end{definition}
For a given generalized marked surface $(\Sigma,\mathcal{P})$, we often let
$$
\frak{S} = \Sigma\setminus\mathcal{P},
$$
and identify it with the data $(\Sigma,\mathcal{P})$, e.g. we refer to $\frak{S}$ as a generalized marked surface. Let
$$
\partial \frak{S} = (\partial \Sigma) \setminus \mathcal{P}, \qquad
\mathring{\frak{S}} = \frak{S} \setminus \partial \frak{S}.
$$
A basic ingredient is an {\em ideal triangulation} of a surface $\frak{S}$.
\begin{definition}[\cite{Le17} \cite{Le18}]
\label{def:ideal_triangulation}
Let $(\Sigma,\mathcal{P})$ be a generalized marked surface, and $\frak{S} = \Sigma\setminus\mathcal{P}$.
\begin{itemize}\itemsep0em
\item An \ul{\em ideal arc} in $\frak{S}$ is the image of an immersion $\alpha : [0,1] \to \Sigma$ such that $\alpha(\{0,1\}) \subset \mathcal{P}$ and $\alpha|_{(0,1)}$ is an embedding into $\frak{S}$. Call $\alpha((0,1))$ the \ul{\em interior} of this ideal arc. Two ideal arcs are \ul{\em isotopic} if they are isotopic within the class of ideal arcs. An ideal arc is called a \ul{\em boundary arc} if it lies in $\partial \Sigma$. An ideal arc is called an \ul{\em internal arc} if its interior lies in $\mathring{\frak{S}}$.

\item The generalized marked surface $\frak{S}$ is said to be \ul{\em triangulable} if it is none of the following:

-- \ul{\em monogon}, i.e. a closed disc with a single marked point on the boundary,

-- \ul{\em biangle}, i.e. a closed disc with two marked points on the boundary,

-- sphere with less than three punctures.

\item An \ul{\em ideal triangulation} of a triangulable generalized marked surface $\frak{S}$ is a collection $\Delta$ of ideal arcs of $\frak{S}$ such that

-- no arc of $\Delta$ bounds a disc whose interior is in $\frak{S}$;

-- no two arcs of $\Delta$ are isotopic or intersect each other in $\frak{S}$;

-- $\Delta$ is maximal among the collections satisfying the above two conditions.

We often identify two ideal triangulations if their members are simultaneously isotopic.

We assume that each constituent arc isotopic to a boundary arc is a boundary arc.
\end{itemize}

\end{definition}
An ideal triangulation $\Delta$ of $\frak{S}$ divides $\frak{S}$ into regions called \ul{\em (ideal) triangles} of $\Delta$, each of which is bounded by three ideal arcs, called the \ul{\em sides} of this triangle, counted with multiplicity. 
\begin{definition}
\label{def:regular_triangulation}
An ideal triangulation $\Delta$ of a triangulable generalized marked surface $\frak{S} = \Sigma\setminus\mathcal{P}$ is \ul{\em self-folded} if there exists a puncture $p$ of $(\Sigma,\mathcal{P})$ such that the valence of $\Delta$ at $p$ is $1$.
\end{definition}

\vspace{-5mm}

$$
\mbox{In the present paper, by an ideal triangulation we always mean a non-self-folded one.}
$$
That is, we assume that the valence of an ideal triangulation at each puncture is at least $2$. Also, by a \ul{\em triangulable} generalized marked surface $\frak{S}$ we mean a surface that admits a non-self-folded ideal triangulation. This means that we should further exclude the case of a monogon with one puncture. When there is a puncture with valence $1$ which occurs exactly at the self-folded side (of multiplicity $2$) of a self-folded triangle, the formulas, as well as proofs, need to be modified, as in the case of ${\rm SL}_2$ which was dealt with in \cite{BW}. A treatment of self-folded triangulations for ${\rm SL}_n$, especially ${\rm SL}_3$, will be established in an upcoming joint work \cite{JK} with Seung-Jo Jung.

\vs

Basic constructions will depend on the choice of an ideal triangulation of a surface, and the heart of the matter is to keep track of what happens if we use a different ideal triangulation. One standard approach in the literature is to deal with the `generators' of all possible changes of ideal triangulations, i.e. the following elementary changes.
\begin{definition}
Two ideal triangulations of a same generalized marked surface are said to be related by a \ul{\em flip at an arc} if, considered up to simultaneous isotopy, they differ precisely by one internal arc.
\end{definition}
We note that in the above definition, we still assume that both of the two ideal triangulations are non-self-folded. When $\Delta$ and $\Delta'$ are related by a flip at an arc, there is a natural bijection between $\Delta$ and $\Delta'$ as sets; each arc of $\Delta$ and the corresponding arc of $\Delta'$ are then denoted by the same symbol. In particular, we would use the same symbol for the flipped arc of $\Delta$ and that for $\Delta'$, although they are actually different as ideal arcs. Say, if the flipped arc is denoted by $i$, then we say $\Delta$ and $\Delta'$ are related by the flip at the arc $i$. We also say that $\Delta'$ is obtained from $\Delta$ by the flip move $\Phi_i$, and write
$$
\Delta' = \Phi_i (\Delta).
$$
A change of ideal triangulations is an ordered pair $(\Delta,\Delta')$ of ideal triangulations, which we often denote by $\Delta \leadsto \Delta'$. In case $\Delta$ and $\Delta'$ are related by a flip at an arc $i$, we denote this change by $\Phi_i$.
\begin{lemma}[{\cite{FST}, \cite[Cor.6.7]{LF}}]
\label{lem:flip_connectedness}
Any two ideal triangulations $\Delta$ and $\Delta'$ are connected by a finite sequence of flips. That is, $\Delta' = \Phi_{i_r} \cdots \Phi_{i_2} \Phi_{i_1} (\Delta)$.
\end{lemma}
In case we allow ideal triangulations to be self-folded, the statement of Lem.\ref{lem:flip_connectedness} is well known; see e.g. \cite{FST} and references therein. The statement when we only allow ideal triangulations that are not self-folded, as we are doing in the present paper, is somewhat less well known and is proved in \cite[Cor.6.7]{LF}.

\vs

The flips satisfy some algebraic relations; that is, sometimes when one applies a certain sequence of flips to a certain ideal triangulation, one gets back the same ideal triangulation. We find it convenient to first recall a well-known {\em signed adjacency matrix} for $\Delta$, which is a $|\Delta| \times |\Delta|$ integer matrix that encodes certain combinatorics of $\Delta$. For each ideal triangle $t$ of $\Delta$, if its sides are $e_1,e_2,e_3$ appearing clockwise in this order, then we say that $e_{i+1}$ is the \ul{\em clockwise next} one to $e_i$ (with $e_4:=e_1$).
\begin{definition}[see e.g. \cite{FST} and references therein]
\label{def:b_ij}
Let $\Delta$ be an ideal triangulation of a triangulable generalized marked surface. The \ul{\em signed adjacency matrix} $(b_{ij})_{i,j}$ of the triangulation $\Delta$ is the integer $\Delta\times \Delta$ matrix defined as
$$
b_{ij} = {\textstyle \sum}_t \, b_{ij}(t)
$$
where the sum is over all ideal triangles $t$ of $\Delta$, and
$$
b_{ij}(t) = \left\{
\begin{array}{ll}
1 & \mbox{if $i,j$ appear as sides of $t$ and $i$ is the clockwise next one to $j$}, \\
-1 & \mbox{if $i,j$ appear as sides of $t$ and $j$ is the clockwise next one to $i$}, \\
0 & \mbox{if at least one of $i,j$ is not a side of $t$}.
\end{array}
\right.
$$
\end{definition}
\begin{remark}
In fact, $(b_{ij})_{i,j}$ is the signed adjacency matrix for the $2$-triangulation quiver $Q_\Delta^{[2]}$ which appeared in the introduction.
\end{remark}

\begin{lemma}[classical consistency relations for flips of ideal triangulations]
\label{lem:classical_flip_relations_for_triangulations}
Fix a triangulable generalized marked surface $\frak{S}$. The flips $\Phi_i$ of ideal triangulations of $\frak{S}$ satisfy the following relations:
\begin{enumerate}\itemsep0em
\item[\rm (1)] $\Phi_i \Phi_i = {\rm id}$ when applied to any ideal triangulation;

\item[\rm (2)] $\Phi_i \Phi_j \Phi_i \Phi_j = {\rm id}$ when applied to an ideal triangulation $\Delta$ with $b_{ij}=0$;

\item[\rm (3)] $\Phi_i \Phi_j \Phi_i \Phi_j \Phi_i = P_{(ij)}$ when applied to an ideal triangulation $\Delta$ with $b_{ij}=\pm 1$, where $P_{(ij)}$ is the label exchange $i\leftrightarrow j$. \qed
\end{enumerate}
\end{lemma}
\begin{proposition}[the completeness of the flip relations]
\label{prop:completeness_of_classical_flip_relations_for_triangulations}
Any algebraic relation among flips is a consequence of the above. That is, any sequence of flips that starts and ends at the same ideal triangulation can be transformed to the empty sequence of flips by applying a finite number of the above three types of relations.
\end{proposition}
Both of these statements are well known; see e.g. \cite{FST} and references therein. We note that, if we use `tagged' ideal triangulations as in \cite{FST} which in particular include self-folded ideal triangulations, instead of just ideal triangulations that are not self-folded as stipulated by Def.\ref{def:regular_triangulation}, then there exist relations of flips that are not consequences of the above; see \cite{FST} \cite{KY}.

\subsection{${\rm SL}_3$-skein algebras and ${\rm SL}_3$-laminations}

\begin{definition}[\cite{S05} \cite{FS} \cite{Higgins} {\cite[Def.5.1]{Kim}}]
\label{def:SL3-web}
Let $(\Sigma,\mathcal{P})$ be a generalized marked surface, and $\frak{S} = \Sigma\setminus\mathcal{P}$. Let
$$
{\bf I} := (-1,1)
$$
be the open interval in $\mathbb{R}$, and let
$
\frak{S} \times {\bf I}
$
be the \ul{\em thickening} of $\frak{S}$, or a \ul{\em thickened surface}. For a point $(x,t) \in \frak{S} \times {\bf I}$, the ${\bf I}$-coordinate $t$ is called the \ul{\em elevation} of $(x,t)$. If $(x,t) \in A\times {\bf I}$ for some subset $A\subset \frak{S}$, we say $(x,t)$ \ul{\em lies over} $A$. For each boundary arc $b$ of $\frak{S}$, the corresponding boundary component $b\times {\bf I}$ of $\frak{S} \times {\bf I}$ is called a \ul{\em boundary wall}.

\vs

An \ul{\em ${\rm SL}_3$-web $W$ in $\frak{S} \times {\bf I}$} consists of

\begin{itemize}\itemsep0em
\item a finite subset of $(\partial \frak{S}) \times {\bf I}$ which we denote by $\partial W$, whose elements are called the \ul{\em external vertices} or the \ul{\em endpoints} of $W$; 

\item a finite subset of $\mathring{\frak{S}} \times {\bf I}$, whose elements are called the \ul{\em internal vertices} of $W$;

\item a finite set of oriented smooth simple non-closed compact curves in $\frak{S}\times {\bf I}$ whose interiors lie in $\mathring{\frak{S}}\times {\bf I}$ and whose ends are at external or internal vertices of $W$, whose elements are called the \ul{\em (oriented) edges} of $W$;

\item a finite set of oriented smooth simple closed curves in $\mathring{\frak{S}} \times {\bf I}$, whose elements are called the \ul{\em (oriented) loops} of $W$;

\item a \ul{\em framing} on $W$, when $W$ is regarded as the union of the constituent edges and loops, i.e. a continuous choice of an element of $T_x(\frak{S}\times {\bf I})\setminus T_x W$ for each point $x\in W$,
\end{itemize}
subject to the following conditions:
\begin{itemize}\itemsep0em
\item each external vertex is $1$-valent, and $W$ meets a boundary wall transversally at an external vertex; 

\item each internal vertex is either a 3-valent sink or a 3-valent source, i.e. for each internal vertex, the orientations of the three incident edges are either all incoming toward the vertex or all outgoing;

\item there is no self-intersection of $W$ except possibly at the 3-valent internal vertices;

\item the framing at each external vertex is \ul{\em upward vertical}, i.e. is parallel to the ${\bf I}$ factor and points toward $1$;

\item for each internal vertex $x$ of $W$, there is a diffeomorphism from a neighborhood $N$ of $x$ in $\frak{S}\times{\bf I}$ to $\mathbb{D}\times {\bf I}$ (where $\mathbb{D}$ is an open disc in $\mathbb{R}^2$) such that the image of $W \cap N$ lies in $\mathbb{D}\times \{0\}$ with an upward vertical framing;

\item for each boundary wall $b\times {\bf I}$, the endpoints of $W$ lying in $b\times {\bf I}$ have mutually distinct elevations.
\end{itemize}
An \ul{\em isotopy} of ${\rm SL}_3$-webs in $\frak{S} \times {\bf I}$ is an isotopy within the class of ${\rm SL}_3$-webs in $\frak{S} \times {\bf I}$.
\end{definition}

\begin{definition}[\cite{S05} \cite{FS} \cite{Higgins} {\cite[Def.5.3]{Kim}}]
\label{def:stated_SL3-skein_algebra}
Let $\frak{S}$ be a generalized marked surface. Let $\mathcal{R}$ be a commutative ring with unity.
\begin{itemize}\itemsep0em
\item The \ul{\em ${\rm SL}_3$-skein algebra} $\mathcal{S}^\omega(\frak{S};\mathcal{R})$ is the free $\mathcal{R}[\omega^{\pm 1/2}]$-module with the set of all isotopy classes of ${\rm SL}_3$-webs in $\frak{S} \times {\bf I}$ as a free basis, mod out by the ${\rm SL}_3$-skein relations in Fig.\ref{fig:A2-skein_relations_quantum}, with eq.\eqref{eq:intro_q_and_omega} in mind. The element of $\mathcal{S}^\omega(\frak{S};\mathcal{R})$ represented by an ${\rm SL}_3$-web $W$ in $\frak{S}\times {\bf I}$ is denoted by $[W]$.

\item A \ul{\em state} of an ${\rm SL}_3$-web $W$ in $\frak{S}\times {\bf I}$ is a map $s:\partial W \to \{1,2,3\}$. A pair $(W,s)$ is called a \ul{\em stated ${\rm SL}_3$-web} in $\frak{S} \times {\bf I}$.

\item The \ul{\em stated ${\rm SL}_3$-skein algebra} $\mathcal{S}^\omega_{\rm s}(\frak{S};\mathcal{R})$ is the free $\mathcal{R}[\omega^{\pm 1/2}]$-module with the set of all isotopy classes of stated ${\rm SL}_3$-webs in $\frak{S}\times {\bf I}$ as a free basis, mod out by the ${\rm SL}_3$-skein relations in Fig.\ref{fig:A2-skein_relations_quantum}, with eq.\eqref{eq:intro_q_and_omega} in mind.

\item The \ul{\em reduced stated ${\rm SL}_3$-skein algebra} $\mathcal{S}^\omega_{\rm s}(\frak{S};\mathcal{R})_{\rm red}$ is the quotient of $\mathcal{S}^\omega_{\rm s}(\frak{S};\mathcal{R})$ by the boundary relations in Fig.\ref{fig:stated_boundary_relations}, where the index-inversion $(r_1(\varepsilon),r_2(\varepsilon))$ for $\varepsilon \in \{1,2,3\}$ is given by
\begin{align}
\label{eq:intro_index_inversion}
(r_1(1),r_2(1)) = (1,2), \quad
(r_1(2),r_2(2)) = (1,3), \quad
(r_1(3),r_2(3)) = (2,3);
\end{align}
in the pictures, $x$ and $x_i$ are labels of endpoints, each picture is assumed to carry a respective state which is usually written as $s$. The element of $\mathcal{S}^\omega_{\rm s}(\frak{S};\mathcal{R})_{\rm red}$ (and that of $\mathcal{S}^\omega_{\rm s}(\frak{S};\mathcal{R})$) represented by a stated ${\rm SL}_3$-web $(W,s)$ in $\frak{S} \times {\bf I}$ is denoted by $[W,s]$.

\item  The multiplication in $\mathcal{S}^\omega_{\rm s}(\frak{S};\mathcal{R})_{\rm red}$ (and that in $\mathcal{S}^\omega_{\rm s}(\frak{S};\mathcal{R})$ or $\mathcal{S}^\omega(\frak{S};\mathcal{R})$) is given by superposition; i.e. $[W,s] \cdot [W',s'] = [W\cup W', s\cup s']$ (and $[W]\cdot[W'] =[W\cup W']$) when $W\subset \frak{S} \times (0,1)$ and $W' \subset \frak{S} \times (-1,0)$, where $s\cup s' : \partial W \cup \partial W' \to \{1,2,3\}$ is defined in an obvious manner.
\end{itemize}
\end{definition}

\begin{figure}[htbp!]
\begin{center}
\hspace*{0mm}\begin{tabular}{c|c}
\raisebox{-0.6\height}{
\begingroup%
  \makeatletter%
  \providecommand\color[2][]{%
    \errmessage{(Inkscape) Color is used for the text in Inkscape, but the package 'color.sty' is not loaded}%
    \renewcommand\color[2][]{}%
  }%
  \providecommand\transparent[1]{%
    \errmessage{(Inkscape) Transparency is used (non-zero) for the text in Inkscape, but the package 'transparent.sty' is not loaded}%
    \renewcommand\transparent[1]{}%
  }%
  \providecommand\rotatebox[2]{#2}%
  \newcommand*\fsize{\dimexpr\f@size pt\relax}%
  \newcommand*\lineheight[1]{\fontsize{\fsize}{#1\fsize}\selectfont}%
  \ifx\svgwidth\undefined%
    \setlength{\unitlength}{42.51968504bp}%
    \ifx\svgscale\undefined%
      \relax%
    \else%
      \setlength{\unitlength}{\unitlength * \real{\svgscale}}%
    \fi%
  \else%
    \setlength{\unitlength}{\svgwidth}%
  \fi%
  \global\let\svgwidth\undefined%
  \global\let\svgscale\undefined%
  \makeatother%
  \begin{picture}(1,1.16666667)%
    \lineheight{1}%
    \setlength\tabcolsep{0pt}%
    \put(0,0){\includegraphics[width=\unitlength,page=1]{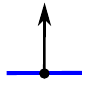}}%
    \put(0.42722329,0.10283054){\color[rgb]{0,0,0}\makebox(0,0)[lt]{\lineheight{1.25}\smash{\begin{tabular}[t]{l}$x$\end{tabular}}}}%
  \end{picture}%
\endgroup%
} \hspace{-3mm}  $=-q^{ - \frac{7}{6}}$ \hspace{-5mm} \raisebox{-0.6\height} {
\begingroup%
  \makeatletter%
  \providecommand\color[2][]{%
    \errmessage{(Inkscape) Color is used for the text in Inkscape, but the package 'color.sty' is not loaded}%
    \renewcommand\color[2][]{}%
  }%
  \providecommand\transparent[1]{%
    \errmessage{(Inkscape) Transparency is used (non-zero) for the text in Inkscape, but the package 'transparent.sty' is not loaded}%
    \renewcommand\transparent[1]{}%
  }%
  \providecommand\rotatebox[2]{#2}%
  \newcommand*\fsize{\dimexpr\f@size pt\relax}%
  \newcommand*\lineheight[1]{\fontsize{\fsize}{#1\fsize}\selectfont}%
  \ifx\svgwidth\undefined%
    \setlength{\unitlength}{42.51968504bp}%
    \ifx\svgscale\undefined%
      \relax%
    \else%
      \setlength{\unitlength}{\unitlength * \real{\svgscale}}%
    \fi%
  \else%
    \setlength{\unitlength}{\svgwidth}%
  \fi%
  \global\let\svgwidth\undefined%
  \global\let\svgscale\undefined%
  \makeatother%
  \begin{picture}(1,1.16666667)%
    \lineheight{1}%
    \setlength\tabcolsep{0pt}%
    \put(0,0){\includegraphics[width=\unitlength,page=1]{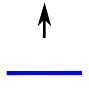}}%
    \put(0.16942518,0.13006975){\color[rgb]{0,0,0}\makebox(0,0)[lt]{\lineheight{1.25}\smash{\begin{tabular}[t]{l}$x_1$\end{tabular}}}}%
    \put(0,0){\includegraphics[width=\unitlength,page=2]{boundary_rel2_quantum.pdf}}%
    \put(0.74016697,0.12462211){\color[rgb]{0,0,0}\makebox(0,0)[lt]{\lineheight{1.25}\smash{\begin{tabular}[t]{l}$x_2$\end{tabular}}}}%
    \put(0.4330944,0.12306698){\color[rgb]{0,0,0}\makebox(0,0)[lt]{\lineheight{1.25}\smash{\begin{tabular}[t]{l}$\prec$\end{tabular}}}}%
  \end{picture}%
\endgroup%
} \hspace{0,5mm} & \raisebox{-0.6\height}{
\begingroup%
  \makeatletter%
  \providecommand\color[2][]{%
    \errmessage{(Inkscape) Color is used for the text in Inkscape, but the package 'color.sty' is not loaded}%
    \renewcommand\color[2][]{}%
  }%
  \providecommand\transparent[1]{%
    \errmessage{(Inkscape) Transparency is used (non-zero) for the text in Inkscape, but the package 'transparent.sty' is not loaded}%
    \renewcommand\transparent[1]{}%
  }%
  \providecommand\rotatebox[2]{#2}%
  \newcommand*\fsize{\dimexpr\f@size pt\relax}%
  \newcommand*\lineheight[1]{\fontsize{\fsize}{#1\fsize}\selectfont}%
  \ifx\svgwidth\undefined%
    \setlength{\unitlength}{42.51968504bp}%
    \ifx\svgscale\undefined%
      \relax%
    \else%
      \setlength{\unitlength}{\unitlength * \real{\svgscale}}%
    \fi%
  \else%
    \setlength{\unitlength}{\svgwidth}%
  \fi%
  \global\let\svgwidth\undefined%
  \global\let\svgscale\undefined%
  \makeatother%
  \begin{picture}(1,1.16666667)%
    \lineheight{1}%
    \setlength\tabcolsep{0pt}%
    \put(0,0){\includegraphics[width=\unitlength,page=1]{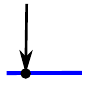}}%
    \put(0.18027884,0.10283054){\color[rgb]{0,0,0}\makebox(0,0)[lt]{\lineheight{1.25}\smash{\begin{tabular}[t]{l}$x_1$\end{tabular}}}}%
    \put(0,0){\includegraphics[width=\unitlength,page=2]{boundary_rel3_quantum.pdf}}%
    \put(0.67416761,0.10283054){\color[rgb]{0,0,0}\makebox(0,0)[lt]{\lineheight{1.25}\smash{\begin{tabular}[t]{l}$x_2$\end{tabular}}}}%
    \put(0.46250103,0.10283054){\color[rgb]{0,0,0}\makebox(0,0)[lt]{\lineheight{1.25}\smash{\begin{tabular}[t]{l}$\prec$\end{tabular}}}}%
  \end{picture}%
\endgroup%
} \hspace{-2mm} $=$ $q$ \hspace{-2mm}
 \raisebox{-0.6\height}{
\begingroup%
  \makeatletter%
  \providecommand\color[2][]{%
    \errmessage{(Inkscape) Color is used for the text in Inkscape, but the package 'color.sty' is not loaded}%
    \renewcommand\color[2][]{}%
  }%
  \providecommand\transparent[1]{%
    \errmessage{(Inkscape) Transparency is used (non-zero) for the text in Inkscape, but the package 'transparent.sty' is not loaded}%
    \renewcommand\transparent[1]{}%
  }%
  \providecommand\rotatebox[2]{#2}%
  \newcommand*\fsize{\dimexpr\f@size pt\relax}%
  \newcommand*\lineheight[1]{\fontsize{\fsize}{#1\fsize}\selectfont}%
  \ifx\svgwidth\undefined%
    \setlength{\unitlength}{42.51968504bp}%
    \ifx\svgscale\undefined%
      \relax%
    \else%
      \setlength{\unitlength}{\unitlength * \real{\svgscale}}%
    \fi%
  \else%
    \setlength{\unitlength}{\svgwidth}%
  \fi%
  \global\let\svgwidth\undefined%
  \global\let\svgscale\undefined%
  \makeatother%
  \begin{picture}(1,1.16666667)%
    \lineheight{1}%
    \setlength\tabcolsep{0pt}%
    \put(0,0){\includegraphics[width=\unitlength,page=1]{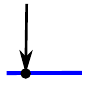}}%
    \put(0.18027884,0.10283054){\color[rgb]{0,0,0}\makebox(0,0)[lt]{\lineheight{1.25}\smash{\begin{tabular}[t]{l}$x_2$\end{tabular}}}}%
    \put(0,0){\includegraphics[width=\unitlength,page=2]{boundary_rel4_quantum.pdf}}%
    \put(0.67416761,0.10283054){\color[rgb]{0,0,0}\makebox(0,0)[lt]{\lineheight{1.25}\smash{\begin{tabular}[t]{l}$x_1$\end{tabular}}}}%
    \put(0.46250103,0.10283054){\color[rgb]{0,0,0}\makebox(0,0)[lt]{\lineheight{1.25}\smash{\begin{tabular}[t]{l}$\prec$\end{tabular}}}}%
  \end{picture}%
\endgroup%
}
 \hspace{-3mm} $+$ \hspace{-2mm}
\raisebox{-0.6\height}{
\begingroup%
  \makeatletter%
  \providecommand\color[2][]{%
    \errmessage{(Inkscape) Color is used for the text in Inkscape, but the package 'color.sty' is not loaded}%
    \renewcommand\color[2][]{}%
  }%
  \providecommand\transparent[1]{%
    \errmessage{(Inkscape) Transparency is used (non-zero) for the text in Inkscape, but the package 'transparent.sty' is not loaded}%
    \renewcommand\transparent[1]{}%
  }%
  \providecommand\rotatebox[2]{#2}%
  \newcommand*\fsize{\dimexpr\f@size pt\relax}%
  \newcommand*\lineheight[1]{\fontsize{\fsize}{#1\fsize}\selectfont}%
  \ifx\svgwidth\undefined%
    \setlength{\unitlength}{42.51968504bp}%
    \ifx\svgscale\undefined%
      \relax%
    \else%
      \setlength{\unitlength}{\unitlength * \real{\svgscale}}%
    \fi%
  \else%
    \setlength{\unitlength}{\svgwidth}%
  \fi%
  \global\let\svgwidth\undefined%
  \global\let\svgscale\undefined%
  \makeatother%
  \begin{picture}(1,1.16666667)%
    \lineheight{1}%
    \setlength\tabcolsep{0pt}%
    \put(0,0){\includegraphics[width=\unitlength,page=1]{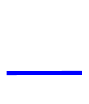}}%
    \put(0.13088999,0.10283052){\color[rgb]{0,0,0}\makebox(0,0)[lt]{\lineheight{1.25}\smash{\begin{tabular}[t]{l}$x_2$\end{tabular}}}}%
    \put(0,0){\includegraphics[width=\unitlength,page=2]{boundary_rel5_quantum.pdf}}%
    \put(0.72355657,0.10283052){\color[rgb]{0,0,0}\makebox(0,0)[lt]{\lineheight{1.25}\smash{\begin{tabular}[t]{l}$x_1$\end{tabular}}}}%
    \put(0,0){\includegraphics[width=\unitlength,page=3]{boundary_rel5_quantum.pdf}}%
    \put(0.44133449,0.10283052){\color[rgb]{0,0,0}\makebox(0,0)[lt]{\lineheight{1.25}\smash{\begin{tabular}[t]{l}$\prec$\end{tabular}}}}%
  \end{picture}%
\endgroup%
}  \\
{\rm (B1)} $s(x)=\varepsilon$, $s(x_1)=r_1(\varepsilon)$, $s(x_2)=r_2(\varepsilon)$ 
& {\rm (B2)} $s(x_1)=\varepsilon_1$, $s(x_2)= \varepsilon_2$, with  $\varepsilon_1>\varepsilon_2$  \\ \hline
\raisebox{-0.6\height}{
\begingroup%
  \makeatletter%
  \providecommand\color[2][]{%
    \errmessage{(Inkscape) Color is used for the text in Inkscape, but the package 'color.sty' is not loaded}%
    \renewcommand\color[2][]{}%
  }%
  \providecommand\transparent[1]{%
    \errmessage{(Inkscape) Transparency is used (non-zero) for the text in Inkscape, but the package 'transparent.sty' is not loaded}%
    \renewcommand\transparent[1]{}%
  }%
  \providecommand\rotatebox[2]{#2}%
  \newcommand*\fsize{\dimexpr\f@size pt\relax}%
  \newcommand*\lineheight[1]{\fontsize{\fsize}{#1\fsize}\selectfont}%
  \ifx\svgwidth\undefined%
    \setlength{\unitlength}{42.51968504bp}%
    \ifx\svgscale\undefined%
      \relax%
    \else%
      \setlength{\unitlength}{\unitlength * \real{\svgscale}}%
    \fi%
  \else%
    \setlength{\unitlength}{\svgwidth}%
  \fi%
  \global\let\svgwidth\undefined%
  \global\let\svgscale\undefined%
  \makeatother%
  \begin{picture}(1,1.16666667)%
    \lineheight{1}%
    \setlength\tabcolsep{0pt}%
    \put(0,0){\includegraphics[width=\unitlength,page=1]{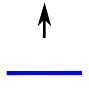}}%
    \put(0.16616776,0.10283052){\color[rgb]{0,0,0}\makebox(0,0)[lt]{\lineheight{1.25}\smash{\begin{tabular}[t]{l}$x$\end{tabular}}}}%
    \put(0,0){\includegraphics[width=\unitlength,page=2]{boundary_rel6_quantum.pdf}}%
    \put(0.72355657,0.10283052){\color[rgb]{0,0,0}\makebox(0,0)[lt]{\lineheight{1.25}\smash{\begin{tabular}[t]{l}$y$\end{tabular}}}}%
    \put(0.44133449,0.10283052){\color[rgb]{0,0,0}\makebox(0,0)[lt]{\lineheight{1.25}\smash{\begin{tabular}[t]{l}$\prec$\end{tabular}}}}%
  \end{picture}%
\endgroup%
} \hspace{-2mm} $=0$
& \raisebox{-0.6\height}{
\begingroup%
  \makeatletter%
  \providecommand\color[2][]{%
    \errmessage{(Inkscape) Color is used for the text in Inkscape, but the package 'color.sty' is not loaded}%
    \renewcommand\color[2][]{}%
  }%
  \providecommand\transparent[1]{%
    \errmessage{(Inkscape) Transparency is used (non-zero) for the text in Inkscape, but the package 'transparent.sty' is not loaded}%
    \renewcommand\transparent[1]{}%
  }%
  \providecommand\rotatebox[2]{#2}%
  \newcommand*\fsize{\dimexpr\f@size pt\relax}%
  \newcommand*\lineheight[1]{\fontsize{\fsize}{#1\fsize}\selectfont}%
  \ifx\svgwidth\undefined%
    \setlength{\unitlength}{42.51968504bp}%
    \ifx\svgscale\undefined%
      \relax%
    \else%
      \setlength{\unitlength}{\unitlength * \real{\svgscale}}%
    \fi%
  \else%
    \setlength{\unitlength}{\svgwidth}%
  \fi%
  \global\let\svgwidth\undefined%
  \global\let\svgscale\undefined%
  \makeatother%
  \begin{picture}(1,1.16666667)%
    \lineheight{1}%
    \setlength\tabcolsep{0pt}%
    \put(0,0){\includegraphics[width=\unitlength,page=1]{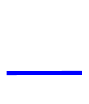}}%
    \put(0.04567213,0.10283052){\color[rgb]{0,0,0}\makebox(0,0)[lt]{\lineheight{1.25}\smash{\begin{tabular}[t]{l}$x_1$\end{tabular}}}}%
    \put(0,0){\includegraphics[width=\unitlength,page=2]{boundary_rel7_quantum.pdf}}%
    \put(0.44839,0.10283052){\color[rgb]{0,0,0}\makebox(0,0)[lt]{\lineheight{1.25}\smash{\begin{tabular}[t]{l}$x_2$\end{tabular}}}}%
    \put(0.83644529,0.10283052){\color[rgb]{0,0,0}\makebox(0,0)[lt]{\lineheight{1.25}\smash{\begin{tabular}[t]{l}$x_3$\end{tabular}}}}%
    \put(0.23672332,0.10283052){\color[rgb]{0,0,0}\makebox(0,0)[lt]{\lineheight{1.25}\smash{\begin{tabular}[t]{l}$\prec$\end{tabular}}}}%
    \put(0.66005655,0.10283052){\color[rgb]{0,0,0}\makebox(0,0)[lt]{\lineheight{1.25}\smash{\begin{tabular}[t]{l}$\prec$\end{tabular}}}}%
  \end{picture}%
\endgroup%
} \hspace{-2mm} $=$ $-q^{\frac{7}{2}}$ \hspace{-5mm} \raisebox{-0.6\height}{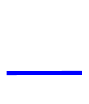} \\ 
{\rm (B3)} $s(x)=s(y)$ & {\rm (B4)} $s(x_1)=1$, $s(x_2)=2$, $s(x_3)=3$
\end{tabular}
\end{center}
\vspace{-2mm}
\caption{Boundary relations for stated ${\rm SL}_3$-skeins (horizontal blue line is boundary); the endpoints in the figure are consecutive in the elevation ordering for that boundary component (i.e. $\nexists$ other endpoint with elevation in between these), and $x\prec y$ means $y$ has a higher elevation than $x$; see eq.\eqref{eq:intro_index_inversion} for the definition of $r_1,r_2$ appearing in (B1).}
\vspace{-3mm}
\label{fig:stated_boundary_relations}
\end{figure}

The main object of study of the present paper is the ${\rm SL}_3$ quantum trace map which is to be reviewed at the end of the present section; the domain of this map is the reduced stated ${\rm SL}_3$-skein algebra $\mathcal{S}^\omega_{\rm s}(\frak{S};\mathbb{Z})_{\rm red}$. It is useful to notice that in case $\frak{S}$ is a triangulable punctured surface, in particular has no boundary, then $\mathcal{S}^\omega_{\rm s}(\frak{S};\mathcal{R})_{\rm red}$ coincides with $\mathcal{S}^\omega(\frak{S};\mathcal{R})$.

\vs

The above boundary relations (B1)--(B4) are used in \cite{Kim} as a modified version of the ones in \cite{Higgins}. Meanwhile, Frohman and Sikora \cite{FS} used somewhat different boundary relations to define their version of the `reduced' (non-stated) skein algebra denoted by $\mathcal{RS}(\frak{S})$, and they found a basis of this algebra consisting of {\em reduced non-elliptic} ${\rm SL}_3$-webs, and constructed a coordinate system for such ${\rm SL}_3$-webs. These basic ${\rm SL}_3$-webs have no {\em crossings} and have upward vertical framing everywhere; hence, they can be projected down to $\frak{S}$ via the projection map
\begin{align}
\label{eq:pi}
\pi ~:~\frak{S} \times {\bf I} \to \frak{S}
\end{align}
and viewed as objects living in the surface $\frak{S}$; see Def.\ref{def:web_in_surface}(BE2) below for a precise meaning of having no crossings. Generalizing these reduced non-elliptic ${\rm SL}_3$-webs living in a surface $\frak{S}$, the notion of {\em${\rm SL}_3$-laminations} in $\frak{S}$ is defined and studied in \cite{Kim}, and a coordinate system on them is established in \cite{Kim} based on Douglas and Sun's coordinates \cite{DS1} which are certain modification of Frohman and Sikora's coordinates \cite{FS}. We review these in a concise manner.
\begin{definition}[modified from \cite{Kuperberg} \cite{SW} \cite{FS}] 
\label{def:web_in_surface}
Let $\frak{S}$ be a generalized marked surface. Let $W$ be an ${\rm SL}_3$-web in $\frak{S}\times {\bf I}$ (without a chosen state), such that 
\begin{enumerate}\itemsep0em
\item[\rm (NE1)] the framing is upward vertical everywhere;

\item[\rm (NE2)] $W$ has no crossing, in the sense that the restriction $\pi|_W : W \to \pi(W)$ of the projection $\pi$ to $W$ is one-to-one;
\end{enumerate}
The projection $\pi(W)$ in $\frak{S}$ is called an \ul{\em ${\rm SL}_3$-web in the surface $\frak{S}$}. If furthermore the following condition is satisfied:
\begin{enumerate}
\item[\rm (NE3)] in $\pi(W)$ there is no contractible region bounded by a loop, a $2$-gon or a $4$-gon, as appearing in (S1)--(S3) of Fig.\ref{fig:A2-skein_relations_quantum},
\end{enumerate}
then $\pi(W)$ is said to be \ul{\em non-elliptic}. In addition to (NE1), (NE2) and (NE3), if furthermore the following condition is satisfied:
\begin{enumerate}
\item[\rm (NE4)] in $\pi(W)$ there is no boundary $2$-gon \raisebox{-0.3\height}{\scalebox{0.7}{
\begingroup%
  \makeatletter%
  \providecommand\color[2][]{%
    \errmessage{(Inkscape) Color is used for the text in Inkscape, but the package 'color.sty' is not loaded}%
    \renewcommand\color[2][]{}%
  }%
  \providecommand\transparent[1]{%
    \errmessage{(Inkscape) Transparency is used (non-zero) for the text in Inkscape, but the package 'transparent.sty' is not loaded}%
    \renewcommand\transparent[1]{}%
  }%
  \providecommand\rotatebox[2]{#2}%
  \newcommand*\fsize{\dimexpr\f@size pt\relax}%
  \newcommand*\lineheight[1]{\fontsize{\fsize}{#1\fsize}\selectfont}%
  \ifx\svgwidth\undefined%
    \setlength{\unitlength}{42.51968504bp}%
    \ifx\svgscale\undefined%
      \relax%
    \else%
      \setlength{\unitlength}{\unitlength * \real{\svgscale}}%
    \fi%
  \else%
    \setlength{\unitlength}{\svgwidth}%
  \fi%
  \global\let\svgwidth\undefined%
  \global\let\svgscale\undefined%
  \makeatother%
  \begin{picture}(1,0.56666667)%
    \lineheight{1}%
    \setlength\tabcolsep{0pt}%
    \put(0,0){\includegraphics[width=\unitlength,page=1]{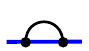}}%
  \end{picture}%
\endgroup%
}}, $3$-gon \raisebox{-0.3\height}{\scalebox{0.7}{
\begingroup%
  \makeatletter%
  \providecommand\color[2][]{%
    \errmessage{(Inkscape) Color is used for the text in Inkscape, but the package 'color.sty' is not loaded}%
    \renewcommand\color[2][]{}%
  }%
  \providecommand\transparent[1]{%
    \errmessage{(Inkscape) Transparency is used (non-zero) for the text in Inkscape, but the package 'transparent.sty' is not loaded}%
    \renewcommand\transparent[1]{}%
  }%
  \providecommand\rotatebox[2]{#2}%
  \newcommand*\fsize{\dimexpr\f@size pt\relax}%
  \newcommand*\lineheight[1]{\fontsize{\fsize}{#1\fsize}\selectfont}%
  \ifx\svgwidth\undefined%
    \setlength{\unitlength}{42.51968504bp}%
    \ifx\svgscale\undefined%
      \relax%
    \else%
      \setlength{\unitlength}{\unitlength * \real{\svgscale}}%
    \fi%
  \else%
    \setlength{\unitlength}{\svgwidth}%
  \fi%
  \global\let\svgwidth\undefined%
  \global\let\svgscale\undefined%
  \makeatother%
  \begin{picture}(1,0.56666667)%
    \lineheight{1}%
    \setlength\tabcolsep{0pt}%
    \put(0,0){\includegraphics[width=\unitlength,page=1]{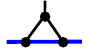}}%
  \end{picture}%
\endgroup%
}}, or $4$-gon \raisebox{-0.3\height}{\scalebox{0.7}{
\begingroup%
  \makeatletter%
  \providecommand\color[2][]{%
    \errmessage{(Inkscape) Color is used for the text in Inkscape, but the package 'color.sty' is not loaded}%
    \renewcommand\color[2][]{}%
  }%
  \providecommand\transparent[1]{%
    \errmessage{(Inkscape) Transparency is used (non-zero) for the text in Inkscape, but the package 'transparent.sty' is not loaded}%
    \renewcommand\transparent[1]{}%
  }%
  \providecommand\rotatebox[2]{#2}%
  \newcommand*\fsize{\dimexpr\f@size pt\relax}%
  \newcommand*\lineheight[1]{\fontsize{\fsize}{#1\fsize}\selectfont}%
  \ifx\svgwidth\undefined%
    \setlength{\unitlength}{42.51968504bp}%
    \ifx\svgscale\undefined%
      \relax%
    \else%
      \setlength{\unitlength}{\unitlength * \real{\svgscale}}%
    \fi%
  \else%
    \setlength{\unitlength}{\svgwidth}%
  \fi%
  \global\let\svgwidth\undefined%
  \global\let\svgscale\undefined%
  \makeatother%
  \begin{picture}(1,0.56666667)%
    \lineheight{1}%
    \setlength\tabcolsep{0pt}%
    \put(0,0){\includegraphics[width=\unitlength,page=1]{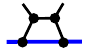}}%
  \end{picture}%
\endgroup%
}},
\end{enumerate}
then $\pi(W)$ is said to be \ul{\em reduced}. An isotopy of non-elliptic ${\rm SL}_3$-webs in $\frak{S}$ means an isotopy within the class of non-elliptic ${\rm SL}_3$-webs in $\frak{S}$.
\end{definition}
\begin{definition}[\cite{Kim}]
\label{def:SL3-lamination}
Let $\frak{S}$ be a generalized marked surface. 

$\bullet$ A simple loop in $\frak{S}$ is called a \ul{\em peripheral loop} if it bounds a region in $\frak{S}$ diffeomorphic to a disc with one puncture in the interior; if the corresponding puncture is $p\in \mathcal{P}$, we say that this peripheral loop \ul{\em surrounds} $p$. A \ul{\em peripheral arc} in $\frak{S}$ is a simple curve in $\frak{S}$ that ends at points of $\partial \frak{S}$ and bounds a region in $\frak{S}$ diffeomorphic to an upper half-disc with one puncture on the boundary. Peripheral loops and peripheral arcs are called \ul{\em peripheral curves}.

\vs

An \ul{\em ${\rm SL}_3$-lamination} $\ell$ in $\frak{S}$ is the isotopy class of a reduced non-elliptic web $W=W(\ell)$ in $\frak{S}$ equipped with integer weights on the components, subject to the following conditions and the equivalence relation: 
\begin{enumerate}\itemsep0em
\item[\rm (L1)] the weight of each component of $W$ containing an internal $3$-valent vertex is $1$;

\item[\rm (L2)] the weight of each component of $W$ that is not a peripheral curve is non-negative;

\item[\rm (L3)] an ${\rm SL}_3$-lamination containing a component of weight zero is equivalent to the ${\rm SL}_3$-lamination with this component removed;

\item[\rm (L4)] an ${\rm SL}_3$-lamination with two of its components being homotopic, respecting orientations, with weights $a$ and $b$ is equivalent to the ${\rm SL}_3$-lamination with one of these components removed and the other having weight $a+b$.
\end{enumerate}
Let $\mathscr{A}_{\rm L}(\frak{S};\mathbb{Z})$ be the set of all ${\rm SL}_3$-laminations in $\frak{S}$.
\end{definition}

A statement about a coordinate system on $\mathscr{A}_{\rm L}(\frak{S};\mathbb{Z})$ is postponed until \S\ref{subsec:the_balanced_algebras_and_quantum_coordinate_change_maps_for_them}.

\subsection{${\rm PGL}_3$ Fock-Goncharov algebras for surfaces}

A \ul{\em quiver} $Q$ consists of a set $\mathcal{V}(Q)$ of \ul{\em nodes} and a set $\mathcal{E}(Q)$ of \ul{\em arrows} between the nodes, where an arrow is an ordered pair $(v,w)$ of nodes, depicted in pictures as $\overset{v}{\circ} \hspace{-1,3mm} \longrightarrow \hspace{-1,3mm}  \overset{w}{\circ}$.  The \ul{\em signed adjacency matrix} of a quiver $Q$ is the $\mathcal{V}(Q) \times \mathcal{V}(Q)$ matrix $\varepsilon_Q = \varepsilon$ whose $(v,w)$-th entry is
$$
\varepsilon_{vw} = \varepsilon_{v,w} = (\mbox{number of arrows from $v$ to $w$}) - (\mbox{number of arrows from $w$ to $v$}).
$$
If a quiver $Q'$ can be obtained from a quiver $Q$ by deleting a cycle of length $1$ or $2$ (i.e. either an arrow of the form $(v,v)$ or the pair of arrows of the form $(v,w)$ and $(w,v)$), we say $Q$ and $Q'$ are \ul{\em equivalent}; this generates an equivalence relation on the set of all quivers. The set of equivalence classes of all quivers for a fixed set of nodes $\mathcal{V}$ is in bijection with the set of all skew-symmetric $\mathcal{V}\times \mathcal{V}$ integer matrices.

\vs

Let's consider a \ul{\em generalized quiver} $Q$ based on a set $\mathcal{V}$ of nodes, which corresponds to a skew-symmetric $\mathcal{V} \times \mathcal{V}$ matrix with entries in $\frac{1}{2} \mathbb{Z}$. This can be thought of as a collection of \ul{\em half-arrows} $\overset{v}{\circ} \hspace{-0,7mm} \dashrightarrow \hspace{-0,7mm}  \overset{w}{\circ}$, so that the signed adjacency matrix is given by
$$
\varepsilon_{vw} = {\textstyle \frac{1}{2}} (\mbox{number of half-arrows from $v$ to $w$}) - {\textstyle \frac{1}{2}}  (\mbox{number of half-arrows from $w$ to $v$}).
$$
In practice, one can define a generalized quiver as a collection of half-arrows and (usual solid) arrows, and consider an equivalence relation generated by the move deleting a cycle of half-arrows of length $1$ or $2$, and the move replacing two half-arrows from $v$ to $w$ by an (usual solid) arrow from $v$ to $w$. In particular, quivers are generalized quivers.  We will identify two generalized quivers if they are equivalent, unless there is a confusion.

\vs

Throughout the paper, when we say a quiver we will mean a generalized quiver.

\begin{definition}
\label{def:3-triangulation_quiver}
Let $\Delta$ be an ideal triangulation of a triangulable generalized marked surface $\frak{S}$. The \ul{\em (extended) $3$-triangulation quiver $Q_\Delta$} for $\Delta$ is a quiver defined as follows. The set of nodes $\mathcal{V}(Q_\Delta)$ is realized as a subset of $\frak{S}$, consisting of one point in the interior of each ideal triangle of $\Delta$ and two points lying in the interior of each ideal arc of $\Delta$. The quiver $Q_\Delta = Q_\Delta^{[3]}$ is obtained by gluing (i.e. taking the union of) all the quivers defined for the ideal triangles of $\Delta$ as in Fig.\ref{fig:n-triangulation}. Denote by $\varepsilon_\Delta = \varepsilon$ the signed adjacency matrix for $Q_\Delta$.
\end{definition}

We now present the quantum algebra of the Fock-Goncharov-Shen (cluster) Poisson moduli space $\mathscr{P}_{{\rm PGL}_3,\frak{S}}$ \cite{FG06} \cite{GS19} considered in \cite{FG09a} \cite{FG09b} \cite{GS19}, and its cube-root version considered in \cite{Douglas1} \cite{Kim}. For later use, we find it convenient to define these algebras for more general setting as follows, not just for ideal triangulations of surfaces, and also to present their standard $*$-structures at this point.

\begin{definition}
\label{def:FG_algebra}
Let $Q$ be a quiver. Denote by $\mathcal{V} = \mathcal{V}(Q)$ its set of nodes, and by $\varepsilon = (\varepsilon_{vw})_{v,w\in \mathcal{V}}$ its signed adjacency matrix.

\vs

Define the \ul{\em Fock-Goncharov algebra} $\mathcal{X}^q_Q$ for $Q$ as the free associative $*$-algebra over $\mathbb{Z}[q^{\pm 1/18}]$ generated by $\{{\bf X}_v^{\pm 1} \, | \, v\in \mathcal{V}\}$ mod out by the relations
$$
{\bf X}_v {\bf X}_w = q^{2\varepsilon_{vw}} {\bf X}_w {\bf X}_v, \quad \forall v,w\in \mathcal{V}, \qquad {\bf X}_v {\bf X}_v^{-1} = {\bf X}_v^{-1} {\bf X}_v = 1, \quad \forall v\in \mathcal{V},
$$
where the $*$-structure on $\mathcal{X}^q_Q$ is defined to be the unique anti-ring-homomorphism $* : \mathcal{X}^q_Q \to \mathcal{X}^q_Q$, ${\bf U} \mapsto {\bf U}^*$, sending each ${\bf X}_v^\epsilon$ to itself, $\forall v\in \mathcal{V}$, $\forall \epsilon \in \{1,-1\}$, and $q^{\pm 1/18}$ to $q^{\mp 1/18}$.

\vs

Define the \ul{\em cube-root Fock-Goncharov algebra} $\mathcal{Z}^\omega_Q$ for $Q$ by the free associative $*$-algebra over $\mathbb{Z}[\omega^{\pm 1/2}]$ generated by $\{{\bf Z}_v^{\pm 1} \, | \, v\in \mathcal{V}\}$ mod out by the relations
$$
{\bf Z}_v {\bf Z}_w = \omega^{2\varepsilon_{vw}} {\bf Z}_w {\bf Z}_v, \quad \forall v,w\in \mathcal{V}, \qquad
{\bf Z}_v {\bf Z}_v^{-1} = {\bf Z}_v^{-1} {\bf Z}_v = 1, \quad \forall v\in \mathcal{V},
$$
where the $*$-structure on $\mathcal{Z}^\omega_Q$ is the unique anti-ring-homomorphism $*:\mathcal{Z}^\omega_Q \to \mathcal{Z}^\omega_Q$, ${\bf U} \mapsto {\bf U}^*$, sending each ${\bf Z}_v^\epsilon$ to itself, $\forall v\in \mathcal{V}$, $\forall \epsilon \in\{1,-1\}$, and $\omega^{\pm 1/2}$ to $\omega^{\mp 1/2}$. 

\vs

The Fock-Goncharov algebra $\mathcal{X}^q_Q$ is regarded as being embedded into $\mathcal{Z}^\omega_Q$ as
$$
\mathcal{X}^q_Q \hookrightarrow \mathcal{Z}^\omega_Q ~:~ {\bf X}_v^\epsilon \mapsto {\bf Z}_v^{3\epsilon}, ~ \forall v\in \mathcal{V}, ~ \forall \epsilon \in\{1,-1\}, \qquad
q^{\pm 1/18} \mapsto \omega^{\pm 1/2}.
$$
We will denote ${\bf Z}_v^{\pm 1}$ by ${\bf X}_v^{\pm 1/3}$. For $a\in \frac{1}{3}\mathbb{Z}$, the symbol ${\bf X}_v^a$ means the element ${\bf Z}_v^{3a} \in \mathcal{Z}^\omega_Q$.

\vs

For each of $\mathcal{X}^q_Q$ and $\mathcal{Z}^\omega_Q$, the $*$-structure is also called the $*$-map, and an element ${\bf U}$ (of $\mathcal{X}^q_Q$ or $\mathcal{Z}^\omega_Q$) is said to be \ul{\em $*$-invariant} if it is fixed by the $*$-map, that is, if ${\bf U}^* = {\bf U}$.

\vs

For an ideal triangulation $\Delta$ of a triangulable generalized marked surface $\frak{S}$, we write
$$
\mathcal{X}^q_\Delta := \mathcal{X}^q_{Q_\Delta} \quad\mbox{and}\quad \mathcal{Z}^\omega_\Delta := \mathcal{Z}^\omega_{Q_\Delta}
$$
for convenience.
\end{definition}

\begin{lemma}
\label{lem:relations_of_cube-root_algebra_by_X}
Under the notation convention of the above definition, ${\bf X}_v^{a_v} {\bf X}_w^{a_w} = q^{2\varepsilon_{vw} a_v a_w} {\bf X}_w^{a_w} {\bf X}_v^{a_v}$ holds for all $v,w\in \mathcal{V}$ and $a_v,a_w\in \frac{1}{3}\mathbb{Z}$, and ${\bf Z}_v^{\alpha_v} {\bf Z}_w^{\alpha_w} = \omega^{2\varepsilon_{vw} \alpha_v \alpha_w} {\bf Z}_w^{\alpha_w} {\bf Z}_v^{\alpha_v}$ holds for all $v,w\in \mathcal{V}$ and $\alpha_v,\alpha_w\in\mathbb{Z}$. \qed
\end{lemma}

Some versions of these algebras appeared in the literature before the works \cite{FG09a} \cite{FG09b} \cite{GS19}, and are often called quantum torus algebras. The Fock-Goncharov algebras $\mathcal{X}^q_Q$ will serve as the quantum algebras that quantize the moduli space $\mathscr{P}_{{\rm PGL}_3,\frak{S}}$ in the end, but the quantization method suggested in \cite{Kim} makes a crucial use of the cube-root Fock-Goncharov algebras $\mathcal{Z}^\omega_Q$ as an intermediate step. The present paper deals mostly with $\mathcal{Z}^\omega_Q$ and a certain subalgebra of it that properly contains $\mathcal{X}^q_Q$, so the terminology we use focuses more on $\mathcal{Z}^\omega_Q$ rather than $\mathcal{X}^q_Q$. For example, by a `(quantum) Laurent polynomial' for $Q$ or $\Delta$ we would mean an element of $\mathcal{Z}^\omega_Q$ or $\mathcal{Z}^\omega_\Delta$, instead of an element of $\mathcal{X}^q_Q$ or $\mathcal{X}^q_\Delta$. 
\begin{remark}
\label{rem:coefficient_ring}
It is expected in \cite[\S5]{Kim} that the `correct' quantum algebras $\mathcal{X}^q_Q$ that one would want in the end should be the quantum torus algebras defined over $\mathbb{Z}[q^{\pm 1/2}]$ instead of over $\mathbb{Z}[q^{\pm 1/18}]$ (in fact, over $\mathbb{Z}[q^{\pm 1}]$ for punctured surfaces $\frak{S}$, without boundary), but at the moment the arguments work only over the coefficient ring $\mathbb{Z}[q^{\pm 1/18}]$. See Prop.\ref{prop:value_of_the_SL3_quantum_trace_is_Delta-balanced} which is a result of \cite{Kim} written only with respect to the coefficient ring $\mathbb{Z}[\omega^{\pm 1/2}]=\mathbb{Z}[q^{\pm 1/18}]$, and also the proof of Lem.\ref{lem:nu_omega_k_preserves_star} which requires the usage of this coefficient ring as opposed to $\mathbb{Z}[q^{\pm 1/2}]$.
\end{remark}

\vs

\begin{definition}
\label{def:Laurent}
Let $Q$, $\mathcal{V} = \mathcal{V}(Q)$, $\mathcal{X}^q_Q$ and $\mathcal{Z}^\omega_Q$ be as in Def.\ref{def:FG_algebra}. Enumerate the elements of $\mathcal{V}$ as $v_1,v_2,\ldots,v_N$, according to any chosen ordering on $\mathcal{V}$. 

\vs

A \ul{\em Laurent polynomial for $Q$} is an element of $\mathcal{Z}^\omega_Q$. A \ul{\em Laurent monomial for $Q$} is an element of $\mathcal{Z}^\omega_Q$ of the form
$$
\epsilon \, \omega^m \, {\bf X}_{v_1}^{\alpha_{v_1}/3} {\bf X}_{v_2}^{\alpha_{v_2}/3} \cdots {\bf X}_{v_N}^{\alpha_{v_N}/3}
= \epsilon \, \omega^m \, {\bf Z}_{v_1}^{\alpha_{v_1}} {\bf Z}_{v_2}^{\alpha_{v_2}} \cdots {\bf Z}_{v_N}^{\alpha_{v_N}}
$$
for some $\epsilon\in\{1,-1\}$, $m\in \frac{1}{2}\mathbb{Z}$ and $\alpha_{v_1},\ldots,\alpha_{v_N} \in \mathbb{Z}$.

\vs

For an ideal triangulation $\Delta$ of a triangulable generalized marked surface $\frak{S}$, a \ul{\em Laurent polynomial for $\Delta$} is a Laurent polynomial for $Q_\Delta$, and a \ul{\em Laurent monomial for $\Delta$} is a Laurent monomial for $Q_\Delta$.
\end{definition}
It is easy to see that the above definition does not depend on the choice of an ordering on $\mathcal{V}$, and that each Laurent polynomial for $Q$ can be expressed as a finite sum of Laurent monomials for $Q$.

\vs

One convenient technical tool used in the present paper is the following well-known notion.
\begin{definition}
\label{def:Weyl-ordered}
Let $Q$, $\mathcal{V} = \mathcal{V}(Q)$ and $\mathcal{Z}^\omega_Q$ be as  in Def.\ref{def:FG_algebra}. Let $\varepsilon=(\varepsilon_{vw})_{v,w\in\mathcal{V}}$ be the signed adjacency matrix for $Q$. Enumerate the elements of $\mathcal{V}$ as $v_1,v_2,\ldots,v_N$, according to any chosen ordering on $\mathcal{V}$. 

\vs

A Laurent monomial $\epsilon \, \omega^m {\bf Z}_{v_1}^{\alpha_{v_1}} \cdots {\bf Z}_{v_N}^{\alpha_{v_N}}$ for $Q$, with $\epsilon\in\{1,-1\}$, $m\in \frac{1}{2}\mathbb{Z}$, $\alpha_{v_1},\ldots,\alpha_{v_N} \in \mathbb{Z}$, is said to be \ul{\em Weyl-ordered} if $\epsilon=1$ and $m = -\sum_{1\le i<j\le N} \varepsilon_{v_iv_j} \alpha_{v_i} \alpha_{v_j}$. We denote such a Laurent monomial by
$$
\left[ {\textstyle \prod}_{v\in \mathcal{V}} {\bf X}_v^{\alpha_v/3} \right]_{\rm Weyl}
= \left[ {\textstyle \prod}_{v\in \mathcal{V}} {\bf Z}_v^{\alpha_v} \right]_{\rm Weyl}
:= \omega^{-\sum_{i<j} \varepsilon_{v_iv_j} \alpha_{v_i} \alpha_{v_j}} {\bf Z}_{v_1}^{\alpha_{v_1}} \cdots {\bf Z}_{v_N}^{\alpha_{v_N}},
$$
associated to $(\alpha_v)_{v\in \mathcal{V}} \in \mathbb{Z}^\mathcal{V}$.

\vs

Let ${\bf f} \in \mathcal{Z}^\omega_Q$ be a Laurent polynomial for $Q$. Express ${\bf f}$ as a finite sum of Laurent monomials $\epsilon \, \omega^m {\bf Z}_{v_1}^{\alpha_{v_1}} \ldots {\bf Z}_{v_N}^{\alpha_{v_N}}$ for $Q$. Denote by $[{\bf f}]_{\rm Weyl}\in \mathcal{Z}^\omega_Q$ what is obtained from this sum expression of ${\bf f}$ by replacing each summand term $\epsilon \, \omega^m {\bf Z}_{v_1}^{\alpha_{v_1}} \ldots {\bf Z}_{v_N}^{\alpha_{v_N}}$ by $\epsilon \, [{\bf Z}_{v_1}^{\alpha_{v_1}} \ldots {\bf Z}_{v_N}^{\alpha_{v_N}}]_{\rm Weyl}$. Call $[{\bf f}]_{\rm Weyl}$ the \ul{\em (term-by-term) Weyl-ordering} of ${\bf f}$. Such an element $[{\bf f}]_{\rm Weyl}$ of $\mathcal{Z}^\omega_Q$ is called a \ul{\em (term-by-term) Weyl-ordered Laurent polynomial for $Q$}.

\vs

For a matrix ${\bf M}$ with entries in $\mathcal{Z}^\omega_Q$, denote by $[{\bf M}]_{\rm Weyl}$ the matrix obtained by replacing each entry of ${\bf M}$ by its (term-by-term) Weyl-ordering.
\end{definition}

It is a well-known straightforward exercise to show that the Weyl-ordered Laurent monomial denoted by $\left[ {\textstyle \prod}_{v\in \mathcal{V}} {\bf X}_v^{\alpha_v/3} \right]_{\rm Weyl}
= \left[ {\textstyle \prod}_{v\in \mathcal{V}} {\bf Z}_v^{\alpha_v} \right]_{\rm Weyl}$ is a well-defined element of $\mathcal{Z}^\omega_Q$ independent of the choice of ordering on $\mathcal{V} = \mathcal{V}(Q_\Delta)$. It is also easy to see that $[{\bf f}]_{\rm Weyl}$ depends only on ${\bf f}$ but not on the choice of an expression of ${\bf f}$ as a sum of Laurent monomials.

\vs

The following lemma is easily verifiable as well, and will become handy. In particular, it gives one characterization of Weyl-ordered Laurent monomials for $Q$.

\begin{lemma}
\label{lem:star-invariance_and_Weyl-ordering_new}
Let $Q$, $\mathcal{V} = \mathcal{V}(Q)$ and $\mathcal{Z}^\omega_Q$ be as in Def.\ref{def:FG_algebra}. Let $(\alpha_v)_{v\in\mathcal{V}} \in \mathbb{Z}^\mathcal{V}$.
\begin{enumerate}
\item[\rm (1)] The Weyl-ordered Laurent monomial $[{\textstyle \prod}_{v\in \mathcal{V}} {\bf Z}_v^{\alpha_v}]_{\rm Weyl}$ for $Q$ is $*$-invariant.

\item[\rm (2)] For $m\in \frac{1}{2}\mathbb{Z}$, the element $\omega^m [{\textstyle \prod}_{v\in \mathcal{V}} {\bf Z}_v^{\alpha_v}]_{\rm Weyl} \in \mathcal{Z}^\omega_Q$ is $*$-invariant if and only if $m=0$. 

\item[\rm (3)] A \ul{\em multiplicity-free} Laurent polynomial for $Q$, i.e. a Laurent polynomial for $Q$ that can be written as a sum of Laurent monomials so that no two appearing Laurent monomials $\epsilon \, \omega^m [\prod_{v\in\mathcal{V}} {\bf Z}_v^{\alpha_v}]_{\rm Weyl}$ and $\epsilon' \, \omega^{m'} [\prod_{v\in\mathcal{V}} {\bf Z}_v^{\alpha_v'}]_{\rm Weyl}$ have the same degrees $\alpha_v=\alpha_v'$, $\forall v$, is fixed by the $*$-map if and only if each of its terms is fixed by the $*$-map, i.e. it is a term-by-term Weyl-ordered Laurent polynomial. \qed
\end{enumerate}
\end{lemma}
As a consequence, for $\alpha_{v_1},\ldots,\alpha_{v_N}\in \mathbb{Z}$, the unique $m\in \frac{1}{2}\mathbb{Z}$ that makes $\omega^m {\bf Z}_{v_1}^{\alpha_{v_1}}\cdots {\bf Z}_{v_N}^{\alpha_{v_N}}$ $*$-invariant is $m = -\sum_{i<j} \varepsilon_{v_iv_j} \alpha_{v_i} \alpha_{v_j}$, in which case $\omega^m {\bf Z}_{v_1}^{\alpha_{v_1}}\cdots {\bf Z}_{v_N}^{\alpha_{v_N}} = [{\textstyle \prod}_{v\in \mathcal{V}} {\bf Z}_v^{\alpha_v}]_{\rm Weyl}$.

\subsection{${\rm SL}_3$ quantum trace maps}

Both the ${\rm SL}_3$-skein algebra and the ${\rm PGL}_3$ Fock-Goncharov algebra could be viewed as certain versions of the quantum algebras for the Poisson moduli space $\mathscr{X}_{{\rm PGL}_3,\frak{S}}$ or $\mathscr{P}_{{\rm PGL}_3,\frak{S}}$. The main object of study of the present paper, the {\em ${\rm SL}_3$ quantum trace map} \cite{Kim}, is a map that connects these two algebras. Its characteristic property is the compatibility under the cutting and gluing of the surfaces. We first recall the process of cutting.

\begin{definition}[cutting process; see e.g. {\cite[Lem.5.6]{Kim}}]
\label{def:cutting_process}
Let $\frak{S}$ be a generalized marked surface. Let $e$ be an ideal arc of $\frak{S}$ whose interior lies in the interior of $\frak{S}$, such that $e$ is not isotopic to a boundary arc of $\frak{S}$ (Def.\ref{def:ideal_triangulation}). Denote by $\frak{S}_e$ the generalized marked surface obtained from $\frak{S}$ by cutting along $e$, which is uniquely determined up to isomorphism. Let $g_e : \frak{S}_e \to \frak{S}$ be the gluing map.

\vs

If $\Delta$ is an ideal triangulation of $\frak{S}$ containing $e$ as one of the members, then the cutting process yields an ideal triangulation $\Delta_e := g_e^{-1}(\Delta_e)$ of $\frak{S}_e$.

\vs

Let $W$ be an ${\rm SL}_3$-web in $\frak{S}\times {\bf I}$ that meets $e\times {\bf I}$ transversally, where the intersection points of $W\cap(e\times {\bf I})$ are not internal $3$-valent vertices of $W$ and have mutually distinct elevations, and the framing at these intersections are upward vertical. Then let $W_e := (g_e \times {\rm id})^{-1}(W)$ be the ${\rm SL}_3$-web in $\frak{S}_e\times {\bf I}$ obtained from $W$ by the cutting process along $e$. A state $s_e$ of $W_e$ is said to be \ul{\em compatible} with a state $s$ of $W$ if $s_e(x) = s((g_e \times {\rm id})(x))$ for all $x\in \partial W_e \cap (g_e \times {\rm id})^{-1}(\partial W)$ and $s_e(x_1)=s_e(x_2)$ for all $x_1,x_2 \in \partial W_e \cap (g_e \times {\rm id})^{-1}(e\times {\bf I})$ such that $(g_e \times {\rm id})(x_1)=(g_e \times {\rm id})(x_2)$.

\vs

Denote by $g_e : \mathcal{V}(Q_{\Delta_e}) \to \mathcal{V}(Q_\Delta)$ the corresponding map between the nodes of the $3$-triangulation quivers. Define the induced \ul{\em cutting map} between the cube-root Fock-Goncharov algebras
$$
i_{\Delta,\Delta_e} ~:~ \mathcal{Z}^\omega_{\Delta} \to \mathcal{Z}^\omega_{\Delta_e}
$$
to be the $\mathbb{Z}[\omega^{\pm 1/2}]$-algebra map given on the generators as
$$
i_{\Delta,\Delta_e}({\bf Z}_v) = {\textstyle \prod}_{w\in g_e^{-1}(v)} {\bf Z}_w, \qquad \forall v\in \mathcal{V}(Q_\Delta),
$$
and likewise on their inverses. 
\end{definition}
Note that $g_e^{-1}(v)\subset \mathcal{V}(Q_{\Delta_e})$ has a single element unless $v\in \mathcal{V}(Q_\Delta)$ is one of the two nodes lying in $e$, in which case $g^{-1}_e(v)$ has two elements; in this case, for $v_1,v_2 \in g_e^{-1}(v)$, we have ${\bf Z}_{v_1} {\bf Z}_{v_2} = {\bf Z}_{v_2} {\bf Z}_{v_1}$ in $\mathcal{Z}^\omega_{\Delta_e}$, so the product expression ${\textstyle \prod}_{w\in g_e^{-1}(v)} {\bf Z}_w$ makes sense.

\vs

The following basic observation on the cutting map $i_{\Delta,\Delta_e}$ will be used later.
\begin{lemma}
\label{lem:i_injective}
The cutting map $i_{\Delta,\Delta_e} : \mathcal{Z}^\omega_\Delta \to \mathcal{Z}^\omega_{\Delta_e}$ defined in Def.\ref{def:cutting_process} is injective.
\end{lemma}

{\it Proof.} Let $\mathcal{V} = \mathcal{V}(Q_\Delta)$ and $\mathcal{V}_e = \mathcal{V}(Q_{\Delta_e})$. Let $\mathcal{W} = \{ [\prod_{v\in\mathcal{V}} {\bf Z}^{\alpha_v}_v ]_{\rm Weyl} \, | \, (\alpha_v)_{v\in\mathcal{V}} \in \mathbb{Z}^\mathcal{V}\} \subset\mathcal{Z}^\omega_\Delta$ and $\mathcal{W}_e = \{[\prod_{w\in \mathcal{V}_e} {\bf Z}_w^{\alpha_w}]_{\rm Weyl} \, | \, (\alpha_w)_{w\in\mathcal{V}_e}\in \mathbb{Z}^{\mathcal{V}_e}\}\subset\mathcal{Z}^\omega_{\Delta_e}$. A standard fact about Laurent polynomial rings says that $\mathcal{W}$ forms a $\mathbb{Z}[\omega^{\pm 1/2}]$-basis of $\mathcal{Z}^\omega_\Delta$ and that $\mathcal{W}_e$ forms a $\mathbb{Z}[\omega^{\pm 1/2}]$-basis of $\mathcal{Z}^\omega_{\Delta_e}$. Let $\mathcal{W}_e' := \{ [\prod_{w\in\mathcal{V}_e} {\bf Z}_w^{\alpha_w}]_{\rm Weyl} \, | \, (\alpha_w)_{w\in\mathcal{V}_e}\in\mathbb{Z}^{\mathcal{V}_e}, \, \mbox{$\alpha_w=\alpha_{w'}$ whenever $g_e(w)=g_e(w')$}\} \subset \mathcal{Z}^\omega_{\Delta_e}$. It is straightforward to see that $i_{\Delta,\Delta_e}$ induces a set bijection from $\mathcal{W}$ to $\mathcal{W}_e'$. Meanwhile, since $\mathcal{W}'_e \subset\mathcal{W}_e$ and $\mathcal{W}_e$ is $\mathbb{Z}[\omega^{\pm 1/2}]$-linearly independent, one also observes that $\mathcal{W}_e$ is $\mathbb{Z}[\omega^{\pm 1/2}]$-linearly independent. Since $i_{\Delta,\Delta_e}$ is a $\mathbb{Z}[\omega^{\pm 1/2}]$-linear map that sends a $\mathbb{Z}[\omega^{\pm 1/2}]$-basis to a $\mathbb{Z}[\omega^{\pm 1/2}]$-linearly independent set, it is injective. \qed

\vs

\begin{theorem}[{\cite[Thm.1.27, Thm.5.8]{Kim}}; see also \cite{Douglas1} \cite{LY23}]
\label{thm:SL3_quantum_trace}
There exists a family of $\mathbb{Z}[\omega^{\pm 1/2}]$-algebra homomorphisms
$$
{\rm Tr}^\omega_\Delta = {\rm Tr}^\omega_{\Delta;\frak{S}} ~:~ \mathcal{S}^\omega_{\rm s}(\frak{S};\mathbb{Z})_{\rm red} \to \mathcal{Z}^\omega_\Delta,
$$
called the \ul{\em ${\rm SL}_3$ quantum trace maps}, defined for each triangulable generalized marked surface $\frak{S}$ and its ideal triangulation $\Delta$, such that the following hold.

\begin{enumerate}\itemsep0em
\item[\rm (QT1)] (cutting/gluing) Let $(W,s)$ be a stated ${\rm SL}_3$-web in $\frak{S}\times {\bf I}$, and $e$ an internal arc of $\Delta$. Let $\frak{S}_e$, $\Delta_e$ and $W_e \subset \frak{S}_e \times {\bf I}$ be obtained from $\frak{S}$, $\Delta$ and $W\subset \frak{S}\times {\bf I}$ by cutting along $e$ as in Def.\ref{def:cutting_process}. Then 
\begin{align}
\label{eq:cutting_axiom}
i_{\Delta,\Delta_e} {\rm Tr}^\omega_{\Delta;\frak{S}}([W,s]) = \underset{s_e}{\textstyle \sum} \, {\rm Tr}^\omega_{\Delta_e;\frak{S}_e} ([W_e,s_e]), 
\end{align}
where the sum is over all states $s_e$ of $W_e$ that are compatible with $s$ in the sense as in Def.\ref{def:cutting_process}, and the cutting map $i_{\Delta,\Delta_e}$ is as in Def.\ref{def:cutting_process}.

\item[\rm (QT2)] (values of oriented edges over a triangle) Let $(W,s)$ be a stated ${\rm SL}_3$-web in $t\times {\bf I}$, where $t$ is a triangle, viewed as a generalized marked surface with a unique ideal triangulation $\Delta$. Denote the sides of $t$ by $e_1,e_2,e_3$ (with $e_4=e_1$), and the nodes of $Q_\Delta$ by $v_{e_\alpha,1}$, $v_{e_\alpha,2}$, $v_t$ (for $\alpha=1,2,3$) as in Fig.\ref{fig:3-triangulation_node_names}.  
\begin{enumerate}
\item[\rm (QT2-1)] If $W$ consists of a single left-turn oriented edge in $t\times {\bf I}$, i.e. a crossingless ${\rm SL}_3$-web with upward vertical framing consisting of a single oriented edge, with the initial point $x$ lying over $e_\alpha$ and the terminal point lying over $e_{\alpha+1}$, then ${\rm Tr}^\omega_{\Delta;t}([W,s])$ is the $(s(x),s(y))$-th entry\footnote{We use the usual convention that the $(a,b)$-th entry of a matrix means the entry at the $a$-th row and the $b$-th column.} of the following $3\times 3$ matrix with entries in $\mathcal{Z}^\omega_\Delta$:
\begin{align}
\label{eq:QT2-1}
\hspace{-2mm} \left[ \smallmatthree{ {\bf Z}_{v_{e_\alpha,2}} {\bf Z}_{v_{e_\alpha,1}}^2 }{0}{0}{0}{\hspace{-2mm} {\bf Z}_{v_{e_\alpha,2}} {\bf Z}_{v_{e_\alpha,1}}^{-1} \hspace{-2mm}}{0}{0}{0}{ {\bf Z}_{v_{e_\alpha,2}}^{-2} {\bf Z}_{v_{e_\alpha,1}}^{-1} }
\hspace{-2mm}
\smallmatthree{ {\bf Z}_{v_t}^2 }{ {\bf Z}_{v_t}^2 + {\bf Z}_{v_t}^{-1} }{ {\bf Z}_{v_t}^{-1} }{0}{ {\bf Z}_{v_t}^{-1}}{ {\bf Z}_{v_t}^{-1}}{0}{0}{ {\bf Z}_{v_t}^{-1}}
\hspace{-2mm}
\smallmatthree{ {\bf Z}_{v_{e_{\alpha+1},1}} {\bf Z}_{v_{e_{\alpha+1},2}}^2 }{0}{0}{0}{\hspace{-2mm} {\bf Z}_{v_{e_{\alpha+1},1}} {\bf Z}_{v_{e_{\alpha+1},2}}^{-1} \hspace{-2mm}}{0}{0}{0}{ {\bf Z}_{v_{e_{\alpha+1},1}}^{-2} {\bf Z}_{v_{e_{\alpha+1},2}}^{-1} } \right]_{\rm Weyl}
\end{align}

\item[\rm (QT2-2)] If $W$ consists of a single right-turn oriented edge in $t\times {\bf I}$, i.e. a single crossingless oriented edge from $x\in e_{\alpha+1}\times {\bf I}$ to $y\in e_\alpha \times {\bf I}$, then  ${\rm Tr}^\omega_{\Delta;t}([W,s])$ is the $(s(x),s(y))$-th entry of
\begin{align}
\nonumber
\hspace{-2mm} \left[ \smallmatthree{ {\bf Z}_{v_{e_{\alpha+1},2}} {\bf Z}_{v_{e_{\alpha+1},1}}^2 }{0}{0}{0}{\hspace{-2mm} {\bf Z}_{v_{e_{\alpha+1},2}} {\bf Z}_{v_{e_{\alpha+1},1}}^{-1} \hspace{-2mm}}{0}{0}{0}{ {\bf Z}_{v_{e_{\alpha+1},2}}^{-2} {\bf Z}_{v_{e_{\alpha+1},1}}^{-1} } 
\hspace{-2mm}
\smallmatthree{ {\bf Z}_{v_t} }{ 0 }{ 0} { {\bf Z}_{v_t} }{ {\bf Z}_{v_t} }{ 0 }{ {\bf Z}_{v_t} }{ {\bf Z}_{v_t} + {\bf Z}_{v_t}^{-2} }{ {\bf Z}_{v_t}^{-2} }
\hspace{-2mm} \smallmatthree{ {\bf Z}_{v_{e_\alpha,1}} {\bf Z}_{v_{e_\alpha,2}}^2 }{0}{0}{0}{\hspace{-2mm} {\bf Z}_{v_{e_\alpha,1}} {\bf Z}_{v_{e_\alpha,2}}^{-1} \hspace{-2mm}}{0}{0}{0}{ {\bf Z}_{v_{e_\alpha,1}}^{-2} {\bf Z}_{v_{e_\alpha,2}}^{-1} }
\right]_{\rm Weyl}
\end{align}
\end{enumerate}
\end{enumerate}
\end{theorem}

The ${\rm SL}_3$ quantum trace map, which can be viewed as the ${\rm SL}_3$ version of Bonahon and Wong's ${\rm SL}_2$ quantum trace map \cite{BW} and which later generalized to the ${\rm SL}_n$ quantum trace map by L\^e and Yu \cite{LY23}, is supposed to be a quantum deformed version of the ${\rm SL}_3$ classical trace map, whose value at an oriented loop is the `trace-of-monodromy function' on the space $\mathscr{X}_{{\rm PGL}_3,\frak{S}}$ along that loop. The values on the basic cases (QT2-1)--(QT2-2), which already appeared in \cite{Douglas1}, are natural candidates for these cases, deforming the entries of the corresponding classical monodromy matrices of Fock and Goncharov \cite{FG06}, or more precisely, the suitably normalized versions. By the cutting/gluing property (QT1), the values of the ${\rm SL}_3$ quantum trace maps for a surface $\frak{S}$ are completely determined by the values for the triangle $t$, viewed as a standalone generalized marked surface. However, the above version of Thm.\ref{thm:SL3_quantum_trace} doesn't tell us how to compute the values of ${\rm Tr}^\omega_{\Delta;t}([W,s])$ for all ${\rm SL}_3$-webs $W$ in $t\times {\bf I}$. In \cite{Kim}, what is crucially used in the proof of Thm.\ref{thm:SL3_quantum_trace} above as well as in the computation of the values is the biangle analog of the ${\rm SL}_3$ quantum trace map.
\begin{proposition}[{\cite[Prop.5.26 and \S5]{Kim}}]
\label{prop:biangle_quantum_trace}
Let $B$ be a biangle, viewed as a generalized marked surface, diffeomorphic to a closed disc with two marked points on the boundary and no puncture in the interior. There exists a $\mathbb{Z}[\omega^{\pm 1/2}]$-algebra homomorphism
$$
{\rm Tr}^\omega_B ~:~ \mathcal{S}^\omega_{\rm s}(B;\mathbb{Z})_{\rm red} \to \mathbb{Z}[\omega^{\pm 1/2}],
$$
called  the \ul{\em biangle ${\rm SL}_3$ quantum trace map}, satisfying the following.
\begin{enumerate}\itemsep0em
\item[\rm (BQT1)] The cutting/gluing property for each internal ideal arc $e$ in $B$ connecting the two marked points of $B$ holds:
$$
{\rm Tr}^\omega_B([W,s]) = {\textstyle \sum}_{s_1,s_2} {\rm Tr}^\omega_{B_1} ([W_1,s_1]) \, {\rm Tr}^\omega_{B_2} ([W_2,s_2]),
$$
where cutting $B$ along $e$ yields $B_e = B_1 \sqcup B_2$ which is a disjoint union of two biangles, with $W$ cut into $W_e = W_1 \sqcup W_2$, and the sum is over all states $s_1,s_2$ such that the state $s_e := s_1 \sqcup s_2$ of $W_e = W_1 \sqcup W_2$ is compatible with $s$ in the sense of Def.\ref{def:cutting_process}.

\item[\rm (BQT2)] When $W$ consists of a single crossingless oriented edge connecting the two boundary walls of $B \times {\bf I}$,
$$
{\rm Tr}^\omega_B([W,s]) = \left\{
\begin{array}{ll}
1 & \mbox{if $s$ assigns the same state values to the two endpoints of $W$,} \\
0 & \mbox{otherwise}.
\end{array}
\right.
$$

\item[\rm (BQT3)] Let $\frak{S}$ be a triangulable generalized marked surface, and $e$ an internal ideal arc of $\frak{S}$ isotopic to a boundary arc $b$, so that cutting $\frak{S}$ along $e$ yields $\frak{S}_e = \frak{S}_0 \sqcup B$ with $B$ being a biangle and $\frak{S}_0$ being isomorphic to $\frak{S}$. Let $\Delta$ be an ideal triangulation of $\frak{S}$ not meeting the interior of $e$. Let $\Delta_0$ be the ideal triangulation of $\frak{S}_0$ obtained from $\frak{S}$ by replacing $b$ by $e$. For an ${\rm SL}_3$-web $W$ in $\frak{S}\times {\bf I}$ such that the cutting process along $e$ yields a well-defined ${\rm SL}_3$-web $W_e = W_0 \sqcup W_B$ in $\frak{S}_e \times {\bf I} = (\frak{S}_0 \times {\bf I}) \sqcup (B\times {\bf I})$, one has
$$
{\rm Tr}^\omega_{\Delta;\frak{S}}([W,s])
= {\textstyle \sum}_{s_0,s_B} {\rm Tr}^\omega_{\Delta_0;\frak{S}_0}([W_0,s_0]) \, {\rm Tr}^\omega_B([W_B,s_B]),
$$
where the sum is over all states $s_0$ and $s_B$ of $W_0$ and $W_B$ constituting a state $s_e := s_0 \sqcup s_B$ of $W_e$ that is compatible with $s$ in the sense of Def.\ref{def:cutting_process}, and the algebras $\mathcal{Z}^\omega_\Delta$ and $\mathcal{Z}^\omega_{\Delta_0}$ are naturally being identified.
\end{enumerate}
\end{proposition}
It is the property (BQT3) that yields a `state-sum formula' \cite[\S5.3]{Kim} for the ${\rm SL}_3$ quantum trace ${\rm Tr}^\omega_{\Delta;\frak{S}}$ for a triangulable generalized marked surface, which we briefly review. Consider a \ul{\em split ideal triangulation} $\wh{\Delta}$ of $\Delta$, obtained by adding one ideal arc $e'$ per each arc $e$ of $\Delta$ so that $e'$ is isotopic to $e$ and $\wh{\Delta}$ is still a collection of arcs that do not meet each other in their interiors. So, an arc of $\Delta$ now becomes two parallel arcs, forming a biangle. Cutting along all internal arcs of $\wh{\Delta}$ yields a disjoint union of triangles $\wh{t}$ and biangles $B$, where each $\wh{t}$ corresponds to an ideal triangle $t$ of $\Delta$, and each $B$ corresponds to an ideal arc of $\Delta$. Let $(W,s)$ be a stated ${\rm SL}_3$-web in $\frak{S}\times {\bf I}$, and assume that $W$ meets $\wh{\Delta} \times {\bf I}$ transversally, and that $W\cap (\wh{t}\times {\bf I})$ and $W\cap (B\times I)$ are ${\rm SL}_3$-webs in $\wh{t}\times {\bf I}$ and in $B\times I$, for each $\wh{t}$ and $B$. We call elements of $W\cap (\wh{\Delta}\times {\bf I})$ the {\it junctures} of $W$ with respect to $\wh{\Delta}$, and call a function $J : W\cap (\wh{\Delta} \times {\bf I}) \to \{1,2,3\}$ a {\it juncture-state} of $W$ with respect to $\wh{\Delta}$. By using Thm.\ref{thm:SL3_quantum_trace}(QT1) and Prop.\ref{prop:biangle_quantum_trace}(BQT3), one obtains the {\it state-sum formula}
\begin{align}
\label{eq:full_state-sum_formula}
{\rm Tr}^\omega_\Delta([W,s]) = {\textstyle \sum}_J ({\textstyle \prod}_B {\rm Tr}^\omega_B([W\cap (B\times {\bf I}), J])) \, ({\textstyle \prod}_{\wh{t}} {\rm Tr}^\omega_{\wh{t}}([W\cap (\wh{t}\times {\bf I}),J])),
\end{align}
where the sum $\sum_J$ is over all juncture states $J$ of $W$ with respect to $\wh{\Delta}$ that restricts to $s$ at $\partial W$, the product $\prod_B$ is over all biangles $B$ of $\wh{\Delta}$, and the product $\prod_{\wh{t}}$ is over all triangles $\wh{t}$ of $\wh{\Delta}$. Note that $J$ can be thought of as yielding a state for $W\cap (B\times {\bf I})$ and that for $W\cap (\wh{t}\times {\bf I})$ via restriction to their endpoints. Eq.\eqref{eq:full_state-sum_formula} is written without the cutting maps (Def.\ref{def:cutting_process}) between different cube-root Fock-Goncharov algebras, so we need to say a few words about it. 
A value of ${\rm Tr}^\omega_B$ lies in $\mathbb{Z}[\omega^{\pm 1/2}]$, and a value of ${\rm Tr}^\omega_{\wh{t}}$ lies in the cube-root Fock-Goncharov algebra $\wh{Z}^\omega_t$ for the triangle $t$ of $\Delta$ (Def.\ref{def:FG_algebra}); here we naturally identified $\mathcal{Z}^\omega_{\wh{t}}$ and $\mathcal{Z}^\omega_t$, which are defined using the unique triangulations of the surfaces $\wh{t}$ and $t$. So the right-hand side of eq.\eqref{eq:full_state-sum_formula} can be viewed as an element of $\bigotimes_t \mathcal{Z}^\omega_t$, where the tensor product is over all triangles $t$ of $\Delta$. The algebra $\mathcal{Z}^\omega_\Delta$ embeds into this tensor product algebra (see \cite[Def.5.5]{Kim}), by the embedding map which sends each generator ${\bf Z}_v^\epsilon$, $v\in \mathcal{V}(Q_\Delta)$, $\epsilon \in \{1,-1\}$, to the product of all generators ${\bf Z}_{w,t}^\epsilon$ of $\mathcal{Z}^\omega_t$ such that $v$ corresponds to a node $w$ of the 3-triangulation quiver $Q_t$ for the triangle $t$ (each viewed as a surface, with a unique triangulation). For each $v$, there are either one or two such pairs $(w,t)$, where the latter happens only when $v$ is a node lying in an internal arc of $\Delta$. It is shown in \cite{Kim} (see e.g. Prop.5.47 there) that the value of the sum in the right-hand side of eq.\eqref{eq:full_state-sum_formula} lies in $\mathcal{Z}^\omega_\Delta \subset \bigotimes_t \mathcal{Z}^\omega_t$.

\vs

If $(W,s)$ is a general stated ${\rm SL}_3$-web in $\frak{S}\times {\bf I}$, then one would apply an isotopy to $(W,s)$ so that $W$ satisfies the above-mentioned hypothesis for the state-sum formula in eq.\eqref{eq:full_state-sum_formula}, and then consider the sum in eq.\eqref{eq:full_state-sum_formula}. It is shown in \cite{Kim} (see e.g. Prop.5.50--5.51 there) that the value of the sum does not depend on the isotopy, and hence depends only on the isotopy class of $(W,s)$.
In practice, one would want to isotope $W$ so that the complexities, e.g. the 3-valent vertices, are pushed to be located over biangles. Then the cutting/gluing properties let us compute ${\rm Tr}^\omega_{\Delta;\frak{S}}$ in terms of ${\rm Tr}^\omega_t$ for triangles $t$ and ${\rm Tr}^\omega_B$ for biangles $B$. One could have isotoped $W$ so that the pieces of $W$ living over each triangle are as in (QT2-1)--(QT2-2), and that each piece lives at a constant elevation; then one says that $W$ is in a {\it gool position}\footnote{This terminology is coined in \cite{Kim}, as stronger version of the `good position' which is used in \cite{Kim} as an ${\rm SL}_3$-analog of Bonahon and Wong's `good position' appearing in \cite{BW} for the ${\rm SL}_2$ quantum trace maps.}, in which case each value of ${\rm Tr}^\omega_t$ is given by product of values presented in (QT2-1)--(QT2-2). Then the hard computation should be done over the biangles (i.e. for the values of ${\rm Tr}^\omega_B$), which are relatively easier than a similar computation for triangles. Note that the biangle ${\rm SL}_3$ quantum trace can be viewed either as an incarnation of the Reshetikhin-Turaev invariant for tangles associated to the standard 3-dimensional representation of the quantum group $\mathcal{U}_q(\frak{sl}_3)$ \cite{RT}, or as the counit of the quantum group $\mathcal{O}_q({\rm SL}_3)$ \cite{Higgins}; see \cite{Kim} for more details. The state-sum formula provides one algorithmic way of computing the values of the ${\rm SL}_3$ quantum trace, and at the same time a way of proving the very existence of the ${\rm SL}_3$ quantum trace maps. However, we will only make a relatively mild use of this state-sum formula for $\wh{\Delta}$ in the present paper.

\vs

Some nice favorable properties of the ${\rm SL}_3$ quantum trace maps ${\rm Tr}^\omega_\Delta$ are proved in \cite{Kim}, but there are still more to be proved. Among the remaining, perhaps the most important property is the `naturality' under the mapping class group action, i.e. the independence on the choice of ideal triangulations $\Delta$. The present paper undertakes the task of properly formulating and proving this property.

\section{Quantum coordinate change for flips of ideal triangulations}
\label{sec:quantum_coordinate_change_for_flips_of_ideal_triangulations}

Per change of ideal triangulations $\Delta\leadsto \Delta'$ of a triangulable generalized marked surface, we investigate the quantum coordinate change maps between various versions of the Fock-Goncharov quantum algebras associated to $\Delta$ and $\Delta'$. This is the first necessary step toward the main result of the present paper, and can be viewed as the ${\rm SL}_3$ analog of Hiatt's result on the square-root version of the quantum mutation maps for ${\rm SL}_2$ \cite{Hiatt} \cite{Miri}.

\subsection{Classical cluster $\mathscr{X}$-mutations}
\label{subsec:classical_mutations}

We begin by reviewing the classical setting of \cite{FG06} \cite{FG09a}. Let $\mathcal{V}$ be any fixed nonempty finite set, and let
$$
\mathcal{F} = \mathbb{Q}(\{X^{\circ}_v \, | \, v\in \mathcal{V}\})
$$
be the field of rational functions on algebraically independent variables enumerated by $\mathcal{V}$; we refer to $\mathcal{F}$ as the \ul{\em ambient field}. The set $\mathcal{V}$ will play a role of the set of nodes of the generalized quivers to be considered, so the elements of $\mathcal{V}$ are called nodes. Choose any subset $\mathcal{V}_{\rm fr}$ of $\mathcal{V}$; the elements of $\mathcal{V}_{\rm fr}$ are called the \ul{\em frozen} nodes, and those of $\mathcal{V}\setminus\mathcal{V}_{\rm fr}$ are called the \ul{\em unfrozen} nodes. With these choices, a \ul{\em cluster $\mathscr{X}$-seed} (or, just a \ul{\em seed}) is defined as a pair $\Gamma = (Q,(X_v)_{v\in \mathcal{V}})$, where $Q$ is a quiver whose set of nodes $\mathcal{V}(Q)$ is $\mathcal{V}$, whose signed adjacency matrix is denoted by $\varepsilon=(\varepsilon_{vw})_{v,w\in\mathcal{V}}$, sometimes called the \ul{\em exchange matrix} of the seed, and $X_v$'s are elements of $\mathcal{F}$ such that $\{X_v\,|\, v\in \mathcal{V}\}$ is a transcendence basis of $\mathcal{F}$ over $\mathbb{Q}$, called the \ul{\em cluster $\mathscr{X}$-variables} of the seed. We require $\varepsilon_{vw}$ to be integers unless both $v$ and $w$ are frozen. For any unfrozen node $k$ of $Q$, i.e. $k \in \mathcal{V} \setminus \mathcal{V}_{\rm fr}$, one defines a process of \ul{\em mutation} $\mu_k$ at the node $k$, which transforms the seed $\Gamma$ into another seed $\mu_k(\Gamma)$. Denoting by $\mu_k(\Gamma) = \Gamma' = (Q', (X'_v)_{v\in \mathcal{V}})$, the quiver $Q'=\mu_k(Q)$ is defined by the following \ul{\em quiver mutation} formula for its signed adjacency matrix $\varepsilon'$:
\begin{align}
\label{eq:quiver_mutation_formula}
\varepsilon'_{vw} = \left\{
\begin{array}{ll}
-\varepsilon_{vw} & \mbox{if $k \in \{v,w\}$}, \\
\varepsilon_{vw} + \frac{1}{2}(\varepsilon_{vk}|\varepsilon_{kw}| + |\varepsilon_{vk}| \varepsilon_{kw}) & \mbox{if $k \not\in \{v,w\}$},
\end{array}
\right.
\end{align}
and the variables $X'_v$ for $\Gamma'$ are defined as elements of $\mathcal{F}$ given by the following \ul{\em cluster $\mathscr{X}$-mutation} formula
$$
X'_v = \left\{
\begin{array}{ll}
X_k^{-1} & \mbox{if $v =k$,} \\
X_v(1+X_k^{-{\rm sgn}(\varepsilon_{vk})})^{-\varepsilon_{vk}} & \mbox{if $v\neq k$},
\end{array}
\right.
$$
where ${\rm sgn}(\sim) \in \{1,-1\}$ denotes the sign, i.e. ${\rm sgn}(a)=1$ if $a>0$ and ${\rm sgn}(a)=-1$ if $a<0$. Another way of transforming a seed $\Gamma = (Q, (X_v)_{v\in \mathcal{V}})$  into a new seed is the \ul{\em seed automorphism} $P_\sigma$ associated to a permutation $\sigma$ of the set $\mathcal{V}$. The new seed $P_\sigma (\Gamma) = \Gamma' = (Q',(X'_v)_{v\in \mathcal{V}})$ is given by
$$
\varepsilon'_{\sigma(v)\, \sigma(w)} = \varepsilon_{vw}, \quad X'_{\sigma(v)} = X_v.
$$
One can apply the mutations and seed automorphisms repeatedly. In general, one begins with one cluster $\mathscr{X}$-seed, referred to as an {\em initial} cluster $\mathscr{X}$-seed, and considers only those cluster $\mathscr{X}$-seeds connected to the initial one by (finite) sequences of mutations and seed automorphisms. The quivers appearing in these seeds are said to be \ul{\em mutation-equivalent} to each other. 

\vs

Let $\frak{S}$ be a triangulable generalized marked surface. To each ideal triangulation $\Delta$ of $\frak{S}$ is associated the seed $\Gamma_\Delta = (Q_\Delta, (X_v)_{v\in \mathcal{V}(Q_\Delta)})$, where $Q_\Delta$ is the 3-triangulation quiver for $\Delta$ defined in Def.\ref{def:3-triangulation_quiver}, whose signed adjacency matrix is denoted by $\varepsilon=\varepsilon_\Delta$. The set $\mathcal{V}(Q_\Delta)$ of all nodes of the quiver $Q_\Delta$ plays the role of $\mathcal{V}$, and the set $\mathcal{V}_{\rm fr}$ of frozen nodes is defined to be the subset of $\mathcal{V}(Q_\Delta)$ consisting of the nodes of $Q_\Delta$ lying in the boundary arcs of $\frak{S}$. A crucial aspect is the relationship between the seeds associated to different ideal triangulations. Suppose that $\Delta \leadsto \Delta'$ is a flip at an arc. It is known \cite{FG06} that the corresponding $3$-triangulation quivers $Q_\Delta$ and $Q_{\Delta'}$ are related by a sequence of four mutations; namely, starting from $Q_\Delta$, first mutate at the two nodes lying in the arc being flipped (in an arbitrary order), then mutate at the two nodes lying in the interiors of the two triangles having this flipped arc as a side (in an arbitrary order). If we denote the nodes of $Q_\Delta$ appearing in the two triangles having this flipped arc as a side by $v_1,v_2,\ldots,v_{12}$ as in Fig.\ref{fig:mutations_for_a_flip}, where some of these nodes may be identical nodes depending on the situation, then one can write
$$
Q_{\Delta'} = \mu_{v_{12}} \mu_{v_7} \mu_{v_4} \mu_{v_3} Q_\Delta,
$$
as seen in Fig.\ref{fig:mutations_for_a_flip}. This mutation sequence also naturally yields an identification between $\mathcal{V}(Q_\Delta)$ and $\mathcal{V}(Q_{\Delta'})$. Now, not only the quivers, but also the variables should be related under this mutation sequence. That is to say, one could view the situation as starting from the seed $\Gamma_\Delta$, and defining a new seed $\Gamma_{\Delta'}$ by
$$
\Gamma_{\Delta'} := \mu_{v_{12}} \mu_{v_7} \mu_{v_4} \mu_{v_3} \Gamma_{\Delta}.
$$

The original formulation of \cite{FG06} is to construct a rational coordinate system for the moduli space $\mathscr{X}_{{\rm PGL}_3,\frak{S}}$ per each ideal triangulation $\Delta$, so that a coordinate function is associated to each node of $Q_\Delta$, and to show that the coordinate systems for ideal triangulations $\Delta$ and $\Delta'$ differing by a flip are related by the coordinate change formula given by the composition of the above particular sequence of cluster $\mathscr{X}$-mutations. Here we are being more abstract, just using the concept of cluster $\mathscr{X}$-seeds, not having a geometric moduli space at hand. One thing to keep in mind in the abstract setting is that if two seeds are connected by a sequence of mutations and seed automorphisms, and if the composition of the corresponding coordinate change maps for the variables is the identity map, then we identify the two seeds. A consequence of the above original geometric formulation of \cite{FG06} is that the consistency relations for flips of triangulations in Lem.\ref{lem:classical_flip_relations_for_triangulations} also hold for the above abstract setting. 

\begin{lemma}[classical consistency relations for flips, for $3$-triangulation quivers, and for seeds]
\label{lem:classical_consistency_relations_for_flips}
For each flip $\Delta \leadsto \Delta'$ of ideal triangulations of a triangulable generalized marked surface $\frak{S}$ at an internal arc $i$ of $\Delta$, where the nodes are denoted as in Fig.\ref{fig:mutations_for_a_flip}, denote by
$$
\Phi_{\Delta \Delta'} = \Phi_i := \mu_{v_{12}} \mu_{v_7} \mu_{v_4} \mu_{v_3},
$$
which can be applied to quivers or to cluster $\mathscr{X}$-seeds, so that in particular, $Q_{\Delta'} = \Phi_i (Q_\Delta)$ and $\Gamma_{\Delta'} = \Phi_i(\Gamma_\Delta)$. Then, $\Phi_i$'s satisfy the following relations, when applied to the $3$-triangulation quiver $Q_\Delta$ or to the cluster $\mathscr{X}$-seed $\Gamma_\Delta$ for an initial triangulation $\Delta$ satisfying the respective conditions:
\begin{enumerate}
\item[\rm (1)] $\Phi_i \Phi_i = {\rm id}$, for any internal arc $i$ of any triangulation $\Delta$;

\item[\rm (2)] $\Phi_i \Phi_j \Phi_i \Phi_j = {\rm id}$, if the internal arcs $i$ and $j$ of $\Delta$ satisfy $b_{ij}=0$ (see Def.\ref{def:b_ij} for $b_{ij}$);

\item[\rm (3)] $\Phi_i \Phi_j \Phi_i \Phi_j \Phi_i = P_{\sigma_{ij}^{[3]}}$, if the internal arcs $i$ and $j$ of $\Delta$ satisfy $b_{ij}=\pm 1$, where $\sigma_{ij}^{[3]}$ is a suitable permutation which permutes the seven nodes involved in the mutations in the left-hand side and fixes all other nodes. \qed
\end{enumerate}

\end{lemma}

It is not hard to write down the permutation $\sigma_{ij}^{[3]}$ explicitly, once one chooses node labels; we leave it as an exercise. A more basic well-known result is about the consistency relations for mutations.
\begin{lemma}[classical consistency relations for mutations of $\mathscr{X}$-seeds]
\label{lem:classical_consistency_relations_for_mutations}
The mutations $\mu_k$'s of quivers and cluster $\mathscr{X}$-seeds satisfy the following relations:
\begin{enumerate}
\item[\rm (1)] $\mu_v \mu_v = {\rm id}$, when applied to any seed $\Gamma$, for any non-frozen node $v$;

\item[\rm (2)] $\mu_v \mu_w \mu_v \mu_w = {\rm id}$, when applied to a seed $\Gamma$ such that $\varepsilon_{vw}=0$;

\item[\rm (3)] $\mu_v \mu_w \mu_v \mu_w \mu_v = P_{(vw)}$, when applied to a seed $\Gamma$ such that $\varepsilon_{vw} = \pm 1$, where $(vw)$ stands for the permutation of the nodes that exchanges $v$ and $w$ and fixes all other nodes.  \qed
\end{enumerate}
\end{lemma}

As of now, the proof of Lem.\ref{lem:classical_consistency_relations_for_flips} relies on the geometry of the moduli space $\mathscr{X}_{{\rm PGL}_3,\frak{S}}$ \cite{FG06} (or $\mathscr{P}_{{\rm PGL}_3,\frak{S}}$ \cite{GS19}). One can try to prove it directly using the more basic algebraic lemma, i.e. Lem.\ref{lem:classical_consistency_relations_for_mutations}. For example, the left-hand side of the item (1) of Lem.\ref{lem:classical_consistency_relations_for_flips}, when applied to $\Delta$, can be written as
$$
(\mu_{v_3} \mu_{v_4} \mu_{v_7} \mu_{v_{12}})
(\mu_{v_{12}} \mu_{v_7} \mu_{v_4} \mu_{v_3})
$$
which one can prove to equal ${\rm id}$, with the help of the items (1) and (2) of Lem.\ref{lem:classical_consistency_relations_for_mutations}. Similarly, the item (2) of Lem.\ref{lem:classical_consistency_relations_for_flips} can be shown using the item (2) of Lem.\ref{lem:classical_consistency_relations_for_mutations}. Showing the item (3) of Lem.\ref{lem:classical_consistency_relations_for_flips} using Lem.\ref{lem:classical_consistency_relations_for_mutations} would be more involved, and we leave this as an exercise to the readers.

\subsection{Quantum mutations for $\mathscr{X}$-seeds and for Fock-Goncharov algebras}
\label{subsec:quantum_X-mutations_for_FG_algebras}

For the quantum setting, here we first review known constructions established in \cite{FG06} \cite{FG09b} \cite{BZ}, restricted and adapted to the setting of the present paper. 

\vs

First, for a general cluster $\mathscr{X}$-variety setting: to a cluster $\mathscr{X}$-seed $\Gamma=(Q,(X_v)_{v\in\mathcal{V}})$ whose underlying quiver is $Q$, we associate the Fock-Goncharov algebra $\mathcal{X}^q_Q$ defined in Def.\ref{def:FG_algebra} and write it as
\begin{align}
\label{eq:FG_algebra_for_Gamma}
\mathcal{X}^q_\Gamma := \mathcal{X}^q_Q,
\end{align}
which is to serve as a quantum algebra for the seed $\Gamma$. As mentioned before, this algebra is an example of the so-called quantum torus algebras, which are known to satisfy the (right) Ore conditions \cite{Cohn}, i.e. ${\bf P}(\mathcal{X}^q_\Gamma\setminus \{0\}) \cap {\bf Q} \mathcal{X}^q_\Gamma \neq {\O}$ for each ${\bf P},{\bf Q} \in \mathcal{X}^q_\Gamma$ with ${\bf Q}\neq 0$, hence the skew field of (right) fractions ${\rm Frac}(\mathcal{X}^q_\Gamma)$, makes sense. An element of ${\rm Frac}(\mathcal{X}^q_\Gamma)$ is represented by a formal expression of the form ${\bf P}{\bf Q}^{-1}$, with ${\bf P},{\bf Q}\in \mathcal{X}^q_\Gamma$, ${\bf Q}\neq 0$, where two such expressions ${\bf P}_1 {\bf Q}_1^{-1}$ and ${\bf P}_2 {\bf Q}_2^{-1}$ represent the same element of ${\rm Frac}(\mathcal{X}^q_\Gamma)$ if there exist nonzero ${\bf S}_1,{\bf S}_2 \in \mathcal{X}^q_\Gamma$ such that ${\bf P}_1 {\bf S}_1 = {\bf P}_2 {\bf S}_2$ and ${\bf Q}_1 {\bf S}_1 = {\bf Q}_2 {\bf S}_2$. The product of two such expressions can be expressed again in the form ${\bf P}{\bf Q}^{-1}$, by algebraic manipulations using the defining commutation relations of $\mathcal{X}^q_\Gamma$; see the proof given right after Def.\ref{def:balanced_fraction_algebra} and the proof of Lem.\ref{lem:nu_omega_k_preserves_star} for examples of such manipulations. We now recall the quantum mutation map associated to a mutation of classical $\mathscr{X}$-seeds. Before giving the formula for this map, it is useful to recall the following crucial ingredient.
\begin{definition}[compact quantum dilogarithm \cite{FaKa}]
\label{def:Psi_q}
The \ul{\em quantum dilogarithm} for a quantum parameter $q$ is the function
$$
\Psi^q(x) = {\textstyle \prod}_{r=0}^\infty (1+q^{2r+1} x)^{-1}
$$
\end{definition}
One way of understanding this function is to view it as a formal series. One can also view it as a meromorphic function on the complex plane, when $q$ is a complex number such that $|q|<1$. Its characteristic property is the difference equation
$$
\Psi^q(q^2 x) = (1+qx) \Psi^q(x),
$$
which is clear at least in a formal sense. 
This equation inspires the following definition.
\begin{definition}[ratio of quantum dilogarithm]
\label{def:F_q}
For $\alpha \in \mathbb{Z}$, define
\begin{align}
\label{eq:F_q_formal}
F^q(x; \alpha) := \Psi^q(q^{2\alpha} x) \, \Psi^q(x)^{-1}
\end{align}
formally. More precisely, $F^q(x; \alpha)$ is defined as the rational expression in $x$ and $q$ given by 
\begin{align}
\label{eq:F_q}
F^q(x;\alpha)
:= {\textstyle \prod}_{r=1}^{|\alpha|} (1+q^{(2r-1){\rm sgn}(\alpha)} x)^{{\rm sgn}(\alpha)},
\end{align}
if $\alpha \neq 0$. When $\alpha=0$, set $F^q(x;0)=1$.
\end{definition}
One can use the formal definition of the quantum dilogarithm $\Psi^q$ and eq.\eqref{eq:F_q_formal} to obtain heuristic ideas of proofs; however, we will make our actual proofs depend solely on eq.\eqref{eq:F_q}. 

\vspace{1mm}

We now describe Fock and Goncharov's quantum mutation formula.

\begin{definition}[quantum $\mathscr{X}$-mutation for Fock-Goncharov algebras; \cite{FG09a} \cite{FG09b}]
\label{def:FG_quantum_mutation}
For a mutation $\Gamma \leadsto \Gamma' = \mu_k(\Gamma)$ of cluster $\mathscr{X}$-seeds (whose underlying set of nodes are $\mathcal{V}$) at an unfrozen node $k\in \mathcal{V}\setminus\mathcal{V}_{\rm fr}$, define the \ul{\em quantum mutation map} as the algebra homomorphism between the skew fields of fractions
$$
\mu^q_{\Gamma\Gamma'} = \mu^q_k ~:~ {\rm Frac}(\mathcal{X}^q_{\Gamma'}) \to {\rm Frac}(\mathcal{X}^q_\Gamma)
$$
given by the composition
$$
\mu^q_k = \mu^{\sharp q}_k \circ \mu'_k,
$$
where the \ul{\em monomial-transformation part}
$$
\mu'_k ~:~ {\rm Frac}(\mathcal{X}^q_{\Gamma'}) \to {\rm Frac}(\mathcal{X}^q_\Gamma)
$$
is the skew field homomorphism (over $\mathbb{Z}[q^{\pm 1/18}]$) given on the generators by
\begin{align}
\label{eq:mu_prime_formula}
\mu'_k( ({\bf X}'_v)^\epsilon ) = \left\{
{\renewcommand{\arraystretch}{1.4} \begin{array}{ll}
{\bf X}_k^{-\epsilon} & \mbox{if $v=k$}, \\
([{\bf X}_v {\bf X}_k^{[\varepsilon_{vk}]_+}]_{\rm Weyl})^\epsilon
& \mbox{if $v\neq k$},
\end{array}}
\right. \qquad \forall \epsilon\in \{1,-1\}
\end{align}
where $\varepsilon=(\varepsilon_{vw})_{v,w\in\mathcal{V}}$ is the exchange matrix for $\Gamma$, and the \ul{\em automorphism part}
$$
\mu^{\sharp q}_k ~:~ {\rm Frac}(\mathcal{X}^q_{\Gamma}) \to {\rm Frac}(\mathcal{X}^q_\Gamma)
$$
is given formally as the conjugation by the expression $\Psi^q({\bf X}_k)$
$$
\mu^{\sharp q}_k = {\rm Ad}_{\Psi^q({\bf X}_k)};
$$
more precisely, $\mu^{\sharp q}_k$ is the skew field homomorphism (over $\mathbb{Z}[q^{\pm 1/18}]$) given on the generators by
\begin{align}
\label{eq:mu_sharp_formula}
\mu^{\sharp q}_k({\bf X}_v^\epsilon) = {\bf X}_v^\epsilon \cdot F^q({\bf X}_k; \epsilon\cdot\varepsilon_{kv}),
\end{align}
for each $v\in\mathcal{V}$ and $\epsilon\in\{1,-1\}$, where $F^q$ is as in eq.\eqref{eq:F_q}.

\vs

For a seed automorphism $\Gamma \leadsto \Gamma' = P_\sigma(\Gamma)$, define the \ul{\em quantum seed automorphism} as the skew field homomorphism (over $\mathbb{Z}[q^{\pm 1/18}]$) given on the generators by
$$
P_\sigma ~:~ {\rm Frac}(\mathcal{X}^q_{\Gamma'}) \to {\rm Frac}(\mathcal{X}^q_\Gamma), \quad ({\bf X}'_{\sigma(v)})^\epsilon \mapsto {\bf X}_v^\epsilon, \quad \forall v\in \mathcal{V}, \quad \forall \epsilon \in \{1,-1\}.
$$
\end{definition}

The most basic property of $\mu^q_k$ and $P_\sigma$ is that they recover the classical mutation and seed automorphism formulas $\mu_k$ and $P_\sigma$ when we put $q^{1/18}=1$. More importantly, they satisfy the quantum counterparts of the consistency relations of the classical mutations as in Lem.\ref{lem:classical_consistency_relations_for_mutations}:
\begin{proposition}[\cite{BZ} {\cite[\S3.3]{FG09a}}]
\label{prop:quantum_mutation_relations}
The quantum mutations $\mu^q_k$'s satisfy the following, when applied to ${\rm Frac}(\mathcal{X}^q_\Gamma)$ for an initial seed $\Gamma$ satisfying the respective conditions:
\begin{enumerate}
\item[\rm (1)] $\mu^q_v \mu^q_v = {\rm id}$ for any non-frozen node $v$, for any initial seed $\Gamma$;

\item[\rm (2)] $\mu^q_v \mu^q_w \mu^q_v \mu^q_w = {\rm id}$, when $\varepsilon_{vw}=0$ holds for the initial seed $\Gamma$;

\item[\rm (3)] $\mu^q_v \mu^q_w \mu^q_v \mu^q_w \mu^q_v = P_{(vw)}$, when $\varepsilon_{vw} = \pm 1$ holds for the initial seed $\Gamma$.  \qed
\end{enumerate}

\end{proposition}

We now apply the quantum mutation construction to our setting, namely for the cluster $\mathscr{X}$-seeds $\Gamma_\Delta$ for $\mathscr{X}_{{\rm PGL}_3,\frak{S}}$, or $\mathscr{P}_{{\rm PGL}_3,\frak{S}}$, associated to an ideal triangulation $\Delta$ of a surface $\frak{S}$. 

\begin{definition}[quantum coordinate change map for cluster $\mathscr{X}$-charts for a flip]
\label{def:Phi_q_i}
For each flip $\Delta \leadsto \Delta'$ of ideal triangulations of a triangulable generalized marked surface $\frak{S}$ at an internal arc $i$ of $\Delta$, where the nodes of the $3$-triangulation quivers $Q_\Delta$ and $Q_{\Delta'}$ are denoted as in Fig.\ref{fig:mutations_for_a_flip}, define the quantum coordinate change map
$$
\Phi^q_{\Delta\Delta'} = \Phi^q_i ~:~ {\rm Frac}(\mathcal{X}^q_{\Delta'}) \to {\rm Frac}(\mathcal{X}^q_\Delta)
$$
between the skew fields of fractions of the Fock-Goncharov algebras (Def.\ref{def:FG_algebra} and eq.\eqref{eq:FG_algebra_for_Gamma}) as
\begin{align}
\label{eq:mu_q_i}
\Phi^q_i := \mu^q_{v_3} \mu^q_{v_4} \mu^q_{v_7} \mu^q_{v_{12}}.
\end{align}
\end{definition}
The reason why the order of composition of the four quantum mutations looks reversed from that in the classical setting is that each quantum mutation $\mu^q_k$ is written in a contravariant manner, for it is a deformation of the pullback map $\mu^*_k$ of the classical mutation birational map $\mu_k$ between the split algebraic tori $(\mathbb{G}_m)^\mathcal{V}$ associated to seeds. That is, the classical mutation sequence $\mu_{v_{12}}\mu_{v_7} \mu_{v_4} \mu_{v_3}$ can be applied to a quiver, a seed, or the split algebraic torus for a seed, while the pullback maps on the (coordinate) functions should be written as $\mu_{v_3}^* \mu_{v_4}^* \mu_{v_7}^* \mu_{v_{12}}^*$, and the above $\Phi^q_i$ is a deformation of this last composition.

\vs

One can naturally expect that the quantum counterparts of the consistency relations of the flips, i.e. Lem.\ref{lem:classical_consistency_relations_for_flips}, should hold.
\begin{proposition}[quantum consistency relations for flips for $\mathscr{X}_{{\rm PGL}_3,\frak{S}}$]
\label{prop:quantum_consistency_of_flips}
For a triangulable generalized marked surface $\frak{S}$, the quantum coordinate change maps $\Phi^q_i$ associated to flips at arcs $i$ of triangulations satisfy the following relations, when applied to the initial seed $\Gamma_\Delta$ for a triangulation $\Delta$ of $\frak{S}$ satisfying the respective conditions:
\begin{enumerate}
\item[\rm (1)] $\Phi^q_i \Phi^q_i = {\rm id}$, for any internal arc $i$ of any triangulation $\Delta$;

\item[\rm (2)] $\Phi^q_i \Phi^q_j \Phi^q_i \Phi^q_j = {\rm id}$, if the internal arcs $i$ and $j$ of $\Delta$ satisfy $b_{ij}=0$ (see Def.\ref{def:b_ij} for $b_{ij}$);

\item[\rm (3)] $\Phi^q_i \Phi^q_j \Phi^q_i \Phi^q_j \Phi^q_i = P_{\sigma_{ij}^{[3]}}$, if the internal arcs $i$ and $j$ of $\Delta$ satisfy $b_{ij}=\pm 1$, where $\sigma_{ij}^{[3]}$ is as in Lem.\ref{lem:classical_consistency_relations_for_flips}(3).
\end{enumerate}
\end{proposition}

{\it Proof.} It is proved in \cite{KN} that a relation satisfied by classical cluster $\mathscr{X}$-mutations is also satisfied by the corresponding quantum cluster $\mathscr{X}$-mutations. Since the classical relations hold by Lem.\ref{lem:classical_consistency_relations_for_mutations}, we are done. \qed

\subsection{The balanced algebras, and the quantum coordinate change maps for them}
\label{subsec:the_balanced_algebras_and_quantum_coordinate_change_maps_for_them}

In this subsection we introduce main technical definitions of the present paper.

\vs

Let $\Delta$ be an ideal triangulation of a triangulable generalized marked surface $\frak{S}$, and $(W,s)$ a stated ${\rm SL}_3$-web in $\frak{S}\times {\bf I}$ (Def.\ref{def:SL3-web}--\ref{def:stated_SL3-skein_algebra}). In general, the value ${\rm Tr}^\omega_\Delta([W,s])$ of the ${\rm SL}_3$ quantum trace map (Thm.\ref{thm:SL3_quantum_trace}) lies in $\mathcal{Z}^\omega_\Delta$ (Def.\ref{def:FG_algebra}), i.e. is a Laurent polynomial in the variables ${\bf Z}_v$'s, $v\in \mathcal{V}(Q_\Delta)$, instead of lying in $\mathcal{X}^q_\Delta$, i.e. being a Laurent polynomial in ${\bf X}_v = {\bf Z}_v^3$'s. Suppose that $\Delta \leadsto \Delta'$ is a flip at an arc. The main purpose of the present paper is to show that the ${\rm SL}_3$ quantum trace values ${\rm Tr}^\omega_\Delta([W,s])$ and ${\rm Tr}^\omega_{\Delta'}([W,s])$ are related by a suitable quantum mutation map. So far, the only known quantum mutation map is $\Phi^q_{\Delta \Delta'} : {\rm Frac}(\mathcal{X}^q_{\Delta'}) \to {\rm Frac}(\mathcal{X}^q_{\Delta})$, which tells us how the variables ${\bf X}_v'$, $v\in \mathcal{V}(Q_{\Delta'})$, are related to the variables ${\bf X}_v$, $v\in \mathcal{V}(Q_\Delta)$. Hence, one needs first to establish a quantum mutation map for the cube-root variables ${\bf Z}_v'$ and ${\bf Z}_v$. Similarly as in the ${\rm SL}_2$ case \cite{Hiatt} \cite{BW} \cite{Miri}, one can find a (quantum) rational formula only between certain subalgebras of the skew fields of fractions ${\rm Frac}(\mathcal{Z}^\omega_{\Delta'})$ and ${\rm Frac}(\mathcal{Z}^\omega_\Delta)$. So, the very first step is to identify these special subalgebras. Following the terminology for the known constructions for the ${\rm SL}_2$ case \cite{Hiatt} \cite{BW}, we call these subalgebras the {\em balanced} subalgebras. The ${\rm SL}_3$ version of the balancedness condition comes from that of the {\em tropical coordinate} systems \cite{Kim} \cite{DS1} on the set $\mathscr{A}_{\rm L}(\frak{S};\mathbb{Z})$ of all ${\rm SL}_3$-laminations in $\frak{S}$. 
\begin{definition}[from {\cite[Prop.3.30]{Kim}}]
\label{def:Delta-balanced_elements}
Let $\Delta$ be an ideal triangulation of a triangulable generalized marked surface $\frak{S}$. An element $(a_v)_{v\in \mathcal{V}(Q_\Delta)} \in (\frac{1}{3}\mathbb{Z})^{\mathcal{V}(Q_\Delta)}$ is said to be \ul{\em $\Delta$-balanced} if for each triangle $t$ of $\Delta$, the following hold: denoting by $e_1,e_2,e_3$ the sides of $t$ in the clockwise order (with $e_4=e_1$), by $v_{e_\alpha,1}$, $v_{e_\alpha,2}$ the nodes of $Q_\Delta$ lying in $e_\alpha$ so that $v_{e_\alpha,1} \to v_{e_\alpha,2}$ matches the clockwise orientation of the boundary of $t$, and by $v_t$ the node of $Q_\Delta$ lying in the interior of $t$,  as in Fig.\ref{fig:3-triangulation_node_names},
\begin{enumerate}\itemsep0em
\item[\rm (BE1)] the numbers $\sum_{\alpha=1}^3 a_{v_{e_\alpha,1}}$ and $\sum_{\alpha=1}^3 a_{v_{e_\alpha,2}}$ belong to $\mathbb{Z}$;

\item[\rm (BE2)] $\forall \alpha = 1,2,3$, the number $a_{v_{e_\alpha,1}} + a_{v_{e_\alpha,2}}$ belongs to $\mathbb{Z}$;

\item[\rm (BE3)] $\forall \alpha=1,2,3$, the number $-a_{v_t}+a_{v_{e_\alpha,2}} + a_{v_{e_{\alpha+1},1}}$ belongs to $\mathbb{Z}$. \\ \hspace*{18mm} (or equivalently, the number $a_{v_t}+a_{v_{e_\alpha,1}} + a_{v_{e_{\alpha+1},2}}$ belongs to $\mathbb{Z}$.)
\end{enumerate}
\end{definition}

As mentioned in the introduction, see \cite{DS2} (and \cite{DS1}) for a combinatorial and representation-theoretic formulation of the $\Delta$-balancedness condition in terms of the Knutson-Tao rhombi \cite{KT99} \cite{GS15}.

\begin{proposition}[{\cite[Prop.3.30]{Kim}}, \cite{DS1}]
\label{prop:coordinates_are_balanced}
Let $\Delta$ be an ideal triangulation of a triangulable generalized marked surface $\frak{S}$. There exists an injective map
$$
{\bf a}_\Delta ~:~ \mathscr{A}_{\rm L}(\frak{S};\mathbb{Z}) \to ({\textstyle \frac{1}{3}}\mathbb{Z})^{\mathcal{V}(Q_\Delta)}, \quad \ell \mapsto ({\rm a}_v(\ell))_{v\in \mathcal{V}(Q_\Delta)},
$$
called the \ul{\em tropical coordinate system} (associated to $\Delta$) on the set $\mathscr{A}_{\rm L}(\frak{S};\mathbb{Z})$ of all ${\rm SL}_3$-laminations in $\frak{S}$ (Def.\ref{def:SL3-lamination}), satisfying favorable properties. Moreover, ${\bf a}_\Delta$ is a bijection onto the set of all $\Delta$-balanced elements of $({\textstyle \frac{1}{3}}\mathbb{Z})^{\mathcal{V}(Q_\Delta)}$.
\end{proposition}
We say that ${\rm a}_v(\ell) \in \frac{1}{3}\mathbb{Z}$ is the {\it tropical coordinate} of the ${\rm SL}_3$-lamination $\ell$ at the node $v$ of the $3$-triangulation quiver $Q_\Delta$. For the reason why this is called a `tropical' coordinate, see \S\ref{subsec:compatibility_of_quantum_duality_maps}.

\vs

The following definition and the next one constitute the first set of main technical definitions introduced in the present paper.
\begin{definition}
\label{def:balanced_subalgebras}
Let $\Delta$ be an ideal triangulation of a triangulable generalized marked surface $\frak{S}$. Let $\mathcal{V} = \mathcal{V}(Q_\Delta)$.

\vs

A Laurent monomial $\epsilon \, \omega^m [\prod_{v\in \mathcal{V}} {\bf X}_v^{a_v}]_{\rm Weyl}=\epsilon \, \omega^m [\prod_{v\in \mathcal{V}} {\bf Z}_v^{3a_v}]_{\rm Weyl} \in \mathcal{Z}^\omega_\Delta$ for $\Delta$ (see Def.\ref{def:FG_algebra} and Def.\ref{def:Laurent}), with $\epsilon \in \{1,-1\}$, $m \in \frac{1}{2}\mathbb{Z}$, $(a_v)_{v\in \mathcal{V}} \in (\frac{1}{3}\mathbb{Z})^\mathcal{V}$, is said to be \ul{\em $\Delta$-balanced} if $(a_v)_{v\in \mathcal{V}} \in (\frac{1}{3}\mathbb{Z})^\mathcal{V}$ is $\Delta$-balanced in the sense of Def.\ref{def:Delta-balanced_elements}.

\vs

A Laurent polynomial for $\Delta$, i.e. an element of $\mathcal{Z}^\omega_\Delta$, is said to be \ul{\em $\Delta$-balanced} if it can be expressed as a sum of $\Delta$-balanced Laurent monomials for $\Delta$.

\vs

Let the \ul{\em $\Delta$-balanced cube-root Fock-Goncharov algebra} $\wh{\mathcal{Z}}^\omega_\Delta$ for $\Delta$ be the subset of $\mathcal{Z}^\omega_\Delta$ consisting of all $\Delta$-balanced Laurent polynomials for $\Delta$.

\end{definition}
It is easy to observe that $\wh{\mathcal{Z}}^\omega_\Delta$ is indeed a subalgebra of $\mathcal{Z}^\omega_\Delta$, and that
\begin{align}
\label{eq:X_q_in_wh_Z_omega}
\mathcal{X}^q_\Delta \subset \wh{\mathcal{Z}}^\omega_\Delta.
\end{align}

\vs

For each cluster $\mathscr{X}$-seed $\Gamma = (Q,(X_v)_{v\in\mathcal{V}})$, we associate the cube-root Fock-Goncharov algebra $\mathcal{Z}^\omega_Q$ defined in Def.\ref{def:FG_algebra}, and write it as
$$
\mathcal{Z}^\omega_\Gamma := \mathcal{Z}^\omega_Q.
$$
A \ul{\em Laurent monomial for $\Gamma$} and a \ul{\em Laurent polynomial for $\Gamma$} mean a Laurent monomial for $Q$ and a Laurent polynomial for $Q$ respectively (Def.\ref{def:Laurent}).

\vs

Just as for the case of $\mathcal{X}^q_\Gamma$, since the algebra $\mathcal{Z}^\omega_\Gamma$ is an example of quantum torus algebras, it satisfies the (right) Ore condition, and hence the skew field of (right) fractions ${\rm Frac}(\mathcal{Z}^\omega_\Gamma)$ makes sense. As a special example, ${\rm Frac}(\mathcal{Z}^\omega_\Delta)$ makes sense, for an ideal triangulation $\Delta$ of a surface $\frak{S}$.

\begin{definition}
\label{def:balanced_fraction_algebra}
Let $\Delta$ be an ideal triangulation of a triangulable generalized marked surface $\frak{S}$. Define the \ul{\em $\Delta$-balanced fraction (cube-root Fock-Goncharov) algebra} for $\Delta$ as the subset $\wh{\rm Frac}(\mathcal{Z}^\omega_\Delta)$ of ${\rm Frac}(\mathcal{Z}^\omega_\Delta)$ consisting of all elements that can be expressed as ${\bf P}{\bf Q}^{-1}$ with ${\bf P} \in \wh{\mathcal{Z}}^\omega_\Delta\subset\mathcal{Z}^\omega_\Delta$ and $0\neq {\bf Q}\in \mathcal{X}^q_\Delta \subset \mathcal{Z}^\omega_\Delta$.
\end{definition}

One can observe that $\wh{\rm Frac}(\mathcal{Z}^\omega_\Delta)$ is a subalgebra of ${\rm Frac}(\mathcal{Z}^\omega_\Delta)$. Here we provide a proof of this observation, which might serve as a preparation for some of the upcoming arguments. First, we show that it is closed under addition. Consider two elements ${\bf P}_1 {\bf Q}_1^{-1}$ and ${\bf P}_2 {\bf Q}_2^{-1}$, with ${\bf P}_1,{\bf P}_2 \in \wh{\mathcal{Z}}^\omega_\Delta$ and ${\bf Q}_1,{\bf Q}_2 \in \mathcal{X}^q_\Delta\setminus\{0\}$. By the Ore condition of $\mathcal{X}^q_\Delta$ (see \S\ref{subsec:quantum_X-mutations_for_FG_algebras}) applied to ${\bf Q}_1$ and ${\bf Q}_2$ we have ${\bf Q}_1(\mathcal{X}^q_\Delta\setminus\{0\}) \cap {\bf Q}_2 \mathcal{X}^q_\Delta \neq {\O}$, hence there exist ${\bf Q}_3,{\bf Q}_4 \in \mathcal{X}^q_\Delta$ with ${\bf Q}_3\neq 0$, satisfying ${\bf Q}_1 {\bf Q}_3 = {\bf Q}_2 {\bf Q}_4$. Let ${\bf Q}_5 := {\bf Q}_1 {\bf Q}_3 = {\bf Q}_2 {\bf Q}_4$; we see ${\bf Q}_5 \neq 0$. Note that ${\bf P}_1{\bf Q}_1^{-1} + {\bf P}_2{\bf Q}_2^{-1} = {\bf P}_1 {\bf Q}_3 ({\bf Q}_1{\bf Q}_3)^{-1} + {\bf P}_2 {\bf Q}_4 ({\bf Q}_2 {\bf Q}_4)^{-1} = ({\bf P}_1 {\bf Q}_3 + {\bf P}_2 {\bf Q}_4){\bf Q}_5^{-1}$. By eq.\eqref{eq:X_q_in_wh_Z_omega} we see that ${\bf P}_3:= {\bf P}_1{\bf Q}_3 + {\bf P}_2 {\bf Q}_4$ belongs to $\wh{\mathcal{Z}}^\omega_\Delta$. Thus we have shown that ${\bf P}_1 {\bf Q}_1^{-1} + {\bf P}_2 {\bf Q}_2^{-1} = {\bf P}_3 {\bf Q}_5^{-1}$ belongs to $\wh{\rm Frac}(\mathcal{Z}^\omega_\Delta)$, as desired. Next, let's show that $\wh{\rm Frac}(\mathcal{Z}^\omega_\Delta)$ is closed under multiplication; we should show ${\bf P}_1 {\bf Q}_1^{-1} {\bf P}_2 {\bf Q}_2^{-1} \in \wh{\rm Frac}(\mathcal{Z}^\omega_\Delta)$. By additivity, it suffices to deal with the case when ${\bf P}_2$ is a Laurent monomial in $\wh{\mathcal{Z}}^\omega_\Delta \subset \mathcal{Z}^\omega_\Delta$; then ${\bf P}_2$ is invertible, and its inverse is also a Laurent monomial. Then, one can easily show that ${\bf Q}_6 := {\bf P}_2^{-1} {\bf Q}_1 {\bf P}_2$ belongs to $\mathcal{X}^q_\Delta\setminus\{0\}$, by using the basic commutation relations in Lem.\ref{lem:relations_of_cube-root_algebra_by_X}. Note that ${\bf P}_2 {\bf Q}_6 = {\bf Q}_1 {\bf P}_2$, hence ${\bf P}_1 {\bf Q}_1^{-1} {\bf P}_2 {\bf Q}_2^{-1} = {\bf P}_1 {\bf P}_2 {\bf Q}_6^{-1} {\bf Q}_2^{-1} = ({\bf P}_1{\bf P}_2)({\bf Q}_2 {\bf Q}_6)^{-1}$; this is an example of the standard algebraic manipulations done in the skew field of fractions as mentioned in the beginning of \S\ref{subsec:quantum_X-mutations_for_FG_algebras}. As ${\bf P}_1 {\bf P}_2 \in \wh{\mathcal{Z}}^\omega_\Delta$ and ${\bf Q}_2{\bf Q}_6\in \mathcal{X}^q_\Delta\setminus\{0\}$, we have shown that ${\bf P}_1{\bf Q}_1^{-1} {\bf P}_2{\bf Q}_2^{-1} \in \wh{\rm Frac}(\mathcal{Z}^\omega_\Delta)$, as desired.

\vspace{1mm}

Another useful observation which is easy to see is that $\wh{\rm Frac}(\mathcal{Z}^\omega_\Delta)$ contains $\wh{\mathcal{Z}}^\omega_\Delta$ and ${\rm Frac}(\mathcal{X}^q_\Delta)$ as subalgebras.

\vspace{1mm}

As a matter of fact, almost verbatim argument of \cite[Rem.3.11]{KLS} shows the following lemma.
\begin{lemma}
\label{lem:balanced_fraction_skew-field}
$\wh{\rm Frac}(\mathcal{Z}^\omega_\Delta)$ is a skew field, and coincides with the skew field of fractions ${\rm Frac}(\wh{\mathcal{Z}}^\omega_\Delta)$ of $\wh{\mathcal{Z}}^\omega_\Delta$, where ${\rm Frac}(\wh{\mathcal{Z}}^\omega_\Delta)$ is naturally viewed as a subalgebra of ${\rm Frac}(\mathcal{Z}^\omega_\Delta)$. \qed
\end{lemma}
Note that the elements of $\wh{\rm Frac}(\mathcal{Z}^\omega_\Delta)$ are of the form ${\bf P}{\bf Q}^{-1}$ with ${\bf P}\in \wh{\mathcal{Z}}^\omega_\Delta \subset \mathcal{Z}^\omega_\Delta$ and $0 \neq {\bf Q}\in \mathcal{X}^q_\Delta \subset \wh{\mathcal{Z}}^\omega_\Delta \subset \mathcal{Z}^\omega_\Delta$, whereas those of ${\rm Frac}(\wh{\mathcal{Z}}^\omega_\Delta)$ are of the form ${\bf P}{\bf Q}^{-1}$ with ${\bf P}\in \wh{\mathcal{Z}}^\omega_\Delta$ and $0\neq {\bf Q} \in \wh{\mathcal{Z}}^\omega_\Delta$.

\vs

The definition of balancedness is inspired by the following important basic statement.
\begin{proposition}[values of the ${\rm SL}_3$ quantum trace are $\Delta$-balanced]
\label{prop:value_of_the_SL3_quantum_trace_is_Delta-balanced}
Let $\Delta$ be an ideal triangulation of a triangulable generalized marked surface $\frak{S}$, and $(W,s)$ be a stated ${\rm SL}_3$-web in $\frak{S} \times {\bf I}$. Then the value ${\rm Tr}^\omega_\Delta([W,s])$ of the ${\rm SL}_3$ quantum trace is a $\Delta$-balanced Laurent polynomial for $\Delta$, i.e.
$$
{\rm Tr}^\omega_\Delta([W,s]) ~\in~ \wh{\mathcal{Z}}^\omega_\Delta ~~ \subset \mathcal{Z}^\omega_\Delta.
$$
\end{proposition}

Prop.\ref{prop:value_of_the_SL3_quantum_trace_is_Delta-balanced} is not stated in \cite{Kim}, but follows from the following statements about the highest-term degrees and the congruence of degrees of terms. The highest-term statement is proved in \cite[Prop.5.80]{Kim} for non-elliptic ${\rm SL}_3$-webs without endpoints and played a crucial role in that whole paper. Almost verbatim  proof of \cite[Prop.5.80]{Kim} yields the following version for a generalized marked surface, possibly with boundary.
\begin{proposition}[the highest term of the ${\rm SL}_3$ quantum trace value; {\cite[Prop.5.80]{Kim}}]
\label{prop:highest_term}
Let $\Delta$ be an ideal triangulation of a triangulable generalized marked surface $\frak{S}$, and $(W,s)$ a stated ${\rm SL}_3$-web in $\frak{S}\times {\bf I}$, such that $W$ has upward vertical framing everywhere, $\pi(W)$ is a reduced non-elliptic ${\rm SL}_3$-web in $\frak{S}$ (with $\pi$ in eq.\eqref{eq:pi}), and $s$ assigns $1$ to all endpoints. Let $\mathcal{V} = \mathcal{V}(Q_\Delta)$. Then ${\rm Tr}^\omega_\Delta([W,s]) \in \mathcal{Z}^\omega_\Delta$ can be written as a sum of Laurent monomials for $\Delta$, such that $\omega^m \, [\prod_{v\in\mathcal{V}} {\bf X}_v^{{\rm a}_v(\pi(W))}]_{\rm Weyl}$, for some $m\in \frac{1}{2}\mathbb{Z}$, is the unique Laurent monomial of the highest preorder induced by the powers of the Laurent monomials, where ${\rm a}_v(\pi(W))\in\frac{1}{3}\mathbb{Z}$ is the tropical coordinate at $v\in\mathcal{V}$ of $\pi(W)$ viewed as an ${\rm SL}_3$-lamination equipped with weight $1$ on all components. That is, for any other Laurent monomial $\epsilon' \omega^{m'} [\prod_{v\in\mathcal{V}} {\bf X}_v^{b_v}]_{\rm Weyl}$ appearing in ${\rm Tr}^\omega_\Delta([W,s])$ (with $\epsilon'\in\{1,-1\}$, $m' \in \frac{1}{2}\mathbb{Z}$, $(b_v)_{v\in \mathcal{V}} \in (\frac{1}{3}\mathbb{Z})^{\mathcal{V}}$), we have ${\rm a}_v(\pi(W)) \ge b_v$ for all $v\in \mathcal{V}$, with at least one of these inequalities being strict. Moreover, if $W$ has no endpoints, then $m=0$. \qed
\end{proposition}
It is worthwhile to remark that the reason why the highest term of ${\rm Tr}^\omega_\Delta([W,s])$ is possibly not equal the Weyl-ordered Laurent monomial $[\prod_{v\in\mathcal{V}} {\bf X}_v^{{\rm a}_v(\pi(W))}]_{\rm Weyl}$ but is equal only up to a factor $\omega^m$ is that when $W$ has endpoints, one can find that an ${\rm SL}_3$-web $W'$ in $\frak{S} \times {\bf I}$ with upward vertical framing whose projection $\pi(W')$ in $\frak{S}$ coincides with $\pi(W)$ but such that $W$ and $W'$ are not isotopic as ${\rm SL}_3$-webs in $\frak{S} \times {\bf I}$. We refer the readers to Rem.\ref{rem:highest_term_using_elevation-reversing} for more detailed arguments and a more refined statement. For now, let's move on to the investigation of the terms of ${\rm Tr}^\omega_\Delta([W,s])$ other than the highest term.
\begin{proposition}[congruence of terms of the ${\rm SL}_3$ quantum trace value; {\cite[Prop.5.76]{Kim}}]
\label{prop:congruence_of_SL3_quantum_trace}
Any pair of two terms $\epsilon' \omega^{m'} [\prod_{v\in\mathcal{V}} {\bf X}_v^{b_v}]_{\rm Weyl}$ and $\epsilon'' \omega^{m''} [\prod_{v\in\mathcal{V}} {\bf X}_v^{b'_v}]_{\rm Weyl}$ of ${\rm Tr}^\omega_\Delta([W,s])$ appearing in Prop.\ref{prop:highest_term} satisfies \\$b_v - b_v' \in \mathbb{Z}$, $\forall v\in \mathcal{V} =\mathcal{V}(Q_\Delta)$; this holds for any state $s$ of $W$, not just for the special state $s$ used in Prop.\ref{prop:highest_term}. \qed
\end{proposition}
\begin{proposition}
\label{prop:congruence_of_SL3_quantum_trace2}
$b_v - {\rm a}_v(\pi(W)) \in \mathbb{Z}$, $\forall v\in \mathcal{V}$, for {\it any} state $s$ for the above $W$. \qed
\end{proposition}

{\it Proof of Prop.\ref{prop:congruence_of_SL3_quantum_trace2}.} For the special state $s$ that assigns $1$ to all endpoints of $W$, the statement follows from Prop.\ref{prop:highest_term} and Prop.\ref{prop:congruence_of_SL3_quantum_trace}. So it suffices to prove a claim saying that, if we change the values of the state $s$, then the degrees of the Laurent monomial terms for ${\rm Tr}^\omega_\Delta([W,s])$ change by integers. This claim can be checked from Thm.\ref{thm:SL3_quantum_trace}(QT2) by inspection in the case when the surface $\frak{S}$ is an ideal triangle. For a general surface, this claim can be proved by using the state-sum formula in eq.\eqref{eq:full_state-sum_formula}, as follows. Put $W$ into a gool position (see the discussion after eq.\eqref{eq:full_state-sum_formula}) by isotopy, so that the values of ${\rm Tr}^\omega_t$ in the right-hand side of eq.\eqref{eq:full_state-sum_formula} are given by products of the values presented in Thm.\ref{thm:SL3_quantum_trace}(QT2). Then one can observe that the degrees of the Laurent monomials appearing in the values of ${\rm Tr}^\omega_t$ change by integers as we change $s$. This leads to the desired claim for the value of ${\rm Tr}^\omega_\Delta([W,s])$, hence the proposition. \qed

\vs

{\it Proof of Prop.\ref{prop:value_of_the_SL3_quantum_trace_is_Delta-balanced}.} With the help of the ${\rm SL}_3$ skein relations and isotopy, we can observe that ${\rm Tr}^\omega_\Delta([W,s])$ can be expressed as a $\mathbb{Z}[\omega^{\pm1/2}]$-linear combination of ${\rm Tr}^\omega_\Delta([W',s'])$ with $W'$ (but not necessarily $s'$) satisfying the conditions of Prop.\ref{prop:highest_term}. By Prop.\ref{prop:congruence_of_SL3_quantum_trace2}, the degree $(b_v)_{v\in\mathcal{V}} \in (\frac{1}{3}\mathbb{Z})^\mathcal{V}$ of each Laurent monomial term $\epsilon' \omega^{m'} [\prod_{v\in\mathcal{V}} {\bf X}_v^{b_v}]_{\rm Weyl}$ of ${\rm Tr}^\omega_\Delta([W',s'])$ belongs to $({\rm a}_v(\pi(W')))_{v\in \mathcal{V}} + \mathbb{Z}^\mathcal{V}$. By Prop.\ref{prop:coordinates_are_balanced} we know that $({\rm a}_v(\pi(W')))_{v\in \mathcal{V}}$ is $\Delta$-balanced, and it is easy to see from definition that the $\Delta$-balancedness of an element of $(\frac{1}{3}\mathbb{Z})^\mathcal{V}$ is preserved by shift by an element of $\mathbb{Z}^\mathcal{V}$. Hence $(b_v)_{v\in\mathcal{V}}$ is $\Delta$-balanced. It follows that all Laurent monomial terms of ${\rm Tr}^\omega_\Delta([W',s'])$, hence also those of ${\rm Tr}^\omega_\Delta([W,s])$, are $\Delta$-balanced. \qed

\vs

We will now extend the quantum coordinate change map $\Phi^q_i = \mu^q_{v_3} \mu^q_{v_4} \mu^q_{v_7} \mu^q_{v_{12}}$ in eq.\eqref{eq:mu_q_i} of Def.\ref{def:Phi_q_i} to the balanced fraction algebras. Note that $\Phi_i = \mu_{v_{12}} \mu_{v_7} \mu_{v_4} \mu_{v_3}$ connects the seed $\Gamma_\Delta$ for the triangulation $\Delta$ to the seed $\Gamma_{\Delta'}$ for the triangulation $\Delta'$ that is related to $\Delta$ by the flip at an internal arc $i$. For convenience, we name the intermediate seeds as follows
$$
\Gamma_\Delta=\Gamma_{\Delta^{(0)}} \overset{\mu_{v_3}}{\leadsto} \Gamma_{\Delta^{(1)}} \overset{\mu_{v_4}}{\leadsto} \Gamma_{\Delta^{(2)}} \overset{\mu_{v_7}}{\leadsto} \Gamma_{\Delta^{(3)}} \overset{\mu_{v_{12}}}{\leadsto} \Gamma_{\Delta^{(4)}} = \Gamma_{\Delta'}
$$
So, for $r=0,1,2,3,4$, $\Delta^{(r)}$ is just a formal symbol for the seed denoted by $\Gamma_{\Delta^{(r)}}$, not necessarily representing an ideal triangulation; we may view $\Delta^{(r)}$ as an `imaginary' ideal triangulation, to which a generalized quiver $Q_{\Delta^{(r)}}$ is associated, whose signed adjacency matrix is denoted by $\varepsilon^{(r)}$:
\begin{align}
\label{eq:Delta_r}
\Delta=\Delta^{(0)} \overset{\mu_{v_3}}{\leadsto} \Delta^{(1)} \overset{\mu_{v_4}}{\leadsto} \Delta^{(2)} \overset{\mu_{v_7}}{\leadsto} \Delta^{(3)} \overset{\mu_{v_{12}}}{\leadsto} \Delta^{(4)} = \Delta'
\end{align}
This notation is reflected already in Fig.\ref{fig:mutations_for_a_flip}. Note that this sequence of four mutations yields natural identifications among the sets of nodes $\mathcal{V}(Q_{\Delta^{(r)}})$, $r=0,1,2,3,4$, so that we may write
\begin{align}
\label{eq:V_for_four_mutations}
\mathcal{V}(Q_\Delta) = \mathcal{V}(Q_{\Delta^{(0)}}) = \mathcal{V}(Q_{\Delta^{(1)}}) = \mathcal{V}(Q_{\Delta^{(2)}}) = \mathcal{V}(Q_{\Delta^{(3)}}) = \mathcal{V}(Q_{\Delta^{(4)}}) = \mathcal{V}(Q_{\Delta'}).
\end{align}
For convenience, we write the four nodes being mutated at as
\begin{align}
\label{eq:four_mutated_nodes}
v^{(1)} = v_3, \quad v^{(2)} = v_4, \quad v^{(3)} = v_7, \quad v^{(4)} = v_{12}.
\end{align}

\vs

We first define the cube-root versions $\nu'_{v^{(r)}}$ of the monomial parts $\mu'_{v^{(r)}}$ of the quantum mutation maps $\mu^q_{v^{(r)}}$ (see Def.\ref{def:FG_quantum_mutation}), for $r=1,2,3,4$. We actually do it in a more general cluster $\mathscr{X}$-seed setting.

\begin{definition}
\label{def:cube-root_monomial_transformation}
For a mutation $\Gamma \leadsto \Gamma' = \mu_k(\Gamma)$ of cluster $\mathscr{X}$-seeds (whose underlying set of nodes being $\mathcal{V}$) at an unfrozen node $k\in \mathcal{V}\setminus\mathcal{V}_{\rm fr}$, define the \ul{cube-root monomial transformation}
$$
\nu'_k ~:~ \mathcal{Z}^\omega_{\Gamma'} \to \mathcal{Z}^\omega_\Gamma
$$
as the $\mathbb{Z}[\omega^{\pm 1/2}]$-algebra homomorphism given on the generators by 
\begin{align}
\label{eq:nu_prime_formula}
\nu'_k (({\bf Z}'_v)^\epsilon) = \left\{
\begin{array}{ll}
{\bf Z}_k^{-\epsilon} & \mbox{if $v=k$}, \\
( [ {\bf Z}_v {\bf Z}_k^{[\varepsilon_{vk}]_+} ]_{\rm Weyl} )^\epsilon & \mbox{if $v \neq k$},
\end{array}
\right. \qquad \forall \epsilon\in\{1,-1\},
\end{align}
where $\varepsilon = (\varepsilon_{vw})_{v,w\in \mathcal{V}}$ is the exchange matrix for $\Gamma$.

\vspace{1mm}

Let
$$
\nu'_k ~:~ {\rm Frac}(\mathcal{Z}^\omega_{\Gamma'}) \to {\rm Frac}(\mathcal{Z}^\omega_\Gamma)
$$
be the unique extension to a homomorphism of skew fields.
\end{definition}

The following basic observation will become handy.

\begin{lemma}[basic properties of the cube-root monomial transformation]
\label{lem:nu_prime_k_basic_properties}
One has:
\begin{enumerate}
\item[\rm (1)] The map $\nu'_k : \mathcal{Z}^\omega_{\Gamma'} \to \mathcal{Z}^\omega_\Gamma$ described in Def.\ref{def:cube-root_monomial_transformation} is well defined. 

\item[\rm (2)] The map $\nu'_k ~:~ {\rm Frac}(\mathcal{Z}^\omega_{\Gamma'}) \to {\rm Frac}(\mathcal{Z}^\omega_\Gamma)$ extends $\mu'_k : {\rm Frac}(\mathcal{X}^q_{\Gamma'}) \to {\rm Frac}(\mathcal{X}^q_\Gamma)$ defined in Def.\ref{def:FG_quantum_mutation}.

\item[\rm (3)]  The map $\nu'_k : \mathcal{Z}^\omega_{\Gamma'} \to \mathcal{Z}^\omega_\Gamma$ defined in Def.\ref{def:cube-root_monomial_transformation} preserves the $*$-structures, i.e.
$$
\nu'_k({\bf U}^*) = ( \nu'_k({\bf U}) )^*, \qquad \forall {\bf U} \in \mathcal{Z}^\omega_{\Gamma'}.
$$

\item[\rm (4)] The map $\nu'_k : \mathcal{Z}^\omega_{\Gamma'} \to \mathcal{Z}^\omega_\Gamma$ sends a $*$-invariant element to a $*$-invariant element.

\item[\rm (5)] The map $\nu'_k : \mathcal{Z}^\omega_{\Gamma'} \to \mathcal{Z}^\omega_\Gamma$ sends a Weyl-ordered Laurent monomial to a Weyl-ordered Laurent monomial.

\item[\rm (6)] The map $\nu'_k : \mathcal{Z}^\omega_{\Gamma'} \to \mathcal{Z}^\omega_\Gamma$, as well as $\nu'_k : {\rm Frac}(\mathcal{Z}^\omega_{\Gamma'}) \to {\rm Frac}(\mathcal{Z}^\omega_\Gamma)$, is invertible, and the inverse is given on the generators by
\begin{align}
\label{eq:nu_prime_inverse_formula}
(\nu'_k)^{-1} ({\bf Z}_v^\epsilon) = \left\{
\begin{array}{ll}
{{\bf Z}_k'}^{-\epsilon} & \mbox{if $v=k$}, \\
( [ {\bf Z}'_v {{\bf Z}'}_k^{[\varepsilon_{vk}]_+} ]_{\rm Weyl} )^\epsilon & \mbox{if $v \neq k$},
\end{array}
\right. \qquad \forall \epsilon\in\{1,-1\},
\end{align}
and sends a Laurent monomial to a Laurent monomial.
\end{enumerate}
\end{lemma}
{\it Proof.} (1) Let $\check{\mathcal{X}}^q_\Gamma$ be the subring of $\mathcal{X}^q_\Gamma$ generated by $\{{\bf X}_v^\epsilon \, | \, v\in \mathcal{V}\} \cup \{q^{\pm 1/2}\}$. Note that $\mathcal{Z}^\omega_\Gamma$ is isomorphic to $\check{\mathcal{X}}^q_\Gamma$ as rings via an isomorphism sending each ${\bf Z}_v^\epsilon$ to ${\bf X}_v^\epsilon$ and $\omega^{\pm 1/2}$ to $q^{\pm 1/2}$, and likewise for $\Gamma'$. Therefore, the statement to be proved follows from the corresponding statement for the monomial transformation $\mu'_k : {\rm Frac}(\mathcal{X}^q_{\Gamma'}) \to {\rm Frac}(\mathcal{X}^q_\Gamma)$ defined in Def.\ref{def:FG_quantum_mutation}, or more precisely that for a restriction $\mu_k' : {\rm Frac}(\check{\mathcal{X}}^q_{\Gamma'}) \to {\rm Frac}(\check{\mathcal{X}}^q_\Gamma)$, which in turn is proved e.g. in \cite[Cor.2.9, Remark 2 after Def.3.1]{FG09b}. One could also verify it directly, as follows. It suffices to check whether the defining relations are preserved. That is, we should check that
$$
\nu'_k({\bf Z}'_v) \nu'_k({\bf Z}'_w) = \omega^{2\varepsilon'_{vw}} \nu'_k({\bf Z}'_w) \nu'_k({\bf Z}'_v), \quad \forall v,w\in \mathcal{V}, \qquad
\nu'_k({\bf Z}'_v) \nu'_k(({\bf Z}'_v)^{-1}) = \nu'_k(({\bf Z}'_v)^{-1}) \nu'_k({\bf Z}'_v), \quad \forall v\in \mathcal{V},
$$
where $\varepsilon' = (\varepsilon'_{vw})_{v,w\in\mathcal{V}}$ is the exchange matrix for $\Gamma'$. To present here one example of such a checking, suppose that $k\notin\{v,w\}$; we should then check whether $[{\bf Z}_v {\bf Z}_k^{[\varepsilon_{vk}]_+}]_{\rm Weyl} [{\bf Z}_w {\bf Z}_k^{[\varepsilon_{wk}]_+}]_{\rm Weyl} = \omega^{2\varepsilon_{vw}'} [{\bf Z}_w {\bf Z}_k^{[\varepsilon_{wk}]_+}]_{\rm Weyl} [{\bf Z}_v {\bf Z}_k^{[\varepsilon_{vk}]_+}]_{\rm Weyl}$ holds. In view of the definition of the Weyl-ordering of a Laurent monomial (Def.\ref{def:Weyl-ordered}), it suffices to check the version of the equation with the Weyl-ordering symbols $[\sim]_{\rm Weyl}$ removed, i.e. $({\bf Z}_v {\bf Z}_k^{[\varepsilon_{vk}]_+}) ({\bf Z}_w {\bf Z}_k^{[\varepsilon_{wk}]_+}) = \omega^{2\varepsilon'_{vw}} ({\bf Z}_w {\bf Z}_k^{[\varepsilon_{wk}]_+}) ({\bf Z}_v {\bf Z}_k^{[\varepsilon_{vk}]_+})$. Using the commutation relations as in Lem.\ref{lem:relations_of_cube-root_algebra_by_X}, observe that
\begin{align*}
({\bf Z}_v \myul{ {\bf Z}_k^{[\varepsilon_{vk}]_+}) ({\bf Z}_w } {\bf Z}_k^{[\varepsilon_{wk}]_+}) & = \omega^{2\varepsilon_{kw} [\varepsilon_{vk}]_+} \myul{ {\bf Z}_v {\bf Z}_w} \, \myul{ {\bf Z}_k^{[\varepsilon_{vk}]_+} {\bf Z}_k^{[\varepsilon_{wk}]_+}}  \\
& = \omega^{2\varepsilon_{kw} [\varepsilon_{vk}]_+} \omega^{2\varepsilon_{vw}} {\bf Z}_w \myul{ {\bf Z}_v {\bf Z}_k^{[\varepsilon_{wk}]_+} } {\bf Z}_k^{[\varepsilon_{vk}]_+} \\
& = \omega^{2\varepsilon_{kw} [\varepsilon_{vk}]_+} \omega^{2\varepsilon_{vw}} \omega^{2\varepsilon_{vk} [\varepsilon_{wk}]_+} ({\bf Z}_w {\bf Z}_k^{[\varepsilon_{wk}]_+}) ({\bf Z}_v {\bf Z}_k^{[\varepsilon_{vk}]_+}).
\end{align*}
This indeed equals $\omega^{2\varepsilon'_{vw}} ({\bf Z}_w {\bf Z}_k^{[\varepsilon_{wk}]_+}) ({\bf Z}_v {\bf Z}_k^{[\varepsilon_{vk}]_+})$ because
$$
\varepsilon_{kw}[\varepsilon_{vk}]_+ + \varepsilon_{vw} + \varepsilon_{vk} [\varepsilon_{wk}]_+
= \varepsilon_{vw} + \varepsilon_{kw} \frac{\varepsilon_{vk}+|\varepsilon_{vk}|}{2} + \varepsilon_{vk} \frac{\varepsilon_{wk}+|\varepsilon_{wk}|}{2}
= \varepsilon_{vw} +\frac{\varepsilon_{vk}|\varepsilon_{kw}|+|\varepsilon_{vk}|\varepsilon_{kw}|}{2} = \varepsilon'_{vw},
$$
where the last equality is from eq.\eqref{eq:quiver_mutation_formula}. The remaining computational checks are similar but easier, and are left to the readers. One might find Def.\ref{def:Heisenberg} and Lem.\ref{lem:BCH}, which are to appear later in the present paper, useful for such computations; this idea of using `log variables' is also taken in the proof in \cite{FG09b} implicitly.

\vs

(2) As each generator ${\bf X}_v^\epsilon$ of $\mathcal{X}^q_\Gamma$ is the third power $({\bf Z}_v^\epsilon)^3$ of the corresponding generator ${\bf Z}_v^\epsilon$ of $\mathcal{Z}^\omega_\Gamma$, and likewise for $\Gamma'$, it suffices to use the basic fact about Weyl-ordering, namely $[{\bf U}^3]_{\rm Weyl} = ([{\bf U}]_{\rm Weyl})^3$. It is a special case of the following well-known lemma, which we formulate slightly more generally than necessary now, for later use:
\begin{lemma} 
\label{lem:Weyl-ordering_basic}
Let $\Gamma$ be a cluster $\mathscr{X}$-seed.
\item[\rm (A)] If ${\bf U}_1,\ldots,{\bf U}_r$ are elements of $\mathcal{Z}^\omega_\Gamma$ that mutually commute. Then $[{\bf U}_1 \cdots {\bf U}_r]_{\rm Weyl} = [{\bf U}_1]_{\rm Weyl} \cdots [{\bf U}_r]_{\rm Weyl}$. 

\item[\rm (B)] For each ${\bf U} \in \mathcal{Z}^\omega_\Gamma$ and $m \in \mathbb{Z}_{\ge 0}$, one has $[{\bf U}^m]_{\rm Weyl} = ([{\bf U}]_{\rm Weyl})^m$.

\item[\rm (C)] If ${\bf U} \in \mathcal{Z}^\omega_\Gamma$ is invertible, then $[{\bf U}]_{\rm Weyl}$ is invertible, and $[{\bf U}^m]_{\rm Weyl} = ([{\bf U}]_{\rm Weyl})^m$ holds for all $m\in \mathbb{Z}$.
\end{lemma}

We remark that the Weyl-ordering of a product of elements of $\mathcal{Z}^\omega_\Gamma$ does not in general coincide with the product of the Weyl-orderings of those elements. 

\vs

{\it Proof of Lem.\ref{lem:Weyl-ordering_basic}.} (A) Write each ${\bf U}_i$ as a sum of Laurent monomials. By the definition of Laurent monomials (Def.\ref{def:Laurent}) and the defining relations of $\mathcal{Z}^\omega_\Gamma$ (Def.\ref{def:FG_algebra}), one can see that the product of Laurent monomials is again a Laurent monomial. Keeping in mind that $[{\bf U}]_{\rm Weyl}$ for ${\bf U} \in \mathcal{Z}^\omega_\Gamma$ is defined as the term-by-term Weyl-ordering (Def.\ref{def:Weyl-ordered}), one observes that it suffices to show the statement for the case when each ${\bf U}_i$ is a Laurent monomial. Then we observe that both $[{\bf U}_1 \cdots {\bf U}_r]_{\rm Weyl}$ and $[{\bf U}_1]_{\rm Weyl} \cdots [{\bf U}_r]_{\rm Weyl}$ are Laurent monomials, and they differ by multiplication by $\omega^m$ for some $m\in \frac{1}{2}\mathbb{Z}$. Note that $[{\bf U}_1 \cdots {\bf U}_r]_{\rm Weyl}$ is $*$-invariant by Lem.\ref{lem:star-invariance_and_Weyl-ordering_new}(1). Meanwhile, note that $[{\bf U}_i]_{\rm Weyl} = \omega^{m_i} {\bf U}_i$ for some $m_i \in \frac{1}{2}\mathbb{Z}$, so $[{\bf U}_1]_{\rm Weyl}$, \ldots, $[{\bf U}_r]_{\rm Weyl}$ commute with each other (because ${\bf U}_1,\ldots,{\bf U}_r$ commute with each other). Since $[{\bf U}_i]_{\rm Weyl}$ is $*$-invariant by Lem.\ref{lem:star-invariance_and_Weyl-ordering_new}(1), we observe that
$$
([{\bf U}_1]_{\rm Weyl} \cdots [{\bf U}_r]_{\rm Weyl})^*
= ([{\bf U}_r]_{\rm Weyl})^* \cdots ([{\bf U}_1]_{\rm Weyl})^*
= [{\bf U}_r]_{\rm Weyl} \cdots [{\bf U}_1]_{\rm Weyl}
= [{\bf U}_1]_{\rm Weyl} \cdots [{\bf U}_r]_{\rm Weyl},
$$
so $[{\bf U}_1]_{\rm Weyl} \cdots [{\bf U}_r]_{\rm Weyl}$ is $*$-invariant. Thus, by using Lem.\ref{lem:star-invariance_and_Weyl-ordering_new}(2) we can deduce that $[{\bf U}_1 \cdots {\bf U}_r]_{\rm Weyl}$ and $[{\bf U}_1]_{\rm Weyl} \cdots [{\bf U}_r]_{\rm Weyl}$ coincide with each other. We note that one could also obtain this statement (A) as an easy consequence of Def.\ref{def:Heisenberg} and Lem.\ref{lem:BCH} which are to appear later.

\vs

(B) This follows immediately from part (A). 

\vs

(C) By part (B), we should prove the statement only when $m<0$. If ${\bf U} \in \mathcal{Z}^\omega_\Gamma$ is invertible, then $[{\bf U}]_{\rm Weyl} [{\bf U}^{-1}]_{\rm Weyl} = [{\bf U} {\bf U}^{-1}]_{\rm Weyl} = [1]_{\rm Weyl}=1$ by part (A), and likewise we get $[{\bf U}^{-1}]_{\rm Weyl} [{\bf U}]_{\rm Weyl}=1$. Thus $[{\bf U}]_{\rm Weyl}$ is invertible, and $[{\bf U}]_{\rm Weyl}^{-1} = [{\bf U}^{-1}]_{\rm Weyl}$. So, if $m<0$, then
$$
[{\bf U}^m]_{\rm Weyl} = [({\bf U}^{-1})^{-m}]_{\rm Weyl} \, \stackrel{{\rm part~(B)}}{=} \, ([{\bf U}^{-1}]_{\rm Weyl})^{-m}
= ([{\bf U}]_{\rm Weyl}^{-1})^{-m} = [{\bf U}]_{\rm Weyl}^m
$$
as desired. \quad {\it [End of proof of Lem.\ref{lem:Weyl-ordering_basic}]}

\vs

(3) First, when ${\bf U}$ is a generator $({\bf Z}_v')^\epsilon$, note that $\nu'_k({\bf U}^*) = ( \nu'_k({\bf U}) )^*$ holds because each generator is $*$-invariant by the definition of the $*$-structure on $\mathcal{Z}^\omega_{\Gamma'}$ (see Def.\ref{def:FG_algebra}), and the image of each generator is also $*$-invariant by the definition of the $*$-structure on $\mathcal{Z}^\omega_\Gamma$ (see Def.\ref{def:FG_algebra}) and by Lem.\ref{lem:star-invariance_and_Weyl-ordering_new}(1). Now, consider an element of $\mathcal{Z}^\omega_{\Gamma'}$ of the form ${\bf U} = \epsilon \, \omega^m {\bf U}_1 {\bf U}_2 \cdots {\bf U}_r$, with $\epsilon\in\{1,-1\}$, $m\in\frac{1}{2}\mathbb{Z}$, and each of ${\bf U}_1,\ldots,{\bf U}_r$ being one of the generators of $\mathcal{Z}^\omega_{\Gamma'}$. Note that
\begin{align*}
\nu'_k( (\epsilon\, \omega^m {\bf U}_1 \cdots {\bf U}_r)^* ) 
& = \nu'_k( \epsilon \, \omega^{-m} {\bf U}_r^* \cdots {\bf U}_1^* ) \\
& = \nu'_k( \epsilon \, \omega^{-m} {\bf U}_r \cdots {\bf U}_1) \\
& = \epsilon \, \omega^{-m} \nu'_k(  {\bf U}_r ) \cdots \nu'_k({\bf U}_1) \\
& = \epsilon \, \omega^{-m} ( \nu'_k(  {\bf U}_r ) )^* \cdots ( \nu'_k({\bf U}_1) )^* \\
& = ( \epsilon \, \omega^m ( \nu'_k({\bf U}_1) \cdots \nu'_k ({\bf U}_r) ) )^* \\
& = (\nu'_k( \epsilon \, \omega^m {\bf U}_1 \cdots {\bf U}_r ) )^*
\end{align*}
So $\nu'_k({\bf U}^*) = ( \nu'_k({\bf U}) )^*$ holds for such an element ${\bf U}$. As an element of $\mathcal{Z}^\omega_{\Gamma'}$ is a sum of such elements, and since $\nu'_k$ and the $*$-maps of $\mathcal{Z}^\omega_{\Gamma'}$ and $\mathcal{Z}^\omega_\Gamma$ all preserve the addition operation, we get that $\nu'_k({\bf U}^*) = ( \nu'_k({\bf U}) )^*$ holds for every element ${\bf U}$ of $\mathcal{Z}^\omega_{\Gamma'}$.

\vs

(4) This immediately follows from the item (3).

\vs

(5) This follows from the item (4) in view of Lem.\ref{lem:star-invariance_and_Weyl-ordering_new}(1)--(2). 

\vs

(6) Let's first show that $\nu'_k$ is invertible. Since $\nu'_k$ sends generators to Laurent monomials, the problem boils down to a simple linear algebra. Namely, once we package the exponents of the generators in the right-hand side of eq.\eqref{eq:nu_prime_formula} as the $\mathcal{V} \times \mathcal{V}$ matrix $(a_{vw})_{v,w\in \mathcal{V}}$, by setting $a_{kk} = -1$, $a_{vv} = 1$ for all $v\neq k$, $a_{kv} = [\varepsilon_{vk}]_+$ for all $v\neq k$, and $a_{wv} =0$ otherwise, so that the entries of the $v$-th column represent the powers of the generators ${\bf Z}_w$ appearing in the image of the $v$-th generator ${\bf Z}_v'$, then all we need to check is whether this matrix has an inverse with integer entries. It is easy to observe that this matrix is the inverse of itself, which implies that $\nu'_k$ is invertible. In fact, this also implies that the inverse map $(\nu'_k)^{-1}$ also sends generators to Laurent monomials, and that the formula for $(\nu'_k)^{-1}({\bf Z}_v^\epsilon)$ for each generator ${\bf Z}_v^\epsilon$ is given as in the right-hand side of eq.\eqref{eq:nu_prime_inverse_formula}, perhaps up to some integer power of $\omega^{1/2}$, which may differ for each generator. Note that item (4) implies that $(\nu'_k)^{-1}$ preserves the $*$-structures, and hence $(\nu'_k)^{-1}({\bf Z}_v^\epsilon)$ is $*$-invariant, because ${\bf Z}_v^\epsilon$ is. By Lem.\ref{lem:star-invariance_and_Weyl-ordering_new}(1)--(2), we can conclude that $(\nu'_k)^{-1}({\bf Z}_v^\epsilon)$ is given precisely as in the right-hand side of eq.\eqref{eq:nu_prime_inverse_formula}, without any powers of $\omega^{1/2}$. \qed

\vs

For convenience, for $r=0,1,2,3,4$ we shall write
$$
\mathcal{Z}^\omega_{\Delta^{(r)}} := \mathcal{Z}^\omega_{\Gamma_{\Delta^{(r)}}}.
$$
The generators ${\bf Z}_v^{\pm 1}$ of $\mathcal{Z}^\omega_{\Delta^{(r)}}$ may be denoted by $({\bf Z}_v^{(r)})^{\pm 1}$, to emphasize the dependence on $r$. Likewise for ${\bf X}_v^{\pm 1} = ({\bf X}_v^{(r)})^{\pm 1}$, and more generally ${\bf X}_v^a = ({\bf X}_v^{(r)})^a$ for $a\in \frac{1}{3}\mathbb{Z}$.

\vs

The goal is to define the cube-root version $\nu^\omega_{v^{(r)}}$ of the quantum mutation $\mu^q_{v^{(r)}}$ as the composition
$$
\nu^\omega_{v^{(r)}} = \nu^{\sharp\omega}_{v^{(r)}} \circ \nu'_{v^{(r)}} ~:~ \mbox{a subset of}~{\rm Frac}(\mathcal{Z}^\omega_{\Delta^{(r)}}) \to \mbox{a subset of}~{\rm Frac}(\mathcal{Z}^\omega_{\Delta^{(r-1)}}),
$$
where the monomial-transformation part $\nu'_{v^{(r)}}$ is the restriction of the map $\nu'_{v^{(r)}} : {\rm Frac}(\mathcal{Z}^\omega_{\Delta^{(r)}}) \to {\rm Frac}(\mathcal{Z}^\omega_{\Delta^{(r-1)}})$ defined in Def.\ref{def:cube-root_monomial_transformation}, and the automorphism part
$$
\nu^{\sharp\omega}_{v^{(r)}} =  {\rm Ad}_{\Psi^q({\bf X}_{v^{(r)}}^{(r-1)})} ~:~ \mbox{a subset of}~{\rm Frac}(\mathcal{Z}^\omega_{\Delta^{(r-1)}}) \to \mbox{a subset of}~{\rm Frac}(\mathcal{Z}^\omega_{\Delta^{(r-1)}})
$$
is given by conjugation by the formal expression $\Psi^q({\bf X}_{v^{(r)}}^{(r-1)})$ in terms of the quantum dilogarithm (Def.\ref{def:Psi_q}), just like $\mu^q_{v^{(r)}}$. One of the major tasks to be done is to find a natural subset of ${\rm Frac}(\mathcal{Z}^\omega_{\Delta^{(r)}})$, $r=0,1,2,3,4$, such that $\nu^\omega_{v^{(r)}}$ is well defined. For $r=0$ and $r=4$, we already have candidates, namely the balanced fraction algebras $\wh{\rm Frac}(\mathcal{Z}^\omega_\Delta)$ and $\wh{\rm Frac}(\mathcal{Z}^\omega_{\Delta'})$. Instead of finding and justifying the best candidates for the intermediate seeds $\Delta^{(r)}$, $r=1,2,3$, we will be content with a model where the maps $\nu^\omega_{v^{(r)}}$ are well defined, which is a bare minimum for the purpose of the present paper. 

\vspace{1mm}

The following technical lemma is a crucial part of this bare minimum condition.

\begin{lemma}
\label{lem:alpha}
Let $v^{(r)}$ be as in eq.\eqref{eq:four_mutated_nodes}, for $r=1,2,3,4$. Denote the set of nodes appearing in eq.\eqref{eq:V_for_four_mutations} by $\mathcal{V}$.
Let $(a_v')_{v\in \mathcal{V}} = (a_v^{(4)})_{v\in \mathcal{V}} \in (\frac{1}{3}\mathbb{Z})^\mathcal{V}$ be $\Delta'$-balanced in the sense of Def.\ref{def:Delta-balanced_elements}.

\vspace{1mm}

Recursively define $(a_v^{(r-1)})_{v\in \mathcal{V}} \in (\frac{1}{3}\mathbb{Z})^\mathcal{V}$ for $r=4,3,2,1$ as 
\begin{align}
\label{eq:a_r_v_recursive}
a^{(r-1)}_v = \left\{
\begin{array}{ll}
-a^{(r)}_{v^{(r)}} + \sum_{w\in \mathcal{V}} [\varepsilon^{(r-1)}_{w,v^{(r)}}]_+ \, a^{(r)}_w, & \mbox{if $v=v^{(r)}$,} \\
a^{(r)}_v & \mbox{if $v\neq v^{(r)}$},
\end{array}
\right.
\end{align}
Then the number
\begin{align}
\label{eq:alpha_for_nu_sharp}
\alpha^{(r)}:= 
\underset{v\in \mathcal{V}}{\textstyle \sum} \,\,\varepsilon^{(r-1)}_{v^{(r)},v} \, a^{(r-1)}_v ~ \in ~ {\textstyle \frac{1}{3}}\mathbb{Z}
\end{align}
belongs to $\mathbb{Z}$, for $r=1,2,3,4$.
\end{lemma}

\vs

{\it Proof.} With the node names as in Fig.\ref{fig:mutations_for_a_flip}, denote $a^{(r)}_{v_j}$ by $a^{(r)}_j$, and $\varepsilon^{(r)}_{v_j v_k}$ by $\varepsilon^{(r)}_{jk}$.  Let's compute $(a_v^{(r-1)})_{v\in \mathcal{V}}$ for $r=4,3,2,1$ in terms of $(a'_v)_{v\in \mathcal{V}} =(a^{(4)}_v)_{v\in \mathcal{V}}$.  For $r=4$, note that $v^{(4)}=v_{12}$, and from Fig.\ref{fig:mutations_for_a_flip} note that
\begin{align}
\label{eq:varepsilon_3}
\varepsilon^{(3)}_{3,12}=\varepsilon^{(3)}_{10,12}=-1, \quad
\varepsilon^{(3)}_{4,12}=\varepsilon^{(3)}_{9,12}=1, \quad
\varepsilon^{(3)}_{j,12}=0, \quad \forall j\not\in\{3,4,9,10\},
\end{align}
hence from eq.\eqref{eq:a_r_v_recursive} we get
$$
a^{(3)}_{12} = - a'_{12} + a'_4 + a'_9, \quad a^{(3)}_j = a'_j, \quad \forall j \neq 12.
$$
For $r=3$, $v^{(3)}=v_7$, and from Fig.\ref{fig:mutations_for_a_flip} we have
\begin{align}
\label{eq:varepsilon_2}
\varepsilon^{(2)}_{1,7}=\varepsilon^{(2)}_{4,7}=-1, \quad
\varepsilon^{(2)}_{3,7}=\varepsilon^{(2)}_{6,7}=1, \quad
\varepsilon^{(2)}_{j,7}=0, \quad \forall j \not\in\{1,3,4,6\},
\end{align}
hence from eq.\eqref{eq:a_r_v_recursive} we get
\begin{align}
\label{eq:a_2_as_a_prime}
\left\{
{\renewcommand{\arraystretch}{1.4} \begin{array}{l}
a^{(2)}_7 = -a^{(3)}_7 + a^{(3)}_3 + a^{(3)}_6
= -a'_7 + a'_3 + a'_6, \quad
a^{(2)}_{12} = a^{(3)}_{12} = -a'_{12}+a'_4+a'_9, \\
a^{(2)}_j = a^{(3)}_j = a_j', \quad \forall j \not\in\{7,12\}.
\end{array}} \right.
\end{align}
For $r=2$, $v^{(2)}=v_4$, and from Fig.\ref{fig:mutations_for_a_flip} we have
$$
\varepsilon^{(1)}_{5,4}=\varepsilon^{(1)}_{12,4}=1, \quad
\varepsilon^{(1)}_{7,4} = \varepsilon^{(1)}_{11,4}=-1, \quad
\varepsilon^{(1)}_{j,4}=0, \quad \forall j \not\in\{5,7,11,12\},
$$
hence from eq.\eqref{eq:a_r_v_recursive} we get
\begin{align*}
& a^{(1)}_4 = -a^{(2)}_4 + a^{(2)}_5 + a^{(2)}_{12} = -a'_4 + a'_5 + (-a'_{12}+a'_4+a'_9) = a'_5 + a'_9 - a'_{12}, \\
& a^{(1)}_7 = a^{(2)}_7 = -a'_7+a'_3+a'_6, \quad
a^{(1)}_{12} = a^{(2)}_{12} = -a'_{12}+a'_4+a'_9, \quad
a^{(1)}_j = a^{(2)}_j = a'_j, \quad \forall j \not\in\{4,7,12\}.
\end{align*}
For $r=1$, $v^{(1)} = v_3$, and from Fig.\ref{fig:mutations_for_a_flip} we have
$$
\varepsilon^{(0)}_{2,3}=\varepsilon^{(0)}_{12,3}=-1, \quad
\varepsilon^{(0)}_{7,3}=\varepsilon^{(0)}_{8,3}=1, \quad
\varepsilon^{(0)}_{j,3}=0, \quad \forall j \not\in\{2,7,8,12\}.
$$
hence from eq.\eqref{eq:a_r_v_recursive} we get
\begin{align}
\label{eq:a_0_as_a_prime}
\left\{
{\renewcommand{\arraystretch}{1.4} \begin{array}{l}
a^{(0)}_3 = - a^{(1)}_3 + a^{(1)}_7 + a^{(1)}_8 = -a'_3 +(-a'_7 + a'_3+a'_6)+a'_8 = - a'_7 + a'_6  + a'_8, \\
 a^{(0)}_4 = a^{(1)}_4 = a'_5 + a'_9 - a'_{12}, \quad
a^{(0)}_7 = a^{(1)}_7 = -a'_7 + a'_3 + a'_6, \quad
a^{(0)}_{12} = a^{(1)}_{12} = -a'_{12}+a'_4+a'_9, \\
 a^{(0)}_j = a^{(1)}_j = a'_j, \quad \forall j \not\in\{3,4,7,12\}.
 \end{array}} \right.
\end{align}

\vs

In view of Fig.\ref{fig:mutations_for_a_flip}, the $\Delta'$-balancedness condition of $(a'_v)_{v\in \mathcal{V}}$ for the two triangles of $\Delta'$ having the flipped arc as a side says that the following numbers are integers:
\begin{align}
\label{eq:b_i}
\left\{ 
{\renewcommand{\arraystretch}{1.4} \begin{array}{rl}
\mbox{(BE1)} : &  b_1 := a_1'+a_8'+a'_{12}, \quad
b_2 := a_2'+a_9'+a_7', \quad
b_3 := a_5' + a_7'+a_{10}', \quad
b_4 := a_6'+a'_{12}+a'_{11}, \\
\mbox{(BE2)} : & b_5 := a'_1+a'_2, \quad b_6 := a'_8+a'_9, \quad b_7 := a'_7+a'_{12}, \quad b_8 := a'_{10}+a'_{11}, \quad b_9 := a'_5 + a'_6, \\
\mbox{(BE3)} : & b_{10} := -a_3'+a'_2+a'_8, \quad b_{11} := -a'_3 + a'_9+a'_{12}, \quad b_{12} := -a'_3+a'_7+a'_1, \\
&  b_{13} := -a'_4+a'_6+a'_7, \quad b_{14} := -a'_4+a'_{12}+a'_{10}, \quad b_{15} := -a'_4+a_{11}'+a'_5.
\end{array} }
\right.
\end{align}
We now show that the numbers $\alpha^{(r)} \in \frac{1}{3}\mathbb{Z}$ defined in eq.\eqref{eq:alpha_for_nu_sharp} are integers, for $r=1,2,3,4$. Observe that
\begin{align}
\label{eq:alpha_r_in_a_prime}
\left\{
{\renewcommand{\arraystretch}{1.4} \begin{array}{rcl}
\alpha^{(4)} & = & a^{(3)}_3 + a^{(3)}_{10} - a^{(3)}_4 - a^{(3)}_9
= a'_3 + a'_{10} - a'_4 - a'_9
= - b_{11} + b_{14}, \\ 
\alpha^{(3)} & = & a^{(2)}_1 + a^{(2)}_4 - a^{(2)}_3 - a^{(2)}_6
= a'_1 + a'_4 - a'_3 - a'_6 = b_{12} - b_{13}, \\ 
\alpha^{(2)} & = & -a^{(1)}_5-a^{(1)}_{12}+a^{(1)}_7+a^{(1)}_{11}
= -a'_5-(-a'_{12}+a'_4+a'_9)+(-a'_7+a'_3+a'_6)+a'_{11} \\
& = & -b_3 + b_4 - b_{11} + b_{14}, \\ 
\alpha^{(1)} & = & a^{(0)}_2 + a^{(0)}_{12} - a^{(0)}_7 - a^{(0)}_8
= a'_2 + (-a'_{12}+a'_4+a'_9)
- (-a'_7+a'_3+a'_6) - a'_8 \\
& = & -b_1 + b_2 + b_{12} - b_{13}. 
\end{array}} \right.
\end{align}
So, in view of eq.\eqref{eq:b_i}, one can see that $\alpha^{(r)}$ are all integers, for $r=1,2,3,4$. \qed

\vs

Consider applying the conjugation ${\rm Ad}_{\Psi^q({\bf X}_{v^{(r)}}^{(r-1)})}$ to a Weyl-ordered Laurent monomial for $\Delta^{(r-1)}$
$$
[{\textstyle \prod}_{v\in \mathcal{V}} ({\bf X}_v^{(r-1)})^{a^{(r-1)}_v}]_{\rm Weyl}
~\in~\mathcal{Z}^\omega_{\Delta^{(r-1)}} ~ \subset {\rm Frac}(\mathcal{Z}^\omega_{\Delta^{(r-1)}})
$$
for some $(a_v^{(r-1)})_{v\in \mathcal{V}} \in (\frac{1}{3}\mathbb{Z})^\mathcal{V}$, not necessarily the specific one defined via the recursive formulas in Lem.\ref{lem:alpha}; here $\mathcal{V}$ should be regarded as $\mathcal{V}(Q_{\Delta^{(r-1)}})$. Define $\alpha^{(r)} \in \frac{1}{3}\mathbb{Z}$ by the formula eq.\eqref{eq:alpha_for_nu_sharp}. Using the relations among the generators of $\mathcal{Z}^\omega_{\Delta^{(r-1)}}$ as written in Lem.\ref{lem:relations_of_cube-root_algebra_by_X}, one can observe that
$$
{\bf X}_{v^{(r)}}^{(r-1)} \, [{\textstyle \prod}_{v\in\mathcal{V}} ({\bf X}_v^{(r-1)})^{a^{(r-1)}_v}]_{\rm Weyl}
= q^{2\alpha^{(r)}} \, [{\textstyle \prod}_{v\in\mathcal{V}} ({\bf X}_v^{(r-1)})^{a^{(r-1)}_v}]_{\rm Weyl} \, {\bf X}_{v^{(r)}}^{(r-1)};
$$
this is in fact where the formula for $\alpha^{(r)}$ in eq.\eqref{eq:alpha_for_nu_sharp} came from. Then we have, at least formally (using eq.\eqref{eq:F_q_formal}),
\begin{align*}
{\rm Ad}_{\Psi^q({\bf X}_{v^{(r)}}^{(r-1)})} ( [{\textstyle \prod}_{v\in\mathcal{V}} ({\bf X}_v^{(r-1)})^{a^{(r-1)}_v}]_{\rm Weyl} )
& = \Psi^q({\bf X}_{v^{(r)}}^{(r-1)}) \, [{\textstyle \prod}_{v\in\mathcal{V}} ({\bf X}_v^{(r-1)})^{a^{(r-1)}_v}]_{\rm Weyl} \, (\Psi^q({\bf X}_{v^{(r)}}^{(r-1)}))^{-1} \\
& = [{\textstyle \prod}_{v\in\mathcal{V}} ({\bf X}_v^{(r-1)})^{a^{(r-1)}_v}]_{\rm Weyl}\, \Psi^q(q^{2\alpha^{(r)}} {\bf X}_{v^{(r)}}^{(r-1)}) \, (\Psi^q({\bf X}_{v^{(r)}}^{(r-1)}))^{-1} \\
& = [{\textstyle \prod}_{v\in\mathcal{V}} ({\bf X}_v^{(r-1)})^{a^{(r-1)}_v}]_{\rm Weyl} \, F^q({\bf X}_{v^{(r)}}^{(r-1)};\alpha^{(r)}).
\end{align*}
In order for the last resulting expression to make sense, we must have $\alpha^{(r)} \in \mathbb{Z}$, which motivated Lem.\ref{lem:alpha}. The above formal computation inspires the following definitions. We will rely only on these rigorous definitions, and not on the formal heuristics above.

\begin{definition}
\label{def:Gamma_k_balanced}
Let $\Gamma$ be a cluster $\mathscr{X}$-seed, whose set of nodes of the underlying quiver and the exchange matrix are denoted by $\mathcal{V}$ and $\varepsilon = (\varepsilon_{vw})_{v,w\in\mathcal{V}}$. Let $k \in \mathcal{V} \setminus\mathcal{V}_{\rm fr}$ be an unfrozen node.

\vs

A Laurent monomial $\epsilon \, \omega^m [\prod_{v\in \mathcal{V}} {\bf X}_v^{a_v}]_{\rm Weyl} \in \mathcal{Z}^\omega_\Gamma$ for the seed $\Gamma$, with $\epsilon \in \{1,-1\}$, $m\in \frac{1}{2}\mathbb{Z}$ and $(a_v)_{v\in\mathcal{V}} \in (\frac{1}{3}\mathbb{Z})^\mathcal{V}$, is said to be \ul{\em $(\Gamma,k)$-balanced} if the number $\alpha := \sum_{v\in \mathcal{V}} \varepsilon_{kv} a_v\in \frac{1}{3}\mathbb{Z}$ belongs to $\mathbb{Z}$.

\vs

A Laurent polynomial for $\Gamma$, i.e. an element of $\mathcal{Z}^\omega_\Gamma$, is said to be \ul{\em $(\Gamma,k)$-balanced} if it is a sum of $(\Gamma,k)$-balanced Laurent monomials.

\vs

The \ul{\em $(\Gamma,k)$-balanced fraction algebra}, denoted by $\wh{\rm Frac}_k(\mathcal{Z}^\omega_\Gamma)$, is defined as the subset of ${\rm Frac}(\mathcal{Z}^\omega_\Gamma)$ consisting of all elements of the form ${\bf P} {\bf Q}^{-1}$ with ${\bf P}$ being a $(\Gamma,k)$-balanced Laurent polynomial for $\Gamma$ and $0\neq {\bf Q} \in \mathcal{X}^q_\Gamma \subset \mathcal{Z}^\omega_\Gamma$.
\end{definition}
An argument similar to the proof given right after Def.\ref{def:balanced_fraction_algebra}, together with an easy observation that the product of $(\Gamma,k)$-balanced elements of $\mathcal{Z}^\omega_\Gamma$ is again $(\Gamma,k)$-balanced, proves that $\wh{\rm Frac}_k(\mathcal{Z}^\omega_\Gamma)$ indeed forms a subalgebra.

\begin{definition}
\label{def:nu_omega_k}
Define the \ul{\em balanced (cube-root) quantum mutation map} associated to the mutation $\Gamma \overset{\mu_k}{\leadsto} \mu_k(\Gamma) = \Gamma'$ of a cluster $\mathscr{X}$-seed $\Gamma$ (whose underlying set of nodes and the exchange matrix are $\mathcal{V}$ and $\varepsilon=(\varepsilon_{vw})_{v,w\in\mathcal{V}}$) at an unfrozen node $k\in \mathcal{V} \setminus \mathcal{V}_{\rm fr}$ as the map
$$
\nu^\omega_k ~:~ \wh{\rm Frac}_k(\mathcal{Z}^\omega_{\Gamma'}) \to \wh{\rm Frac}_k(\mathcal{Z}^\omega_\Gamma)
$$
defined as the composition 
$$
\xymatrix{
\nu^\omega_k := \nu^{\sharp \omega}_k \circ \nu'_k ~:~ \wh{\rm Frac}_k (\mathcal{Z}^\omega_{\Gamma'}) \ar[r]^-{\nu'_k} & \wh{\rm Frac}_k(\mathcal{Z}^\omega_{\Gamma}) \ar[r]^-{\nu^{\sharp \omega}_k} & \wh{\rm Frac}_k(\mathcal{Z}^\omega_{\Gamma}),
}
$$
where $\nu'_k : \wh{\rm Frac}_k(\mathcal{Z}^\omega_{\Gamma'}) \to \wh{\rm Frac}_k(\mathcal{Z}^\omega_\Gamma)$ is the restriction of the map $\nu'_k : {\rm Frac}(\mathcal{Z}^\omega_{\Gamma'}) \to {\rm Frac}(\mathcal{Z}^\omega_\Gamma)$ defined in Def.\ref{def:cube-root_monomial_transformation}, and $\nu^{\sharp\omega}_k : \wh{\rm Frac}_k(\mathcal{Z}^\omega_\Gamma) \to \wh{\rm Frac}_k(\mathcal{Z}^\omega_\Gamma)$ is given formally by the conjugation action
$$
\nu^{\sharp \omega}_k := {\rm Ad}_{\Psi^q({\bf X}_k)},
$$
or more precisely, given on the $(\Gamma,k)$-balanced Laurent monomials as
\begin{align}
\label{eq:nu_sharp_formula}
\nu^{\sharp \omega}_k(\epsilon \, \omega^m [{\textstyle \prod}_{v\in \mathcal{V}} {\bf X}_v^{a_v}]_{\rm Weyl}) = (\epsilon \, \omega^m [{\textstyle \prod}_{v\in \mathcal{V}} {\bf X}_v^{a_v}]_{\rm Weyl}) \cdot F^q({\bf X}_k;\alpha),
\end{align}
where $\alpha = \sum_{v\in \mathcal{V}} \varepsilon_{kv} a_v$, and $F^q$ is as in eq.\eqref{eq:F_q}.
\end{definition}

\begin{lemma}
\label{lem:nu_omega_u_well-defined}
Above $\nu^\omega_k$ is well defined.
\end{lemma}

{\it Proof.} Let's first check whether $\nu'_k : {\rm Frac}(\mathcal{Z}^\omega_{\Gamma'}) \to {\rm Frac}(\mathcal{Z}^\omega_\Gamma)$ of Def.\ref{def:cube-root_monomial_transformation} restricts to $\nu'_k : \wh{\rm Frac}_k(\mathcal{Z}^\omega_{\Gamma'}) \to \wh{\rm Frac}_k(\mathcal{Z}^\omega_\Gamma)$; that is, we should check that the $(\Gamma',k)$-balanced fraction algebra $\wh{\rm Frac}_k(\mathcal{Z}^\omega_{\Gamma'})$ is sent by $\nu'_k$ to the $(\Gamma,k)$-balanced fraction algebra $\wh{\rm Frac}_k(\mathcal{Z}^\omega_\Gamma)$. First, recall from Lem.\ref{lem:nu_prime_k_basic_properties}(2) that $\nu'_k : {\rm Frac}(\mathcal{Z}^\omega_{\Gamma'}) \to {\rm Frac}(\mathcal{Z}^\omega_\Gamma)$ of Def.\ref{def:cube-root_monomial_transformation} extends $\mu'_k : {\rm Frac}(\mathcal{X}^q_{\Gamma'}) \to {\rm Frac}(\mathcal{X}^q_\Gamma)$ of Def.\ref{def:FG_quantum_mutation}. Hence $\nu'_k$ sends $\mathcal{X}^q_{\Gamma'} \subset \mathcal{Z}^\omega_{\Gamma'} \subset {\rm Frac}(\mathcal{Z}^\omega_{\Gamma'})$ to $\mathcal{X}^q_\Gamma \subset \mathcal{Z}^\omega_\Gamma \subset {\rm Frac}(\mathcal{Z}^\omega_\Gamma)$. Therefore it remains to show that each element of $\mathcal{Z}^\omega_{\Gamma'}$ that is $(\Gamma',k)$-balanced is sent by $\nu'_k$ to an element of $\mathcal{Z}^\omega_\Gamma$ that is $(\Gamma,k)$-balanced.

\vs

We show this at the level of Laurent monomials. Consider a (Weyl-ordered) $(\Gamma',k)$-balanced Laurent monomial $[\prod_{v\in \mathcal{V}}({\bf X}'_v)^{a'_v}]_{\rm Weyl} \in \mathcal{Z}^\omega_{\Gamma'}$ for $\Gamma'$, with $(a'_v)_{v\in \mathcal{V}} \in (\frac{1}{3}\mathbb{Z})^\mathcal{V}$. By Lem.\ref{lem:nu_prime_k_basic_properties}(5), we have
\begin{align}
\label{eq:image_under_nu_prime_of_Laurent_monomial}
\nu'_k([{\textstyle \prod}_{v\in \mathcal{V}} ({\bf X}'_v)^{a_v'}]_{\rm Weyl}) = [{\textstyle \prod}_{v\in \mathcal{V}} {\bf X}_v^{a_v}]_{\rm Weyl}
\end{align}
for some $(a_v)_{v\in \mathcal{V}} \in (\frac{1}{3}\mathbb{Z})^\mathcal{V}$. In view of the formula for the image under $\nu'_k$ of the generators as in eq.\eqref{eq:nu_prime_formula}, one can observe that
\begin{align}
\label{eq:image_under_nu_prime_of_Laurent_monomial2}
a_v = a_v', \quad \forall v \neq k, \qquad
a_k = -a'_k + \underset{v\in\mathcal{V}}{\textstyle \sum} [\varepsilon_{vk}]_+ a'_v
\end{align}
For the current lemma we don't need $a_k$ but we recorded it for later use. The $(\Gamma',k)$-balancedness condition of $[{\textstyle \prod}_{v\in \mathcal{V}} ({\bf X}'_v)^{a_v'}]_{\rm Weyl}$ says that $\sum_{v\in \mathcal{V}} \varepsilon'_{kv} a'_v \in \mathbb{Z}$. Since $\varepsilon'_{kv} = -\varepsilon_{kv}$ for all $v\neq k$ (see eq.\eqref{eq:quiver_mutation_formula}) and $\varepsilon'_{kk}=0 = \varepsilon_{kk}$ (by skew-symmetry), it follows that 
\begin{align}
\label{eq:alpha_equals_minus_alpha_prime}
\underset{v\in \mathcal{V}}{\textstyle \sum} \varepsilon_{kv} a_v = \underset{v\neq k}{\textstyle \sum} \varepsilon_{kv} a_v  = \underset{v\neq k}{\textstyle \sum} \varepsilon_{kv} a'_v = - \underset{v\neq k}{\textstyle \sum} \varepsilon_{kv}' a'_v = - \underset{v\in \mathcal{V}}{\textstyle \sum} \varepsilon'_{kv} a'_v,
\end{align}
Hence $\sum_{v\in \mathcal{V}} \varepsilon_{kv} a_v \in \mathbb{Z}$, and the image $[{\textstyle \prod}_{v\in \mathcal{V}} {\bf X}_v^{a_v}]_{\rm Weyl} \in \mathcal{Z}^\omega_\Gamma$ is $(\Gamma,k)$-balanced, as desired.

\vs

For $\nu^{\sharp\omega}_k$, in case the Laurent monomial $\epsilon \, \omega^m [\prod_{v\in \mathcal{V}} {\bf X}_v^{a_v}]_{\rm Weyl} \in \mathcal{Z}^\omega_\Gamma$ for $\Gamma$ is $(\Gamma,k)$-balanced, then the number $\alpha = \sum_{v\in \mathcal{V}} \varepsilon_{kv} a_v \in \frac{1}{3}\mathbb{Z}$ is an integer, hence the expression $F^q({\bf X}_k;\alpha)$ in eq.\eqref{eq:nu_sharp_formula} makes sense as an element of ${\rm Frac}(\mathcal{X}^q_\Gamma)$. Any element of $\mathcal{X}^q_\Gamma$ is a sum of terms of the form $\epsilon \, \omega^m [\prod_{v\in \mathcal{V}} {\bf X}_v^{c_v}]_{\rm Weyl}$ with $\epsilon \in \{1,-1\}$, $m\in \frac{1}{2}\mathbb{Z}$, $(c_v)_{v\in\mathcal{V}} \in \mathbb{Z}^\mathcal{V}$, each of which is easily observed to be $(\Gamma,k)$-balanced as we have $\sum_{v\in \mathcal{V}} \varepsilon_{kv} c_v \in \mathbb{Z}$ (keeping in mind that $\varepsilon_{kv}$ is an integer, since $k$ is an unfrozen node); this shows that elements of $\mathcal{Z}^\omega_\Gamma$ that belong to $\mathcal{X}^q_\Gamma$ are $(\Gamma,k)$-balanced. Note also that the product of $(\Gamma,k)$-balanced elements of $\mathcal{Z}^\omega_\Gamma$ is again $(\Gamma,k)$-balanced. It follows that the element in the right-hand side of eq.\eqref{eq:nu_sharp_formula} is an element of the $(\Gamma,k)$-balanced fraction algebra $\wh{\rm Frac}_k(\mathcal{Z}^\omega_\Gamma)$, provided that $\epsilon \, \omega^m [\prod_{v\in \mathcal{V}} {\bf X}_v^{a_v}]_{\rm Weyl} \in \mathcal{Z}^\omega_\Gamma$ is $(\Gamma,k)$-balanced. In case $\epsilon \, \omega^m [\prod_{v\in \mathcal{V}} {\bf X}_v^{a_v}]_{\rm Weyl} \in \mathcal{Z}^\omega_\Gamma$ belongs to $\mathcal{X}^q_\Gamma$, one observes that the right-hand side of eq.\eqref{eq:nu_sharp_formula} is an element of ${\rm Frac}(\mathcal{X}^q_\Gamma)$. Combining these observations together with the fact that $\wh{\rm Frac}_k(\mathcal{Z}^\omega_\Gamma)$ is closed under addition, one observes that $\nu^{\sharp\omega}_k$ sends an element of $\wh{\rm Frac}_k(\mathcal{Z}^\omega_\Gamma)$ to an element of $\wh{\rm Frac}_k(\mathcal{Z}^\omega_\Gamma)$.

\vs

One last thing to check for the well-definedness of $\nu^{\sharp\omega}_k$ is whether the formula of eq.\eqref{eq:nu_sharp_formula} is consistent with the algebraic relations among $(\Gamma,k)$-balanced Laurent monomials. It suffices to check that
$$
\nu^{\sharp\omega}_k(\epsilon \, \omega^m [{\textstyle \prod}_{v\in \mathcal{V}} {\bf X}_v^{a_v}]_{\rm Weyl}) \, \nu^{\sharp\omega}_k(\epsilon' \omega^{m'} [{\textstyle \prod}_{v\in \mathcal{V}} {\bf X}_v^{b_v}]_{\rm Weyl}) = \nu^{\sharp\omega}_k(\epsilon \,\omega^m [{\textstyle \prod}_{v\in \mathcal{V}} {\bf X}_v^{a_v}]_{\rm Weyl} \cdot \epsilon' \omega^{m'} [{\textstyle \prod}_{v\in \mathcal{V}} {\bf X}_v^{b_v}]_{\rm Weyl})
$$
holds whenever $\epsilon \, \omega^m [{\textstyle \prod}_{v\in \mathcal{V}} {\bf X}_v^{a_v}]_{\rm Weyl}$ and $\epsilon' \omega^{m'} [{\textstyle \prod}_{v\in \mathcal{V}} {\bf X}_v^{b_v}]_{\rm Weyl}$ are $(\Gamma,k)$-balanced Laurent monomials. Notice that the argument of $\nu^{\sharp\omega}_k$ in the right-hand side is a $(\Gamma,k)$-balanced monomial $\epsilon'' \omega^{m''} [\prod_{v\in \mathcal{V}} {\bf X}_v^{a_v + b_v}]_{\rm Weyl}$, so the value of the right-hand side is also determined by the formula eq.\eqref{eq:nu_sharp_formula} to be
$$
\nu^{\sharp\omega}_k(\epsilon \,\omega^m [{\textstyle \prod}_{v\in \mathcal{V}} {\bf X}_v^{a_v}]_{\rm Weyl} \cdot \epsilon' \omega^{m'} [{\textstyle \prod}_{v\in \mathcal{V}} {\bf X}_v^{b_v}]_{\rm Weyl})
= (\epsilon \,\omega^m [{\textstyle \prod}_{v\in \mathcal{V}} {\bf X}_v^{a_v}]_{\rm Weyl} \cdot \epsilon' \omega^{m'} [{\textstyle \prod}_{v\in \mathcal{V}} {\bf X}_v^{b_v}]_{\rm Weyl}) \cdot F^q({\bf X}_k; \alpha+\beta),
$$ 
where $\alpha = \sum_{v\in \mathcal{V}} \varepsilon_{kv} a_v$ and $\beta = \sum_{v\in \mathcal{V}} \varepsilon_{kv} b_v$; we used the fact that $\sum_{v\in\mathcal{V}} \varepsilon_{kv}(a_v+b_v) = \alpha+\beta$. On the other hand, the value of the left-hand side is determined by eq.\eqref{eq:nu_sharp_formula} to be
\begin{align*}
& \nu^{\sharp\omega}_k(\epsilon \, \omega^m [{\textstyle \prod}_{v\in \mathcal{V}} {\bf X}_v^{a_v}]_{\rm Weyl}) \, \nu^{\sharp\omega}_k(\epsilon' \omega^{m'} [{\textstyle \prod}_{v\in \mathcal{V}} {\bf X}_v^{b_v}]_{\rm Weyl}) \\
& = (\epsilon \, \omega^m [{\textstyle \prod}_{v\in \mathcal{V}} {\bf X}_v^{a_v}]_{\rm Weyl}) \, \myul{F^q({\bf X}_k;\alpha) \, (\epsilon' \omega^{m'} [{\textstyle \prod}_{v\in \mathcal{V}} {\bf X}_v^{b_v}]_{\rm Weyl}) } \, F^q({\bf X}_k;\beta) \\
& = (\epsilon \, \omega^m [{\textstyle \prod}_{v\in \mathcal{V}} {\bf X}_v^{a_v}]_{\rm Weyl}) \, (\epsilon' \omega^{m'} [{\textstyle \prod}_{v\in \mathcal{V}} {\bf X}_v^{b_v}]_{\rm Weyl}) \, F^q(q^{2\beta} {\bf X}_k; \alpha) \, F^q({\bf X}_k;\beta)
\end{align*}
where the last equality is deduced from ${\bf X}_k [\prod_{v\in \mathcal{V}} {\bf X}_v^{b_v}]_{\rm Weyl} = q^{2\beta} [\prod_{v\in \mathcal{V}} {\bf X}_v^{b_v}]_{\rm Weyl} {\bf X}_k$ which in turn follows from Lem.\ref{lem:relations_of_cube-root_algebra_by_X}. The only thing that remains to be shown is that $F^q(q^{2\beta} {\bf X}_k;\alpha) F^q({\bf X}_k;\beta)$ equals $F^q({\bf X}_k;\alpha+\beta)$, which follows from the corresponding statement for the function $F^q$:
$$
F^q(x;\alpha+\beta) = F^q(q^{2\beta} x;\alpha) F^q(x;\beta) \mbox{ holds for all $\alpha,\beta \in \mathbb{Z}$}.
$$
One can see this using a formal definition of $F^q$ in eq.\eqref{eq:F_q_formal} as $F^q(x;\alpha+\beta) = \Psi^q(q^{2(\alpha+\beta)}x)\Psi^q(x)^{-1} $ \\ $= \Psi^q(q^{2\alpha} (q^{2\beta}x)) \Psi^q(q^{2\beta}x)^{-1} \Psi^q(q^{2\beta}x) \Psi^q(x)^{-1} = F^q(q^{2\beta}x;\alpha) F^q(x;\beta)$, and can also obtain a rigorous proof using eq.\eqref{eq:F_q}. \qed

\begin{remark}
This codomain subset $\wh{\rm Frac}_k(\mathcal{Z}^\omega_\Gamma)$ of ${\rm Frac}(\mathcal{Z}^\omega_\Gamma)$ is sufficient for the purpose of the present paper. However, one may seek for the most natural subset of ${\rm Frac}(\mathcal{Z}^\omega_\Gamma)$, e.g. on which the balanced cube-root quantum mutation can be applied at all possible unfrozen nodes. For the seed $\Gamma = \Gamma_\Delta$ associated to an ideal triangulation $\Delta$, our candidate is the $\Delta$-balanced fraction algebra $\wh{\rm Frac}(\mathcal{Z}^\omega_\Delta)$, based on the $\Delta$-balancedness condition on the powers $(a_v)_{v\in\mathcal{V}}$ of the Laurent monomials. We suggest that, for a general seed $\Gamma$, connected to the triangulation seeds $\Gamma_\Delta$ by mutations, a `$\Gamma$-balancedness' condition on $(a_v)_{v\in\mathcal{V}} \in (\frac{1}{3}\mathbb{Z})^\mathcal{V}$ should be defined by applying the tropical versions of the cluster $\mathscr{A}$-mutations to the $\Delta$-balanced elements of $\Gamma_\Delta$. In the present paper, we shall not try to verify whether this is indeed a good choice of a subset of ${\rm Frac}(\mathcal{Z}^\omega_\Gamma)$, nor try to come up with a description of it in terms of a general seed $\Gamma$.
\end{remark}

From Lem.\ref{lem:nu_prime_k_basic_properties}(2), and the formulas in eq.\eqref{eq:nu_sharp_formula} and eq.\eqref{eq:mu_sharp_formula}, one observes the following lemma.
\begin{lemma}
\label{lem:nu_omega_k_extending}
Above $\nu^\omega_k:\wh{\rm Frac}_k(\mathcal{Z}^\omega_{\Gamma'}) \to \wh{\rm Frac}_k(\mathcal{Z}^\omega_\Gamma)$ of Def.\ref{def:nu_omega_k} extends $\mu^q_k : {\rm Frac}(\mathcal{X}^q_{\Gamma'}) \to {\rm Frac}(\mathcal{X}^q_\Gamma)$ of Def.\ref{def:FG_quantum_mutation}. \qed
\end{lemma}
We establish one useful basic fact, to be used later:
\begin{lemma}
\label{lem:nu_omega_k_preserves_star}
Above $\nu^\omega_k : \wh{\rm Frac}_k(\mathcal{Z}^\omega_{\Gamma'}) \to \wh{\rm Frac}_k(\mathcal{Z}^\omega_\Gamma)$ preserves the $*$-structures. More precisely, $\wh{\rm Frac}_k(\mathcal{Z}^\omega_{\Gamma'})$ is closed under the $*$-structure of ${\rm Frac}(\mathcal{Z}^\omega_{\Gamma'})$ (induced by the $*$-structure of $\mathcal{Z}^\omega_{\Gamma'}$ defined in Def.\ref{def:FG_algebra}), $\wh{\rm Frac}_k(\mathcal{Z}^\omega_\Gamma)$ is closed under the $*$-structure of ${\rm Frac}(\mathcal{Z}^\omega_\Gamma)$, and for each ${\bf F} \in \wh{\rm Frac}_k(\mathcal{Z}^\omega_{\Gamma'})$ one has
$$
\nu^\omega_k({\bf F}^*) = ( \nu^\omega_k({\bf F}))^*
$$
\end{lemma}

{\it Proof.} Let's first show that $\wh{\rm Frac}_k(\mathcal{Z}^\omega_\Gamma)$ is closed under the $*$-map of ${\rm Frac}(\mathcal{Z}^\omega_\Gamma)$. We begin by describing this $*$-map in more detail. An element of ${\rm Frac}(\mathcal{Z}^\omega_\Gamma)$ is of the form ${\bf F} = {\bf P}{\bf Q}^{-1}$ with ${\bf P} \in \mathcal{Z}^\omega_\Gamma$ and $0\neq {\bf Q} \in \mathcal{X}^q_\Gamma \subset \mathcal{Z}^\omega_\Gamma$. By definition of the $*$-structure on the skew field of fractions ${\rm Frac}(\mathcal{Z}^\omega_\Gamma)$ induced by that of $\mathcal{Z}^\omega_\Gamma$, we have ${\bf F}^* = ({\bf Q}^*)^{-1} {\bf P}^*$, which can be written as ${\bf P}_1 {\bf Q}_1^{-1}$ for some ${\bf P}_1 \in \mathcal{Z}^\omega_\Gamma$ and $0\neq {\bf Q}_1 \in \mathcal{X}^q_\Gamma \subset \mathcal{Z}^\omega_\Gamma$ by some algebraic manipulations.

\vs

Now, in view of Def.\ref{def:Gamma_k_balanced}, the element ${\bf F} = {\bf P}{\bf Q}^{-1}$ belongs to $\wh{\rm Frac}_k(\mathcal{Z}^\omega_\Gamma)$ if and only if ${\bf P}$ is $(\Gamma,k)$-balanced in the sense of Def.\ref{def:Gamma_k_balanced}, that is, if ${\bf P}$ is a sum of $(\Gamma,k)$-balanced Laurent monomials. So, one can see that, in order to show that $\wh{\rm Frac}_k(\mathcal{Z}^\omega_\Gamma)$ is closed under the $*$-map, it suffices to deal with the case when ${\bf P}$ is a single $(\Gamma,k)$-balanced Laurent monomial ${\bf P} = [\prod_{v\in \mathcal{V}} {\bf X}_v^{a_v}]_{\rm Weyl}$; here $(a_v)_{v\in\mathcal{V}} \in (\frac{1}{3}\mathbb{Z})^\mathcal{V}$, and $\sum_{v\in\mathcal{V}} \varepsilon_{kv} a_v \in \mathbb{Z}$. Let's write ${\bf Q} = \sum \epsilon \, q^m [\prod_{v\in \mathcal{V}} {\bf X}_v^{b_v}]_{\rm Weyl}$; here $\sum$ is a finite sum, and we wrote a summand in a generic form, with $\epsilon \in \{1,-1\}$, $m \in \frac{1}{18}\mathbb{Z}$, $(b_v)_{v\in\mathcal{V}} \in \mathbb{Z}^\mathcal{V}$, without introducing the index variable for the sum. 

\vs

By Lem.\ref{lem:star-invariance_and_Weyl-ordering_new}(1), a Weyl-ordered Laurent monomial $[\sim]_{\rm Weyl}$ is $*$-invariant. Hence, ${\bf P}^* = {\bf P}$, and ${\bf Q}^* = \sum \epsilon \, q^{-m} [\prod_{v\in \mathcal{V}} {\bf X}_v^{b_v}]_{\rm Weyl}$, in view of the definition of the $*$-map in Def.\ref{def:FG_algebra}. Using Lem.\ref{lem:relations_of_cube-root_algebra_by_X}, one observes that the commutation relation of ${\bf P} = [\prod_{v\in\mathcal{V}} {\bf X}_v^{a_v}]_{\rm Weyl}$ and $[\prod_{v\in\mathcal{V}} {\bf X}_v^{b_v}]_{\rm Weyl}$ appearing in the summand of ${\bf Q}$ or ${\bf Q}^*$ reads
\begin{align}
\label{eq:Weyl-ordered_terms_commutation}
[{\textstyle \prod}_{v\in\mathcal{V}} {\bf X}_v^{b_v}]_{\rm Weyl} \, [{\textstyle \prod}_{v\in\mathcal{V}} {\bf X}_v^{a_v}]_{\rm Weyl} = q^{2c} \, [{\textstyle \prod}_{v\in\mathcal{V}} {\bf X}_v^{a_v}]_{\rm Weyl} \, [{\textstyle \prod}_{v\in\mathcal{V}} {\bf X}_v^{b_v}]_{\rm Weyl}, \quad
c = {\textstyle \sum}_{v,w\in\mathcal{V}} \varepsilon_{vw} b_v a_v.
\end{align}
Note that $c \in \frac{1}{6}\mathbb{Z}$ in the current setting, so $q^{2c} \in \mathbb{Z}[q^{\pm 1/3}] \subset \mathbb{Z}[q^{\pm 1/18}] = \mathbb{Z}[\omega^{\pm 1/2}]$ (but not necessarily $q^{2c} \in \mathbb{Z}[q^{\pm 1/2}]$). One can rewrite eq.\eqref{eq:Weyl-ordered_terms_commutation} in terms of conjugation by ${\bf P}$:
\begin{align}
\label{eq:Weyl-ordered_terms_commutation2}
{\bf P}^{-1} [{\textstyle \prod}_{v\in\mathcal{V}} {\bf X}_v^{b_v}]_{\rm Weyl} {\bf P} = q^{2c} [{\textstyle \prod}_{v\in\mathcal{V}} {\bf X}_v^{b_v}]_{\rm Weyl},
\end{align}
which can be used to move ${\bf P}$ around, as follows:
\begin{align*}
{\bf F}^* = ({\bf Q}^*)^{-1} {\bf P}^* & = ({\textstyle \sum} \epsilon \, q^{-m} [{\textstyle \prod}_{v\in\mathcal{V}} {\bf X}_v^{b_v}]_{\rm Weyl})^{-1} {\bf P} \\
& = {\bf P} {\bf P}^{-1} ({\textstyle \sum} \epsilon \, q^{-m} [{\textstyle \prod}_{v\in\mathcal{V}} {\bf X}_v^{b_v}]_{\rm Weyl})^{-1} {\bf P} \\
& = {\bf P} ({\textstyle \sum} \epsilon \, q^{-m} {\bf P}^{-1} [{\textstyle \prod}_{v\in\mathcal{V}} {\bf X}_v^{b_v}]_{\rm Weyl} {\bf P})^{-1} \\
& = {\bf P} ({\textstyle \sum} \epsilon \, q^{-m} q^{2c} [{\textstyle \prod}_{v\in\mathcal{V}} {\bf X}_v^{b_v}]_{\rm Weyl} )^{-1} \quad (\because\mbox{eq.\eqref{eq:Weyl-ordered_terms_commutation2}})
\end{align*}
Above is a typical example of what we referred to as an `algebraic manipulation' in the skew field of fractions ${\rm Frac}(\mathcal{Z}^\omega_\Gamma)$. So we see that ${\bf F}^*$ is of the form ${\bf P}_1 {\bf Q}_1^{-1}$, with ${\bf P}_1 = {\bf P}$ and ${\bf Q}_1 = {\textstyle \sum} \epsilon \, q^{-m} q^{2c} [{\textstyle \prod}_{v\in\mathcal{V}} {\bf X}_v^{b_v}]_{\rm Weyl} \in \mathcal{X}^q_\Gamma$. One may see here a reason why we used $\mathbb{Z}[q^{\pm 1/18}]=\mathbb{Z}[\omega^{\pm 1/2}]$ as the coefficient ring, instead of $\mathbb{Z}[q^{\pm 1/2}]$, in the definition of $\mathcal{X}^q_\Gamma$ in Def.\ref{def:FG_algebra}; see Rem.\ref{rem:coefficient_ring}. From ${\bf Q}_1 = {\bf P}^{-1} {\bf Q}^* {\bf P}$ we see that ${\bf Q} \neq 0  \Rightarrow {\bf Q}^*\neq 0 \Rightarrow {\bf Q}_1 \neq 0$. Note finally that ${\bf P}_1 = {\bf P}$ is $(\Gamma,k)$-balanced (because we assumed ${\bf P}$ is $(\Gamma,k)$-balanced), and therefore ${\bf F}^* \in \wh{\rm Frac}_k(\mathcal{Z}^\omega_\Gamma)$, as desired. This shows that $\wh{\rm Frac}_k(\mathcal{Z}^\omega_\Gamma)$ is closed under the $*$-map. This statement also applies to the seed $\Gamma'$, yielding that $\wh{\rm Frac}_k(\mathcal{Z}^\omega_{\Gamma'})$ is also closed under the $*$-map.

\vs

Let us first deal with the easiest case when ${\bf F}$ equals a single $(\Gamma',k)$-balanced Weyl-ordered Laurent monomial ${\bf F} = [\prod_{v\in\mathcal{V}} {{\bf X}_v'}^{a'_v}]_{\rm Weyl} \in \mathcal{Z}^\omega_{\Gamma'}$; here $(a'_v)_{v\in\mathcal{V}} \in (\frac{1}{3}\mathbb{Z})^\mathcal{V}$, with $\sum_{v\in\mathcal{V}} \varepsilon'_{kv} a'_v \in \mathbb{Z}$. In this case we have ${\bf F}^* = {\bf F}$ (see Lem.\ref{lem:star-invariance_and_Weyl-ordering_new}(1)). We observed in the proof of Lem.\ref{lem:nu_omega_u_well-defined} (eq.\eqref{eq:image_under_nu_prime_of_Laurent_monomial}--\eqref{eq:image_under_nu_prime_of_Laurent_monomial2}) that $\nu'_k[\prod_{v\in\mathcal{V}} {{\bf X}'_v}^{a'_v} ]_{\rm Weyl} = [\prod_{v\in\mathcal{V}} {\bf X}_v^{a_v}]_{\rm Weyl}$, where $a_v = a_v'$ for all $v\neq k$, and $a_k = -a'_k + \sum_{v\in\mathcal{V}} [\varepsilon_{vk}]_+ a'_v$, and also that $[\prod_{v\in\mathcal{V}} {\bf X}_v^{a_v}]_{\rm Weyl} \in \mathcal{Z}^\omega_\Gamma$ is $(\Gamma,k)$-balanced, i.e. $\alpha = \sum_{v\in\mathcal{V}} \varepsilon_{kv} a_v \in \mathbb{Z}$. Now eq.\eqref{eq:nu_sharp_formula} applies and gives $\nu^{\sharp \omega}_k([\prod_{v\in\mathcal{V}} {\bf X}_v^{a_v}]_{\rm Weyl}) = [\prod_{v\in\mathcal{V}} {\bf X}_v^{a_v}]_{\rm Weyl} \, F^q({\bf X}_k;\alpha)$, where $F^q$ is as in eq.\eqref{eq:F_q}. Thus
$$
\nu^\omega_k({\bf F}^*) = \nu^\omega_k({\bf F}) = \nu^{\sharp \omega}_k(\nu'_k([{\textstyle \prod}_{v\in \mathcal{V}} {{\bf X}_v'}^{a'_v}]_{\rm Weyl})) = [{\textstyle \prod}_{v\in \mathcal{V}} {\bf X}_v^{a_v}]_{\rm Weyl} \, F^q({\bf X}_k;\alpha).
$$
On the other hand, note that
\begin{align*}
(\nu^\omega_k({\bf F}))^* & = ( [{\textstyle \prod}_{v\in \mathcal{V}} {\bf X}_v^{a_v}]_{\rm Weyl} \, F^q({\bf X}_k;\alpha) )^* \\
& = ( F^q({\bf X}_k;\alpha) )^* \, ( [{\textstyle \prod}_{v\in \mathcal{V}} {\bf X}_v^{a_v}]_{\rm Weyl})^*  \\
& = F^{1/q}({\bf X}_k;\alpha) \, [{\textstyle \prod}_{v\in \mathcal{V}} {\bf X}_v^{a_v}]_{\rm Weyl} \qquad (\because {\bf X}_k^*={\bf X}_k, \, q^* = 1/q, \mbox{Lem.\ref{lem:star-invariance_and_Weyl-ordering_new}(1)}).
\end{align*}
In order to move $[{\textstyle \prod}_{v\in \mathcal{V}} {\bf X}_v^{a_v}]_{\rm Weyl}$ to the left, we perform an algebraic manipulation as before. First note from Lem.\ref{lem:relations_of_cube-root_algebra_by_X} the commutation relation ${\bf X}_k [{\textstyle \prod}_{v\in \mathcal{V}} {\bf X}_v^{a_v}]_{\rm Weyl} = q^{2\sum_{v\in\mathcal{V}} \varepsilon_{kv} a_v} [{\textstyle \prod}_{v\in \mathcal{V}} {\bf X}_v^{a_v}]_{\rm Weyl} {\bf X}_k$, thus
\begin{align}
\label{eq:Weyl-ordered_terms_commutation3}
([{\textstyle \prod}_{v\in \mathcal{V}} {\bf X}_v^{a_v}]_{\rm Weyl} )^{-1 } {\bf X}_k [{\textstyle \prod}_{v\in \mathcal{V}} {\bf X}_v^{a_v}]_{\rm Weyl} = q^{2\alpha} {\bf X}_k.
\end{align}
Hence, moving $[{\textstyle \prod}_{v\in\mathcal{V}} {\bf X}_v^{a_v}]_{\rm Weyl}$ to the left, we get
\begin{align*}
(\nu^\omega_k({\bf F}))^* & =[{\textstyle \prod}_{v\in\mathcal{V}} {\bf X}_v^{a_v}]_{\rm Weyl} \, \, F^{1/q}( ([{\textstyle \prod}_{v\in\mathcal{V}} {\bf X}_v^{a_v}]_{\rm Weyl} \,)^{-1} {\bf X}_k [{\textstyle \prod}_{v\in\mathcal{V}} {\bf X}_v^{a_v}]_{\rm Weyl} \,;\alpha) \\
& =
[{\textstyle \prod}_{v\in\mathcal{V}} {\bf X}_v^{a_v}]_{\rm Weyl} \, F^{1/q}(q^{2\alpha} {\bf X}_k;\alpha).
\end{align*}
We claim that
$$
F^{1/q}(q^{2m} x;m) = F^q(x;m)
$$
holds as equality of functions, for any $m\in \mathbb{Z}$. Indeed, from the definition of $F^q$ as in eq.\eqref{eq:F_q}, in case $m\ge 0$ we have
\begin{align*}
F^{1/q}(q^{2m} x;m) & = (1+q^{-1}(q^{2m}x))(1+q^{-3}(q^{2m}x)) \cdots (1+q^{-(2m-1)}(q^{2m}x)) \\
& = (1+q^{2m-1}x)(1+q^{2m-3}x)\cdots(1+qx) \\
& = (1+qx) \cdots (1+q^{2m-3}x)(1+q^{2m-1}x) = F^q(x;m),
\end{align*}
and in case $m<0$ we have
\begin{align*}
F^{1/q}(q^{2m} x;m) &  = (1+ q (q^{2m} x))^{-1} (1+ q^3 (q^{2m} x))^{-1}  \cdots (1+ q^{2(-m)-1} (q^{2m} x))^{-1} \\
& = (1+q^{2m+1} x)^{-1} (1+q^{2m+3} x)^{-1} \cdots (1+q^{-3} x)^{-1} (1+q^{-1} x)^{-1} \\
& = (1+q^{-1} x)^{-1} (1+q^{-3} x)^{-1} \cdots (1+q^{-(2(-m)-3)})^{-1} (1+q^{-(2(-m)-1)})^{-1} = F^q(x;m).
\end{align*}
So we get $(\nu^\omega_k({\bf F}))^* = [\prod_{v\in\mathcal{V}} {\bf X}_v^{a_v}]_{\rm Weyl} F^q({\bf X}_k;\alpha) = \nu^\omega_k({\bf F})$, hence finally $\nu^\omega_k({\bf F}^*) = (\nu^\omega_k({\bf F}))^*$, as desired.

\vs

We just showed that $\nu^\omega_k({\bf F}^*) = (\nu^\omega_k({\bf F}))^*$ holds for every $(\Gamma',k)$-balanced Weyl-ordered Laurent monomial ${\bf F}$ for $\Gamma'$. By the complex-conjugate linearity of the $*$-maps, we see that this equality holds for any $(\Gamma',k)$-balanced Laurent polynomial ${\bf F} \in \mathcal{Z}^\omega_{\Gamma'}$ for $\Gamma'$. Since any element $\mathcal{X}^q_{\Gamma'} \subset \mathcal{Z}^\omega_{\Gamma'}$ is a $(\Gamma',k)$-balanced Laurent polynomial, this equality holds for all ${\bf F} \in \mathcal{X}^q_{\Gamma'}$. A general element ${\bf F}$ of $\wh{\rm Frac}_k(\mathcal{Z}^\omega_{\Gamma'})$ is of the form ${\bf P}{\bf Q}^{-1}$ with ${\bf P} \in \mathcal{Z}^\omega_{\Gamma'}$ being a $(\Gamma',k)$-balanced Laurent polynomial for $\Gamma'$, and $0\neq {\bf Q} \in \mathcal{X}^q_{\Gamma'}$. We just saw that $\nu^\omega_k({\bf P}^*) = (\nu^\omega_k({\bf P}))^*$ and $\nu^\omega_k({\bf Q}^*) = (\nu^\omega_k({\bf Q}))^*$ holds. Let's assume for now that ${\bf P}$ is a Laurent monomial $[\prod_{v\in\mathcal{V}} {{\bf X}_v'}^{a'_v}]_{\rm Weyl}$; in particular, ${\bf P}^* = {\bf P}$ in this case. Using the computation of ${\bf F}^*$ we performed above, we have ${\bf F}^* = {\bf P}{\bf Q}_1^{-1}$, with ${\bf Q}_1 = {\bf P}^{-1} {\bf Q}^* {\bf P}$. Note that $\nu^\omega_k({\bf F}) = \nu^\omega_k({\bf P}) (\nu^\omega_k({\bf Q}))^{-1}$, and similarly $\nu^\omega_k({\bf F}^*) = \nu^\omega_k({\bf P}) (\nu^\omega_k({\bf Q}_1))^{-1}$; here $\nu^\omega_k({\bf Q}), \nu^\omega_k({\bf Q}_1) \in \mathcal{X}^q_\Gamma$ (see e.g. Lem.\ref{lem:nu_omega_k_extending}). Note that
\begin{align*}
(\nu^\omega_k({\bf F}))^* & = ( (\nu^\omega_k({\bf Q}))^* )^{-1} \, (\nu^\omega_k({\bf P}))^* \\
& = (\nu^\omega_k({\bf Q}^*))^{-1} \nu^\omega_k({\bf P}^*) \\
& = (\nu^\omega_k({\bf Q}^*))^{-1} \nu^\omega_k({\bf P}) \\
& = \nu^\omega_k({\bf P}) (\nu^\omega_k({\bf P}))^{-1} (\nu^\omega_k({\bf Q}^*))^{-1} \nu^\omega_k({\bf P}) \\
& = \nu^\omega_k({\bf P})  \left( (\nu^\omega_k({\bf P}))^{-1} \, \nu^\omega_k({\bf Q}^*) \, \nu^\omega_k({\bf P}) \right)^{-1}  \\
& = \nu^\omega_k({\bf P}) (\nu^\omega_k( {\bf P}^{-1} {\bf Q}^* {\bf P}) )^{-1} \quad (\because\mbox{$\nu^\omega_k$ is a homomorphism}) \\
& = \nu^\omega_k({\bf P}) (\nu^\omega_k({\bf Q}_1))^{-1} \\
& = \nu^\omega_k({\bf F}^*).
\end{align*}
So $\nu^\omega_k({\bf F}^*) = (\nu^\omega_k({\bf F}))^*$ holds for such ${\bf F}$. Any general ${\bf F}$ of $\wh{\rm Frac}_k(\mathcal{Z}^\omega_{\Gamma'})$ is a $\mathbb{Z}[\omega^{1/2}]$-linear combination of such kind of ${\bf F}$ we just dealt with. By the complex-conjugate linearity of the $*$-maps, $\nu^\omega_k({\bf F}^*) = (\nu^\omega_k({\bf F}))^*$ holds for such a general ${\bf F}$. \qed

\vs

A useful observation that one can obtain by similar arguments as presented in the beginning of the proof of Lem.\ref{lem:nu_omega_k_preserves_star} is that for any ideal triangulation $\Delta$ of a triangulable punctured surface, the $\Delta$-balanced fraction algebra $\wh{\rm Frac}(\mathcal{Z}^\omega_\Delta)$ is closed under the $*$-structure of ${\rm Frac}(\mathcal{Z}^\omega_\Delta)$.

\vs

The following definition is the second main definition of the present paper.

\begin{definition}[balanced cube-root quantum coordinate change map for a flip]
\label{def:Theta}
Let $\Delta$ and $\Delta'$ be ideal triangulations of a triangulable generalized marked surface $\frak{S}$ that differ by a flip at an arc $i$, with notations as in Fig.\ref{fig:mutations_for_a_flip}. Define the \ul{\em balanced (cube-root) quantum coordinate change map} between the balanced fractions algebras (Def.\ref{def:balanced_fraction_algebra})
$$
\Theta^\omega_{\Delta \Delta'} = \Theta^\omega_i ~:~ \wh{\rm Frac}(\mathcal{Z}^\omega_{\Delta'}) \to \wh{\rm Frac}(\mathcal{Z}^\omega_{\Delta})
$$
as
$$
\Theta^\omega_{\Delta\Delta'} := \nu^\omega_{v_3} \nu^\omega_{v_4} \nu^\omega_{v_7} \nu^\omega_{v_{12}}.
$$
\end{definition}
The following two lemmas establish the well-definedness of $\Theta^\omega_{\Delta\Delta'} = \Theta^\omega_i$.

\begin{lemma}
\label{lem:Theta_well-definedness1}
The above $\Theta^\omega_{\Delta\Delta'}$ is a well-defined map from $\wh{\rm Frac}(\mathcal{Z}^\omega_{\Delta'})$ to ${\rm Frac}(\mathcal{Z}^\omega_\Delta)$.
\end{lemma}

\begin{lemma}
\label{lem:Theta_well-definedness2}
The image of $\Theta^\omega_{\Delta\Delta'}$ lies in the $\Delta$-balanced fraction algebra $\wh{\rm Frac}(\mathcal{Z}^\omega_\Delta)$ (Def.\ref{def:balanced_fraction_algebra}).
\end{lemma}

As a preliminary step for a proof of the above two lemmas, we first establish the following simple observation.

\begin{lemma}
\label{lem:alpha_alternative}
Let $v^{(r)}$ be as in eq.\eqref{eq:four_mutated_nodes}, for $r=1,2,3,4$. Denote the set of nodes appearing in eq.\eqref{eq:V_for_four_mutations} by $\mathcal{V}$. Let $(a'_v)_{v\in \mathcal{V}} = (a^{(4)}_v)_{v\in\mathcal{V}} \in (\frac{1}{3}\mathbb{Z})^\mathcal{V}$, which is not necessarily assumed to be $\Delta'$-balanced.

\vs

Recursively define $(a_v^{(r-1)})_{v\in\mathcal{V}} \in (\frac{1}{3} \mathbb{Z})^\mathcal{V}$, for $r=4,3,2,1$, using the formula in eq.\eqref{eq:a_r_v_recursive} applied to $(a'_v)_{v\in \mathcal{V}} = (a^{(4)}_v)_{v\in\mathcal{V}}$. Define the numbers $\alpha^{(r)} \in \frac{1}{3}\mathbb{Z}$, for $r=1,2,3,4$, by the formula in eq.\eqref{eq:alpha_for_nu_sharp}, i.e. we let $\alpha^{(r)} = \sum_{v\in \mathcal{V}} \varepsilon^{(r-1)}_{v^{(r)},v} a^{(r-1)}_v$. Then 
$$
\alpha^{(r)} = - \underset{v\in \mathcal{V}}{\textstyle \sum} \varepsilon^{(r)}_{v^{(r)},v} a^{(r)}_v, \qquad \forall r=1,2,3,4.
$$
\end{lemma}

{\it Proof of Lem.\ref{lem:alpha_alternative}.} The statement of this lemma is a special case of eq.\eqref{eq:alpha_equals_minus_alpha_prime} which we established in the proof of Lem.\ref{lem:nu_omega_u_well-defined}. To make sure that the current situation is indeed an example of the setting for eq.\eqref{eq:alpha_equals_minus_alpha_prime}, one can verify by inspection that eq.\eqref{eq:a_r_v_recursive} is an example of eq.\eqref{eq:image_under_nu_prime_of_Laurent_monomial2}; in fact, as mentioned in the proof of Lem.\ref{lem:nu_omega_u_well-defined}, what we need from eq.\eqref{eq:a_r_v_recursive} in order to apply eq.\eqref{eq:alpha_equals_minus_alpha_prime} is just $a^{(r-1)}_v = a^{(r)}_v$, $\forall v\neq v^{(r)}$ . \qed

\vs

{\it Proof of Lem.\ref{lem:Theta_well-definedness1}.}  What needs to be shown is that the four maps $\nu^\omega_{v^{(r)}}$, $r=4,3,2,1$, can indeed be applied in this order to each element of $\wh{\rm Frac}(\mathcal{Z}^\omega_{\Delta'})$ (Def.\ref{def:balanced_fraction_algebra}); that is, we should make sure that at each stage, the relevant element lies in the domain of $\nu^\omega_{v^{(r)}}$, namely the $(\Gamma_{\Delta^{(r)}},v^{(r)})$-balanced fraction algebra $\wh{\rm Frac}_{v^{(r)}}(\mathcal{Z}^\omega_{\Delta^{(r)}})$ (Def.\ref{def:Gamma_k_balanced}). If we restrict to ${\rm Frac}(\mathcal{X}^q_{\Delta'}) \subset \wh{\rm Frac}(\mathcal{Z}^\omega_{\Delta'})$, then by Lem.\ref{lem:nu_omega_k_extending} the map $\Theta^\omega_{\Delta\Delta'} = \nu^\omega_{v^{(1)}} \nu^\omega_{v^{(2)}} \nu^\omega_{v^{(3)}} \nu^\omega_{v^{(4)}}$ coincides with $\Phi^q_{\Delta\Delta'} = \mu^q_{v^{(1)}} \mu^q_{v^{(2)}} \mu^q_{v^{(3)}} \mu^q_{v^{(4)}}$; in particular, it is well defined on ${\rm Frac}(\mathcal{X}^q_{\Delta'})$. Therefore, it suffices to show that $\Theta^\omega_{\Delta\Delta'}$ is well defined on the $\Delta'$-balanced Laurent polynomials in $\mathcal{Z}^\omega_{\Delta'}$, or just on the $\Delta'$-balanced Laurent monomials in $\mathcal{Z}^\omega_{\Delta'}$ (see Def.\ref{def:balanced_subalgebras} for $\Delta'$-balancedness).

\vs

Let ${\bf U}'  = [\prod_{v\in \mathcal{V}} ({\bf X}_v')^{a'_v}]_{\rm Weyl} \in \mathcal{Z}^\omega_{\Delta'}$, with $(a'_v)_{v\in\mathcal{V}}  = (a^{(4)}_v)_{v\in \mathcal{V}} \in (\frac{1}{3}\mathbb{Z})^\mathcal{V}$, be a $\Delta'$-balanced Laurent monomial. So $(a'_v)_{v\in\mathcal{V}} = (a^{(4)}_v)_{v\in\mathcal{V}}$ is a $\Delta'$-balanced element of $(\frac{1}{3}\mathbb{Z})^\mathcal{V}$ in the sense of Def.\ref{def:Delta-balanced_elements}. To ease the computation, we establish some auxiliary elements; let ${\bf U}^{(4)} := {\bf U}' = [\prod_{v\in \mathcal{V}} ({\bf X}_v^{(4)})^{a^{(4)}_v}]_{\rm Weyl} \in \mathcal{Z}^\omega_{\Delta^{(4)}}$, and for each $r=4,3,2,1$ define
$$
{\bf U}^{(r-1)} := \nu'_{v^{(r)}} ({\bf U}^{(r)}) \in \mathcal{Z}^\omega_{\Delta^{(r-1)}}.
$$
Since $\nu'_{v^{(r)}} : \mathcal{Z}^\omega_{\Delta^{(r)}} \to \mathcal{Z}^\omega_{\Delta^{(r-1)}}$ sends a Weyl-ordered Laurent monomial to a Weyl-ordered monomial for $r=4,3,2,1$ ($\because$ Lem.\ref{lem:nu_prime_k_basic_properties}(5)), and since ${\bf U}^{(4)}$ is a Weyl-ordered Laurent monomial, it follows that ${\bf U}^{(r-1)}$ is a Weyl-ordered Laurent monomial, for $r=4,3,2,1$. Hence, for each $r=4,3,2,1$, one can write
\begin{align}
\label{eq:U_r-1_as_Laurent_monomial}
{\bf U}^{(r-1)} = [{\textstyle \prod}_{v\in\mathcal{V}} ({\bf X}_v^{(r-1)})^{a^{(r-1)}_v}]_{\rm Weyl}
\end{align}
for some $(a^{(r-1)}_v)_{v\in \mathcal{V}} \in (\frac{1}{3}\mathbb{Z})^\mathcal{V}$. For each $r=4,3,2,1$, one can compute $(a_v^{(r-1)})_{v\in\mathcal{V}}$ in terms of $(a_v^{(r)})_{v\in\mathcal{V}}$ using the formula in eq.\eqref{eq:nu_prime_formula}; namely, we can deduce that $a_v^{(r-1)} = a_v^{(r)}$ for all $v\neq v^{(r)}$, and that $a_{v^{(r)}}^{(r-1)} = - a_{v^{(r)}}^{(r)}+ \sum_{v\in\mathcal{V}} [\varepsilon_{v,v^{(r)}}^{(r-1)}]_+ a_v^{(r)}$, which exactly matches eq.\eqref{eq:a_r_v_recursive}. Therefore, Lem.\ref{lem:alpha} applies to the current situation. 

\vs

\ul{Claim} : For each $r=4,3,2,1$, the Laurent monomial ${\bf U}^{(r)} \in \mathcal{Z}^\omega_{\Delta^{(r)}}$ is $(\Gamma_{\Delta^{(r)}},v^{(r)})$-balanced (Def.\ref{def:Gamma_k_balanced}).

\vs

In view of Def.\ref{def:Gamma_k_balanced}, to prove this Claim we need to check whether the number $\sum_{v\in \mathcal{V}} \varepsilon^{(r)}_{v^{(r)},v} a^{(r)}_v \in \frac{1}{3}\mathbb{Z}$ belongs to $\mathbb{Z}$ for $r=4,3,2,1$. This indeed holds by Lem.\ref{lem:alpha} and Lem.\ref{lem:alpha_alternative}; we note that this integer is denoted by $-\alpha^{(r)}$ in those lemmas. Now that we proved the Claim, we can apply $\nu^\omega_{v^{(r)}}$ to ${\bf U}^{(r)}$ for each $r=4,3,2,1$, in view of Def.\ref{def:nu_omega_k} (and Lem.\ref{lem:nu_omega_u_well-defined}); the result is
\begin{align*}
\nu^\omega_{v^{(r)}}({\bf U}^{(r)})
= \nu^{\sharp\omega}_{v^{(r)}} \nu'_{v^{(r)}}({\bf U}^{(r)})
= \nu^{\sharp\omega}_{v^{(r)}} ({\bf U}^{(r-1)})
& = \nu^{\sharp\omega}_{v^{(r)}} ([{\textstyle \prod}_{v\in \mathcal{V}} ({\bf X}_v^{(r-1)})^{a_v^{(r-1)}}]_{\rm Weyl}) \\
& \hspace{-4mm} \stackrel{{\rm eq}.\eqref{eq:nu_sharp_formula}}{=} [{\textstyle \prod}_{v\in \mathcal{V}} ({\bf X}_v^{(r-1)})^{a_v^{(r-1)}}]_{\rm Weyl} \cdot F^q({\bf X}_{v^{(r)}}^{(r-1)}; \underset{v\in \mathcal{V}}{\textstyle \sum} \varepsilon_{v^{(r)},v}^{(r-1)} a_v^{(r-1)}) \\
& \hspace{-4mm} \stackrel{{\rm eq}.\eqref{eq:alpha_for_nu_sharp}}{=}  {\bf U}^{(r-1)} \cdot F^q({\bf X}_{v^{(r)}}^{(r-1)}; \alpha^{(r)}).
\end{align*}
For convenience, let ${\bf V}^{(r-1)} := F^q({\bf X}^{(r-1)}_{v^{(r)}};\alpha^{(r)})$, for $r=4,3,2,1$. Then we have
\begin{align}
\label{eq:nu_omega_v_r_U_r}
\nu^\omega_{v^{(r)}}({\bf U}^{(r)}) = {\bf U}^{(r-1)} {\bf V}^{(r-1)}, \quad \forall r=4,3,2,1,
\end{align}
where ${\bf V}^{(r-1)} \in {\rm Frac}(\mathcal{X}^q_{\Delta^{(r-1)}})$.

\vs

Let's now start applying the maps $\nu^\omega_{v^{(r)}}$, $r=4,3,2,1$, to the element ${\bf U}' \in \mathcal{Z}^\omega_{\Delta'}$ in this order; we should make sure that these maps are indeed applicable, and then investigate how the final result looks like. First, by the above Claim, we can observe that ${\bf U}' = {\bf U}^{(4)}$ is in the domain of $\nu^\omega_{v^{(4)}} : \wh{\rm Frac}_{v^{(4)}}(\mathcal{Z}^\omega_{\Delta^{(4)}}) \to \wh{\rm Frac}_{v^{(4)}} (\mathcal{Z}^\omega_{\Delta^{(3)}})$. In eq.\eqref{eq:nu_omega_v_r_U_r}, we computed the image as $\nu^\omega_{v^{(4)}}({\bf U}^{(4)}) = {\bf U}^{(3)} {\bf V}^{(3)}$, with ${\bf V}^{(3)} \in {\rm Frac}(\mathcal{X}^q_{\Delta^{(3)}})$. Hence, by the Claim for ${\bf U}^{(3)}$, we see that $\nu^\omega_{v^{(4)}}({\bf U}^{(4)}) = {\bf U}^{(3)} {\bf V}^{(3)}$ belongs to $\wh{\rm Frac}_{v^{(3)}}(\mathcal{Z}^\omega_{\Delta^{(3)}})$, which is the domain of $\nu^\omega_{v^{(3)}} : \wh{\rm Frac}_{v^{(3)}}(\mathcal{Z}^\omega_{\Delta^{(3)}}) \to \wh{\rm Frac}_{v^{(3)}}(\mathcal{Z}^\omega_{\Delta^{(2)}})$. The image is
\begin{align}
\label{eq:two_nu_applied_to_Laurent_monomial}
\nu^\omega_{v^{(3)}} (\nu^\omega_{v^{(4)}}({\bf U}^{(4)}) ) = \nu^\omega_{v^{(3)}} ( {\bf U}^{(3)} {\bf V}^{(3)})
= \nu^\omega_{v^{(3)}} ({\bf U}^{(3)}) \nu^\omega_{v^{(3)}}({\bf V}^{(3)})
= {\bf U}^{(2)} \underbrace{ {\bf V}^{(2)} \mu^q_{v^{(3)}} ({\bf V}^{(3)})},
\end{align}
where for the last equality we used eq.\eqref{eq:nu_omega_v_r_U_r} and Lem.\ref{lem:nu_omega_k_extending}. Since the underbraced part belongs to ${\rm Frac}(\mathcal{X}^q_{\Delta^{(2)}})$, using the Claim for ${\bf U}^{(2)}$ we get that $\nu^\omega_{v^{(3)}} (\nu^\omega_{v^{(4)}}({\bf U}^{(4)}) ) = {\bf U}^{(2)} {\bf V}^{(2)} \mu^q_{v^{(3)}} ({\bf V}^{(3)})$ belongs to $\wh{\rm Frac}_{v^{(2)}}(\mathcal{Z}^\omega_{\Delta^{(2)}})$, which is the domain of $\nu^\omega_{v^{(2)}} : \wh{\rm Frac}_{v^{(2)}} (\mathcal{Z}^\omega_{\Delta^{(2)}}) \to \wh{\rm Frac}_{v^{(2)}}(\mathcal{Z}^\omega_{\Delta^{(1)}})$. The image is
\begin{align*}
\nu^\omega_{v^{(2)}}(\nu^\omega_{v^{(3)}}(\nu^\omega_{v^{(4)}}({\bf U}^{(4)})))
= \nu^\omega_{v^{(2)}}( {\bf U}^{(2)} {\bf V}^{(2)} \mu^q_{v^{(3)}}({\bf V}^{(3)}))
= {\bf U}^{(1)} \underbrace{ {\bf V}^{(1)} \mu^q_{v^{(2)}}({\bf V}^{(2)} \mu^q_{v^{(3)}}({\bf V}^{(3)}))},
\end{align*}
where for the last equality we used eq.\eqref{eq:nu_omega_v_r_U_r} and Lem.\ref{lem:nu_omega_k_extending}. Since the underbraced expression is an element of ${\rm Frac}(\mathcal{X}^q_{\Delta^{(1)}})$, using the Claim for ${\bf U}^{(1)}$ we get that $\nu^\omega_{v^{(2)}}(\nu^\omega_{v^{(3)}}(\nu^\omega_{v^{(4)}}({\bf U}^{(4)})))$ belongs to $\wh{\rm Frac}_{v^{(1)}}(\mathcal{Z}^\omega_{\Delta^{(1)}})$, which is the domain of $\nu^\omega_{v^{(1)}} : \wh{\rm Frac}_{v^{(1)}}(\mathcal{Z}^\omega_{\Delta^{(1)}}) \to \wh{\rm Frac}_{v^{(1)}}(\mathcal{Z}^\omega_{\Delta^{(0)}})$. Therefore, $\nu^\omega_{v^{(1)}} (\nu^\omega_{v^{(2)}}(\nu^\omega_{v^{(3)}}(\nu^\omega_{v^{(4)}}({\bf U}^{(4)}))))$ makes sense as a well-defined element of $\wh{\rm Frac}_{v^{(1)}}(\mathcal{Z}^\omega_{\Delta^{(0)}}) \subset {\rm Frac}(\mathcal{Z}^\omega_{\Delta^{(0)}}) \subset {\rm Frac}(\mathcal{Z}^\omega_\Delta)$. This finishes the proof. \qed

\vs

{\it Proof of Lem.\ref{lem:Theta_well-definedness2}.} In the beginning of the above proof of Lem.\ref{lem:Theta_well-definedness1}, we observed that $\Theta^\omega_{\Delta\Delta'}$ sends an element of ${\rm Frac}(\mathcal{X}^q_{\Delta'}) \subset \wh{\rm Frac}(\mathcal{Z}^\omega_{\Delta'}) \subset {\rm Frac}(\mathcal{Z}^\omega_{\Delta'})$ to an element of ${\rm Frac}(\mathcal{X}^q_\Delta) \subset {\rm Frac}(\mathcal{Z}^\omega_\Delta)$ which is contained in the $\Delta$-balanced fraction algebra $\wh{\rm Frac}(\mathcal{Z}^\omega_\Delta)$ (Def.\ref{def:balanced_fraction_algebra}). So it suffices to show that $\Theta^\omega_{\Delta\Delta'}$ sends a $\Delta'$-balanced Laurent polynomial in $\mathcal{Z}^\omega_{\Delta'}$ to a $\Delta$-balanced Laurent polynomial in $\mathcal{Z}^\omega_\Delta$. We will show that $\Theta^\omega_{\Delta\Delta'}$ sends a $\Delta'$-balanced Laurent monomial in $\mathcal{Z}^\omega_{\Delta'}$ to a $\Delta$-balanced Laurent monomial in $\mathcal{Z}^\omega_\Delta$. 

\vs

We use the same notations as in the above proof of Lem.\ref{lem:Theta_well-definedness1}. In particular ${\bf U}' = {\bf U}^{(4)} = [\prod_{v\in \mathcal{V}} ({\bf X}_v')^{a'_v}]_{\rm Weyl} \in \mathcal{Z}^\omega_{\Delta'}$ is a $\Delta'$-balanced Laurent monomial (Def.\ref{def:balanced_subalgebras}), with $(a'_v)_{v\in \mathcal{V}} \in (\frac{1}{3}\mathbb{Z})^\mathcal{V}$; in view of Def.\ref{def:balanced_subalgebras}, we see that $(a'_v)_{v\in \mathcal{V}}$ is $\Delta'$-balanced in the sense of Def.\ref{def:Delta-balanced_elements}. To continue the last steps of the proof of Lem.\ref{lem:Theta_well-definedness1}, we arrive at the following computational result for the image of ${\bf U}' = {\bf U}^{(4)}$ under $\Theta^\omega_{\Delta\Delta'} = \nu^\omega_{v^{(1)}}\nu^\omega_{v^{(2)}}\nu^\omega_{v^{(3)}}\nu^\omega_{v^{(4)}}$:
\begin{align*}
\Theta^\omega_{\Delta\Delta'}({\bf U}) = \nu^\omega_{v^{(1)}} (\nu^\omega_{v^{(2)}}(\nu^\omega_{v^{(3)}}(\nu^\omega_{v^{(4)}}({\bf U}^{(4)})))) & = \nu^\omega_{v^{(1)}} ( {\bf U}^{(1)} {\bf V}^{(1)} \mu^q_{v^{(2)}}({\bf V}^{(2)} \mu^q_{v^{(3)}}({\bf V}^{(3)})) ) \\
& = {\bf U}^{(0)} \underbrace{ {\bf V}^{(0)} \mu^q_{v^{(1)}}( \underset{\in \, {\rm Frac}(\mathcal{X}^q_{\Delta^{(1)}})}{\myul{{\bf V}^{(1)} \mu^q_{v^{(2)}}({\bf V}^{(2)} \mu^q_{v^{(3)}}({\bf V}^{(3)})) }} )}_{\in \, {\rm Frac}(\mathcal{X}^q_{\Delta^{(0)}})},
\end{align*}
where for the last equality we used eq.\eqref{eq:nu_omega_v_r_U_r} and Lem.\ref{lem:nu_omega_k_extending}. Note that the underbraced part belongs to ${\rm Frac}(\mathcal{X}^q_{\Delta^{(0)}}) = {\rm Frac}(\mathcal{X}^q_\Delta)$ which is contained in $\wh{\rm Frac}(\mathcal{Z}^\omega_\Delta)$. Therefore, to show $\Theta^\omega_{\Delta\Delta'}({\bf U}) \in \wh{\rm Frac}(\mathcal{Z}^\omega_\Delta)$ it suffices to show ${\bf U}^{(0)} \in \wh{\rm Frac}(\mathcal{Z}^\omega_\Delta)$; as ${\bf U}^{(0)} = [\prod_{v\in \mathcal{V}} ({\bf X}^{(0)}_v)^{a_v^{(0)}}]_{\rm Weyl}$ is a Laurent monomial in $\mathcal{Z}^\omega_{\Delta^{(0)}} = \mathcal{Z}^\omega_\Delta$ as written in eq.\eqref{eq:U_r-1_as_Laurent_monomial}, we need to show that ${\bf U}^{(0)}$ is $\Delta$-balanced (Def.\ref{def:balanced_subalgebras}), i.e. that $(a_v^{(0)})_{v\in \mathcal{V}} \in (\frac{1}{3}\mathbb{Z})^\mathcal{V}$ is $\Delta$-balanced (Def.\ref{def:Delta-balanced_elements}). For convenience, write $(a_v^{(0)})_{v\in \mathcal{V}} = (a_v)_{v\in\mathcal{V}}$. In view of Fig.\ref{fig:mutations_for_a_flip}, the $\Delta$-balancedness condition of $(a_v)_{v\in \mathcal{V}}$ for the two triangles of $\Delta$ appearing in Fig.\ref{fig:mutations_for_a_flip} which are the triangles being involved in the flip $\Delta \leadsto \Delta'$ says that the following fifteen numbers are integers:
\begin{align}
\label{eq:fifteen_numbers_Delta}
\left\{
{\renewcommand{\arraystretch}{1.3} \begin{array}{rl}
{\rm (BE1)} : & a_1 +a_3+a_5, \quad  a_2+a_4+a_6,  \quad a_8+a_{10}+a_4, \quad a_9+a_{11}+a_3, \\
{\rm (BE2)} : & a_1+a_2, \quad a_3+a_4, \quad a_5+a_6,  \quad a_8+a_9, \quad a_{10}+a_{11}, \\
{\rm (BE3)}: & -a_7+a_6+a_1, \quad  -a_7+a_2+a_3, \quad -a_7+a_4+a_5, \\
& -a_{12}+a_3+a_8, \quad -a_{12}+a_9+a_{10}, \quad -a_{12}+a_{11}+a_4.
\end{array}}
\right.
\end{align}

In the proof of Lem.\ref{lem:alpha}, we expressed $(a_v)_{v\in\mathcal{V}}$ in terms of $(a'_v)_{v\in\mathcal{V}}$, through the formulas written in eq.\eqref{eq:a_0_as_a_prime}. As we observed in the proof of Lem.\ref{lem:alpha}, the $\Delta'$-balancedness condition of $(a'_v)_{v\in\mathcal{V}}$ for the two triangles of $\Delta'$ appearing in Fig.\ref{fig:mutations_for_a_flip} which are involved in the flip translates to the condition that the fifteen numbers $b_1,b_2,\ldots,b_{15} \in \frac{1}{3}\mathbb{Z}$ in eq.\eqref{eq:b_i} are integers. Below we show that the integrality of the numbers in eq.\eqref{eq:b_i} implies that of the numbers in eq.\eqref{eq:fifteen_numbers_Delta}:
\begin{align*}
a_1+a_3+a_5 & = a'_1+(-a'_7+a'_6+a'_8)+a'_5 = (a'_1+a'_8+a'_{12}) - (a'_7 + a'_{12} ) + (a'_5 + a'_6) \\
& = b_1 - b_7 + b_9 ~\in~\mathbb{Z}, \\
a_2 + a_4 + a_6 & = a'_2 + (a'_5+a'_9-a'_{12}) + a'_6
= (a'_2+a'_9+a'_7) - (a'_7+a'_{12}) + (a'_5+a'_6) \\
& = b_2 - b_7 + b_9 ~\in~\mathbb{Z}, \\
a_8+a_{10}+a_4 & = a'_8+a'_{10}+(a'_5+a'_9-a'_{12})
= (a'_5+a'_7+a'_{10}) - (a'_7+a'_{12}) + (a'_8+a'_9) \\
& = b_3 - b_7 + b_6 ~\in~\mathbb{Z}, \\
a_9+a_{11}+a_3 & = a'_9 + a'_{11} + (-a'_7+a'_6+a'_8)  = (a'_6+a'_{12}+a'_{11}) - (a'_7 + a'_{12}) + (a'_8+a'_9) \\
& = b_4 - b_7 + b_6 ~\in~\mathbb{Z}, \\
a_1+a_2 & = a'_1 + a'_2 = b_5 ~\in~\mathbb{Z}, \\
a_3+a_4 & = (-a'_7+a'_6+a'_8) + (a'_5+a'_9-a'_{12}) = -(a'_7+a'_{12}) + (a'_5+a'_6) + (a'_8+a'_9) \\
& = - b_7 + b_9 + b_6 ~\in~\mathbb{Z}, \\
a_5+a_6 & = a'_5+a'_6 = b_9 ~\in~\mathbb{Z}, \\
a_8+a_9 & = a'_8+a'_9 = b_6 ~\in~\mathbb{Z}, \\
a_{10}+a_{11} & = a'_{10}+a'_{11} = b_8 ~\in~\mathbb{Z}, \\
-a_7 + a_6 + a_1 & = -(-a'_7+a'_3+a'_6) + a'_6 + a'_1 = b_{12} ~\in~\mathbb{Z}, \\
-a_7+a_2+a_3 & = -(-a'_7+a'_3+a'_6)+a'_2 + (-a'_7+a'_6+a'_8) = b_{10} ~\in~\mathbb{Z}, \\
-a_7+a_4+a_5 & = -(-a'_7+a'_3+a'_6) + (a'_5+a'_9-a'_{12}) + a'_5 \\
& = 3a'_5 - 3a'_{12} 
+( - a'_3+ a'_9 +a'_{12}) 
-(a'_5+a'_6)  + (a'_7+a'_{12}) \\
& = 3a'_5-3a'_{12} +b_{11} - b_9 + b_7 ~\in~\mathbb{Z}, \\
-a_{12}+a_3+a_8 & = -(-a'_{12}+a'_4+a'_9) + (-a'_7+a'_6+a'_8) + a'_8 \\
& = 3 a'_8 - 3a'_7
+(-a'_4+a'_6+a'_7) +(a'_7+ a'_{12}) - (a'_8+a'_9)  \\
& = 3a'_8 - 3a'_7 + b_{13} + b_7 - b_6 ~\in~\mathbb{Z}, \\
-a_{12}+a_9+a_{10} & = -(-a'_{12}+a'_4+a'_9) + a'_9 + a'_{10} = - a'_4 + a'_{12}  + a'_{10} = b_{14} ~ \in ~ \mathbb{Z}, \\
-a_{12} + a_{11}+a_4 & = -(-a'_{12}+a'_4+a'_9) + a'_{11} + (a'_5+a'_9-a'_{12}) = - a'_4  + a'_{11} + a'_5 = b_{15} ~ \in ~ \mathbb{Z}.
\end{align*}
This is sufficient for us to deduce that $(a_v)_{v\in\mathcal{V}}$ is $\Delta$-balanced, because the $\Delta$-balancedness condition of $(a_v)_{v\in \mathcal{V}}$ for triangles other than the two triangles of $\Delta$ in Fig.\ref{fig:mutations_for_a_flip} coincides with the $\Delta'$-balancedness condition of $(a_v')_{v\in \mathcal{V}'}$ for triangles other than the two triangles of $\Delta'$ in Fig.\ref{fig:mutations_for_a_flip}. This finishes the proof of $\Theta^\omega_{\Delta\Delta'}({\bf U}) \in \wh{\rm Frac}(\mathcal{Z}^\omega_\Delta)$. \qed
}

\vs

The following is easily observed from Lem.\ref{lem:nu_omega_k_extending}.
\begin{lemma}
\label{lem:Theta_extends_Phi}
Above map $\Theta^\omega_{\Delta\Delta'}$ extends $\Phi^q_{\Delta\Delta'} : {\rm Frac}(\mathcal{X}^q_{\Delta'}) \to {\rm Frac}(\mathcal{X}^q_\Delta)$ of Def.\ref{def:Phi_q_i}. \qed
\end{lemma}

Because of the similarity of the formulas of $\nu^\omega_k$ and those of $\mu^q_k$, the proofs in \cite{BZ} {\cite[\S3.3]{FG09a}} \cite{KN} \cite{Kim21} of the consistency relations for the quantum mutations $\mu^q_k$ (as in Prop.\ref{prop:quantum_mutation_relations} and Prop.\ref{prop:quantum_consistency_of_flips}) apply almost verbatim to the proof of those for the corresponding balanced quantum mutations $\nu^\omega_k$.
\begin{proposition}
\label{prop:consistency_for_Theta}
The consistency relations satisfied by the quantum mutations $\mu^q_k$ and the quantum coordinate change maps $\Phi^q_i$ for flips of ideal triangulations are satisfied by the balanced counterparts $\nu^\omega_k$ and $\Theta^\omega_i$, whenever the relations make sense. \qed
\end{proposition}
We elaborate on the phrase `whenever the relations make sense', which is about the domains and codomains of the maps $\nu^\omega_k$ involved in the relations. For example,
\begin{align}
\label{eq:twice_flip_balanced}
\nu^\omega_k \nu^\omega_k = {\rm id}
\end{align}
holds, when understood on appropriate domains and codomains; namely the composition
\begin{align}
\label{eq:twice_flip_balanced_with_domains}
\xymatrix{
\wh{\rm Frac}_k (\mathcal{Z}^\omega_\Gamma) \ar[r]^-{\nu^\omega_k} & \wh{\rm Frac}_k (\mathcal{Z}^\omega_{\mu_k(\Gamma)}) \ar[r]^-{\nu^\omega_k} & \wh{\rm Frac}_k(\mathcal{Z}^\omega_\Gamma)
}
\end{align}
coincides with the identity map on $\wh{\rm Frac}_k(\mathcal{Z}^\omega_\Gamma)$, for any cluster $\mathscr{X}$-seed $\Gamma$ and any unfrozen node $k$ of $\Gamma$. 

\vs

As a consequence of Prop.\ref{prop:consistency_for_Theta}, one obtains a balanced quantum coordinate change map $\Theta^\omega_{\Delta\Delta'} : \wh{\rm Frac}(\mathcal{Z}^\omega_{\Delta'}) \to \wh{\rm Frac}(\mathcal{Z}^\omega_\Delta)$ for each change of ideal triangulations $\Delta\leadsto \Delta'$, i.e. for each pair of ideal triangulations $\Delta$ and $\Delta'$, not just for flips. Namely, since any two ideal triangulations $\Delta$ and $\Delta'$ are connected by a finite sequence of flips, one can find a sequence of ideal triangulations $\Delta=\Delta_0$, $\Delta_1$, \ldots, $\Delta_n = \Delta'$, so that each $\Delta_{i-1} \leadsto \Delta_i$ is a flip. Then define
\begin{align}
\label{eq:Theta_omega_general}
\Theta^\omega_{\Delta \Delta'} = \Theta^\omega_{\Delta_0 \Delta_n} := \Theta^\omega_{\Delta_0 \Delta_1} \circ \Theta^\omega_{\Delta_1 \Delta_2} \circ \cdots \circ \Theta^\omega_{\Delta_{n-1} \Delta_n}.
\end{align}
Prop.\ref{prop:consistency_for_Theta} and Prop.\ref{prop:completeness_of_classical_flip_relations_for_triangulations} guarantee that the resulting map depends, up to (index) permutation of nodes of the 3-triangulation quivers for ideal triangulations (Def.\ref{def:3-triangulation_quiver}), only on the initial and the terminal triangulations $\Delta,\Delta'$ and not on the choice of the decomposition $\Delta_0,\Delta_1,\ldots,\Delta_n$.

\vs

Another important property of the balanced quantum coordinate change map $\Theta^\omega_{\Delta\Delta'}$ is its compatibility with the cutting map $i_{\Delta,\Delta_e}$ which appeared in Def.\ref{def:cutting_process}.
\begin{proposition}[compatibility of the balanced quantum coordinate change under cutting]
\label{prop:compatibility_of_Theta_under_cutting}
Let $\Delta$ be an ideal triangulation of a triangulable generalized marked surface $\frak{S}$, and $e$ an internal arc of $\Delta$. Let $\frak{S}_e$ be the generalized marked surface obtained from $\frak{S}$ be cutting along $e$, and let $\Delta_e$ be the triangulation of $\frak{S}_e$ obtained from $\Delta$ by cutting along $e$ (Def.\ref{def:cutting_process}). Denote by $i_{\Delta,\Delta_e} : {\rm Frac}(\mathcal{Z}^\omega_\Delta) \to {\rm Frac}(\mathcal{Z}^\omega_{\Delta_e})$ the natural extension of the cutting map $i_{\Delta,\Delta_e} : \mathcal{Z}^\omega_\Delta \to \mathcal{Z}^\omega_{\Delta_e}$ (Def.\ref{def:cutting_process}), defined by the map sending the fraction expression ${\bf  P}{\bf Q}^{-1} \in {\rm Frac}(\mathcal{Z}^\omega_\Delta)$ to $(i_{\Delta,\Delta_e}({\bf P})) (i_{\Delta,\Delta_e}({\bf Q}))^{-1} \in {\rm Frac}(\mathcal{Z}^\omega_{\Delta_e})$, where ${\bf P},{\bf Q} \in \mathcal{Z}^\omega_\Delta$ with ${\bf Q}\neq 0$.
\begin{enumerate}
\label{prop:Theta_compatibility_with_cutting}
\item[\rm (1)] The image of $\wh{\rm Frac}(\mathcal{Z}^\omega_\Delta) \subset {\rm Frac}(\mathcal{Z}^\omega_\Delta)$ (Def.\ref{def:balanced_fraction_algebra}) under the (extended) cutting map $i_{\Delta,\Delta_e} : {\rm Frac}(\mathcal{Z}^\omega_\Delta) \to {\rm Frac}(\mathcal{Z}^\omega_{\Delta_e})$ described above lies in $\wh{\rm Frac}(\mathcal{Z}^\omega_{\Delta_e}) \subset {\rm Frac}(\mathcal{Z}^\omega_\Delta)$.

\item[\rm (2)] Suppose that $\Delta'$ is another ideal triangulation of $\frak{S}$ such that $e$ is an internal arc of $\Delta'$ and that $\Delta$ and $\Delta'$ are connected by a finite sequence of flips at arcs other than $e$. Let $\Delta'_e$ be the ideal triangulation of $\frak{S}_e$ obtained from $\Delta'$ by cutting along $e$. Then the diagram
$$
\xymatrix@C+5mm{
\wh{\rm Frac}(\mathcal{Z}^\omega_{\Delta'}) \ar[r]^-{\Theta^\omega_{\Delta\Delta'}} \ar[d]_{i_{\Delta',\Delta'_e}} & \wh{\rm Frac}(\mathcal{Z}^\omega_\Delta) \ar[d]^{i_{\Delta,\Delta_e}} \\
\wh{\rm Frac}(\mathcal{Z}^\omega_{\Delta'_e}) \ar[r]^-{\Theta^\omega_{\Delta_e \, \Delta'_e}} & \wh{\rm Frac}(\mathcal{Z}^\omega_{\Delta_e})
}
$$
commutes, i.e.
\begin{align}
\label{eq:balanced_mutation_commutes_with_cutting}
i_{\Delta,\Delta_e} \circ \Theta^\omega_{\Delta\Delta'} = \Theta^\omega_{\Delta_e \, \Delta'_e} \circ i_{\Delta',\Delta'_e} ~:~ \wh{\rm Frac}(\mathcal{Z}^\omega_{\Delta'}) \to \wh{\rm Frac}(\mathcal{Z}^\omega_\Delta)
\end{align}

\end{enumerate}
\end{proposition}

{\it Proof.} (1) Let $\mathcal{V}=\mathcal{V}(Q_\Delta)$ and $\mathcal{V}_e = \mathcal{V}(Q_{\Delta_e})$. Consider a Laurent monomial $[\prod_{v\in\mathcal{V}}{\bf X}_v^{a_v}]_{\rm Weyl} \in \mathcal{Z}^\omega_\Delta$ for $\Delta$, with $(a_v)_{v\in\mathcal{V}} \in (\frac{1}{3}\mathbb{Z})^\mathcal{V}$. In view of the definition of the map $i_{\Delta,\Delta_e} : \mathcal{Z}^\omega_\Delta \to \mathcal{Z}^\omega_{\Delta_e}$ as in Def.\ref{def:cutting_process}, one can observe that the image of this Laurent monomial for $\Delta$ is a Laurent monomial for $\Delta_e$, i.e. $i_{\Delta,\Delta_e}[\prod_{v\in\mathcal{V}} {\bf X}_v^{a_v}]_{\rm Weyl} = [\prod_{w \in\mathcal{V}_e} {\bf X}_w^{b_w}]_{\rm Weyl}$ for some $(b_w)_{w\in\mathcal{V}_e}\in(\frac{1}{3}\mathbb{Z})^{\mathcal{V}_e}$, and that $(a_v)_{v\in\mathcal{V}} \in \mathbb{Z}^\mathcal{V}$ implies $(b_w)_{w\in\mathcal{V}_e} \in \mathbb{Z}^{\mathcal{V}_e}$. The latter observation means that $i_{\Delta,\Delta_e}$ sends $\mathcal{X}^q_\Delta \subset \mathcal{Z}^\omega_\Delta$ to $\mathcal{X}^q_{\Delta_e} \subset \mathcal{Z}^\omega_{\Delta_e}$. So, in order to show that $i_{\Delta,\Delta_e} : {\rm Frac}(\mathcal{Z}^\omega_\Delta) \to {\rm Frac}(\mathcal{Z}^\omega_{\Delta_e})$ sends the $\Delta$-balanced fraction algebra $\wh{\rm Frac}(\mathcal{Z}^\omega_\Delta)$ (Def.\ref{def:balanced_fraction_algebra}) to the $\Delta_e$-balanced fraction algebra $\wh{\rm Frac}(\mathcal{Z}^\omega_{\Delta_e})$, it suffices to show that $i_{\Delta,\Delta_e}$ sends a $\Delta$-balanced Laurent monomial to a $\Delta_e$-balanced Laurent monomial.

\vs

Notice from the definition of the balancedness condition as in Def.\ref{def:Delta-balanced_elements} and Def.\ref{def:balanced_subalgebras} that the balancedness condition of a Laurent monomial is described for each ideal triangle, in terms of the powers in the Laurent monomial of the variables for the nodes living in this triangle. From Def.\ref{def:cutting_process}, observe that triangles of $\Delta$ naturally correspond to triangles of $\Delta_e$, and that the power $a_v$ in a Laurent monomial $[\prod_{v\in\mathcal{V}}{\bf X}_v^{a_v}]_{\rm Weyl}$ for $\Delta$ for a node $v$ of $Q_\Delta$ living in some triangle $t$ coincides with the power $b_w$ in the image Laurent monomial $i_{\Delta,\Delta_e} [\prod_{v\in\mathcal{V}} {\bf X}_v^{a_v}]_{\rm Weyl} = [\prod_{w\in \mathcal{V}_e} {\bf X}_w^{b_w}]_{\rm Weyl}$ for $\Delta_e$ for each node $w$ of $Q_{\Delta_e}$ corresponding to $v$ in the cutting process (i.e. $g_e(w)=v$ in the notation of Def.\ref{def:cutting_process}); in particular, $w$ lives in the triangle $t_e$ of $\Delta_e$ corresponding to $t$ of $\Delta$. From this observation one can deduce that the $\Delta$-balancedness of $[\prod_{v\in\mathcal{V}} {\bf X}_v^{a_v}]_{\rm Weyl}$ implies the $\Delta_e$-balancedness of $[\prod_{w\in \mathcal{V}_e} {\bf X}_w^{b_w}]_{\rm Weyl}$, as desired.

\vs

(2) It suffices to show this for a flip (by Lem.\ref{lem:flip_connectedness} and Prop.\ref{prop:consistency_for_Theta}), and it suffices to show the compatibility of each balanced quantum mutation $\nu^\omega_{v^{(r)}}$, $r=1,2,3,4$. For convenience, write the sets of nodes of $Q_{\Delta^{(r)}}$ by $\mathcal{V}$, and those of $Q_{\Delta_e^{(r)}}$ by $\mathcal{V}_e$.

\vs

For each $r=0,1,2,3,4$, denote the quantum cluster $\mathscr{X}$-variables for $\Delta^{(r)}_e$ by ${\bf Y}^{(r)}_w$ per each $w\in \mathcal{V}_e$. Let $v_1,v_2$ be the nodes of $\mathcal{V}$ lying in $e$. Let $w_1,w_1'$ be the nodes of $\mathcal{V}_e$ corresponding to $v_1$, and $w_2,w_2' \in \mathcal{V}_e$ corresponding to $v_2$, so that $w_1$ and $w_2$ lie in the same boundary arc of $\frak{S}_e$, while $w_1'$, $w_2'$ lie on another boundary arc of $\frak{S}_e$. Apart from these nodes, there is a natural bijection between the nodes of 
$\mathcal{V}$ and those of $\mathcal{V}_e$. Let $g: \mathcal{V}_e \to \mathcal{V}$ be the gluing map of the nodes (see Def.\ref{def:cutting_process}, where it is denoted by $g_e$), i.e. $g(w_1) = g(w_1') = v_1$, $g(w_2)=g(w_2') = v_2$, $g(v)=v$, $\forall v\in\mathcal{V}_e\setminus \{w_1,w_1',w_2,w_2'\} = \mathcal{V}\setminus\{v_1,v_2\}$.

\vs

For each $r=4,3,2,1$, let's show
\begin{align}
\label{eq:compatibility_of_cutting_and_balanced_mutation}
i_{\Delta^{(r-1)},\Delta^{(r-1)}_e} \circ \nu^\omega_{v^{(r)}} = \nu^\omega_{v^{(r)}} \circ i_{\Delta^{(r)},\Delta^{(r)}_e} ~:~ \wh{\rm Frac}_{v^{(r)}}(\mathcal{Z}^\omega_{\Delta^{(r)}}) \to \wh{\rm Frac}_{v^{(r)}}(\mathcal{Z}^\omega_{\Delta^{(r-1)}_e})
\end{align}
We will first show this equality when applied to a $(\Delta^{(r)},v^{(r)})$-balanced Laurent monomial $[\prod_{v\in \mathcal{V}} ({\bf X}_v^{(r)})^{a^{(r)}_v}]_{\rm Weyl} \in \wh{\rm Frac}_{v^{(r)}} (\mathcal{Z}^\omega_{\Delta^{(r)}})$, with $(a_v)_{v\in\mathcal{V}} \in (\frac{1}{3}\mathbb{Z})^\mathcal{V}$. That is, we aim to show that
\begin{align}
\label{eq:compatibility_of_cutting_and_balanced_mutation_for_Laurent_monomials}
i_{\Delta^{(r-1)},\Delta^{(r-1)}_e} (\nu^\omega_{v^{(r)}} ([{\textstyle \prod}_{v\in \mathcal{V}} ({\bf X}_v^{(r)})^{a^{(r)}_v}]_{\rm Weyl})) = \nu^\omega_{v^{(r)}} ( i_{\Delta^{(r)},\Delta^{(r)}_e} ([{\textstyle \prod}_{v\in \mathcal{V}} ({\bf X}_v^{(r)})^{a^{(r)}_v}]_{\rm Weyl})).
\end{align}

Note that, in view of Def.\ref{def:cutting_process}, for each $r=4,3,2,1$ the cutting map $i_{\Delta^{(r)},\Delta^{(r)}_e}$ sends the Laurent monomial $[\prod_{v\in \mathcal{V}} ({\bf X}_v^{(r)})^{a^{(r)}_v}]_{\rm Weyl} \in \wh{\rm Frac}_{v^{(r)}} (\mathcal{Z}^\omega_{\Delta^{(r)}})$ to
\begin{align}
\label{eq:cutting_map_using_gluing_of_nodes}
i_{\Delta^{(r)},\Delta_e^{(r)}} [{\textstyle \prod}_{v\in \mathcal{V}} ({\bf X}_v^{(r)})^{a^{(r)}_v}]_{\rm Weyl}
= [{\textstyle \prod}_{w\in\mathcal{V}_e} ({\bf Y}^{(r)}_w)^{a^{(r)}_{g(w)}}]_{\rm Weyl}.
\end{align}

\vs

Consider the left-hand side of eq.\eqref{eq:compatibility_of_cutting_and_balanced_mutation_for_Laurent_monomials} first. By eq.\eqref{eq:nu_omega_v_r_U_r}, we have
$$
\nu^\omega_{v^{(r)}}([{\textstyle \prod}_{v\in\mathcal{V}}({\bf X}_v^{(r)})^{a^{(r)}_v}]_{\rm Weyl}) = 
[{\textstyle \prod}_{v\in \mathcal{V}} ({
\bf X}_v^{(r-1)})^{a^{(r-1)}_v}]_{\rm Weyl} \, F^q({\bf X}^{(r-1)}_{v^{(r)}}; \alpha^{(r)})
$$
where $(a^{(r-1)}_v)_{v\in\mathcal{V}}$ is given by the formula in eq.\eqref{eq:a_r_v_recursive} in terms of $(a_v^{(r)})_{v\in\mathcal{V}}$, and $\alpha^{(r)}$ by eq.\eqref{eq:alpha_for_nu_sharp}; in turn, in view of Def.\ref{def:cutting_process}, this is sent via $i_{\Delta^{(r-1)},\Delta^{(r-1)}_e}$ to
$$
i_{\Delta^{(r-1)},\Delta_e^{(r-1)}} (\nu^\omega_{v^{(r)}}([{\textstyle \prod}_{v\in\mathcal{V}}({\bf X}_v^{(r)})^{a^{(r)}_v}]_{\rm Weyl})) = 
[{\textstyle \prod}_{w\in \mathcal{V}_e} ({\bf Y}_w^{(r-1)})^{a^{(r-1)}_{g(w)}}]_{\rm Weyl} \, F^q({\bf Y}^{(r-1)}_{v^{(r)}};\alpha^{(r)}).
$$

\vs

For the right-hand side of eq.\eqref{eq:compatibility_of_cutting_and_balanced_mutation_for_Laurent_monomials}, first, note that $[\prod_{v\in \mathcal{V}} ({\bf X}_v^{(r)})^{a^{(r)}_v}]_{\rm Weyl} \in \wh{\rm Frac}_{v^{(r)}} (\mathcal{Z}^\omega_{\Delta^{(r)}})$ is sent via $i_{\Delta^{(r)},\Delta_e^{(r)}}$ to the element $[\prod_{w\in \mathcal{V}_e} ({\bf Y}_w^{(r)})^{a^{(r)}_{g(w)}}]_{\rm Weyl}$ as written in eq.\eqref{eq:cutting_map_using_gluing_of_nodes}. In the meantime, the signed adjacency matrices $\varepsilon^{(r)}_e$ for $Q_{\Delta^{(r)}_e}$ and $\varepsilon^{(r)}$ for $Q_{\Delta^{(r)}}$ are related by
\begin{align}
\label{eq:varepsilon_gluing}
\varepsilon^{(r)}_{vu} = \underset{
v' \in g^{-1}(v), \, u' \in g^{-1}(u)
}{\textstyle \sum} (\varepsilon^{(r)}_e)_{v'u'}, \qquad \forall v,u\in \mathcal{V}.
\end{align}
One observation is that, no cancellation occurs in the right-hand side of eq.\eqref{eq:varepsilon_gluing} except for two possible cases. One is when $\{v,u\} = \{v_1,v_2\}$ holds. In the other case, one of $v$ and $u$ is a node lying in the edge $e$, while the other is a node lying in the edge $f$ which is the edge being flipped in the current proof; moreover, the common endpoint puncture $p$ shared by $e$ and $f$ is of valence $2$, and each of $v$ and $u$ is the node closer to $p$ among the two nodes lying in the respective edges $e$ and $f$. However, the latter case does not happen, because after flipping at an edge incident to a puncture of valence $2$ one obtains a puncture of valence 1, which we excluded in \S\ref{subsec:surfaces_and_triangulations}. Therefore, the only case when a cancellation might occur in the right-hand side of eq.\eqref{eq:varepsilon_gluing} is when $\{v,u\} = \{v_1,v_2\}$. Thus, whenever $\{v,u\} \neq \{v_1,v_2\}$, the terms in the right-hand side are either all non-negative, or all non-positive, and hence we have
\begin{align}
\label{eq:varepsilon_sign_coherence}
[\varepsilon^{(r)}_{v u}]_+ = \underset{v' \in g^{-1}(v), \, u' \in g^{-1}(u)}{\textstyle \sum} [(\varepsilon^{(r)}_e)_{v'u'}]_+, \qquad \forall v,u\in \mathcal{V}.
\end{align}

Now, observe from eq.\eqref{eq:nu_omega_v_r_U_r} applied to the case $\Theta^\omega_{\Delta_e \Delta_e'}$ of the cut surfaces (instead of original $\Theta^\omega_{\Delta \Delta'}$) that
\begin{align*}
\nu^\omega_{v^{(r)}} [{\textstyle \prod}_{v\in \mathcal{V}_e} ({\bf Y}_v^{(r)})^{a^{(r)}_{g(v)}}]_{\rm Weyl}
= [{\textstyle \prod}_{v\in \mathcal{V}_e} ({\bf Y}_v^{(r-1)})^{b^{(r-1)}_v}]_{\rm Weyl} F^q({\bf Y}^{(r-1)}_{v^{(r)}};\alpha_e^{(r)}),
\end{align*}
where $(b^{(r-1)}_v)_{v\in\mathcal{V}_e}$ is given similarly as in eq.\eqref{eq:a_r_v_recursive} by
\begin{align}
\label{eq:b_as_a}
b^{(r-1)}_v = \left\{
\begin{array}{ll}
-a^{(r)}_{v^{(r)}} + \sum_{w\in \mathcal{V}_e} [(\varepsilon^{(r-1)}_e)_{w,v^{(r)}}]_+ \, a^{(r)}_{g(w)}, & \mbox{if $v=v^{(r)}$,} \\
a^{(r)}_{g(v)} & \mbox{if $v\neq v^{(r)}$},
\end{array}
\right.
\end{align}
and $\alpha^{(r)}_e$ is given by
$$
\alpha^{(r)}_e = {\textstyle \sum}_{w\in\mathcal{V}_e} (\varepsilon^{(r-1)}_e)_{v^{(r)},w} a^{(r-1)}_{g(w)}.
$$

We claim that
$$
b^{(r-1)}_v= a^{(r-1)}_{g(v)}.
$$
From eq.\eqref{eq:b_as_a} and eq.\eqref{eq:a_r_v_recursive}, this is clear for the case $v\neq v^{(r)}$. For the case $v=v^{(r)}$ this follows from eq.\eqref{eq:varepsilon_gluing} and eq.\eqref{eq:varepsilon_sign_coherence}. 
Moreover, observe from eq.\eqref{eq:varepsilon_gluing} that $\alpha_e^{(r)}$ equals $\alpha^{(r)}$ of eq.\eqref{eq:alpha_for_nu_sharp}. Hence we showed that the left-hand side and the right-hand side of eq.\eqref{eq:compatibility_of_cutting_and_balanced_mutation_for_Laurent_monomials} are equal, as desired. This holds for any $(a_v^{(r)})_{v\in\mathcal{V}} \in (\frac{1}{3}\mathbb{Z})^{\mathcal{V}}$ that is $(\Delta^{(r)},v^{(r)})$-balanced, and in particular when $(a_v^{(r)})_{v\in\mathcal{V}} \in \mathbb{Z}^\mathcal{V}$. So, in view of the definitions of the $(\Delta^{(r)},v^{(r)})$-balanced fraction algebra $\wh{\rm Frac}_{v^{(r)}}(\mathcal{Z}^\omega_{\Delta^{(r)}})$ and that of the $(\Delta^{(r-1)},v^{(r)})$-balanced fraction algebra $\wh{\rm Frac}_{v^{(r)}}(\mathcal{Z}^\omega_{\Delta^{(r-1)}})$ as in Def.\ref{def:Gamma_k_balanced}, it follows that eq.\eqref{eq:compatibility_of_cutting_and_balanced_mutation} holds. The sought-for eq.\eqref{eq:balanced_mutation_commutes_with_cutting} then follows by combining eq.\eqref{eq:compatibility_of_cutting_and_balanced_mutation} for $r=4,3,2,1$, since each of $\Theta^\omega_{\Delta\Delta'}$ and $\Theta^\omega_{\Delta\Delta'}$ is given by the composition $\nu^\omega_{v^{(1)}}\nu^\omega_{v^{(2)}}\nu^\omega_{v^{(3)}}\nu^\omega_{v^{(4)}}$, in an appropriate sense. This finishes the proof. \qed

\section{Naturality of ${\rm SL}_3$ quantum trace maps
under changes of triangulations}

Here comes the main statement of the present paper.
\begin{theorem}[main theorem; naturality of ${\rm SL}_3$ quantum trace maps under changes of triangulations]
\label{thm:main}
Let $\Delta$ and $\Delta'$ be ideal triangulations (without self-folded triangles) of a triangulable generalized marked surface $\frak{S}$.  Let $(W,s)$ be a stated ${\rm SL}_3$-web in $\frak{S} \times {\bf I}$. Then the values under the ${\rm SL}_3$ quantum trace maps of $[W,s] \in \mathcal{S}^\omega_{\rm s}(\frak{S};\mathbb{Z})_{\rm red}$ are compatible under the balanced quantum coordinate change map $\Theta^\omega_{\Delta\Delta'}$ defined in Def.\ref{def:Theta} and eq.\eqref{eq:Theta_omega_general}, i.e.
$$
\Theta^\omega_{\Delta\Delta'}( {\rm Tr}^\omega_{\Delta'}([W,s])) = {\rm Tr}^\omega_\Delta([W,s]).
$$
\end{theorem}
We note that the left-hand side makes sense because ${\rm Tr}^\omega_{\Delta'}([W,s])$ lies in $\wh{\mathcal{Z}}^\omega_{\Delta'}$ due to Prop.\ref{prop:value_of_the_SL3_quantum_trace_is_Delta-balanced}, and hence in $\wh{\rm Frac}(\mathcal{Z}^\omega_{\Delta'})$ (Def.\ref{def:balanced_fraction_algebra}) which is the domain of $\Theta^\omega_{\Delta\Delta'}$ as seen in Def.\ref{def:Theta} and eq.\eqref{eq:Theta_omega_general}.

\vs

The present section is devoted to the proof of this theorem. In view of Lem.\ref{lem:flip_connectedness} and Prop.\ref{prop:consistency_for_Theta}, it suffices to prove this in the case when $\Delta$ and $\Delta'$ are related by a flip at an arc.

\subsection{The base case: crossingless arcs over an ideal quadrilateral}
\label{subsec:the_base_case}

Let $e$ be the internal arc of $\Delta$ that is being flipped, i.e. the only arc of $\Delta$ that is not an arc of $\Delta'$. Let the two triangles of $\Delta$ having $e$ as a side be $t$ and $u$. Collect all the sides of $t$ and $u$ that are not $e$ or boundary arcs of $\frak{S}$, and let this collection be $E$. So the number of members of $E$ is one of $0$, $1$, $2$, $3$ and $4$. The surface obtained from $\frak{S}$ by cutting along all arcs of $E$ is a disjoint union of a quadrilateral $\mathcal{Q}$ and a surface $\frak{S}_E$. Note that $\frak{S}_E$ may also be disconnected, and even be empty in case $\frak{S}$ was already a quadrilateral. Here, by a quadrilateral we mean a generalized marked surface of genus zero with one boundary component, with four marked points on the boundary and no marked point in the interior. By the cutting process (Def.\ref{def:cutting_process}), $\Delta$ yields  triangulations $\Delta_\mathcal{Q}$ and $\Delta_E$ for $\mathcal{Q}$ and $\frak{S}_E$ respectively, while $\Delta'$ yields $\Delta'_{\mathcal{Q}}$ and $\Delta'_E$ for $\mathcal{Q}$ and $\frak{S}_E$. Note that $\Delta_E$ and $\Delta'_E$ coincide with each other, while $\Delta_\mathcal{Q}$ and $\Delta'_\mathcal{Q}$ differ by a flip. Suppose that an ${\rm SL}_3$-web $(W,s)$ in $\frak{S}\times {\bf I}$ meets $E\times {\bf I}$ transversally, and moreover that $W\cap (\mathcal{Q}\times {\bf I})$ and $W\cap (\frak{S}_E\times {\bf I})$ are ${\rm SL}_3$-webs in $\mathcal{Q} \times {\bf I}$ and in $\frak{S}_E \times {\bf I}$ respectively; if not, $(W,s)$ can be isotoped in $\frak{S}\times {\bf I}$ to satisfy this condition. By the cutting/gluing property of the ${\rm SL}_3$ quantum trace maps (Thm.\ref{thm:SL3_quantum_trace}(QT1)), one has
\begin{align*}
i_E ( {\rm Tr}^\omega_{\Delta;\frak{S}}([W, s]) ) & = 
\underset{s_\mathcal{Q}, s_E}{\textstyle \sum} {\rm Tr}^\omega_{\Delta_\mathcal{Q};\mathcal{Q}} ([W\cap (\mathcal{Q} \times {\bf I}),s_\mathcal{Q}]) 
\otimes {\rm Tr}^\omega_{\Delta_E;\frak{S}_E} ([W\cap (\frak{S}_E \times {\bf I}), s_E]), \\
i'_E ( {\rm Tr}^\omega_{\Delta'; \frak{S}}([W, s]) ) & = 
\underset{s_\mathcal{Q}, s_E}{\textstyle \sum} {\rm Tr}^\omega_{\Delta_\mathcal{Q}' ; \mathcal{Q}} ([W\cap (\mathcal{Q} \times {\bf I}),s_\mathcal{Q}]) 
\otimes {\rm Tr}^\omega_{\Delta_E;\frak{S}_E} ([W\cap (\frak{S}_E \times {\bf I}), s_E]),
\end{align*}
where $i_E$ and $i'_E$ are natural embeddings (see Def.\ref{def:cutting_process})
\begin{align*}
i_E : \mathcal{Z}^\omega_\Delta \to \mathcal{Z}^\omega_{\Delta_\mathcal{Q} \sqcup \Delta_E} \cong \mathcal{Z}^\omega_{\Delta_\mathcal{Q}} \otimes \mathcal{Z}^\omega_{\Delta_E}, \\
i'_E : \mathcal{Z}^\omega_{\Delta'} \to \mathcal{Z}^\omega_{\Delta'_\mathcal{Q} \sqcup \Delta_E} \cong \mathcal{Z}^\omega_{\Delta'_\mathcal{Q}} \otimes \mathcal{Z}^\omega_{\Delta_E}, 
\end{align*}
and the sums are over all states $s_\mathcal{Q}$ and $s_E$ of the ${\rm SL}_3$-webs $W\cap (\mathcal{Q}\times {\bf I})$ in $\mathcal{Q} \times {\bf I}$ and $W\cap (\frak{S}_E\times {\bf I})$ in $\frak{S}_E \times {\bf I}$ that constitute states $s_\mathcal{Q} \sqcup s_E$ that are compatible with $s$ in the sense of Def.\ref{def:cutting_process}. By Prop.\ref{prop:compatibility_of_Theta_under_cutting}, which is the compatibility of the quantum mutation maps $\Theta^\omega_{\Delta,\Delta'}$ under the cutting maps along ideal arcs, one observes that
$$
i_E \circ \Theta^\omega_{\Delta,\Delta'} = (\Theta^\omega_{\Delta_\mathcal{Q},\Delta'_\mathcal{Q}} \otimes \underset{={\rm id}}{\myul{ \Theta^\omega_{\Delta_E,\Delta_E} }}) \circ i'_E ~:~ \wh{\mathcal{Z}}^\omega_{\Delta'} \to \wh{\mathcal{Z}}^\omega_{\Delta_\mathcal{Q}} \otimes \wh{\mathcal{Z}}^\omega_{\Delta_E}.
$$
It then suffices to show that
$$
{\rm Tr}^\omega_{\Delta_\mathcal{Q};\mathcal{Q}} ([W\cap (\mathcal{Q} \times {\bf I}),s_\mathcal{Q}]) 
= \Theta^\omega_{\Delta_\mathcal{Q},\Delta'_{\mathcal{Q}}} {\rm Tr}^\omega_{\Delta'_\mathcal{Q};\mathcal{Q}} ([W\cap (\mathcal{Q} \times {\bf I}),s_\mathcal{Q}]).
$$
Thus, we could assume from the beginning that the entire surface $\frak{S}$ is just a quadrilateral $\mathcal{Q}$.

\vs

Fatten each of the four boundary arcs of the quadrilateral $\frak{S} = \mathcal{Q}$ to a biangle, by choosing an ideal arc isotopic to each boundary arc. Let $B_1,B_2,B_3,B_4$ be these biangles. Cutting $\mathcal{Q}$ along these four ideal arcs yields a disjoint union $B_1,B_2,B_3,B_4$, and a quadrilateral, which we denote by $\mathcal{Q}_0$. Suppose that $(W,s)$ is a stated ${\rm SL}_3$-web in $\mathcal{Q}\times {\bf I}$, and suppose that $W\cap (B_i \times {\bf I})$, $i=1,2,3,4$, and $W\cap (\mathcal{Q}_0 \times {\bf I})$ are ${\rm SL}_3$-webs in $B_i \times {\bf I}$, $i=1,2,3,4$, and in $\mathcal{Q}_0 \times {\bf I}$, respectively. Let $\Delta$ and $\Delta'$ be two distinct ideal triangulations of $\mathcal{Q}$, which are related by a flip. Let $\Delta_0$ and $\Delta_0'$ be the ideal triangulation of $\mathcal{Q}_0$ induced by the cutting process (Def.\ref{def:cutting_process}). By the cutting/gluing property of the ${\rm SL}_3$ quantum trace maps (Prop.\ref{prop:biangle_quantum_trace}(BQT3)), we have
\begin{align*}
i_{\Delta,\Delta_0} {\rm Tr}^\omega_{\Delta;\mathcal{Q}} ([W,s]) & = \underset{s_1,s_2,s_3,s_4,s_0}{\textstyle \sum} ({\textstyle \prod}_{i=1}^4 {\rm Tr}^\omega_{B_i}([W\cap (B_i \times {\bf I}), s_i]) ) \, {\rm Tr}^\omega_{\Delta_0;\mathcal{Q}_0}([W\cap (\mathcal{Q}_0 \times {\bf I}),s_0]), \\
i_{\Delta',\Delta'_0} {\rm Tr}^\omega_{\Delta';\mathcal{Q}} ([W,s]) & = \underset{s_1,s_2,s_3,s_4,s_0}{\textstyle \sum} ({\textstyle \prod}_{i=1}^4 {\rm Tr}^\omega_{B_i}([W\cap (B_i \times {\bf I}), s_i]) ) \, {\rm Tr}^\omega_{\Delta_0';\mathcal{Q}_0}([W\cap (\mathcal{Q}_0 \times {\bf I}),s_0])
\end{align*}
where $i_{\Delta,\Delta_0}$ and $i_{\Delta',\Delta'_0}$ are natural isomorphisms
$$
i_{\Delta,\Delta_0} ~:~ \mathcal{Z}^\omega_\Delta \to \mathcal{Z}^\omega_{\Delta_0}, \qquad
i_{\Delta',\Delta'_0} ~:~ \mathcal{Z}^\omega_{\Delta'} \to \mathcal{Z}^\omega_{\Delta_0'},
$$
and the sums are over all states that constitute states compatible with $s$ in the sense of Def.\ref{def:cutting_process}. Observing 
$$
i_{\Delta,\Delta_0} \circ \Theta^\omega_{\Delta,\Delta'}
= \Theta^\omega_{\Delta_0,\Delta'_0} \circ i_{\Delta',\Delta'_0} ~:~ \wh{\mathcal{Z}}^\omega_{\Delta'} \to \wh{\mathcal{Z}}^\omega_{\Delta_0}
$$
from Prop.\ref{prop:compatibility_of_Theta_under_cutting}, it suffices to show that
$$
{\rm Tr}^\omega_{\Delta_0;\mathcal{Q}_0}([W\cap (\mathcal{Q}_0 \times {\bf I}),s_0])
= \Theta^\omega_{\Delta_0,\Delta_0'} {\rm Tr}^\omega_{\Delta_0';\mathcal{Q}_0}([W\cap (\mathcal{Q}_0 \times {\bf I}),s_0]).
$$
The advantage we get by this cutting process is as follows. In the beginning, we do not assume anything about the ${\rm SL}_3$-web $W$ in $\mathcal{Q}\times {\bf I}$. Then, by isotopy, one can push complicated parts of $W$, e.g. all 3-valent vertices, to the biangles $B_1,B_2,B_3,B_4$, and also apply `vertical' isotopy, so that $W$ is `nice' over $\mathcal{Q}_0$ in the following sense:
\begin{lemma}
\label{lem:isotoped_shape}
One can isotope $W$ so that $W_0 := W\cap (\mathcal{Q}_0 \times {\bf I})$ satisfies the following:
\begin{enumerate}\itemsep0em
\item[\rm (1)] $W_0$ has no crossing or $3$-valent vertex;

\item[\rm (2)] each component of $W_0$ is at a constant elevation equipped with upward vertical framing, and the elevations of the components of $W_0$ are mutually distinct;

\item[\rm (3)] each component of the part of $W_0$ lying over each of the ideal triangles of $\Delta_0$ and $\Delta_0'$ is a left-turn or a right-turn oriented edge (as in Thm.\ref{thm:SL3_quantum_trace}(QT2)). 

\item[\rm (4)] each component of $W_0$ is a constant-elevation lift in $\mathcal{Q}_0\times {\bf I}$, with upward vertical framing, of a corner arc of $\mathcal{Q}_0$, i.e. the projection $\pi(W_0)$ in $\mathcal{Q}_0$ is a simple oriented edge connecting two adjacent boundary arcs of $\mathcal{Q}_0$. \qed
\end{enumerate}
\end{lemma}
In the language of \cite{Kim}, the items (1)--(3) mean that $W_0$ can be put into a `gool' position (which is stronger than a `good' position) (\cite[Lem.5.45]{Kim}). Then there could be some component of $W_0$ whose projection in $\mathcal{Q}_0$ is a simple edge connecting two non-adjacent boundary arcs of $\mathcal{Q}_0$; one can avoid having these and get (4) by using isotopies of the form appearing in the right picture of \cite[Fig.17]{Kim} (before applying the vertical isotopy to make sure that the distinct-elevations  condition in (2) holds).

\vs

Since ${\rm Tr}^\omega_{\Delta_0;\mathcal{Q}_0}$, $\Theta^\omega_{\Delta_0,\Delta_0'}$ and ${\rm Tr}^\omega_{\Delta_0';\mathcal{Q}_0}$ all preserve the product structures, the problem boils down to the following. By dropping all the subscripts $0$ for convenience,
what remains to check is the equality 
\begin{align}
\label{eq:boiled_down}
{\rm Tr}^\omega_{\Delta;\mathcal{Q}}([W,s]) = \Theta^\omega_{\Delta \Delta'} {\rm Tr}^\omega_{\Delta';\mathcal{Q}}([W,s]),
\end{align}
when $\Delta$ and $\Delta'$ are ideal triangulations of a quadrilateral surface $\mathcal{Q}$, and $W$ satisfies the following:
$$
\mbox{$W$ is a constant-elevation lift in $\mathcal{Q}\times {\bf I}$, with upward vertical framing, of a corner arc of $\mathcal{Q}$}.
$$
Some of the possible cases are visualized in Fig.\ref{fig:four_cases}.

\vs

One can check eq.\eqref{eq:boiled_down} for these cases by direct computation, which is doable by hand. But we will try an approach that minimizes the amount of computations, which is also more enlightening. Here we perform a preliminary step. Note from Def.\ref{def:Theta} that $\Theta^\omega_{\Delta\Delta'} = \nu^\omega_{v_3} \nu^\omega_{v_4} \nu^\omega_{v_7} \nu^\omega_{v_{12}}$. Thus, by applying $\nu^\omega_{v_4} \nu^\omega_{v_3}$ to both sides of eq.\eqref{eq:boiled_down} and using the involutivity of $\nu^\omega_k$ as in eq.\eqref{eq:twice_flip_balanced}, one can turn eq.\eqref{eq:boiled_down} into
\begin{align}
\label{eq:boiled_down2}
\nu^\omega_{v_4} \nu^\omega_{v_3}( {\rm Tr}^\omega_{\Delta;\mathcal{Q}}([W,s])) = \nu^\omega_{v_7} \nu^\omega_{v_{12}}({\rm Tr}^\omega_{\Delta';\mathcal{Q}}([W,s])),
\end{align}
which is much more symmetric than eq.\eqref{eq:boiled_down}. Indeed, what is done in either side is, first to apply the ${\rm SL}_3$ quantum trace map with respect to an ideal triangulation of the quadrilateral $\mathcal{Q}$, and then apply the balanced cube-root quantum mutation maps $\nu^\omega_k$ (Def.\ref{def:nu_omega_k}) for the two nodes $k$ living in the diagonal arc of this ideal triangulation.

\vs

However, one must be careful about the domains and codomains when applying the balanced cube-root quantum mutation map $\nu^\omega_k$, as well as when applying the relation in eq.\eqref{eq:boiled_down}. So, a priori, eq.\eqref{eq:boiled_down2} and eq.\eqref{eq:boiled_down} should be understood as separate statements. The statement of eq.\eqref{eq:boiled_down2} asserts that both sides are well defined (with respect to the domains of $\nu^\omega_k$) and that they are equal to each other. The remainder of the present section will be devoted to a proof of this eq.\eqref{eq:boiled_down2}. For the purpose of our proof of Thm.\ref{thm:main}, we should also make sure that eq.\eqref{eq:boiled_down2} implies eq.\eqref{eq:boiled_down}, which we show below.

\vs

By Prop.\ref{prop:value_of_the_SL3_quantum_trace_is_Delta-balanced} applied to $\Delta'$ we have ${\rm Tr}^\omega_{\Delta';\mathcal{Q}}([W,s]) \in \wh{\mathcal{Z}}^\omega_{\Delta'} \subset \wh{\rm Frac}(\mathcal{Z}^\omega_{\Delta'})$. In the proof of Lem.\ref{lem:Theta_well-definedness1}, we showed that the four maps $\nu^\omega_{v_{12}}$, $\nu^\omega_{v_7}$, $\nu^\omega_{v_4}$ and $\nu^\omega_{v_3}$ can be applied, in this order, to any element of $\wh{\rm Frac}(\mathcal{Z}^\omega_{\Delta'})$, and that it results in an element of $\wh{\rm Frac}(\mathcal{Z}^\omega_\Delta)$. For the moment, denote $\nu^\omega_{v_3} \nu^\omega_{v_4}\nu^\omega_{v_7}\nu^\omega_{v_{12}} {\rm Tr}^\omega_{\Delta';\mathcal{Q}}[W,s]$ by ${\bf W} \in \wh{\rm Frac}(\mathcal{Z}^\omega_\Delta)$. The last step to reach ${\bf W}$ is the application of $\nu^\omega_{v_3} : \wh{\rm Frac}_{v_3}(\mathcal{Z}^\omega_{\Delta^{(1)}}) \to \wh{\rm Frac}_{v_3}(\mathcal{Z}^\omega_{\Delta})$ to the element $\nu^\omega_{v_4}\nu^\omega_{v_7}\nu^\omega_{v_{12}} {\rm Tr}^\omega_{\Delta';\mathcal{Q}}[W,s]$, so by what we checked in the proof of Lem.\ref{lem:Theta_well-definedness1} we have $\nu^\omega_{v_4}\nu^\omega_{v_7}\nu^\omega_{v_{12}} {\rm Tr}^\omega_{\Delta';\mathcal{Q}}[W,s] \in \wh{\rm Frac}_{v_3}(\mathcal{Z}^\omega_{\Delta^{(1)}})$ and ${\bf W} \in \wh{\rm Frac}_{v_3}(\mathcal{Z}^\omega_\Delta)$. Applying the involutivity of $\nu^\omega_{v_3}$ as in eq.\eqref{eq:twice_flip_balanced} or more precisely as written in eq.\eqref{eq:twice_flip_balanced_with_domains} applied to $\nu^\omega_{v_3} : \wh{\rm Frac}_{v_3}(\mathcal{Z}^\omega_{\Delta^{(1)}}) \to \wh{\rm Frac}_{v_3}(\mathcal{Z}^\omega_\Delta)$ and $\nu^\omega_{v_3} : \wh{\rm Frac}_{v_3}(\mathcal{Z}^\omega_\Delta) \to \wh{\rm Frac}_{v_3}(\mathcal{Z}^\omega_{\Delta^{(1)}})$, we get $\nu^\omega_{v_4} \nu^\omega_{v_7} \nu^\omega_{v_{12}} {\rm Tr}^\omega_{\Delta';\mathcal{Q}}[W,s] = \nu^\omega_{v_3} ({\bf W}) =: {\bf W}'$. The next to last step to reach ${\bf W}$ is the application of $\nu^\omega_{v_4} : \wh{\rm Frac}_{v_4}(\mathcal{Z}^\omega_{\Delta^{(2)}}) \to \wh{\rm Frac}_{v_4}(\mathcal{Z}^\omega_{\Delta^{(1)}})$ to the element $\nu^\omega_{v_7} \nu^\omega_{v_{12}} {\rm Tr}^\omega_{\Delta';\mathcal{Q}}[W,s]$, so by what we checked in the proof of Lem.\ref{lem:Theta_well-definedness1} we have $\nu^\omega_{v_7} \nu^\omega_{v_{12}} {\rm Tr}^\omega_{\Delta';\mathcal{Q}}[W,s] \in \wh{\rm Frac}_{v_4}(\mathcal{Z}^\omega_{\Delta^{(2)}})$ and $\nu^\omega_{v_4}( \nu^\omega_{v_7} \nu^\omega_{v_{12}} {\rm Tr}^\omega_{\Delta';\mathcal{Q}}[W,s]) = {\bf W}' \in \wh{\rm Frac}_{v_4}(\mathcal{Z}^\omega_{\Delta^{(1)}})$. Applying the involutivity of $\nu^\omega_k$ as in eq.\eqref{eq:twice_flip_balanced}--eq.\eqref{eq:twice_flip_balanced_with_domains} applied to $\nu^\omega_{v_4} : \wh{\rm Frac}_{v_4}(\mathcal{Z}^\omega_{\Delta^{(2)}}) \to \wh{\rm Frac}_{v_4}(\mathcal{Z}^\omega_{\Delta^{(1)}})$ and $\nu^\omega_{v_4} : \wh{\rm Frac}_{v_4}(\mathcal{Z}^\omega_{\Delta^{(1)}}) \to \wh{\rm Frac}_{v_4}(\mathcal{Z}^\omega_{\Delta^{(2)}})$, we get $\nu^\omega_{v_7} \nu^\omega_{v_{12}} {\rm Tr}^\omega_{\Delta';\mathcal{Q}}[W,s] = \nu^\omega_{v_4}({\bf W}')$. Thus, so far we proved that
\begin{align}
\label{eq:liberation0}
\nu^\omega_{v_7} \nu^\omega_{v_{12}} {\rm Tr}^\omega_{\Delta';\mathcal{Q}}[W,s] = \nu^\omega_{v_4} \nu^\omega_{v_3}({\bf W}).
\end{align}
As a side remark, one can now see that eq.\eqref{eq:boiled_down} implies eq.\eqref{eq:boiled_down2}; indeed, if eq.\eqref{eq:boiled_down} holds then we have ${\bf W} = {\rm Tr}^\omega_{\Delta;\mathcal{Q}}[W,s]$, hence eq.\eqref{eq:liberation0} becomes eq.\eqref{eq:boiled_down2} in this case.

\vs

However, we wanted to show that eq.\eqref{eq:boiled_down2} implies eq.\eqref{eq:boiled_down}. If we assume that eq.\eqref{eq:boiled_down2} holds, then from eq.\eqref{eq:liberation0} it follows that
\begin{align}
\label{eq:liberation1}
\nu^\omega_{v_4} \nu^\omega_{v_3} {\rm Tr}^\omega_{\Delta;\mathcal{Q}}([W,s]) = \nu^\omega_{v_4} \nu^\omega_{v_3} ({\bf W})
\end{align}
holds. Note that each $\nu^\omega_k : \wh{\rm Frac}_k(\mathcal{Z}^\omega_{\mu_k(\Gamma)}) \to \wh{\rm Frac}_k(\mathcal{Z}^\omega_\Gamma)$ is an isomorphism, with its inverse being $\nu^\omega_k : \wh{\rm Frac}_k(\mathcal{Z}^\omega_{\Gamma}) \to \wh{\rm Frac}_k(\mathcal{Z}^\omega_{\mu_k(\Gamma)})$, due to the involutivity in eq.\eqref{eq:twice_flip_balanced}--eq.\eqref{eq:twice_flip_balanced_with_domains}. So $\nu^\omega_{v_4} \nu^\omega_{v_3}$ is an isomorphism, hence from eq.\eqref{eq:liberation1} we get ${\rm Tr}^\omega_{\Delta;\mathcal{Q}}[W,s] = {\bf W}$; since ${\bf W} = \nu^\omega_{v_3} \nu^\omega_{v_4}\nu^\omega_{v_7}\nu^\omega_{v_{12}} {\rm Tr}^\omega_{\Delta';\mathcal{Q}}[W,s] = \Theta^\omega_{\Delta\Delta'} {\rm Tr}^\omega_{\Delta';\mathcal{Q}}[W,s]$, eq.\eqref{eq:boiled_down} holds. So we indeed showed that eq.\eqref{eq:boiled_down2} implies eq.\eqref{eq:boiled_down}, as desired.

\subsection{The Weyl-ordering and the classical compatibility statement}
\label{subsec:classical_compatibility_and_Weyl-ordering}

It remains to show that eq.\eqref{eq:boiled_down2} holds. Our strategy to show eq.\eqref{eq:boiled_down2} is to use Lem.\ref{lem:star-invariance_and_Weyl-ordering_new}, especially its item (3), which says that a multiplicity-free Laurent polynomial in a cube-root Fock-Goncharov algebra $\mathcal{Z}^\omega_Q$ (i.e. a (cube-root) quantum torus algebra) (Def.\ref{def:FG_algebra}), for a general quiver $Q$, is $*$-invariant if and only if it is a (term-by-term) Weyl-ordered Laurent polynomial (Def.\ref{def:Weyl-ordered}). We find it convenient to establish and use an easy observation that a Weyl-ordered Laurent polynomial is completely determined by its {\it classicalization}.

\vs

\begin{definition}
\label{def:cl_and_Wl} 
Let $Q$ be a quiver, and denote by $\mathcal{V} = \mathcal{V}(Q)$ the set of its nodes, and by $\varepsilon=(\varepsilon_{vw})_{v,w\in\mathcal{V}}$ its signed adjacency matrix. Let $\mathcal{Z}^\omega_Q$ be the cube-root Fock-Goncharov algebra for $Q$, as defined in Def.\ref{def:FG_algebra}, with generators ${\bf Z}_v^{\pm 1} = {\bf X}_v^{\pm 1/3}$, $v\in \mathcal{V}$. Denote by $\mathcal{Z}^1_Q$ the cube-root Fock-Goncharov algebra for $Q$ with the quantum parameter $\omega^{1/2}$ set to be $1$, where the generators are written as $Z_v^{\pm 1} = X_v^{\pm 1/3}$, $v\in \mathcal{V}$; that is, $\mathcal{Z}^1_\Delta$ can be viewed just as the Laurent polynomial ring $\mathbb{Z}[\{Z_v^{\pm 1} \, | \, v\in \mathcal{V}\}]$.

\vs

The \ul{\em classicalization} map
$$
{\rm cl}^\omega_Q ~:~ \mathcal{Z}^\omega_Q \to \mathcal{Z}^1_Q
$$
is defined as the unique ring homomorphism sending $\omega^{\pm 1/2}$ to $1$ and each ${\bf Z}_v^{\pm 1}$ to $Z_v^{\pm 1}$. Define the \\ \ul{\em Weyl-ordering quantization map}
$$
{\rm Wl}^\omega_Q ~:~ \mathcal{Z}^1_Q \to \mathcal{Z}^\omega_Q
$$
as the unique $\mathbb{Z}$-linear map sending each Laurent monomial to its corresponding Weyl-ordered Laurent monomial (Def.\ref{def:Weyl-ordered})
$$
{\rm Wl}^\omega_Q ( {\textstyle \prod}_{v\in\mathcal{V}} X_v^{a_v})  = [ {\textstyle \prod}_{v\in\mathcal{V}} {\bf X}_v^{a_v} ]_{\rm Weyl}, \qquad \forall (a_v)_{v\in\mathcal{V}} \in ({\textstyle \frac{1}{3}}\mathbb{Z})^{\mathcal{V}}.
$$

\vs

When $\Delta$ is an ideal triangulation of a triangulable generalized marked surface $\frak{S}$, we denote $\mathcal{Z}^1_{Q_\Delta}$, ${\rm cl}^\omega_{Q_\Delta}$ and ${\rm Wl}^\omega_{Q_\Delta}$ by $\mathcal{Z}^1_\Delta$, ${\rm cl}^\omega_\Delta$ and ${\rm Wl}^\omega_\Delta$.
\end{definition}

So, for any quiver $Q$, an element of $\mathcal{Z}^\omega_Q$ is a Weyl-ordered Laurent polynomial if and only if it is in the image of ${\rm Wl}^\omega_Q$. Meanwhile, ${\rm cl}^\omega_Q \circ {\rm Wl}^\omega_Q={\rm id}$ obviously holds. Hence ${\rm Wl}^\omega_Q \circ {\rm cl}^\omega_Q \circ {\rm Wl}^\omega_Q={\rm Wl}^\omega_Q$, and therefore ${\rm Wl}^\omega_Q \circ {\rm cl}^\omega_Q={\rm id}$ holds when applied to Weyl-ordered Laurent polynomials. From this one obtains the following lemma.
\begin{lemma}
\label{lem:classicalization_determines_quantum}
A Weyl-ordered Laurent polynomial ${\bf U} \in \mathcal{Z}^\omega_Q$ for a quiver $Q$ is completely determined by its classicalization ${\rm cl}^\omega_Q({\bf U})\in \mathcal{Z}^1_Q$. That is, if ${\bf U}$ and ${\bf V}$ are Weyl-ordered Laurent polynomials in $\mathcal{Z}^\omega_Q$ such that ${\rm cl}^\omega_Q({\bf U}) = {\rm cl}^\omega_Q({\bf V})$, then ${\bf U} = {\bf V}$. \qed
\end{lemma}

Keeping Lem.\ref{lem:classicalization_determines_quantum} in mind, we propose the following strategy to prove eq.\eqref{eq:boiled_down2}. We are using the same notations for $\mathcal{Q}$, $\Delta$, $\Delta'$ and $W$ as used in the last subsection \S\ref{subsec:the_base_case} for eq.\eqref{eq:boiled_down}, and those for intermediate cluster $\mathscr{X}$-seeds $\Delta^{(r)}$ connecting the seeds for $\Delta$ and $\Delta'$ as in \S\ref{subsec:the_balanced_algebras_and_quantum_coordinate_change_maps_for_them}; in particular, $\Delta^{(2)} = \mu_{v_4}\mu_{v_3} \Delta = \mu_{v_7} \mu_{v_{12}} \Delta'$.

\vs

\begin{enumerate}\itemsep0,4em
\item[\rm Step 1.] Show that $\nu^\omega_{v_4} \nu^\omega_{v_3} ({\rm Tr}^\omega_{\Delta;\mathcal{Q}}([W,s])) \in {\rm Frac}(\mathcal{Z}^\omega_{\Delta^{(2)}})$ lies in $\mathcal{Z}^\omega_{\Delta^{(2)}}$ (i.e. is Laurent) and is a multiplicity-free Laurent polynomial in $\mathcal{Z}^\omega_{\Delta^{(2)}}$ (in the sense of Lem.\ref{lem:star-invariance_and_Weyl-ordering_new}(3)).

\item[\rm Step 2.] Show that $\nu^\omega_{v_7} \nu^\omega_{v_{12}} ({\rm Tr}^\omega_{\Delta';\mathcal{Q}}([W,s])) \in {\rm Frac}(\mathcal{Z}^\omega_{\Delta^{(2)}})$ lies in $\mathcal{Z}^\omega_{\Delta^{(2)}}$ (i.e. is Laurent) and is a multiplicity-free Laurent polynomial in $\mathcal{Z}^\omega_{\Delta^{(2)}}$.

\item[\rm Step 3.] Show that $\nu^\omega_{v_4} \nu^\omega_{v_3} ({\rm Tr}^\omega_{\Delta;\mathcal{Q}}([W,s]))$ and $\nu^\omega_{v_7} \nu^\omega_{v_{12}} ({\rm Tr}^\omega_{\Delta';\mathcal{Q}}([W,s]))$ are fixed by the $*$-map.

\item[\rm Step 4.] Show that $\nu^\omega_{v_4} \nu^\omega_{v_3} ({\rm Tr}^\omega_{\Delta;\mathcal{Q}}([W,s]))$ and $\nu^\omega_{v_7} \nu^\omega_{v_{12}} ({\rm Tr}^\omega_{\Delta';\mathcal{Q}}([W,s]))$ have the same classicalizations.
\end{enumerate}

\vs

From Steps 1, 2, and 3 it would follow that both $\nu^\omega_{v_4} \nu^\omega_{v_3} ({\rm Tr}^\omega_{\Delta;\mathcal{Q}}([W,s]))$ and $\nu^\omega_{v_7} \nu^\omega_{v_{12}} ({\rm Tr}^\omega_{\Delta';\mathcal{Q}}([W,s]))$ are $*$-invariant multiplicity-free Laurent polynomials in $\mathcal{Z}^\omega_{\Delta^{(2)}}$, hence are Weyl-ordered Laurent polynomials, by Lem.\ref{lem:star-invariance_and_Weyl-ordering_new}(3). Then, from Step 4 and Lem.\ref{lem:classicalization_determines_quantum} it would follow that they are equal, as desired in eq.\eqref{eq:boiled_down2}.

\vs

Steps 3 and 4 are relatively easy, so we do them here now.

\vs

\ul{Step 3.} We first recall the following result from \cite{Kim}.

\begin{proposition}[elevation reversing and $*$-structure {\cite[Prop.5.25]{Kim}}]
\label{prop:elevation_reversing_and_star}
Let $\Delta$ be an ideal triangulation of a triangulable generalized marked surface $\frak{S}$. Then
$$
{\rm Tr}^\omega_\Delta \circ {\bf r} = * \circ {\rm Tr}^\omega_\Delta
$$
holds, where 
$$
{\bf r} : \mathcal{S}^\omega_{\rm s}(\frak{S};\mathbb{Z})_{\rm red} \to \mathcal{S}^\omega_{\rm s}(\frak{S};\mathbb{Z})_{\rm red}
$$
is the \ul{\em elevation-reversing map}, defined as the $\mathbb{Z}$-linear map sending $\omega^{\pm 1/2}$ to $\omega^{\mp 1/2}$ and $[W,s]$ to $[W',s']$, whenever $W$ and $W'$ are ${\rm SL}_3$-webs in $\frak{S}\times {\bf I}$ with upward vertical framing with no crossings such that $W'$ is obtained from $W$ by reversing the elevation of all points, i.e. replacing each point $(x,t) \in \frak{S} \times {\bf I}$ by $(x,-t)$, and $s'(x,-t) = s(x,t)$, and the $*$-map $* : \mathcal{Z}^\omega_\Delta \to \mathcal{Z}^\omega_\Delta$ on $\mathcal{Z}^\omega_\Delta$ is as defined in Def.\ref{def:FG_algebra}.
\end{proposition}

Note that when $W$ is a constant-elevation lift in $\mathcal{Q}\times {\bf I}$, with upward vertical framing, of a corner arc in $\mathcal{Q}$, the element $[W,s]$ of $\mathcal{S}^\omega_{\rm s}(\mathcal{Q};\mathbb{Z})_{\rm red}$ is fixed by the elevation-reversing map ${\bf r}$, and hence ${\rm Tr}^\omega_{\Delta;\mathcal{Q}}([W,s]) \in \mathcal{Z}^\omega_\Delta$ and ${\rm Tr}^\omega_{\Delta';\mathcal{Q}}([W,s]) \in \mathcal{Z}^\omega_{\Delta'}$ are fixed by the $*$-maps; indeed, using Prop.\ref{prop:elevation_reversing_and_star} note that $({\rm Tr}^\omega_{\Delta;\mathcal{Q}}([W,s]))^* = {\rm Tr}^\omega_{\Delta;\mathcal{Q}}({\bf r}([W,s])) = {\rm Tr}^\omega_{\Delta;\mathcal{Q}}([W,s])$, and likewise for ${\rm Tr}^\omega_{\Delta';\mathcal{Q}}([W,s])$. Now, to finish proving Step 3, it is enough to recall from Lem.\ref{lem:nu_omega_k_preserves_star} that the maps $\nu^\omega_{v_3}$, $\nu^\omega_{v_4}$, $\nu^\omega_{v_7}$ and $\nu^\omega_{v_{12}}$ all preserve the $*$-structures. Indeed, note $(\nu^\omega_{v_4}(\nu^\omega_{v_3}({\rm Tr}^\omega_{\Delta;\mathcal{Q}}([W,s]))))^* = \nu^\omega_{v_4}((\nu^\omega_{v_3}({\rm Tr}^\omega_{\Delta;\mathcal{Q}}([W,s])))^*) = \nu^\omega_{v_4}(\nu^\omega_{v_3}(({\rm Tr}^\omega_{\Delta;\mathcal{Q}}([W,s]))^*)) = \nu^\omega_{v_4}(\nu^\omega_{v_3}({\rm Tr}^\omega_{\Delta;\mathcal{Q}}([W,s])))$, and similarly for $\nu^\omega_{v_7}\nu^\omega_{v_{12}}({\rm Tr}^\omega_{\Delta';\mathcal{Q}}([W,s]))$. So Step 3 is done.

\begin{remark}
\label{rem:highest_term_using_elevation-reversing}
Before going on to Step 4, we take a slight digression which utilizes Prop.\ref{prop:elevation_reversing_and_star}, which was promised in the discussion immediately following Prop.\ref{prop:highest_term}. Namely, let $\frak{S}$, $\Delta$ and $(W,s)$ be as in Prop.\ref{prop:highest_term}. By Prop.\ref{prop:highest_term}, the unique highest term of ${\rm Tr}^\omega_\Delta([W,s]) \in \mathcal{Z}^\omega_\Delta$ equals $\omega^m \, [\prod_{v\in\mathcal{V}} {\bf X}_v^{{\rm a}_v(\pi(W))}]_{\rm Weyl}$ for some $m \in \frac{1}{2}\mathbb{Z}$. The question was why there should be $\omega^m$, in case when $W$ has endpoints. Since $W$ has upward vertical framing and has no crossings, we have ${\bf r}([W,s]) = [W',s']$, where $W'$ is an ${\rm SL}_3$-web in $\frak{S}\times {\bf I}$ with upward vertical framing with no crossings such that the projection $\pi(W')$ coincides with $\pi(W)$, and $s'$ is the state of $W'$ that assigns $1$ to all endpoints. So $W$ and $W'$ can be thought of as being related by changing the orderings of elevations of endpoints. It is a straightforward exercise using \cite[Prop.5.50, Prop.5.27]{Kim} to show that ${\rm Tr}^\omega_\Delta([W',s']) = \omega^{3N} {\rm Tr}^\omega_\Delta([W,s])$ holds for some $N\in \mathbb{Z}$. In the meantime, from ${\bf r}([W,s]) = [W',s']$ one obtains ${\rm Tr}^\omega_\Delta({\bf r}([W,s])) = {\rm Tr}^\omega_\Delta([W',s'])$, whose left-hand side equals $*({\rm Tr}^\omega_\Delta([W,s]))$ by Prop.\ref{prop:elevation_reversing_and_star}, and whose right-hand side equals $\omega^{3N} {\rm Tr}^\omega_\Delta([W,s])$. So we have $*({\rm Tr}^\omega_\Delta([W,s])) = \omega^{3N} {\rm Tr}^\omega_\Delta([W,s])$. It is easy to see that the highest term of the left-hand side is $*(\omega^m \, [\prod_{v\in\mathcal{V}} {\bf X}_v^{{\rm a}_v(\pi(W))}]_{\rm Weyl}) = \omega^{-m} \, [\prod_{v\in\mathcal{V}} {\bf X}_v^{{\rm a}_v(\pi(W))}]_{\rm Weyl}$, whereas the highest term of the right-hand side is $\omega^{3N} \omega^m \, [\prod_{v\in\mathcal{V}} {\bf X}_v^{{\rm a}_v(\pi(W))}]_{\rm Weyl}$. This shows that $\omega^{-m} = \omega^{3N+m}$, hence $m = - \frac{3}{2}N$. In particular, the highest term of ${\rm Tr}^\omega_\Delta([W,s]) \in \mathcal{Z}^\omega_\Delta$ equals $\omega^m \, [\prod_{v\in\mathcal{V}} {\bf X}_v^{{\rm a}_v(\pi(W))}]_{\rm Weyl}$ for some $m \in \frac{3}{2} \mathbb{Z}$, which slightly refines the statement of Prop.\ref{prop:highest_term}. The arguments we gave here, which depend on \cite[Prop.5.50, Prop.5.27]{Kim}, reveal that $m$ is not in general $0$, although one could get a bound of its absolute value in terms of the number of endpoints of $W$.
\end{remark}

\ul{Step 4.} This statement is the classical version, i.e. for $\omega^{1/2}=1$, of  eq.\eqref{eq:boiled_down2}. Recall that we showed at the end of \S\ref{subsec:the_base_case} that eq.\eqref{eq:boiled_down2} is equivalent to eq.\eqref{eq:boiled_down} for a general parameter $\omega^{1/2}$. Each step of the proof works also when $\omega^{1/2}=1$, so one can deduce that the classical version of eq.\eqref{eq:boiled_down2} is equivalent to the classical version of eq.\eqref{eq:boiled_down}. So it is enough to just show the latter. We formulate this classical statement for more general setting as follows, as a classical counterpart of Thm.\ref{thm:main}.

\begin{proposition}[naturality of ${\rm SL}_3$ classical trace maps under changes of triangulations] 
\label{prop:classical_compatibility}
Let $\Delta$ and $\Delta'$ be ideal triangulations of a triangulable generalized marked surface $\frak{S}$. Let $(W,s)$ be a stated ${\rm SL}_3$-web in $\frak{S}\times {\bf I}$. Then the values under the ${\rm SL}_3$ quantum trace maps of $[W,s] \in \mathcal{S}^\omega_{\rm s}(\frak{S};\mathbb{Z})_{\rm red}$ in the case when $\omega^{1/2}=1$ are compatible under the coordinate change map, i.e.
$$
\Theta^1_{\Delta\Delta'}({\rm Tr}^1_{\Delta'}([W,s])) = {\rm Tr}^1_\Delta([W,s]).
$$
\end{proposition}
Note that, by $\Theta^1_{\Delta\Delta'}$ we mean the limit of $\Theta^\omega_{\Delta\Delta'}$ as $\omega^{1/2} \to 1$. This Prop.\ref{prop:classical_compatibility} is proven in \cite[Cor.5.70]{Kim} in the case when $\frak{S}$ is a punctured surface. Essentially the same proof applies to the case when $\frak{S}$ is a generalized marked surface.

\vs

{\it Proof of Prop.\ref{prop:classical_compatibility}.} By the argument of \S\ref{subsec:the_base_case}, it suffices to prove this when $\frak{S}$ is a quadrilateral $\mathcal{Q}$ and $W$ is a constant-elevation lift, with upward vertical framing, of a corner arc of $\mathcal{Q}$. For such a case, one can verify the equality by direct computations. Here we provide a proof using arguments in \cite{Kim}. Suppose the state $s$ assigns $\varepsilon_1,\varepsilon_2\in\{1,2,3\}$ to the initial and the terminal endpoints of $W$; then ${\rm Tr}^1_\Delta([W,s])$ is the $(\varepsilon_1,\varepsilon_2)$-th entry of the normalized monodromy matrix ${\rm M}_{W;\Delta}$ associated to $W$, which is a product of normalized basic monodromy matrices associated to small pieces of $W$, as described in \cite[\S4.2]{Kim}; see (MM1)--(MM3) of \cite[\S4.2]{Kim} for the basic monodromy matrices. Note that ${\rm M}_{W;\Delta}$ is a $3 \times 3$ matrix whose entries are in $\mathcal{Z}^1_\Delta$, having determinant $1$; i.e. ${\rm M}_{W;\Delta} \in {\rm SL}_3(\mathcal{Z}^1_\Delta)$. Likewise, ${\rm Tr}^1_{\Delta'}([W,s])$ is the $(\varepsilon_1,\varepsilon_2)$-th entry of the normalized monodromy matrix ${\rm M}_{W;\Delta'} \in {\rm SL}_3(\mathcal{Z}^1_{\Delta'})$. By inspection of the entries, especially of the $(1,1)$-th entries which are Laurent monomials of the highest preorder (induced by the degrees) whose degrees are given by the tropical coordinates of the ${\rm SL}_3$-laminations $\pi(W)$ (see Prop.\ref{prop:highest_term}, and also \cite[Fig.5]{Kim}), one observes that ${\rm M}_{W;\Delta} \in {\rm SL}_3(\wh{\mathcal{Z}}^1_\Delta)$ and ${\rm M}_{W;\Delta'} \in {\rm SL}_3(\wh{\mathcal{Z}}^1_{\Delta'})$. Meanwhile, as noted in \cite[\S4.2]{Kim}, observe that these monodromy matrices ${\rm M}_{W;\Delta}$ and ${\rm M}_{W;\Delta'}$ are normalized versions of Fock and Goncharov's unnormalized monodromy matrices $\til{\rm M}_{W;\Delta} \in {\rm GL}_3(\mathcal{X}^1_\Delta)$ and $\til{\rm M}_{W;\Delta'} \in {\rm GL}_3(\mathcal{X}^1_{\Delta'})$, appearing in \cite{FG06}. Observe also that in \cite{FG06}, it is shown that $\til{\rm M}_{W;\Delta}$ and $\til{\rm M}_{W;\Delta'}$, viewed as elements of the projective transformations ${\rm PGL}_3(\mathcal{X}^1_\Delta)$ and ${\rm PGL}_3(\mathcal{X}^1_{\Delta'})$, are related to each other by the composition $\mu^1_{v_3} \mu^1_{v_4} \mu^1_{v_7} \mu^1_{v_{12}}$ of the sequence of four (usual) classical $\mathscr{X}$-mutations which appeared in \S\ref{sec:quantum_coordinate_change_for_flips_of_ideal_triangulations}. Viewing the $\mathscr{X}$-coordinates $X_v$'s and $X_v'$'s as real-valued functions (on the set $\mathscr{X}_{{\rm PGL}_3,\frak{S}}(\mathbb{R})$ of $\mathbb{R}$-points of the moduli stack $\mathscr{X}_{{\rm PGL}_3,\frak{S}}$), one can view $Z_v$'s and $Z'_v$'s as unique real-valued cube-roots of $X_v$'s and $X_v'$'s. Now, as in \cite[\S4.2]{Kim}, from the fact that the projection ${\rm SL}_3(\mathbb{R}) \to {\rm PGL}_3(\mathbb{R})$ is bijective, one can deduce that the normalized matrices ${\rm M}_{W;\Delta}$ and ${\rm M}_{W;\Delta'}$ are related by the composition $\mu^1_{v_3} \mu^1_{v_4} \mu^1_{v_7} \mu^1_{v_{12}}$ of the four cluster $\mathscr{X}$-mutations. More precisely, this last statement is about the evaluations at $\mathbb{R}$. One can finish the proof by observing that the mutation formulas for the balanced subalgebras, in case $\omega^{1/2}=1$, when evaluated at $\mathbb{R}$, are compatible with the evaluation of the usual cluster $\mathscr{X}$-mutation. \qed

\vs

This finishes Step 4.

\vs

\begin{figure}[htbp!]
\vspace{-1mm}
\hspace*{-5mm}
\raisebox{-0.5\height}{\scalebox{0.85}{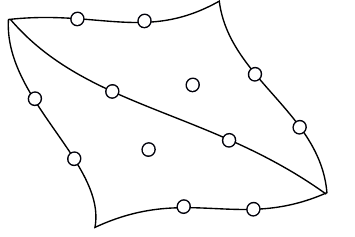}}
\vspace{-4mm}
\caption{Four cases of ${\rm SL}_3$-webs $W$ over a quadrilateral $\mathcal{Q}$, with triangulation $\Delta'$}
\vspace{-2mm}
\label{fig:four_cases}
\end{figure}

Observe that (as mentioned already), for the remaining Steps 1 and 2, the situation for $\nu^\omega_{v_4} \nu^\omega_{v_3} {\rm Tr}^\omega_{\Delta;\mathcal{Q}}([W,s])$ and that for $\nu^\omega_{v_7} \nu^\omega_{v_{12}} {\rm Tr}^\omega_{\Delta';\mathcal{Q}}([W,s])$ are symmetric. Namely, each of these two situations can be interpreted as first applying the ${\rm SL}_3$ quantum trace for an ideal triangulation of $\mathcal{Q}$ and then applying the mutations at the two nodes lying in the unique internal arc of this ideal triangulation. The two mutations can be taken in any order, as they commute with each other, by the relation $\nu^\omega_k \nu^\omega_\ell = \nu^\omega_\ell \nu^\omega_k$ for $\varepsilon_{k\ell}=0$ in Prop.\ref{prop:consistency_for_Theta}. In general, when dealing with relations among $\nu^\omega_k$'s one has to be careful about the domains; here we specifically mean the equations $\nu^\omega_{v_4} \nu^\omega_{v_3} {\rm Tr}^\omega_{\Delta;\mathcal{Q}}([W,s]) = \nu^\omega_{v_3} \nu^\omega_{v_4} {\rm Tr}^\omega_{\Delta;\mathcal{Q}}([W,s])$ and $\nu^\omega_{v_7} \nu^\omega_{v_{12}} {\rm Tr}^\omega_{\Delta';\mathcal{Q}}([W,s]) = \nu^\omega_{v_{12}} \nu^\omega_{v_7} {\rm Tr}^\omega_{\Delta';\mathcal{Q}}([W,s])$, which indeed make sense and hence fall into the case of Prop.\ref{prop:consistency_for_Theta}; one can use Lem.\ref{lem:Theta_well-definedness2} and Prop.\ref{prop:value_of_the_SL3_quantum_trace_is_Delta-balanced} appropriately to check that both sides of both equations do make sense. Because of this symmetry, checking Step 1 for all possible cases when $W$ is a constant-elevation lift of a corner arc of $\mathcal{Q}$
is equivalent to checking Step 2 for all those cases of $W$, and also equivalent to checking Step 2 for the four cases of $W$ depicted in Fig.\ref{fig:four_cases}. We perform the computational check of Step 2 for $W$ in Fig.\ref{fig:four_cases} in the following subsections.

\vspace{0mm}

\subsection{Reducing the amount of computation}
\label{subsec:reducing_the_amount_of_computation}

It remains to do Step 2 
of the previous subsection, for all cases of $W$ as in Fig.\ref{fig:four_cases} . For this task, we make use of the following convenient way of dealing with the Weyl-ordered Laurent monomials by introducing the log variables.
\begin{definition}
\label{def:Heisenberg}
Let $\Delta$ be an ideal triangulation of a triangulable generalized marked surface. Let $\mathcal{V} = \mathcal{V}(Q_\Delta)$. Let $\mathcal{H}_\Delta$ be a free $\mathbb{Z}$-module with the set of symbols $\{ {\bf z}_v \, | \, v\in \mathcal{V}\} \cup \{\frac{1}{18}  {\bf c}\}$ as a free basis, equipped with a skew-symmetric bilinear map $[\cdot ,\cdot]: \mathcal{H}_\Delta\times\mathcal{H}_\Delta\to (\frac{1}{18}{\bf c})\mathbb{Z}$ such that
$$
[{\bf z}_v,{\bf z}_w] = 2\varepsilon_{vw} \cdot 2\cdot {\textstyle \frac{1}{18}} {\bf c}, \qquad
[{\textstyle \frac{1}{18}}{\bf c}, {\bf z}_v] = 0, \qquad \forall v,w\in \mathcal{V}.
$$
Define
$$
{\bf x}_v := 3{\bf z}_v, \quad \forall v\in \mathcal{V},
$$
so that $\frac{1}{3} {\bf x}_v$ means ${\bf z}_v$. Meanwhile, for any $a \in \frac{1}{18}\mathbb{Z}$, $a{\bf c}$ denotes $(18a)\cdot(\frac{1}{18}{\bf c})$. 
Define the {\em exponential} map
$$
\exp ~:~ \mathcal{H}_\Delta \to \mathcal{Z}^\omega_\Delta, \quad {\bf z} \mapsto \exp({\bf z}) =: e^{{\bf z}}
$$
as
\begin{align*}
\exp( \alpha \cdot {\textstyle \frac{1}{18}}{\bf c} + {\textstyle \sum}_{v\in\mathcal{V}} a_v {\bf x}_v) := \omega^{\alpha/2} \, [{\textstyle \prod}_{v\in\mathcal{V}} {\bf X}_v^{a_v}]_{\rm Weyl},
\end{align*}
where $\alpha \in \mathbb{Z}$ and $(a_v)_{v\in\mathcal{V}} \in (\frac{1}{3}\mathbb{Z})^\mathcal{V}$.
\end{definition}
For example, $e^{\frac{1}{18}{\bf c}} = \omega^{1/2}$, $e^{\frac{1}{9}{\bf c}} = e^{2\cdot \frac{1}{18}{\bf c}} = \omega$, $e^{{\bf c}}=q$, $e^{{\bf z}_v} = {\bf Z}_v=e^{\frac{1}{3}{\bf x}_v}$, $e^{{\bf x}_v} = {\bf X}_v$. The bracket $[\cdot, \cdot]$ can be written in terms of ${\bf x}_v$'s as follows, which will become a handy formula to have:
\begin{align}
\label{eq:bracket_in_x}
[a_v {\bf x}_v, \, a_w {\bf x}_w] = 2 a_v a_w \varepsilon_{vw} \, {\bf c}, \qquad \forall v,w\in\mathcal{V}, \quad \forall a_v,a_w \in \textstyle \frac{1}{3}\mathbb{Z}.
\end{align}
We find it useful to have the following well-known fact. We provide a proof, for completeness.
\begin{lemma}[Baker-Campbell-Hausdorff (BCH) formula]
\label{lem:BCH}
For ${\bf x},{\bf y} \in \mathcal{H}_\Delta$,
$$
\exp({\bf x}) \exp({\bf y}) = e^{\frac{1}{2}[{\bf x},{\bf y}]} \exp({\bf x} + {\bf y}). 
$$
\end{lemma}
{\it Proof.} The statement is easily seen to hold if one of ${\bf x}$ and ${\bf y}$ is an integer multiple of $\frac{1}{18}{\bf c}$. So it suffices to show the statement in the case when ${\bf x} = \sum_{v\in \mathcal{V}} a_v {\bf x}_v$ and ${\bf y} = \sum_{v\in\mathcal{V}} b_v {\bf x}_v$, with $(a_v)_{v\in\mathcal{V}}, (b_v)_{v\in\mathcal{V}} \in (\frac{1}{3}\mathbb{Z})^\mathcal{V}$. Then ${\bf x}+ {\bf y} = \sum_{v\in\mathcal{V}}(a_v+b_v){\bf x}_v$. Choose any ordering on $\mathcal{V}$ to write its elements as $v_1,v_2,\ldots,v_N$. For notational convenience, write $a_{v_r} = a_r$ and $b_{v_r} = b_r$ for each $r=1,\ldots,N$, and $\varepsilon_{v_r v_s} =\varepsilon_{rs}$ for $r,s = 1,\ldots,N$. Note that
\begin{align*}
[{\textstyle \prod}_{v\in\mathcal{V}} {\bf X}_v^{a_v}]_{\rm Weyl} [{\textstyle \prod}_{v\in\mathcal{V}} {\bf X}_v^{b_v}]_{\rm Weyl}
& = q^{-\sum_{r<s} \varepsilon_{rs} a_r a_s} ({\bf X}_{v_1}^{a_1} \cdots {\bf X}_{v_N}^{a_N}) \cdot q^{-\sum_{r<s} \varepsilon_{rs} b_r b_s} ({\bf X}_{v_1}^{b_1} \cdots {\bf X}_{v_N}^{b_N}) \quad (\because \mbox{Def.\ref{def:Weyl-ordered}}) \\
& = q^{-\sum_{r<s}\varepsilon_{rs} (a_ra_s + b_rb_s)}  {\bf X}_{v_1}^{a_1} \cdots {\bf X}_{v_N}^{a_N} \, {\bf X}_{v_1}^{b_1} \cdots {\bf X}_{v_N}^{b_N} 
\end{align*}
Now, we use Lem.\ref{lem:relations_of_cube-root_algebra_by_X} to move the factor ${\bf X}_{v_1}^{b_1}$ to the left till it meets ${\bf X}_{v_1}^{a_1}$, then move ${\bf X}_{v_2}^{b_2}$ to the left till it meets ${\bf X}_{v_2}^{a_2}$, etc.
\begin{align*}
& [{\textstyle \prod}_{v\in\mathcal{V}} {\bf X}_v^{a_v}]_{\rm Weyl} [{\textstyle \prod}_{v\in\mathcal{V}} {\bf X}_v^{b_v}]_{\rm Weyl} \\
& = q^{-\sum_{r<s}\varepsilon_{rs} (a_ra_s + b_rb_s)} q^{-\sum_{r=2}^N \varepsilon_{r1} a_r b_1} ({\bf X}_{v_1}^{a_1}{\bf X}_{v_1}^{b_1}) \, {\bf X}_{v_2}^{a_2} \cdots {\bf X}_{v_N}^{a_N} \, {\bf X}_{v_2}^{b_2} \cdots {\bf X}_{v_N}^{b_N} \\
& = q^{-\sum_{r<s}\varepsilon_{rs} (a_ra_s + b_rb_s)} q^{2\sum_{r=2}^N \varepsilon_{r1} a_r b_1}  q^{2\sum_{r=3}^N \varepsilon_{r2} a_r b_2} ({\bf X}_{v_1}^{a_1}{\bf X}_{v_1}^{b_1}) ({\bf X}_{v_2}^{a_2}{\bf X}_{v_2}^{b_2}) \, {\bf X}_{v_3}^{a_3} \cdots {\bf X}_{v_N}^{a_N} \, {\bf X}_{v_3}^{b_3} \cdots {\bf X}_{v_N}^{b_N} \\
& = \cdots = q^{-\sum_{r<s} \varepsilon_{rs}(a_ra_s+b_rb_s)} q^{2\sum_{r>s} \varepsilon_{rs} a_r b_s} ({\bf X}_{v_1}^{a_1} {\bf X}_{v_1}^{b_1})  ({\bf X}_{v_2}^{a_2} {\bf X}_{v_2}^{b_2}) \cdots ({\bf X}_{v_N}^{a_N} {\bf X}_{v_N}^{b_N}) \\
& = q^{-\sum_{r<s} \varepsilon_{rs}(a_ra_s+b_rb_s)}
q^{2\sum_{r<s} \varepsilon_{sr} a_s b_r} {\bf X}_{v_1}^{a_1+b_1} \cdots {\bf X}_{v_N}^{a_N+b_N}
\end{align*}
where for the last equality we exchanged the indices $r$ and $s$ for the second sum over $r$ and $s$ appearing in the exponent of $q$. On the other hand, note from Def.\ref{def:Weyl-ordered} that
$$
[{\textstyle \prod}_{v\in\mathcal{V}} {\bf X}_v^{a_v+b_v}]_{\rm Weyl}
= q^{-\sum_{r<s} \varepsilon_{rs} (a_r+b_r)(a_s+b_s)} {\bf X}_{v_1}^{a_1+b_1} \cdots {\bf X}_{v_N}^{a_N+b_N},
$$
and note that $\frac{1}{2} [{\bf x},{\bf y}] = \frac{1}{2} \sum_{r,s\in\{1,\ldots,N\}} a_r b_s[{\bf x}_{v_r},{\bf x}_{v_s}] = \sum_{r,s\in\{1,\ldots,N\}}  \varepsilon_{rs} a_r b_s {\bf c}$, where we used eq.\eqref{eq:bracket_in_x}. So, in view of $q^a = e^{a{\bf c}}$ for $a\in\frac{1}{18}\mathbb{Z}$, to show $\exp({\bf x}) \exp({\bf y}) = e^{\frac{1}{2}[{\bf x},{\bf y}]} \exp({\bf x} + {\bf y})$ it remains to show
$$
- {\textstyle \sum}_{r<s} \varepsilon_{rs}(a_ra_s+b_rb_s) + 2{\textstyle \sum}_{r<s} \varepsilon_{sr} a_sb_r
= {\textstyle \sum}_{r,s\in\{1,\ldots,N\}} \varepsilon_{rs} a_r b_s - {\textstyle \sum}_{r<s} \varepsilon_{rs}(a_r+b_r)(a_s+b_s)
$$
Subtracting $- {\textstyle \sum}_{r<s} \varepsilon_{rs}(a_ra_s+b_rb_s)$ from both sides, rewriting $\varepsilon_{sr}$ in the second sum of the left-hand side as $-\varepsilon_{rs}$, and adding $\sum_{r<s} \varepsilon_{rs}(a_r b_s + b_r a_s)$ to both sides, this equation becomes
$$
{\textstyle \sum}_{r<s} \varepsilon_{rs}(-2a_sb_r+a_rb_s+b_ra_s) = {\textstyle \sum}_{r,s \in\{1,2,\ldots,N\}} \varepsilon_{rs} a_r b_s.
$$
The left-hand side equals $\sum_{r<s} \varepsilon_{rs} (-a_s b_r + a_r b_s)$. For the right-hand side, write the sum $\sum_{r,s}$ as the sum of three sums $\sum_{r=s}(\sim) + \sum_{r<s} (\sim)+ \sum_{r>s}(\sim)$; the first sum for $r=s$ vanishes because $\varepsilon_{rr}=0$, the second sum equals $\sum_{r<s} \varepsilon_{rs} a_r b_s$, and the third sum can be written with the indices $r$ and $s$ exchanged as $\sum_{s>r} \varepsilon_{sr} a_s b_r = \sum_{r<s} \varepsilon_{rs}(-a_sb_r)$. So one observes that the left-hand side equals the right-hand side. \qed

\vs

Step 2 involves twelve nodes $v_1,\ldots,v_{12}$ appearing in Fig.\ref{fig:mutations_for_a_flip}. We denote the generators $({\bf Z}_{v_j}^{(r)})^{\pm 1} = ({\bf X}_{v_j}^{(r)})^{\pm 1/3}$, $j=1,\ldots,12$, of the cube-root Fock-Goncharov algebra $\mathcal{Z}^\omega_{\Delta^{(r)}}$, $r=0,1,2,3,4$, by $({\bf Z}_j^{(r)})^{\pm 1} = ({\bf X}_j^{(r)})^{\pm 1/3}$, and the `log' generators $\frac{1}{3} {\bf x}^{(r)}_{v_j} = {\bf z}^{(r)}_{v_j}$ of $\mathcal{H}_{\Delta^{(r)}}$ (Def.\ref{def:Heisenberg}) by $\frac{1}{3}{\bf x}^{(r)}_j = {\bf z}^{(r)}_j$. Here $\Delta^{(r)}$, $r=0,1,2,3,4$, is defined as in the discussion above eq.\eqref{eq:Delta_r} in \S\ref{subsec:the_balanced_algebras_and_quantum_coordinate_change_maps_for_them}. The entries $\varepsilon^{(r)}_{v_j,v_k}$ of the exchange matrix $\varepsilon^{(r)}$ for $\Delta^{(r)}$ are written as $\varepsilon^{(r)}_{jk}$. An element $(a^{(r)}_{v_j})_{j=1}^{12}$ of $(\frac{1}{3}\mathbb{Z})^{\mathcal{V}}$ would be denoted by $(a^{(r)}_j)_{j=1}^{12} \in (\frac{1}{3}\mathbb{Z})^{12}$. As in \S\ref{sec:quantum_coordinate_change_for_flips_of_ideal_triangulations}, for $r=4$ the superscript $(4)$ may be replaced by the prime symbol $'$, such as ${\bf Z}_j^{(4)} = {\bf Z}_j'$, ${\bf X}_j^{(4)} = {\bf X}_j'$, ${\bf z}_j^{(4)} = {\bf z}_j'$, ${\bf x}_j^{(4)} = {\bf x}_j'$, $\varepsilon^{(4)}_{jk} = \varepsilon'_{jk}$, $a^{(4)}_j = a'_j$.

\vs

For Step 2 (and Step 1), we need to investigate $\nu^\omega_{v_7} \nu^\omega_{v_{12}} ({\rm Tr}^\omega_{\Delta';\mathcal{Q}}([W,s]))$. By Prop.\ref{prop:value_of_the_SL3_quantum_trace_is_Delta-balanced}, we know that ${\rm Tr}^\omega_{\Delta';\mathcal{Q}}([W,s])$ is a sum of $\Delta'$-balanced Laurent monomials in $\mathcal{Z}^\omega_{\Delta'}$ (Def.\ref{def:balanced_subalgebras}). So we shall first investigate $\nu^\omega_{v_7}\nu^\omega_{v_{12}}({\bf U}')$ for an arbitrary $\Delta'$-balanced Laurent monomial ${\bf U}' = [\prod_{j=1}^{12} {{\bf X}_j'}^{a'_j}]_{\rm Weyl} = \exp(\sum_{j=1}^{12} a'_j {\bf x}'_j)$, where $(a'_j)_{j=1}^{12} \in (\frac{1}{3}\mathbb{Z})^{12} = (\frac{1}{3}\mathbb{Z})^\mathcal{V}$ is $\Delta'$-balanced (Def.\ref{def:Delta-balanced_elements}). In fact, in the proof of Lem.\ref{lem:Theta_well-definedness1} we have already performed this computation. In particular, there we checked that ${\bf U}'$, hence ${\rm Tr}^\omega_{\Delta';\mathcal{Q}}([W,s])$ as well, belongs to $\wh{\rm Frac}_{v_{12}}(\mathcal{Z}^\omega_{\Delta'})$ which is the domain of $\nu^\omega_{v_{12}} = \nu^{\sharp\omega}_{v_{12}} \nu'_{v_{12}} : \wh{\rm Frac}_{v_{12}}(\mathcal{Z}^\omega_{\Delta'}) \to \wh{\rm Frac}_{v_{12}}(\mathcal{Z}^\omega_{\Delta^{(3)}}) \subset {\rm Frac}(\mathcal{Z}^\omega_{\Delta^{(3)}})$, and that $\nu^\omega_{v_{12}}({\bf U}')$ belongs to $\wh{\rm Frac}_{v_7}(\mathcal{Z}^\omega_{\Delta^{(3)}}) \subset {\rm Frac}(\mathcal{Z}^\omega_{\Delta^{(3)}})$ which is the domain of $\nu^\omega_{v_7} = \nu^{\sharp \omega}_{v_7} \nu'_{v_7} : \wh{\rm Frac}_{v_7}(\mathcal{Z}^\omega_{\Delta^{(3)}}) \to \wh{\rm Frac}_{v_7}(\mathcal{Z}^\omega_{\Delta^{(2)}}) \subset {\rm Frac}(\mathcal{Z}^\omega_{\Delta^{(2)}})$. We also investigated explicitly how the image $\nu^\omega_{v_7} \nu^\omega_{v_{12}}({\bf U}') \in {\rm Frac}(\mathcal{Z}^\omega_{\Delta^{(2)}})$ is given, which we could recollect here. However, in order to lessen the amount of computation, we perform a slight manipulation, which essentially kills the effect of the maps $\nu'_{v_7}$ and $\nu'_{v_{12}}$. We first justify this process in the following lemma.

\begin{lemma}
\label{lem:OK_to_check_in_Delta_prime}
An element ${\bf Y}$ of ${\rm Frac}(\mathcal{Z}^\omega_{\Delta^{(2)}})$ is Laurent for $\Delta^{(2)}$ if and only if the element $(\nu'_{v_{12}})^{-1}(\nu'_{v_7})^{-1}({\bf Y})$ of ${\rm Frac}(\mathcal{Z}^\omega_{\Delta'})$ is Laurent for $\Delta'$, where $\nu'_{v_{12}} : {\rm Frac}(\mathcal{Z}^\omega_{\Delta'}) \to {\rm Frac}(\mathcal{Z}^\omega_{\Delta^{(3)}})$ and $\nu'_{v_7} : {\rm Frac}(\mathcal{Z}^\omega_{\Delta^{(3)}}) \to {\rm Frac}(\mathcal{Z}^\omega_{\Delta^{(2)}})$ are as defined in Def.\ref{def:cube-root_monomial_transformation}. Moreover, when ${\bf Y}$ is Laurent for $\Delta^{(2)}$, i.e. when ${\bf Y} \in \mathcal{Z}^\omega_{\Delta^{(2)}}$, it is multiplicity-free (in the sense of Lem.\ref{lem:star-invariance_and_Weyl-ordering_new}(3)) if and only if $(\nu'_{v_{12}})^{-1}(\nu'_{v_7})^{-1}({\bf Y}) \in \mathcal{Z}^\omega_{\Delta'}$ is multiplicity-free.
\end{lemma}
{\it Proof.} It is enough to recall from Lem.\ref{lem:nu_prime_k_basic_properties}(6) that as each of $\nu'_{v_{12}}$ and $\nu'_{v_7}$ is invertible and that each of them and their inverses sends Laurent monomials to Laurent monomials. \qed

\vs

We go on to investigate $(\nu'_{v_{12}})^{-1}(\nu'_{v_7})^{-1}( \nu^\omega_{v_7} \nu^\omega_{v_{12}}({\rm Tr}^\omega_{\Delta';\mathcal{Q}}([W,s])))$.

\vs

\begin{lemma}
\label{lem:section5_modified}
Let $(a'_j)_{j=1}^{12} \in (\frac{1}{3}\mathbb{Z})^{12} = (\frac{1}{3}\mathbb{Z})^\mathcal{V}$ be $\Delta'$-balanced (Def.\ref{def:Delta-balanced_elements}). Then, to the element $\exp({\textstyle \sum}_{j=1}^{12} a'_j {\bf x}_j')$ of $\mathcal{Z}^\omega_{\Delta'}$, the composition of maps $(\nu'_{v_{12}})^{-1}(\nu'_{v_7})^{-1}  \nu^\omega_{v_7} \nu^\omega_{v_{12}}$ can be applied, where $\nu^\omega_{v_{12}} : \wh{\rm Frac}_{v_{12}}(\mathcal{Z}^\omega_{\Delta'}) \to \wh{\rm Frac}_{v_{12}}(\mathcal{Z}^\omega_{\Delta^{(3)}})$ and $\nu^\omega_{v_7} : \wh{\rm Frac}_{v_7}(\mathcal{Z}^\omega_{\Delta^{(3)}}) \to \wh{\rm Frac}_{v_7}(\mathcal{Z}^\omega_{\Delta^{(2)}})$ are as in Def.\ref{def:nu_omega_k}, and $\nu'_{v_7} : {\rm Frac}(\mathcal{Z}^\omega_{\Delta^{(3)}}) \to {\rm Frac}(\mathcal{Z}^\omega_{\Delta^{(2)}})$ and $\nu'_{v_{12}} : {\rm Frac}(\mathcal{Z}^\omega_{\Delta'}) \to {\rm Frac}(\mathcal{Z}^\omega_{\Delta^{(3)}})$ are as in Def.\ref{def:cube-root_monomial_transformation}. The image is given by
\begin{align}
\label{eq:result_of_nu_omega_v_7_and_v_12_modified}
(\nu'_{v_{12}})^{-1}(\nu'_{v_7})^{-1} \nu^\omega_{v_7} \nu^\omega_{v_{12}}(\exp({\textstyle \sum}_{j=1}^{12} a'_j {\bf x}_j'))
=\exp({\textstyle \sum}_{j=1}^{12} a'_j {\bf x}_j') \, F^q( {{\bf X}'}_{\hspace{-1mm}7}^{-1}; \alpha^{(3)}) \, F^q( {{\bf X}'}_{\hspace{-1mm}12}^{-1}; \alpha^{(4)}), 
\end{align}
where $F^q$ is as in eq.\eqref{eq:F_q}, and the numbers $\alpha^{(3)}$ and $\alpha^{(4)}$ are given by
$$
\alpha^{(3)} := a'_1+a'_4-a'_3-a'_6, \qquad
\alpha^{(4)} := a'_3+a'_{10}-a'_4-a'_9,
$$
which are integers.
\end{lemma}

\vs

{\it Proof.} First, note that the domain issues are dealt with in the proof of Lem.\ref{lem:Theta_well-definedness1}. Namely, if we follow the notation there and let ${\bf U}' = \exp(\sum_{j=1}^{12} a'_j {\bf x }_j')$, then it is proved there that $\nu^\omega_{v_7} \nu^\omega_{v_{12}}({\bf U}')$ indeed makes sense as a well-defined element of ${\rm Frac}(\mathcal{Z}^\omega_{\Delta^{(2)}})$. There is no problem applying $(\nu'_{v_{12}})^{-1} (\nu'_{v_7})^{-1}$ to $\nu^\omega_{v_7} \nu^\omega_{v_{12}}({\bf U}')$ (one may want to see Lem.\ref{lem:nu_prime_k_basic_properties}(6)). 

\vs

In eq.\eqref{eq:two_nu_applied_to_Laurent_monomial} in the proof of Lem.\ref{lem:Theta_well-definedness1}, we obtained $\nu^\omega_{v_7} \nu^\omega_{v_{12}}({\bf U}')={\bf U}^{(2)} {\bf V}^{(2)} \mu^q_{v_7}({\bf V}^{(3)})$, where ${\bf U}^{(2)} = \nu'_{v_7}({\bf U}^{(3)}) = \nu'_{v_7}(\nu'_{v_{12}}({\bf U}'))$, ${\bf V}^{(2)} = F^q({\bf X}_7^{(2)};\alpha^{(3)})$ and ${\bf V}^{(3)} = F^q({\bf X}_{12}^{(3)}; \alpha^{(4)})$, where $\alpha^{(3)}$ and $\alpha^{(4)}$ are given by eq.\eqref{eq:alpha_for_nu_sharp}. Note that in the proof of Lem.\ref{lem:alpha} we found formulas to express $\alpha^{(3)}$ and $\alpha^{(4)}$ in terms of $(a'_j)_{j=1}^{12}$ in eq.\eqref{eq:alpha_r_in_a_prime}, namely $\alpha^{(3)} := a'_1+a'_4-a'_3-a'_6$ and $\alpha^{(4)} := a'_3+a'_{10}-a'_4-a'_9$, as written in the statement of the current lemma. The integrality of $\alpha^{(3)}$ and $\alpha^{(4)}$ was the content of Lem.\ref{lem:alpha}. To summarize, we have
$$
\nu^\omega_{v_7} \nu^\omega_{v_{12}}({\bf U}') = \left( \nu'_{v_7} \nu'_{v_{12}}({\bf U}') \right) \, F^q({\bf X}^{(2)}_7;\alpha^{(3)}) \, \mu^q_{v_7} (F^q({\bf X}^{(3)}_{12};\alpha^{(4)})).
$$
Thus, applying the skew field isomorphism $(\nu'_{v_{12}})^{-1} (\nu'_{v_7})^{-1} : {\rm Frac}(\mathcal{Z}^\omega_{\Delta^{(2)}}) \to {\rm Frac}(\mathcal{Z}^\omega_{\Delta'})$ to both sides, we get
\begin{align*}
(\nu'_{v_{12}})^{-1} (\nu'_{v_7})^{-1} \nu^\omega_{v_7} \nu^\omega_{v_{12}}({\bf U}') = {\bf U}' \, F^q( \underbrace{(\nu'_{v_{12}})^{-1} (\nu'_{v_7})^{-1}({\bf X}^{(2)}_7)};\alpha^{(3)}) \,  F^q( \underbrace{(\nu'_{v_{12}})^{-1} (\nu'_{v_7})^{-1} \mu^q_{v_7}({\bf X}^{(3)}_{12})};\alpha^{(4)}),
\end{align*}
and we just need to compute the two underbraced arguments of $F^q$.

\vs

From eq.\eqref{eq:nu_prime_inverse_formula} of Lem.\ref{lem:nu_prime_k_basic_properties}(6), we see that 
$$
(\nu'_{v_{12}})^{-1} (\nu'_{v_7})^{-1}({\bf X}^{(2)}_7)
 = (\nu'_{v_{12}})^{-1} (({\bf X}^{(3)}_7)^{-1})
 = ([{\bf X}'_7 ({\bf X}'_{12})^{[\varepsilon^{(3)}_{7,12}]_+}]_{\rm Weyl})^{-1}
 = ({\bf X}'_7)^{-1},
$$
where in the last equality we used $\varepsilon^{(3)}_{7,12}=0$ (see eq.\eqref{eq:varepsilon_3}). Observe that
\begin{align*}
\mu^q_{v_7}({\bf X}_{12}^{(3)})
& = \mu^{\sharp q}_{v_7} \mu'_{v_7} ({\bf X}_{12}^{(3)}) \qquad (\because \mbox{Def.\ref{def:FG_quantum_mutation}}) \\
& = \mu^{\sharp q}_{v_7} ({\bf X}_{12}^{(2)}) \qquad\quad\,\, (\because \mbox{eq.\eqref{eq:mu_prime_formula}, $\varepsilon^{(2)}_{7,12}=0$ (see eq.\eqref{eq:varepsilon_2})}) \\
& = {\bf X}_{12}^{(2)} \cdot F^q({\bf X}_7^{(2)};\varepsilon^{(2)}_{7,12}) \quad (\because \mbox{eq.\eqref{eq:mu_sharp_formula}}) \\
& = {\bf X}_{12}^{(2)} \qquad (\because \mbox{$\varepsilon^{(2)}_{7,12}=0$, Def.\ref{def:F_q}}),
\end{align*}
and that
\begin{align*}
(\nu'_{v_{12}})^{-1} (\nu'_{v_7})^{-1}({\bf X}^{(2)}_{12})
& = (\nu'_{v_{12}})^{-1}({\bf X}^{(3)}_{12}) \quad (\because \mbox{eq.\eqref{eq:nu_prime_inverse_formula}, $\varepsilon^{(2)}_{12,7}=0$ (see eq.\eqref{eq:varepsilon_2})}) \\
& = ({\bf X}'_{12})^{-1} \quad (\because \mbox{eq.\eqref{eq:nu_prime_inverse_formula}}),
\end{align*}
which finishes the computation of the two underbraced parts, hence the desired statement. \qed

\vs

\begin{remark}
Lem.\ref{lem:section5_modified} is all we need for our purposes, but one could also compute $\nu^\omega_{v_7} \nu^\omega_{v_{12}}(\exp({\textstyle \sum}_{j=1}^{12} a'_j {\bf x}_j'))$ explicitly by similar computations, and obtain
\begin{align}
\nonumber
\nu^\omega_{v_7} \nu^\omega_{v_{12}} ( \exp({\textstyle \sum}_{j=1}^{12} a_j' {\bf x}_j') )
=(\exp( {\textstyle \sum}_{j=1}^{12} a_j^{(2)} {\bf x}_j^{(2)}) ) \, F^q({\bf X}_7^{(2)};\alpha^{(3)}) \, F^q({\bf X}_{12}^{(2)};\alpha^{(4)}),
\end{align}
where $(a_j^{(2)})_{j=1}^{12} \in (\frac{1}{3}\mathbb{Z})^{12}$ is given by
$$
a_7^{(2)} = -a_7' + a_3' + a_6', \quad
a_{12}^{(2)} = -a_{12}' + a_4' + a_9', \quad
a_j^{(2)} = a_j', \quad \forall j\neq 7,12,
$$
(see Lem.\ref{lem:alpha} and eq.\eqref{eq:a_2_as_a_prime}), and $\alpha^{(3)}$ and $\alpha^{(4)}$ are the same as in Lem.\ref{lem:section5_modified}. This result is what is used in a previous version of the present paper, which we now replaced by Lem.\ref{lem:section5_modified}, which lessens the subsequent computations.
\end{remark}

Now we apply the above lemmas. In each case of $W$ and each state $s$ of $W$, Step 2 of the previous subsection requires us to check that $\nu^\omega_{v_7} \nu^\omega_{v_{12}} ({\rm Tr}^\omega_{\Delta';\mathcal{Q}}([W,s]))$ is Laurent for $\Delta^{(2)}$ and is multiplicity-free. Due to Lem.\ref{lem:OK_to_check_in_Delta_prime}, it suffices to check that $(\nu'_{v_{12}})^{-1}(\nu'_{v_7})^{-1} \nu^\omega_{v_7} \nu^\omega_{v_{12}}({\rm Tr}^\omega_{\Delta';\mathcal{Q}}([W,s]))$ is Laurent for $\Delta'$ and is multiplicity-free. We perform this checking in the following subsection, with the help of Lem.\ref{lem:section5_modified}.

\subsection{Checking the quantum Laurent property and the multiplicity-freeness}
\label{subsec:checking_quantum_Laurent}

In each of the four cases of $W$ in Fig.\ref{fig:four_cases}, and each state $s$, we will now check that $(\nu'_{v_{12}})^{-1}(\nu'_{v_7})^{-1}\nu^\omega_{v_7} \nu^\omega_{v_{12}}({\rm Tr}^\omega_{\Delta';\mathcal{Q}}([W,s]))$ is Laurent for $\Delta'$ and is multiplicity-free. We will first compute ${\rm Tr}^\omega_{\Delta';\mathcal{Q}}([W,s])$ and observe that it is a multiplicity-free sum of Weyl-ordered Laurent monomials $\exp(\sum_{j=1}^{12} a'_j {\bf x}'_j)$ for some $(a'_j)_{j=1}^{12} \in (\frac{1}{3}\mathbb{Z})^{12}$. We know from Prop.\ref{prop:value_of_the_SL3_quantum_trace_is_Delta-balanced} that these $(a'_j)_{j=1}^{12}$ are $\Delta'$-balanced in the sense of Def.\ref{def:Delta-balanced_elements}. Hence we can apply Lem.\ref{lem:section5_modified}. As seen in eq.\eqref{eq:result_of_nu_omega_v_7_and_v_12_modified}, the application of $(\nu'_{v_{12}})^{-1}(\nu'_{v_7})^{-1} \nu^\omega_{v_7} \nu^\omega_{v_{12}}$ to $\exp(\sum_{j=1}^{12} a'_j {\bf x}'_j)$ has the effect of gaining the factor $F^q( {{\bf X}'}_{\hspace{-1mm}7}^{-1}; \alpha^{(3)}) \, F^q( {{\bf X}'}_{\hspace{-1mm}12}^{-1}; \alpha^{(4)})$. When this factor is merely rational but not Laurent, checking the Laurentness for Step 2 is non-trivial. When ${\rm Tr}^\omega_{\Delta';\mathcal{Q}}([W,s])$ involves some ${\bf X}_7'$ and ${\bf X}_{12}'$, checking the mulitplicity-freeness for Step 2 is nontrivial, even though ${\rm Tr}^\omega_{\Delta';\mathcal{Q}}([W,s])$ is multiplicity-free. For actual checking of Step 2, for each case of $W$ we will make a table of values of $a'_j$, $\alpha^{(3)}$ and $\alpha^{(4)}$.

\vs

Suppose that $s$ assigns the state values $\varepsilon_1,\varepsilon_2 \in \{1,2,3\}$ to the initial and the terminal endpoints of $W$. We would sometimes write $s$ as the pair $(\varepsilon_1,\varepsilon_2)$. We write ${\rm Tr}^\omega_{\Delta'}$ for ${\rm Tr}^\omega_{\Delta';\mathcal{Q}}$.

\vs

\ul{Case $\circled{1}$}. This case falls into Thm.\ref{thm:SL3_quantum_trace}(QT2-1), hence one can observe that ${\rm Tr}^\omega_{\Delta'}([W,s])$ is the image under the map ${\rm Wl}^\omega_{\Delta'}$ (Def.\ref{def:cl_and_Wl}) of the $(\varepsilon_1,\varepsilon_2)$-th entry of the classical matrix
\begin{align*}
& \smallmatthree{{Z_2'}{Z_1'}^2}{0}{0}{0}{{Z_2'}{Z_1'}^{-1}}{0}{0}{0}{{Z_2'}^{-2}{Z_1'}^{-1}}
\smallmatthree{{Z_3'}^2}{{Z_3'}^2+{Z_3'}^{-1}}{{Z_3'}^{-1}}{0}{{Z_3'}^{-1}}{{Z_3'}^{-1}}{0}{0}{{Z_3'}^{-1}}
\smallmatthree{{Z_8'} {Z_9'}^2}{0}{0}{0}{{Z_8'} {Z_9'}^{-1}}{0}{0}{0}{{Z_8'}^{-2} {Z_9'}^{-1}} \\
& = \smallmatthree{{Z_2'}{Z_1'}^2 {Z_3'}^2 {Z_8'} {Z_9'}^2}{~{Z_2'}{Z_1'}^2 {Z_3'}^2 {Z_8'}{Z_9'}^{-1} + {Z_2'}{Z_1'}^2 {Z_3'}^{-1}{Z_8'}{Z_9'}^{-1}~}{{Z_2'}{Z_1'}^2{Z_3'}^{-1}{Z_8'}^{-2}{Z_9'}^{-1}}{0}{{Z_2'}{Z_1'}^{-1}{Z_3'}^{-1}{Z_8'}{Z_9'}^{-1}}{{Z_2'}{Z_1'}^{-1}{Z_3'}^{-1}{Z_8'}^{-2}{Z_9'}^{-1}}{0}{0}{{Z_2'}^{-2}{Z_1'}^{-1}{Z_3'}^{-1}{Z_8'}^{-2}{Z_9'}^{-1}}
\end{align*}
By inspection, ${\rm Tr}^\omega_{\Delta'}([W,s])$ is a Weyl-ordered multiplicity-free Laurent polynomial for $\Delta'$, with each summand Laurent monomial term being of the form 
$$
{\rm Wl}^\omega_\Delta( {X_2'}^{a_2'} {X_1'}^{a_1'} {X_3'}^{a_3'} {X_8'}^{a_8'} {X_9'}^{a_9'} )
= \exp( {\textstyle \sum}_{j=1}^{12} a_j' {\bf x}_j')
$$
with $(a_j')_{j=1}^{12} \in (\frac{1}{3}\mathbb{Z})^{12}$ and $a_j'=0$ if $j\not\in\{2,1,3,8,9\}$. By Prop.\ref{prop:value_of_the_SL3_quantum_trace_is_Delta-balanced}, $(a_j')_{j=1}^{12} \in (\frac{1}{3}\mathbb{Z})^{12}$ is $\Delta'$-balanced, which one can also directly check easily in this case. For each of these Laurent monomials, we record $a_1',a_3',a_9'$, $\alpha^{(3)} = a_1' + a_4' - a_3' - a_6' = a_1' - a_3'$ and $\alpha^{(4)} = a_3' + a_{10}' - a_4' - a_9' = a_3' - a_9'$ in the following table; in the first row, $(\varepsilon_1,\varepsilon_2)$ stands for the $(\varepsilon_1,\varepsilon_2)$-th entry, and $(\varepsilon_1,\varepsilon_2)_k$ for the $k$-th term of the $(\varepsilon_1,\varepsilon_2)$-th entry.

\begin{center}
\begin{tabular}{c|c|c|c|c|c|c|c}
& $(1,1)$ & $(1,2)_1$ & $(1,2)_2$ & $(1,3)$ & $(2,2)$ & $(2,3)$ & $(3,3)$ \\ \hline
$a_1'$ & $2/3$ & $2/3$ & $2/3$ & $2/3$ & $-1/3$ & $-1/3$ & $-1/3$ \\
$a_3'$ & $2/3$ & $2/3$ & $-1/3$ & $-1/3$ & $-1/3$ & $-1/3$ & $-1/3$ \\
$a_9'$ & $2/3$ & $-1/3$ & $-1/3$ & $-1/3$ & $-1/3$ & $-1/3$ & $-1/3$ \\
$\alpha^{(3)}$ & 0 & 0 & 1 & 1 & 0 & 0 & 0  \\
$\alpha^{(4)}$ & 0 & 1 & 0 & 0 & 0 & 0 & 0
\end{tabular}
\end{center}

Therefore, for each column, since $\alpha^{(3)},\alpha^{(4)} \in \{0,1\}$, we see that $F^q( {{\bf X}'}_{\hspace{-1mm}7}^{-1}; \alpha^{(3)}) \, F^q( {{\bf X}'}_{\hspace{-1mm}12}^{-1}; \alpha^{(4)})$ is a multiplicity-free Laurent polynomial in the variables ${\bf X}_7'$ and ${\bf X}_{12}'$, in view of eq.\eqref{eq:F_q}. As ${\rm Tr}^\omega_{\Delta'}([W,s])$ does not involve ${\bf X}_7'$ or ${\bf X}_{12}'$ for any $s$ in this case, using eq.\eqref{eq:result_of_nu_omega_v_7_and_v_12_modified} of Lem.\ref{lem:section5_modified} we see that $(\nu'_{v_{12}})^{-1}(\nu'_{v_7})^{-1}\nu^\omega_{v_7} \nu^\omega_{v_{12}}({\rm Tr}^\omega_{\Delta'}([W,s]))$ is Laurent for $\Delta'$ and is multiplicity-free, as desired.

\vs

\ul{Case $\circled{2}$}. This case falls into Thm.\ref{thm:SL3_quantum_trace}(QT2-2), hence one can observe that ${\rm Tr}^\omega_{\Delta'}([W,s])$ is the image under the map ${\rm Wl}^\omega_{\Delta'}$ of the $(\varepsilon_1,\varepsilon_2)$-th entry of the classical matrix
$$
\smallmatthree{{Z_9'}{Z_8'}^2}{0}{0}{0}{{Z_9'}{Z_8'}^{-1}}{0}{0}{0}{{Z_9'}^{-2}{Z_8'}^{-1}}
\smallmatthree{{Z_3'}}{0}{0}{{Z_3'}}{{Z_3'}}{0}{{Z_3'}}{{Z_3'}+{Z_3'}^{-2}}{{Z_3'}^{-2}}
\smallmatthree{{Z_1'} {Z_2'}^2}{0}{0}{0}{{Z_1'} {Z_2'}^{-1}}{0}{0}{0}{{Z_1'}^{-2} {Z_2'}^{-1}}
$$
So, by inspection, ${\rm Tr}^\omega_{\Delta'}([W,s])$ is a Weyl-ordered multiplicity-free Laurent polynomial for $\Delta'$, with each summand Laurent monomial being of the form $\exp(\sum_{j=1}^{12} a_j' {\bf x}_j')$ with $(a_j')_{j=1}^{12} \in (\frac{1}{3}\mathbb{Z})^{12}$, where $a_j' =0$ when $j\not\in\{9,8,3,1,2\}$. For nonzero Laurent monomials, we record $a_9',a_3',a_1'$, $\alpha^{(3)} = a'_1-a'_3$ and $\alpha^{(4)} = a'_3 - a'_9$.

\begin{center}
\begin{tabular}{c|c|c|c|c|c|c|c}
& $(1,1)$ & $(2,1)$ & $(2,2)$ & $(3,1)$ & $(3,2)_1$ & $(3,2)_2$ & $(3,3)$ \\ \hline
$a_9'$ & $1/3$ & $1/3$ & $1/3$ & $-2/3$ & $-2/3$ & $-2/3$ & $-2/3$ \\
$a_3'$ & $1/3$ & $1/3$ & $1/3$ & $1/3$ & $1/3$ & $-2/3$ & $-2/3$ \\
$a_1'$ & $1/3$ & $1/3$ & $1/3$ & $1/3$ & $1/3$ & $1/3$ & $-2/3$ \\
$\alpha^{(3)}$ & 0 & 0 & 0 & 0 & 0 & 1 & 0 \\
$\alpha^{(4)}$ & 0 & 0 & 0 & 1 & 1 & 0 & 0 
\end{tabular}
\end{center}
Again, since $\alpha^{(3)},\alpha^{(4)} \in \{0,1\}$ and ${\rm Tr}^\omega_{\Delta'}([W,s])$ does not involve ${\bf X}_7'$ or ${\bf X}_{12}'$, using eq.\eqref{eq:result_of_nu_omega_v_7_and_v_12_modified} of Lem.\ref{lem:section5_modified} we see that $(\nu'_{v_{12}})^{-1}(\nu'_{v_7})^{-1}\nu^\omega_{v_7} \nu^\omega_{v_{12}}({\rm Tr}^\omega_{\Delta'}([W,s]))$ is Laurent for $\Delta'$ and is multiplicity-free, as desired.

\vs

\ul{Case $\circled{3}$}. We claim that ${\rm Tr}^\omega_{\Delta'}([W,s])$ is the image under the map ${\rm Wl}^\omega_{\Delta'}$ of the $(\varepsilon_1,\varepsilon_2)$-th entry of the classical matrix
\begin{align*}
& \hspace{-7mm} \smallmatthree{{Z_6'}{Z_5'}^2}{0}{0}{0}{{Z_6'}{Z_5'}^{\hspace{-0,5mm}-1}}{0}{0}{0}{{Z_6'}^{-2}{Z_5'}^{\hspace{-0,5mm}-1}}
\myul{ \smallmatthree{{Z_4'}^2}{{Z_4'}^2+{Z_4'}^{\hspace{-0,5mm}-1}}{{Z_4'}^{\hspace{-0,5mm}-1}}{0}{{Z_4'}^{\hspace{-0,5mm}-1}}{{Z_4'}^{\hspace{-0,5mm}-1}}{0}{0}{{Z_4'}^{\hspace{-0,5mm}-1}}
\smallmatthree{{Z_7'} {Z_{12}'}^{\hspace{-1,5mm} 2}}{0}{0}{0}{{Z_7'} {Z_{12}'}^{\hspace{-1,5mm}-1}}{0}{0}{0}{{Z_7'}^{\hspace{-0,5mm}-2} {Z_{12}'}^{\hspace{-1,5mm}-1}}
\smallmatthree{{Z_3'}^2}{{Z_3'}^2+{Z_3'}^{\hspace{-0,5mm}-1}}{{Z_3'}^{\hspace{-0,5mm}-1}}{0}{{Z_3'}^{\hspace{-0,5mm}-1}}{{Z_3'}^{\hspace{-0,5mm}-1}}{0}{0}{{Z_3'}^{\hspace{-0,5mm}-1}} }
\smallmatthree{{Z_1'}{Z_2'}^2}{0}{0}{0}{{Z_1'}{Z_2'}^{\hspace{-0,5mm}-1}}{0}{0}{0}{{Z_1'}^{\hspace{-0,5mm}-2} {Z_2'}^{\hspace{-0,5mm}-1}}.
\end{align*}
First let's verify this at the classical level. Conceptually, as mentioned in the proof of Prop.\ref{prop:classical_compatibility}, the `${\rm SL}_3$ classical trace' ${\rm Tr}^1_{\Delta'}([W,\cdot])$ corresponds to the monodromy matrix ${\rm M}_{W;\Delta'}$ associated to $W$, or more precisely to the projection $\pi(W)$ of $W$ in the surface $\frak{S}$, where the monodromy matrix is given by the product of basic monodromy matrices given in (MM1)--(MM3) of \cite[\S4.2]{Kim}. The state $s$ of ${\rm Tr}^1_{\Delta'}([W,s])$ indicates which entry we read from the matrix ${\rm M}_{W;\Delta'}$. The above expression presents the product of five basic monodromy matrices. We first divide $W$ into five small pieces, as done in \cite[\S4.2]{Kim}. The first factor is the `edge matrix' (MM1) associated to the initial piece of $W$ passing through the lower-left-side edge of $\mathcal{Q}$ having $v_5$ and $v_6$. The second factor is the `left-turn matrix' (MM2) associated to the second piece of $W$ in the interior of the lower triangle of $\Delta'$ having $v_4$. The third factor is the edge matrix (MM1) for the piece of $W$ passing through the diagonal edge of $\Delta'$. The fourth factor is the left-turn matrix (MM2) for the next piece of $W$ in the interior of the upper triangle of $\Delta'$. The fifth factor is the edge matrix (MM1) for the last piece of $W$ passing through the upper-left-side edge of $\mathcal{Q}$ having $v_1$ and $v_2$. A proof that the values of the ${\rm SL}_3$ classical trace are indeed the entries of the monodromy matrices can be obtained by using Thm.\ref{thm:SL3_quantum_trace} in the case $\omega^{1/2}=1$, especially by the cutting/gluing axiom (QT1) which justifies the interpretation in terms of the product of monodromy matrices for parts of $W$ living in the two triangles of $\Delta'$. A precise argument can be found in the proof of \cite[Prop.5.69]{Kim}. 

\vs

Before proving the analogous statement for the quantum case ${\rm Tr}^\omega_{\Delta'}([W,s])$, note that the product of the underlined middle three matrices in the above product expression is
\begin{align*}
\smallmatthree{%
{Z_4'}^{\hspace{-0,5mm}2} {Z_7'} {Z_{12}'}^{\hspace{-0,5mm}2} {Z_3'}^{\hspace{-0,5mm}2}}%
{~{Z_4'}^{\hspace{-0,5mm}2} {Z_7'}{Z_{12}'}^{\hspace{-1,5mm}2}({Z_3'}^{\hspace{-0,5mm}2}+{Z_3'}^{\hspace{-0,5mm}-1}) 
+ ({Z_4'}^{\hspace{-0,5mm}2}+{Z_4'}^{\hspace{-0,5mm}-1}) {Z_7'} {Z_{12}'}^{\hspace{-1,5mm}-1} {Z_3'}^{\hspace{-0,5mm}-1}~}%
{({Z_4'}^{\hspace{-0,5mm}2} {Z_7'}{Z_{12}'}^{\hspace{-1,5mm}2}
+ ({Z_4'}^{\hspace{-0,5mm}2}+{Z_4'}^{\hspace{-0,5mm}-1}) {Z_7'} {Z_{12}'}^{\hspace{-1,5mm}-1} + {Z_4'}^{\hspace{-0,5mm}-1} {Z_7'}^{\hspace{-0,5mm}-2} {Z_{12}'}^{\hspace{-1,5mm}-1}) {Z_3'}^{\hspace{-0,5mm}-1}}%
{0}%
{{Z_4'}^{\hspace{-0,5mm}-1} {Z_7'} {Z_{12}'}^{\hspace{-1,5mm}-1} {Z_3'}^{\hspace{-0,5mm}-1}}%
{{Z_4'}^{\hspace{-0,5mm}-1}( {Z_7'} {Z_{12}'}^{\hspace{-1,5mm}-1} + {Z_7'}^{\hspace{-0,5mm}-2} {Z_{12}'}^{\hspace{-1,5mm}-1} ) {Z_3'}^{\hspace{-0,5mm}-1}}%
{0}{0}%
{{Z_4'}^{\hspace{-0,5mm}-1} {Z_7'}^{\hspace{-0,5mm}-2} {Z_{12}'}^{\hspace{-1,5mm}-1} {Z_3'}^{\hspace{-0,5mm}-1}}
\end{align*}
So, by inspection, the ${\rm SL}_3$-classical trace ${\rm Tr}^1_{\Delta'}([W,s])$ is a multiplicity-free Laurent polynomial for $\Delta'$.

\vs

Let's investigate the ${\rm SL}_3$ quantum trace ${\rm Tr}^\omega_{\Delta'}([W,s])$. Cutting the quadrilateral $\mathcal{Q}$ into two ideal triangles $t_1$ and $t_2$ by the internal arc $e$ of $\Delta'$, with $t_1$ being the bottom triangle in Fig.\ref{fig:four_cases}, we apply the cutting/gluing axiom Thm.\ref{thm:SL3_quantum_trace}(QT1). Then
\begin{align}
\label{eq:state_sum_case3}
i_{\Delta',\Delta'_e} {\rm Tr}^\omega_{\Delta'}([W,s]) = \underset{s_1,s_2}{\textstyle \sum} \, {\rm Tr}^\omega_{t_1}([W_1, s_1]) {\rm Tr}^\omega_{t_2}([W_2, s_2]) 
\end{align}
with $W_i := W\cap (t_i \times {\bf I})$ and ${\rm Tr}^\omega_{t_i}$ is with respect to the unique ideal triangulation of $t_i$, where the sum is over all pairs of states $s_1,s_2$ of $W_1$ and $W_2$ that are compatible with $s$ in the sense of Def.\ref{def:cutting_process}. Each of $W_1 \subset t_1\times {\bf I}$ and $W_2 \subset t_2 \times {\bf I}$ falls into the case of Thm.\ref{thm:SL3_quantum_trace}(QT2); for this particular setting of Case \circled{3}, both are left turns, i.e. (QT2-1). Hence from Thm.\ref{thm:SL3_quantum_trace}(QT2) we see that ${\rm Tr}^\omega_{t_i}([W_i, s_i])$ is an entry of a $3\times 3$ matrix with entries being multiplicity-free Laurent polynomials in $\mathcal{Z}^\omega_{t_1}$. Hence, from the state-sum formula in eq.\eqref{eq:state_sum_case3} one can deduce that ${\rm Tr}^\omega_{\Delta'}([W,s])$ is a multiplicity-free Laurent polynomial in $\mathcal{Z}^\omega_{\Delta'}$. In the meantime we have observed by using Prop.\ref{prop:elevation_reversing_and_star} that ${\rm Tr}^\omega_{\Delta'}([W,s])$ is $*$-invariant (see the paragraph following  Prop.\ref{prop:elevation_reversing_and_star}).  Therefore, by Lem.\ref{lem:star-invariance_and_Weyl-ordering_new}(3), ${\rm Tr}^\omega_{\Delta'}([W,s])$ is (term-by-term) Weyl-ordered, and hence by Lem.\ref{lem:classicalization_determines_quantum}, or more precisely by the discussion preceding it, ${\rm Tr}^\omega_{\Delta'}([W,s])$ equals the image under ${\rm Wl}^\omega_{\Delta'}$ of its classicalization, namely ${\rm Tr}^1_{\Delta'}([W,s])$, as claimed.

\vs

Now, by inspection on the classical monodromy matrices for ${\rm Tr}^1_{\Delta'}([W,s])$, one finds out that ${\rm Tr}^\omega_{\Delta'}([W,s])$ is a Weyl-ordered multiplicity-free Laurent polynomial in $\mathcal{Z}^\omega_{\Delta'}$ with each summand Laurent monomial being of the form $\exp(\sum_{j=1}^{12} a_j' {\bf x}'_j)$ with $(a_j')_{j=1}^{12} \in (\frac{1}{3}\mathbb{Z})^{12}$, where $a_j'=0$ when $j\not\in \{6,5,4,7,12,3,1,2\}$. We record the nonzero $a_j'$'s, $\alpha^{(3)} = a'_1+a'_4-a'_3-a'_6$ and $\alpha^{(4)} = a'_3+a'_{10}-a'_4-a'_9 = a'_3 - a'_4$.

\begin{center}
{\setlength{\tabcolsep}{3,5pt} \begin{tabular}{c|c|c|c|c|c|c|c|c|c|c|c|c|c}
& $(1,1)$ & $(1,2)_1$ & $(1,2)_2$ & $(1,2)_3$ & $(1,2)_4$ & $(1,3)_1$ & $(1,3)_2$ & $(1,3)_3$ & $(1,3)_4$ & $(2,2)$ & $(2,3)_1$ & $(2,3)_2$ & $(3,3)$ \\ \hline
$a_6'$ & 1/3 & 1/3 & 1/3 & 1/3 & 1/3 & 1/3 & 1/3 & 1/3 & 1/3  & 1/3 & 1/3 & 1/3 & $-2/3$\\
$a_5'$ & 2/3 & 2/3 & 2/3 & 2/3 & 2/3 & 2/3 & 2/3 & 2/3 & 2/3 & $-1/3$ & $-1/3$ & $-1/3$ & $-1/3$ \\
$a_4'$ & 2/3 & 2/3 & 2/3 & 2/3 & $-1/3$ & 2/3 & 2/3 & $-1/3$ & $-1/3$ & $-1/3$ & $-1/3$ & $-1/3$ & $-1/3$ \\
$a_7'$ & 1/3 & 1/3 & 1/3 & 1/3 & 1/3 & 1/3 & 1/3 & 1/3 & $-2/3$ & 1/3 & 1/3 & $-2/3$ & $-2/3$ \\
$a_{12}'$ & 2/3 & 2/3 & 2/3 & $-1/3$ & $-1/3$ & 2/3 & $-1/3$ & $-1/3$ & $-1/3$ & $-1/3$ & $-1/3$ & $-1/3$ & $-1/3$ \\
$a_3'$ & 2/3 & 2/3 & $-1/3$ & $-1/3$ & $-1/3$ & $-1/3$ &  $-1/3$ &  $-1/3$ &  $-1/3$ &  $-1/3$ &  $-1/3$ &  $-1/3$ &  $-1/3$ \\
$a_1'$ & 1/3 & 1/3 & 1/3 & 1/3 & 1/3 & $-2/3$ & $-2/3$ & $-2/3$ & $-2/3$ & 1/3 & $-2/3$ & $-2/3$  & $-2/3$  \\
$a_2'$ & 2/3 & $-1/3$ & $-1/3$ & $-1/3$ & $-1/3$ & $-1/3$ & $-1/3$ & $-1/3$ & $-1/3$ & $-1/3$ & $-1/3$ & $-1/3$ & $-1/3$  \\
$\alpha^{(3)}$ & 0 & 0 & 1 & 1 & 0 & 0 & 0 & $-1$ & $-1$ & 0 & $-1$ & $-1$ & 0 \\
$\alpha^{(4)}$ & 0 & 0 & $-1$ & $-1$ & 0 & $-1$ & $-1$ & 0 & 0 & 0 & 0  &0 & 0\\
\end{tabular}}
\end{center}
For the entries $(1,1)$, $(2,2)$ and $(3,3)$, we have $\alpha^{(3)}=\alpha^{(4)}=0$, hence $F^q( {{\bf X}'}_{\hspace{-1mm}7}^{-1}; \alpha^{(3)}) \, F^q( {{\bf X}'}_{\hspace{-1mm}12}^{-1}; \alpha^{(4)})=1$, and therefore using eq.\eqref{eq:result_of_nu_omega_v_7_and_v_12_modified} of Lem.\ref{lem:section5_modified} we see that $(\nu'_{v_{12}})^{-1}(\nu'_{v_7})^{-1}\nu^\omega_{v_7} \nu^\omega_{v_{12}}({\rm Tr}^\omega_{\Delta'}([W,s]))$ is Laurent for $\Delta'$ and is multiplicity-free as desired, for the cases when the state $s$ is $(1,1)$, $(2,2)$ or $(3,3)$.

\vs

We have $\alpha^{(4)}=-1$ for $(1,2)_2$ and $(1,2)_3$, which represent the second and the third terms $\exp(\sum_{j=1}^{12} a_j' {\bf x}_j')$ of the $(1,2)$-th entry ${\rm Tr}^\omega_{\Delta'}([W,(1,2)])$; in the corresponding two columns in the above table, the only difference is $a_{12}'$. We compute the sum of these two terms:
\begin{align*}
& \textstyle \exp(\frac{1}{3} {\bf x}'_6 + \frac{2}{3} {\bf x}'_5 + \frac{2}{3} {\bf x}'_4 + \frac{1}{3} {\bf x}'_7 + \frac{2}{3} {\bf x}'_{12} - \frac{1}{3} {\bf x}'_3 + \frac{1}{3} {\bf x}'_1 - \frac{1}{3} {\bf x}'_2)  \\
& \textstyle \quad + \exp(\frac{1}{3} {\bf x}'_6 + \frac{2}{3} {\bf x}'_5 + \frac{2}{3} {\bf x}'_4 + \frac{1}{3} {\bf x}'_7 - \frac{1}{3} {\bf x}'_{12} - \frac{1}{3} {\bf x}'_3 + \frac{1}{3} {\bf x}'_1 - \frac{1}{3} {\bf x}'_2)  \\
& \textstyle = \exp(\frac{1}{3} {\bf x}'_6 + \frac{2}{3} {\bf x}'_5 + \frac{2}{3} {\bf x}'_4 + \frac{1}{3} {\bf x}'_7 + \frac{2}{3} {\bf x}'_{12} - \frac{1}{3} {\bf x}'_3 + \frac{1}{3} {\bf x}'_1 - \frac{1}{3} {\bf x}'_2) (1 + e^{-{\bf c}} \exp(-{\bf x}'_{12})),
\end{align*}
where we used the BCH formula (Lem.\ref{lem:BCH}), and
\begin{align*}
& \textstyle [\frac{1}{3} {\bf x}'_6 + \frac{2}{3} {\bf x}'_5 + \frac{2}{3} {\bf x}'_4 + \frac{1}{3} {\bf x}'_7 + \frac{2}{3} {\bf x}'_{12} - \frac{1}{3} {\bf x}'_3 + \frac{1}{3} {\bf x}'_1 - \frac{1}{3} {\bf x}'_2, \, - {\bf x}'_{12}] \\
& \textstyle = 2( - \frac{1}{3} \varepsilon'_{6,12} - \frac{2}{3} \varepsilon'_{5,12} - \frac{2}{3} \varepsilon'_{4,12} - \frac{1}{3} \varepsilon'_{7,12} - \frac{2}{3} \varepsilon'_{12,12} + \frac{1}{3} \varepsilon'_{3,12} - \frac{1}{3} \varepsilon'_{1,12} + \frac{1}{3} \varepsilon'_{2,12} ) {\bf c}  \qquad (\because \mbox{eq.\eqref{eq:bracket_in_x}}) \\
& \textstyle = 2 ( - \frac{1}{3} \cdot 0 - \frac{2}{3} \cdot 0 - \frac{2}{3} \cdot (-1) - \frac{1}{3} \cdot 0 - \frac{2}{3} \cdot 0 + \frac{1}{3} \cdot 1 - \frac{1}{3} \cdot 0 + \frac{1}{3} \cdot 0 ) {\bf c} \qquad (\because\mbox{read from Fig.\ref{fig:mutations_for_a_flip}}) \\
& = 2 {\bf c}
\end{align*}
In more detail, if we denote $\frac{1}{3} {\bf x}'_6 + \frac{2}{3} {\bf x}'_5 + \frac{2}{3} {\bf x}'_4 + \frac{1}{3} {\bf x}'_7 + \frac{2}{3} {\bf x}'_{12} - \frac{1}{3} {\bf x}'_3 + \frac{1}{3} {\bf x}'_1 - \frac{1}{3} {\bf x}'_2$ by ${\bf x}'$, then what is just done is $e^{{\bf x}'} + e^{{\bf x}' - {\bf x}'_{12}} \stackrel{{\rm Lem}.\ref{lem:BCH}}{=} e^{{\bf x}'} + e^{-\frac{1}{2}[{\bf x}', -{\bf x}'_{12}]} e^{{\bf x}'} e^{-{\bf x}'_{12}} = e^{{\bf x}'} + e^{-{\bf c}} e^{{\bf x}'} e^{-{\bf x}'_{12}}= e^{{\bf x}'}(1+e^{-{\bf c}} e^{-{\bf x}'_{12}})$. Note that
$$
1 + e^{-{\bf c}} \exp(-{\bf x}_{12}')
= 1 + q^{-1} {{\bf X}'}_{\hspace{-1mm}12}^{-1} = ( F^q( {{\bf X}'}_{\hspace{-1mm}12}^{-1}; -1) )^{-1} \qquad (\because \mbox{eq.\eqref{eq:F_q}})
$$
This means that
$$
(\mbox{the sum of the second and third terms of ${\rm Tr}^\omega_{\Delta'}([W,(1,2)])$}) \cdot F^q({{\bf X}'}_{\hspace{-1mm}12}^{-1};-1)
$$
is a Laurent monomial; it equals $e^{{\bf x}'}$ in our notation. So, in view of eq.\eqref{eq:result_of_nu_omega_v_7_and_v_12_modified} of Lem.\ref{lem:section5_modified}, and since $F^q( {{\bf X}'}_{\hspace{-1mm}7}^{-1}; \alpha^{(3)})= F^q( {{\bf X}'}_{\hspace{-1mm}7}^{-1}; 1) = 1+q {{\bf X}'}_{\hspace{-1mm}7}^{-1}$, it follows that $(\nu'_{v_{12}})^{-1}(\nu'_{v_7})^{-1}\nu^\omega_{v_7} \nu^\omega_{v_{12}}({\rm Tr}^\omega_{\Delta'}([W,(1,2)]))$ is Laurent for $\Delta'$. We note that the multiplicity-free-ness is not immediate, although we know that ${\rm Tr}^\omega_{\Delta'}([W,(1,2)])$ is multiplicity-free; this is because of the factors $F^q({\bf X}_7^{(2)};\alpha^{(3)})$ and $F^q({\bf X}_{12}^{(2)};\alpha^{(4)})$. For simplicity, let's denote by $(1,2)_r$ the $r$-th Laurent monomial term of ${\rm Tr}^\omega_{\Delta'}([W,(1,2)])$, so that ${\rm Tr}^\omega_{\Delta'}([W,(1,2)]) = \sum_{r=1}^4 (1,2)_r$. In particular, we have $(1,2)_2 = e^{{\bf x}'}$ and $(1,2)_3 = e^{{\bf x}'- {\bf x}'_{12}}$ in our notation. From eq.\eqref{eq:result_of_nu_omega_v_7_and_v_12_modified} and our computation above, we have
$$
(\nu'_{v_{12}})^{-1}(\nu'_{v_7})^{-1}\nu^\omega_{v_7} \nu^\omega_{v_{12}}({\rm Tr}^\omega_{\Delta'}([W,(1,2)]))
= (1,2)_1 + (1,2)_2 (1+q{{\bf X}'}_{\hspace{-1mm}7}^{-1}) + (1,2)_4
$$
So, to check if this is multiplicity-free, we should check in the above table of values of $a'_j$ that no two of the columns labeled by $(1,2)_1$, $(1,2)_2$, $(1,2)_4$ and a new column obtained from the $(1,2)_2$ column by shifting the value of $a'_7$ by $-1$ (hence $a'_7 = -2/3$ for this new column), which respectively represent the terms $(1,2)_1$, $(1,2)_2$, $(1,2)_4$ and $(1,2)_2 \cdot q{{\bf X}'}^{-1}_{\hspace{-1mm}7}$, have identical values of $a'_j$. This can be done by inspection, and hence we can conclude that $(\nu'_{v_{12}})^{-1}(\nu'_{v_7})^{-1}\nu^\omega_{v_7} \nu^\omega_{v_{12}}({\rm Tr}^\omega_{\Delta'}([W,(1,2)]))$ is Laurent for $\Delta'$ and is multiplicity-free, as desired.

\vs

We do likewise for $(1,3)_1$ and $(1,3)_2$, where $\alpha^{(3)}=0$ and $\alpha^{(4)}=-1$; the two columns in the above table differ only at $a_{12}'$. The sum of the corresponding two terms $\exp(\sum_{j=1}^{12} a'_j {\bf x}'_2)$ is
\begin{align*}
& \textstyle \exp(\frac{1}{3} {\bf x}'_6 + \frac{2}{3} {\bf x}'_5 + \frac{2}{3} {\bf x}'_4 + \frac{1}{3} {\bf x}'_7 + \frac{2}{3} {\bf x}'_{12} - \frac{1}{3} {\bf x}'_3 - \frac{2}{3} {\bf x}'_1 - \frac{1}{3} {\bf x}'_2) \\
& \textstyle \quad + \exp(\frac{1}{3} {\bf x}'_6 + \frac{2}{3} {\bf x}'_5 + \frac{2}{3} {\bf x}'_4 + \frac{1}{3} {\bf x}'_7 - \frac{1}{3} {\bf x}'_{12} - \frac{1}{3} {\bf x}'_3 - \frac{2}{3} {\bf x}'_1 - \frac{1}{3} {\bf x}'_2) \\
& \textstyle = \exp(\frac{1}{3} {\bf x}'_6 + \frac{2}{3} {\bf x}'_5 + \frac{2}{3} {\bf x}'_4 + \frac{1}{3} {\bf x}'_7 + \frac{2}{3} {\bf x}'_{12} - \frac{1}{3} {\bf x}'_3 - \frac{2}{3} {\bf x}'_1 - \frac{1}{3} {\bf x}'_2) (1 + e^{-{\bf c}} \exp(-{\bf x}'_{12}))
\end{align*}
by a similar computation as before; note especially that the difference between this situation and that for $(1,2)_2$ and $(1,2)_3$ is just the coefficient of ${\bf x}'_1$. So
$$
(\mbox{the sum of the first and the second terms of ${\rm Tr}^\omega_{\Delta'}([W,(1,3)])$}) \cdot F^q({{\bf X}'}_{\hspace{-1mm}12}^{-1};-1)
$$
is a Laurent monomial, which is same as the Laurent monomial term for $(1,3)_1$. Let's now investigate $(1,3)_3$ and $(1,3)_4$, where $\alpha^{(3)}=-1$ and $\alpha^{(4)}=0$. The only difference of these two columns is at $a_7'$. The sum of the corresponding two terms $\exp(\sum_{j=1}^{12} a_j' {\bf x}_j')$ is
\begin{align*}
& \textstyle \exp(\frac{1}{3} {\bf x}'_6 + \frac{2}{3} {\bf x}'_5 - \frac{1}{3} {\bf x}'_4 + \frac{1}{3} {\bf x}'_7 - \frac{1}{3} {\bf x}'_{12} - \frac{1}{3} {\bf x}'_3 - \frac{2}{3} {\bf x}'_1 - \frac{1}{3} {\bf x}'_2) \\
& \textstyle \quad + \exp(\frac{1}{3} {\bf x}'_6 + \frac{2}{3} {\bf x}'_5 - \frac{1}{3} {\bf x}'_4 - \frac{2}{3} {\bf x}'_7 - \frac{1}{3} {\bf x}'_{12} - \frac{1}{3} {\bf x}'_3 - \frac{2}{3} {\bf x}'_1 - \frac{1}{3} {\bf x}'_2) \\
& \textstyle = \exp(\frac{1}{3} {\bf x}'_6 + \frac{2}{3} {\bf x}'_5 - \frac{1}{3} {\bf x}'_4 + \frac{1}{3} {\bf x}'_7 - \frac{1}{3} {\bf x}'_{12} - \frac{1}{3} {\bf x}'_3 - \frac{2}{3} {\bf x}'_1 - \frac{1}{3} {\bf x}'_2) (1 + e^{-{\bf c}} \exp(-{\bf x}'_7)),
\end{align*}
where we used the BCH formula (Lem.\ref{lem:BCH}) and
\begin{align*}
& \textstyle [\frac{1}{3} {\bf x}'_6 + \frac{2}{3} {\bf x}'_5 - \frac{1}{3} {\bf x}'_4 + \frac{1}{3} {\bf x}'_7 - \frac{1}{3} {\bf x}'_{12} - \frac{1}{3} {\bf x}'_3 - \frac{2}{3} {\bf x}'_1 - \frac{1}{3} {\bf x}'_2, \, - {\bf x}'_7] \\
& \textstyle = 2 ( - \frac{1}{3} \varepsilon'_{6,7} - \frac{2}{3} \varepsilon'_{5,7} + \frac{1}{3} \varepsilon'_{4,7} - \frac{1}{3} \varepsilon'_{7,7} + \frac{1}{3} \varepsilon'_{12,7} + \frac{1}{3} \varepsilon'_{3,7} + \frac{2}{3} \varepsilon'_{1,7} + \frac{1}{3} \varepsilon'_{2,7}) {\bf c} \qquad (\because \mbox{eq.\eqref{eq:bracket_in_x}}) \\
& \textstyle = 2 ( - \frac{1}{3} \cdot (-1) - \frac{2}{3} \cdot 0 + \frac{1}{3} \cdot 1 - \frac{1}{3} \cdot 0 + \frac{1}{3} \cdot 0 + \frac{1}{3} \cdot (-1) + \frac{2}{3} \cdot 1 + \frac{1}{3} \cdot 0) {\bf c} \qquad (\because\mbox{read from Fig.\ref{fig:mutations_for_a_flip}}) \\
& = 2 {\bf c}
\end{align*}
Since $1+ e^{-{\bf c}} \exp(-{\bf x}'_7) = 1 + q^{-1} {{\bf X}'}_{\hspace{-1mm}7}^{-1} = (F^q({{\bf X}'}_{\hspace{-1mm}7}^{-1};-1))^{-1}$ it follows that 
$$
(\mbox{the sum of the third and the fourth terms of ${\rm Tr}^\omega_{\Delta'}([W,(1,3)])$}) \cdot F^q( {{\bf X}'}_{\hspace{-1mm}7}^{-1}; -1)
$$
is a Laurent monomial, which is the Laurent monomial term for $(1,3)_3$. Combining, in view of eq.\eqref{eq:result_of_nu_omega_v_7_and_v_12_modified} of Lem.\ref{lem:section5_modified}, we conclude that $(\nu'_{v_{12}})^{-1}(\nu'_{v_7})^{-1}\nu^\omega_{v_7} \nu^\omega_{v_{12}}({\rm Tr}^\omega_{\Delta'}([W,(1,3)]))$ is a Laurent polynomial, having two terms, which are the Laurent monomial terms for $(1,3)_1$ and $(1,3)_3$ of ${\rm Tr}^\omega_{\Delta'}([W,(1,3)])$. By inspection on the columns for $(1,3)_1$ and $(1,3)_3$ in the above table for $a'_j$, we see that $(\nu'_{v_{12}})^{-1}(\nu'_{v_7})^{-1}\nu^\omega_{v_7} \nu^\omega_{v_{12}}({\rm Tr}^\omega_{\Delta'}([W,(1,3)]))$ is multiplicity-free, as desired.

\vs

Lastly, we investigate $(2,3)_1$ and $(2,3)_2$, where $\alpha^{(3)}=-1$ and $\alpha^{(4)}=0$; the only difference of the two columns in the table is at $a_7'$. The sum of the corresponding two terms $\exp(\sum_{j=1}^{12} a'_j {\bf x}'_j)$ is
\begin{align*}
& \textstyle \exp( \frac{1}{3} {\bf x}'_6 - \frac{1}{3} {\bf x}'_5 - \frac{1}{3} {\bf x}'_4 + \frac{1}{3} {\bf x}'_7 - \frac{1}{3} {\bf x}'_{12} - \frac{1}{3} {\bf x}'_3 - \frac{2}{3} {\bf x}'_1 - \frac{1}{3} {\bf x}'_2) \\
& \textstyle \quad + \exp( \frac{1}{3} {\bf x}'_6 - \frac{1}{3} {\bf x}'_5 - \frac{1}{3} {\bf x}'_4 - \frac{2}{3} {\bf x}'_7 - \frac{1}{3} {\bf x}'_{12} - \frac{1}{3} {\bf x}'_3 - \frac{2}{3} {\bf x}'_1 - \frac{1}{3} {\bf x}'_2) \\
& \textstyle = \exp ( \frac{1}{3} {\bf x}'_6 - \frac{1}{3} {\bf x}'_5 - \frac{1}{3} {\bf x}'_4 + \frac{1}{3} {\bf x}'_7 - \frac{1}{3} {\bf x}'_{12} - \frac{1}{3} {\bf x}'_3 - \frac{2}{3} {\bf x}'_1 - \frac{1}{3} {\bf x}'_2) (1 + e^{-{\bf c}} \exp(-{\bf x}'_7)),
\end{align*}
by a similar computation as in the situation for $(1,3)_3$ and $(1,3)_4$, which differs from the current situation just at the coefficient of ${\bf x}'_5$. So $(\nu'_{v_{12}})^{-1}(\nu'_{v_7})^{-1}\nu^\omega_{v_7} \nu^\omega_{v_{12}}({\rm Tr}^\omega_{\Delta'}([W,(2,3)]))$ is a Laurent monomial, hence is a multiplicity-free Laurent polynomial, as desired.

\vs

\ul{Case $\circled{4}$}. The flow of logic will be similar. By a similar reasoning as in Case \circled{3}, one can note that ${\rm Tr}^\omega_{\Delta'}([W,s])$ is the image under ${\rm Wl}^\omega_{\Delta'}$ of the $(\varepsilon_1,\varepsilon_2)$-th entry of the classical matrix
\begin{align*}
& \hspace{-7mm} \smallmatthree{{Z_2'}{Z_1'}^{\hspace{-0,4mm}2}}{0}{0}{0}{{Z_2'}{Z_1'}^{\hspace{-0,4mm}-1}}{0}{0}{0}{{Z_2'}^{\hspace{-0,5mm} -2}{Z_1'}^{\hspace{-0,5mm} -1}}
\myul{ \smallmatthree{{Z_3'}}{0}{0}{{Z_3'}}{{Z_3'}}{0}{{Z_3'}}{{Z_3'}+{Z_3'}^{\hspace{-0,5mm} -2}}{{Z_3'}^{\hspace{-0,5mm} -2}}
\smallmatthree{{Z_{12}'} {Z_7'}^{\hspace{-0,5mm} 2}}{0}{0}{0}{{Z_{12}'} {Z_7'}^{\hspace{-0,5mm} -1}}{0}{0}{0}{{Z_{12}'}^{\hspace{-1,5mm}-2} {Z_7'}^{\hspace{-0,5mm} -1}}
\smallmatthree{{Z_4'}}{0}{0}{{Z_4'}}{{Z_4'}}{0}{{Z_4'}}{{Z_4'}+{Z_4'}^{\hspace{-0,5mm} -2}}{{Z_4'}^{\hspace{-0,5mm} -2}} }
\smallmatthree{{Z_5'} {Z_6'}^{\hspace{-0,5mm} 2}}{0}{0}{0}{{Z_5'} {Z_6'}^{\hspace{-0,5mm} -1}}{0}{0}{0}{{Z_5'}^{\hspace{-0,5mm} -2} {Z_6'}^{\hspace{-0,5mm} -1}}
\end{align*}
where the product of the underlined middle three matrices is
$$
\smallmatthree{{Z_3'} {Z_{12}'} {Z_7'}^{\hspace{-0,4mm}2} {Z_4'}}{0}{0}{{Z_3'}{Z_{12}'}({Z_7'}^{\hspace{-0,4mm}2}+{Z_7'}^{\hspace{-0,4mm}-1}){Z_4'}}{ {Z_3'} {Z_{12}'} {Z_7'}^{\hspace{-0,4mm}-1} {Z_4'}}{ 0 }{ ( {Z_3'} {Z_{12}'} {Z_7'}^{\hspace{-0,4mm}2} + ({Z_3'}+{Z_3'}^{\hspace{-0,4mm}-2}) {Z_{12}'} {Z_7'}^{\hspace{-0,4mm}-1} + {Z_3'}^{\hspace{-0,4mm}-2} {Z_{12}'}^{\hspace{-1,5mm}-2} {Z_7'}^{\hspace{-0,4mm}-1}  ) {Z_4'} }{~ ({Z_3'}+{Z_3'}^{\hspace{-0,4mm}-2}) {Z_{12}'} {Z_7'}^{\hspace{-0,4mm}-1} {Z_4'} + {Z_3'}^{\hspace{-0,4mm}-2} {Z_{12}'}^{\hspace{-1,5mm}-2} {Z_7'}^{\hspace{-0,4mm}-1} ({Z_4'} + {Z_4'}^{\hspace{-0,4mm}-2})~ }{ {Z_3'}^{\hspace{-0,4mm}-2} {Z_{12}'}^{\hspace{-1,5mm}-2} {Z_7'}^{\hspace{-0,4mm}-1} {Z_4'}^{\hspace{-0,4mm}-2} }
$$
The table of terms $\exp(\sum_{j=1}^{12} a_j' {\bf x}_j')$ for ${\rm Tr}^\omega_{\Delta'}([W,s])$, together with $\alpha^{(3)}=a_1'+a_4'-a_3'-a_6'$ and $\alpha^{(4)}=a_3'+a_{10}'-a_4'-a_9' = a_3' -a_4'$, are recorded as before:
\begin{center}
{\setlength{\tabcolsep}{3pt} \begin{tabular}{c|c|c|c|c|c|c|c|c|c|c|c|c|c}
& $(1,1)$ & $(2,1)_1$ & $(2,1)_2$ & $(2,2)$ & $(3,1)_1$ & $(3,1)_2$ & $(3,1)_3$ & $(3,1)_4$ & $(3,2)_1$ & $(3,2)_2$ & $(3,2)_3$ & $(3,2)_4$ & $(3,3)$ \\ \hline
$a_2'$ & 1/3 & 1/3 & 1/3 & 1/3 & $-2/3$  & $-2/3$ & $-2/3$ & $-2/3$ & $-2/3$ & $-2/3$ & $-2/3$ & $-2/3$ & $-2/3$ \\
$a_1'$ & 2/3 & $-1/3$ & $-1/3$ & $-1/3$ & $-1/3$ & $-1/3$ & $-1/3$ & $-1/3$ & $-1/3$ & $-1/3$ & $-1/3$ & $-1/3$ & $-1/3$ \\
$a_3'$ & 1/3 & 1/3 & 1/3 & 1/3 & 1/3 & 1/3 & $-2/3$ & $-2/3$ & 1/3 & $-2/3$ & $-2/3$ & $-2/3$ & $-2/3$ \\
$a_{12}'$ & 1/3 & 1/3 & 1/3 & 1/3 & 1/3 & 1/3 & 1/3 & $-2/3$ & 1/3 & 1/3 & $-2/3$ & $-2/3$ & $-2/3$ \\
$a_7'$ & 2/3 & 2/3 & $-1/3$ & $-1/3$ & 2/3 & $-1/3$ & $-1/3$ & $-1/3$ & $-1/3$ & $-1/3$ & $-1/3$ & $-1/3$ & $-1/3$ \\
$a_4'$ & 1/3 & 1/3 & 1/3 & 1/3 & 1/3 & 1/3 & 1/3 & 1/3 & 1/3 & 1/3 & 1/3 & $-2/3$ & $-2/3$ \\
$a_5'$ & 1/3 & 1/3 & 1/3 & 1/3 & 1/3 & 1/3 & 1/3 & 1/3 & 1/3 & 1/3 & 1/3 & 1/3 & $-2/3$ \\
$a_6'$ & 2/3 & 2/3 & 2/3 & $-1/3$ & 2/3 & 2/3 & 2/3 & 2/3 & $-1/3$ & $-1/3$ & $-1/3$ & $-1/3$ & $-1/3$ \\
$\alpha^{(3)}$ & 0 & $-1$ & $-1$ & 0 & $-1$ & $-1$ & 0 & 0 & 0 & 1 & 1 & 0 & 0 \\
$\alpha^{(4)}$ & 0 & 0 & 0 & 0 & 0 & 0 & $-1$ & $-1$ & 0 & $-1$ & $-1$ & 0 & 0
\end{tabular}}
\end{center}
We should focus on the cases when $\alpha^{(3)}$ or $\alpha^{(4)}$ is nonzero. The sum of the relevant terms $\exp(\sum_{j=1}^{12} a'_j {\bf x}'_j)$ for such cases are as follows:
\begin{align*}
\mbox{$(2,1)_1$, $(2,1)_2$}: & \textstyle \exp( \frac{1}{3} {\bf x}'_2 - \frac{1}{3} {\bf x}'_1 + \frac{1}{3} {\bf x}'_3 + \frac{1}{3} {\bf x}'_{12} + \frac{2}{3} {\bf x}'_7 + \frac{1}{3} {\bf x}'_4 + \frac{1}{3} {\bf x}'_5 + \frac{2}{3} {\bf x}'_6) \\
& \textstyle \quad + \exp (\frac{1}{3} {\bf x}'_2 - \frac{1}{3} {\bf x}'_1 + \frac{1}{3} {\bf x}'_3 + \frac{1}{3} {\bf x}'_{12} - \frac{1}{3} {\bf x}'_7 + \frac{1}{3} {\bf x}'_4 + \frac{1}{3} {\bf x}'_5 + \frac{2}{3} {\bf x}'_6) \\
& \textstyle = \exp( \frac{1}{3} {\bf x}'_2 - \frac{1}{3} {\bf x}'_1 + \frac{1}{3} {\bf x}'_3 + \frac{1}{3} {\bf x}'_{12} + \frac{2}{3} {\bf x}'_7 + \frac{1}{3} {\bf x}'_4 + \frac{1}{3} {\bf x}'_5 + \frac{2}{3} {\bf x}'_6) (1 + e^{-{\bf c}} \exp(-{\bf x}'_7)), \\
\mbox{$(3,1)_1$, $(3,1)_2$}: & \textstyle \exp( - \frac{2}{3} {\bf x}'_2 - \frac{1}{3} {\bf x}'_1 + \frac{1}{3} {\bf x}'_3 + \frac{1}{3} {\bf x}'_{12} + \frac{2}{3} {\bf x}'_7 + \frac{1}{3} {\bf x}'_4 + \frac{1}{3} {\bf x}'_5 + \frac{2}{3} {\bf x}'_6) \\
& \textstyle \quad + \exp( - \frac{2}{3} {\bf x}'_2 - \frac{1}{3} {\bf x}'_1 + \frac{1}{3} {\bf x}'_3 + \frac{1}{3} {\bf x}'_{12} - \frac{1}{3} {\bf x}'_7 + \frac{1}{3} {\bf x}'_4 + \frac{1}{3} {\bf x}'_5 + \frac{2}{3} {\bf x}'_6) \\
& \textstyle = \exp( - \frac{2}{3} {\bf x}'_2 - \frac{1}{3} {\bf x}'_1 + \frac{1}{3} {\bf x}'_3 + \frac{1}{3} {\bf x}'_{12} + \frac{2}{3} {\bf x}'_7 + \frac{1}{3} {\bf x}'_4 + \frac{1}{3} {\bf x}'_5 + \frac{2}{3} {\bf x}'_6) (1 + e^{-{\bf c}} \exp(-{\bf x}'_7)), \\
\mbox{$(3,1)_3$, $(3,1)_4$}: & \textstyle \exp( -\frac{2}{3} {\bf x}'_2 - \frac{1}{3} {\bf x}'_1 - \frac{2}{3} {\bf x}'_3 + \frac{1}{3} {\bf x}'_{12} - \frac{1}{3} {\bf x}'_7 + \frac{1}{3} {\bf x}'_4 + \frac{1}{3} {\bf x}'_5 + \frac{2}{3} {\bf x}'_6) \\
& \textstyle \quad + \exp( -\frac{2}{3} {\bf x}'_2 - \frac{1}{3} {\bf x}'_1 - \frac{2}{3} {\bf x}'_3 - \frac{2}{3} {\bf x}'_{12} - \frac{1}{3} {\bf x}'_7 + \frac{1}{3} {\bf x}'_4 + \frac{1}{3} {\bf x}'_5 + \frac{2}{3} {\bf x}'_6) \\
& \textstyle = \exp( -\frac{2}{3} {\bf x}'_2 - \frac{1}{3} {\bf x}'_1 - \frac{2}{3} {\bf x}'_3 + \frac{1}{3} {\bf x}'_{12} - \frac{1}{3} {\bf x}'_7 + \frac{1}{3} {\bf x}'_4 + \frac{1}{3} {\bf x}'_5 + \frac{2}{3} {\bf x}'_6) ( 1 + e^{-{\bf c}} \exp(-{\bf x}'_{12})), \\
\mbox{$(3,2)_2$, $(3,2)_3$}: & \textstyle \exp( - \frac{2}{3} {\bf x}'_2 - \frac{1}{3} {\bf x}'_1 - \frac{2}{3} {\bf x}'_3 + \frac{1}{3} {\bf x}'_{12} - \frac{1}{3} {\bf x}'_7 + \frac{1}{3} {\bf x}'_4 + \frac{1}{3} {\bf x}'_5 - \frac{1}{3} {\bf x}'_6 ) \\
& \textstyle \quad + \exp( - \frac{2}{3} {\bf x}'_2 - \frac{1}{3} {\bf x}'_1 - \frac{2}{3} {\bf x}'_3 - \frac{2}{3} {\bf x}'_{12} - \frac{1}{3} {\bf x}'_7 + \frac{1}{3} {\bf x}'_4 + \frac{1}{3} {\bf x}'_5 - \frac{1}{3} {\bf x}'_6 ) \\
& \textstyle =  \exp( - \frac{2}{3} {\bf x}'_2 - \frac{1}{3} {\bf x}'_1 - \frac{2}{3} {\bf x}'_3 + \frac{1}{3} {\bf x}'_{12} - \frac{1}{3} {\bf x}'_7 + \frac{1}{3} {\bf x}'_4 + \frac{1}{3} {\bf x}'_5 - \frac{1}{3} {\bf x}'_6 ) (1+e^{-{\bf c}}\exp(-{\bf x}'_{12})).
\end{align*}
So, in view of eq.\eqref{eq:result_of_nu_omega_v_7_and_v_12_modified} of Lem.\ref{lem:section5_modified}, for each $s=(\varepsilon_1,\varepsilon_2)$, we can verify by a similar reasoning as in Case \circled{3} that  $(\nu'_{v_{12}})^{-1}(\nu'_{v_7})^{-1}\nu^\omega_{v_7} \nu^\omega_{v_{12}}({\rm Tr}^\omega_{\Delta'}([W,(\varepsilon_1,\varepsilon_2)]))$ is Laurent for $\Delta'$ and is multiplicity-free.

\vs

This finishes Step 2 (hence also Step 1) of the previous subsection, and therefore completes our proof of Thm.\ref{thm:main}, the main theorem of the present paper.

\section{Consequences and conjectures}

\subsection{Naturality of quantum ${\rm SL}_3$-${\rm PGL}_3$ duality maps under changes of triangulations}
\label{subsec:compatibility_of_quantum_duality_maps}

The main consequence and motivation of our main theorem, Thm.\ref{thm:main}, is the naturality of the quantum ${\rm SL}_3$-${\rm PGL}_3$ duality maps of \cite{Kim} for triangulable punctured surfaces under a change of ideal triangulation $\Delta \leadsto \Delta'$. We urge the readers to keep in mind that in the present subsection, the surfaces are assumed to have empty boundary.
\begin{theorem}[naturality of the quantum duality maps under changes of triangulations]
\label{thm:main_consequence}
Let $\frak{S}$ be a triangulable punctured surface. For each ideal triangulation $\Delta$ of $\frak{S}$ (without self-folded triangles (Def.\ref{def:regular_triangulation})), let $\mathcal{X}^q_\Delta$ be the Fock-Goncharov algebra as defined in Def.\ref{def:FG_algebra}. 
Then, the family of quantum ${\rm SL}_3$-${\rm PGL}_3$ duality maps
$$
\mathbb{I}^q_\Delta : \mathscr{A}_{{\rm SL}_3,\frak{S}}(\mathbb{Z}^T) \to \mathcal{X}^q_\Delta
$$
constructed in \cite[Thm.1.28, Thm.5.83]{Kim} is compatible under the quantum coordinate change maps associated to changes of ideal triangulations. That is, if $\Delta$ and $\Delta'$ are ideal triangulations of $\frak{S}$ (without self-folded triangles), and $\Phi^q_{\Delta\Delta'} : {\rm Frac}(\mathcal{X}^q_{\Delta'}) \to {\rm Frac}(\mathcal{X}^q_\Delta)$ is the corresponding quantum coordinate change map (defined in Def.\ref{def:Phi_q_i}), then
$$
\mathbb{I}^q_\Delta(\ell) = \Phi^q_{\Delta\Delta'}(\mathbb{I}^q_{\Delta'}(\ell)), \qquad \forall \ell \in \mathscr{A}_{{\rm SL}_3,\frak{S}}(\mathbb{Z}^T).
$$
\end{theorem}
We briefly recall from the construction in \cite{Kim} of this map $\mathbb{I}^q_\Delta$. We first need to discuss the domain $\mathscr{A}_{{\rm SL}_3,\frak{S}}(\mathbb{Z}^T)$, which is the set of $\mathbb{Z}^T$-points of the space $\mathscr{A}_{{\rm SL}_3,\frak{S}}$. Here $\mathbb{Z}^T$ is the semi-field of tropical integers, which is $\mathbb{Z}$ as a set, equipped with two binary operations called the tropical addition $\oplus$ and the tropical multiplication $\odot$, defined as $a\oplus b := \max(a,b)$ and $a\odot b= a+b$. Instead of recalling the definition of $\mathscr{A}_{{\rm SL}_3,\frak{S}}$ and hence of $\mathscr{A}_{{\rm SL}_3,\frak{S}}(\mathbb{Z}^T)$ \cite{FG06}, we will only see how these are understood in our situation. First, recall the set $\mathscr{A}_{\rm L}(\frak{S};\mathbb{Z})$ of all ${\rm SL}_3$-laminations in $\frak{S}$ (Def.\ref{def:SL3-lamination}), and the `tropical' coordinate-system map (Prop.\ref{prop:coordinates_are_balanced})
$$
{\bf a}_\Delta ~:~ \mathscr{A}_{\rm L}(\frak{S};\mathbb{Z}) \to ({\textstyle \frac{1}{3}}\mathbb{Z})^{\mathcal{V}(Q_\Delta)}, \quad \ell \mapsto ({\rm a}_v(\ell))_{v\in \mathcal{V}(Q_\Delta)},
$$
developed in \cite{DS1} \cite{Kim}, which is injective. Upon the change of triangulations $\Delta \leadsto \Delta'$, it is proved in \cite{DS2} \cite{Kim} that the coordinate maps ${\bf a}_\Delta$ and ${\bf a}_{\Delta'}$ are related to each other by a tropicalized version of the coordinate change formula for the sequence of cluster $\mathscr{A}$-mutations associated to the transformations $\Delta \leadsto \Delta'$ as mentioned in \S\ref{sec:introduction}; here tropicalization means that one replaces the usual addition, product, and division by the tropical counterparts. These tropicalized coordinate change formulas are piecewise $\mathbb{Z}$-linear. In particular, the integrality condition ${\bf a}_\Delta(\ell) \in \mathbb{Z}^{\mathcal{V}(Q_\Delta)}$ holds if and only if ${\bf a}_{\Delta'}(\ell) \in \mathbb{Z}^{\mathcal{V}(Q_{\Delta'})}$. In \cite{Kim}, an ${\rm SL}_3$-lamination $\ell \in \mathscr{A}_{\rm L}(\frak{S};\mathbb{Z})$ is said to be {\em congruent} if ${\bf a}_{\Delta}(\ell) \in \mathbb{Z}^{\mathcal{V}(Q_{\Delta})} \subset (\frac{1}{3}\mathbb{Z})^{\mathcal{V}(Q_\Delta)}$ holds for some $\Delta$, hence for all $\Delta$ (thanks to the piecewise $\mathbb{Z}$-linearity), i.e. all tropical coordinate values which a priori live in $\frac{1}{3}\mathbb{Z}$ are actually integers. So, for a congruent ${\rm SL}_3$-lamination $\ell$, the coordinates ${\bf a}_\Delta(\ell) = ({\rm a}_v(\ell))_{v\in \mathcal{V}(Q_\Delta)}$ are integer-valued at each node $v$ of the quiver $Q_\Delta$, and they transform via the tropicalized cluster $\mathscr{A}$-mutations under the change of triangulations, in a suitable manner. In this sense, it was concluded in \cite{Kim} that the set of all congruent ${\rm SL}_3$-laminations in $\frak{S}$ is in natural bijection with $\mathscr{A}_{{\rm SL}_3,\frak{S}}(\mathbb{Z}^T)$, giving a geometric meaning to the algebraically defined set $\mathscr{A}_{{\rm SL}_3,\frak{S}}(\mathbb{Z}^T)$. See also \cite{GS15} \cite{Le16} \cite{Akh} for other geometric models for what can be called the ${\rm SL}_3$-laminations, i.e. the elements of $\mathscr{A}_{{\rm SL}_3,\frak{S}}(\mathbb{Z}^T)$.

\vs

We shall construct a map $\mathbb{I}^q_\Delta$ not just on the set $\mathscr{A}_{{\rm SL}_3,\frak{S}}(\mathbb{Z}^T) \subset \mathscr{A}_{\rm L}(\frak{S};\mathbb{Z})$ of congruent ${\rm SL}_3$-laminations in $\frak{S}$, but on the set $\mathscr{A}_{\rm L}(\frak{S};\mathbb{Z})$ of all ${\rm SL}_3$-laminations, with the codomain of the map replaced by $\wh{\mathcal{Z}}^\omega_\Delta$. That is, we shall give a construction of a map
$$
\mathbb{I}^q_\Delta : \mathscr{A}_{\rm L}(\frak{S};\mathbb{Z}) \to \wh{\mathcal{Z}}^\omega_\Delta,
$$
and then later restrict on the subset $\mathscr{A}_{{\rm SL}_3,\frak{S}}(\mathbb{Z}^T)$. We note that in \cite{Kim}, this map on $\mathscr{A}_{\rm L}(\frak{S};\mathbb{Z})$ is denoted by $\wh{\mathbb{I}}^\omega_\Delta$ while its restriction to $\mathscr{A}_{{\rm SL}_3,\frak{S}}(\mathbb{Z}^T)$ is denoted by $\mathbb{I}^q$; here we write $\wh{\mathbb{I}}^\omega_\Delta$ just as $\mathbb{I}^q_\Delta$, to save notations.

\vs

Let $\ell \in \mathscr{A}_{\rm L}(\frak{S};\mathbb{Z})$, i.e. an ${\rm SL}_3$-lamination in $\frak{S}$. If $\ell = {\O}$, set $\mathbb{I}^q_\Delta(\ell) := 1$. Now assume that $\ell$ is non-empty.
Write $\ell$ as a disjoint union $\ell = \ell_1\cup \ell_2 \cup \cdots \cup \ell_n$, where each $\ell_i$ is an ${\rm SL}_3$-lamination represented by a single component ${\rm SL}_3$-web $W_i$ in $\frak{S}$, with weight $k_i \in \mathbb{Z}\setminus\{0\}$. Here, following \cite[Def.3.31]{Kim}, ${\rm SL}_3$-laminations are said to be disjoint if they can be represented by non-elliptic ${\rm SL}_3$-webs in $\frak{S}$ that are disjoint, and a union of disjoint ${\rm SL}_3$-laminations is called a disjoint union of ${\rm SL}_3$-laminations; in fact, we consider a union of ${\rm SL}_3$-laminations only when they are disjoint. 

\vs

If $W_i$ is not a peripheral loop (Def.\ref{def:SL3-lamination}), then let
\begin{align}
\label{eq:I_q_Delta_ell_i_non-peripheral}
\mathbb{I}^q_\Delta(\ell_i) := ( {\rm Tr}^\omega_\Delta([\til{W}_i,{\O}]) )^{k_i},
\end{align}
where $\til{W}_i$ is a lift of $W_i$ in $\frak{S}\times {\bf I}$ at a constant elevation with upward vertical framing; the $k_i$-th power makes sense because $k_i > 0$. If $W_i$ is a peripheral loop, let
\begin{align}
\label{eq:I_q_Delta_ell_i_peripheral}
\mathbb{I}^q_\Delta(\ell_i) := \left[ {\textstyle \prod}_{v\in\mathcal{V}(Q_\Delta)} {\bf X}_v^{{\rm a}_v(\ell_i)} \right]_{\rm Weyl};
\end{align}
that is, $\mathbb{I}^q_\Delta(\ell_i) = [ (\prod_{v\in\mathcal{V}(Q_\Delta)} {\bf X}_v^{{\rm a}_v(W_i)})^{k_i}]_{\rm Weyl} = ( [ (\prod_{v\in\mathcal{V}(Q_\Delta)} {\bf X}_v^{{\rm a}_v(W_i)})]_{\rm Weyl} )^{k_i}$ in this case (we used Lem.\ref{lem:Weyl-ordering_basic}(C) and the fact that $\prod_{v\in \mathcal{V}(Q_\Delta)} {\bf X}_v^{{\rm a}_v(W_i)}$ is invertible), where $W_i$ denotes the ${\rm SL}_3$-lamination consisting just of $W_i$ with weight $1$. Finally, let
\begin{align}
\label{eq:I_q_Delta_ell_as_product}
\mathbb{I}^q_\Delta(\ell) = \mathbb{I}^q_\Delta(\ell_1) \mathbb{I}^q_\Delta(\ell_2) \cdots \mathbb{I}^q_\Delta(\ell_n),
\end{align}
given by the product.

\vs

If $W_i$ is not a peripheral loop, then $\mathbb{I}^q_\Delta(\ell_i) = ({\rm Tr}^\omega_\Delta([\til{W}_i,{\O}]))^{k_i}$ lies in $\wh{\mathcal{Z}}^\omega_\Delta \subset \mathcal{Z}^\omega_\Delta$, because ${\rm Tr}^\omega_\Delta([\til{W}_i,{\O}])$ lies in $\wh{\mathcal{Z}}^\omega_\Delta$ due to Prop.\ref{prop:value_of_the_SL3_quantum_trace_is_Delta-balanced}, and $\wh{\mathcal{Z}}^\omega_\Delta$ is a subring of $\mathcal{Z}^\omega_\Delta$. If $W_i$ is a peripheral loop, then from eq.\eqref{eq:I_q_Delta_ell_i_peripheral} and eq.\eqref{eq:X_ell} we have $\mathbb{I}^q_\Delta(\ell_i) = {\bf X}^{\ell_i} = [\prod_{v\in\mathcal{V}} {\bf X}_v^{{\rm a}_v(\ell_i)}]_{\rm Weyl}$. From Prop.\ref{prop:coordinates_are_balanced},  Def.\ref{def:Delta-balanced_elements} and Def.\ref{def:balanced_subalgebras}, we see that $\mathbb{I}^q_\Delta(\ell_i) \in \wh{\mathcal{Z}}^\omega_\Delta$. In particular, we can deduce that each factor $\mathbb{I}^q_\Delta(\ell_i)$ in the right-hand side of eq.\eqref{eq:I_q_Delta_ell_as_product} belongs to $\wh{\mathcal{Z}}^\omega_\Delta$.
\begin{proposition}
\label{prop:I_q_Delta_is_well-defined}
Above $\mathbb{I}^q_\Delta(\ell)$ is a well-defined element of $\wh{\mathcal{Z}}^\omega_\Delta\subset \mathcal{Z}^\omega_\Delta$.
\end{proposition}
To prove Prop.\ref{prop:I_q_Delta_is_well-defined}, what needs to be checked is that each $\mathbb{I}^q_\Delta(\ell_i)$ commutes with each $\mathbb{I}^q_\Delta(\ell_j)$. This is addressed in \cite[Lem.5.74]{Kim}, with the treatment of peripheral loops left as an exercise, as being an ${\rm SL}_3$-analog of \cite[Lem.3.9]{AK} which is for the setting of ${\rm SL}_2$. For completeness, we give a proof here. We note that the method presented below (for ${\rm SL}_3$) can be modified to give an alternative proof of \cite[Lem.3.9]{AK} for ${\rm SL}_2$. We first establish a following statement about peripheral loops. For the remainder of the present section, we let
$$
\mathcal{V} = \mathcal{V}(Q_\Delta)
$$
for convenience, unless there is a possible confusion.
\begin{lemma}
\label{lem:commutativity_of_Laurent_monomials_of_peripheral}
Let $\ell_0$ be an ${\rm SL}_3$-lamination in $\frak{S}$ represented by a union of peripheral loops (Def.\ref{def:SL3-lamination}) with arbitrary weights, and let $\ell'$ be any ${\rm SL}_3$-lamination in $\frak{S}$. Then $[\prod_{v\in\mathcal{V}} {\bf X}_v^{{\rm a}_v(\ell_0)}]_{\rm Weyl}$ commutes with $[\prod_{v\in\mathcal{V}} {\bf X}_v^{{\rm a}_v(\ell')}]_{\rm Weyl}$, that is,
\begin{align}
\label{eq:commutativity_statement}
[{\textstyle \prod}_{v\in\mathcal{V}} {\bf X}_v^{{\rm a}_v(\ell_0)}]_{\rm Weyl} [{\textstyle \prod}_{v\in\mathcal{V}} {\bf X}_v^{{\rm a}_v(\ell')}]_{\rm Weyl}
= [{\textstyle \prod}_{v\in\mathcal{V}} {\bf X}_v^{{\rm a}_v(\ell')}]_{\rm Weyl} [{\textstyle \prod}_{v\in\mathcal{V}} {\bf X}_v^{{\rm a}_v(\ell_0)}]_{\rm Weyl}.
\end{align}
\end{lemma}
{\it Proof of Lem.\ref{lem:commutativity_of_Laurent_monomials_of_peripheral}.} For convenience, let us introduce some temporary notations for this proof. For any ${\rm SL}_3$-lamination $\ell$ in $\frak{S}$, write
\begin{align}
\label{eq:X_ell}
{\bf X}^\ell := [{\textstyle \prod}_{v\in\mathcal{V}} {\bf X}_v^{{\rm a}_v(\ell)}]_{\rm Weyl}.
\end{align}
An ${\rm SL}_3$-lamination in $\frak{S}$ will be called a \ul{\em weight-positive} if it can be represented by a non-elliptic ${\rm SL}_3$-web in $\frak{S}$ with all weights being non-negative. An ${\rm SL}_3$-lamination in $\frak{S}$ is called a \ul{\em peripheral} if it can be represented by a non-elliptic ${\rm SL}_3$-web in $\frak{S}$ given by a union of peripheral loops (with arbitrary weights). A peripheral is called a \ul{\em single peripheral} if it can be represented by a single peripheral loop. So, the sought-for equation eq.\eqref{eq:commutativity_statement} can be written as
\begin{align}
\nonumber
{\bf X}^{\ell_0} {\bf X}^{\ell'} = {\bf X}^{\ell'} {\bf X}^{\ell_0},
\end{align}
which should be proved for each peripheral $\ell_0$ and each weight-positive $\ell'$. We do this in several steps.

\vs

\ul{Step 1} : ${\bf X}^{\ell_0} {\bf X}^{\ell'} = {\bf X}^{\ell'} {\bf X}^{\ell_0}$ holds when $\ell_0$ is a single peripheral and $\ell'$ is a weight-positive.

\vs

Let $\ell_0$ and $\ell'$ be as such. Then $\ell_0$ can be represented by a non-elliptic ${\rm SL}_3$-web $W_0$ in $\frak{S}$ consisting of a peripheral loop, with some weight $k_0 \in \mathbb{Z}$. Let $\til{W}_0$ be a constant-elevation lift of $W_0$ in $\frak{S}\times {\bf I}$ with upward vertical framing. Similarly, $\ell'$ can be represented by a non-elliptic ${\rm SL}_3$-web $W'$ in $\frak{S}$, all of whose components are of weight $1$. Let $\til{W}'$ be a constant-elevation lift of $W'$ in $\frak{S}\times {\bf I}$ with upward vertical framing.

\vs

One can isotope $W_0$ so that $W_0$ is contained in a small enough neighborhood of the puncture it surrounds, so that $W_0$ is disjoint from $W'$. Denote still by $\til{W}_0$ the corresponding lift of $W_0$ in $\frak{S} \times {\bf I}$. Since the product structure of $\mathcal{S}^\omega_{\rm s}(\frak{S};\mathbb{Z})_{\rm red}$ is defined by superposition, one can then observe that the elements $[\til{W}_0,{\O}]$ and $[\til{W}',{\O}]$ of $\mathcal{S}^\omega_{\rm s}(\frak{S};\mathbb{Z})_{\rm red}$ commute with each other. Since ${\rm Tr}^\omega_\Delta$ is an algebra homomorphism (Thm.\ref{thm:SL3_quantum_trace}), we thus have
\begin{align}
\label{eq:product_of_two_quantum_traces_equal}
{\rm Tr}^\omega_{\Delta}([\til{W}_0,{\O}]) {\rm Tr}^\omega_{\Delta}([\til{W}',{\O}]) = {\rm Tr}^\omega_{\Delta}([\til{W}',{\O}]) {\rm Tr}^\omega_{\Delta}([\til{W}_0,{\O}])
\end{align}
Note that Prop.\ref{prop:highest_term}, the highest-term statement, applies to both $[\til{W}_0,{\O}]$ and $[\til{W}',{\O}]$. Thus, ${\rm Tr}^\omega_{\Delta}([\til{W}_0,{\O}])$ is a Laurent polynomial with the unique highest Laurent monomial term being ${\bf X}^{W_0}$, where $W_0$ denotes the ${\rm SL}_3$-lamination represented by $W_0$ with weight $1$,  and ${\rm Tr}^\omega_{\Delta}([\til{W}',{\O}])$ is a Laurent polynomial with the unique highest Laurent monomial term being ${\bf X}^{\ell'}$. We can then observe that the left-hand side of eq.\eqref{eq:product_of_two_quantum_traces_equal} is a Laurent polynomial with the unique highest term being ${\bf X}^{W_0} {\bf X}^{\ell'}$, whereas the right-hand side is a Laurent polynomial with the unique highest term being ${\bf X}^{\ell'} {\bf X}^{W_0}$. The equality of eq.\eqref{eq:product_of_two_quantum_traces_equal} says that these highest terms of the left and the right sides agree, i.e. ${\bf X}^{W_0} {\bf X}^{\ell'} = {\bf X}^{\ell'} {\bf X}^{W_0}$. Note from ${\rm a}_v(\ell_0) = k_0\cdot {\rm a}_v(W_0)$ and Lem.\ref{lem:Weyl-ordering_basic}(C) that ${\bf X}^{\ell_0} = ({\bf X}^{W_0})^{k_0}$. From this we obtain the sought-for commutation relation ${\bf X}^{\ell_0} {\bf X}^{\ell'} = {\bf X}^{\ell'} {\bf X}^{\ell_0}$.

\vs

\ul{Step 2} : ${\bf X}^{\ell_0} {\bf X}^{\ell'} = {\bf X}^{\ell'} {\bf X}^{\ell_0}$ holds when $\ell_0$ and $\ell'$ are single peripherals.

\vs

Let $\ell'$ be represented by a single peripheral loop $W'$ in $\frak{S}$ with weight $k' \in \mathbb{Z}$, and view $W'$ as an ${\rm SL}_3$-lamination itself, with weight $1$. Then, $W'$ is a weight-positive, hence Step 1 applies, yielding ${\bf X}^{\ell_0} {\bf X}^{W'} = {\bf X}^{W'} {\bf X}^{\ell_0}$. By a similar reasoning used at the end of Step 1, we have ${\bf X}^{\ell'} = ({\bf X}^{W'})^{k'}$. Hence ${\bf X}^{\ell_0} {\bf X}^{\ell'} = {\bf X}^{\ell'} {\bf X}^{\ell_0}$ follows, as desired.

\vs

\ul{Step 3} : ${\bf X}^{\ell_0} {\bf X}^{\ell'} = {\bf X}^{\ell'} {\bf X}^{\ell_0}$ holds when $\ell_0$ and $\ell'$ are any peripherals.

\vs

One can write $\ell_0$ as a disjoint union $\ell_0 = \ell_1 \cup \ell_2 \cup \cdots \cup \ell_n$, with each $\ell_i$ being a single peripheral for $i=1,\ldots,n$; then ${\rm a}_v(\ell_0) = {\rm a}_v(\ell_1)+\cdots +{\rm a}_v(\ell_n)$ for all $v\in \mathcal{V}$ (see \cite[Lem.3.32]{Kim}). By Step 2, ${\bf X}^{\ell_i} = [\prod_{v\in\mathcal{V}}{\bf X}_v^{{\rm a}_v(\ell_i)}]_{\rm Weyl}$ and ${\bf X}^{\ell_j}=[\prod_{v\in\mathcal{V}}{\bf X}_v^{{\rm a}_v(\ell_j)}]_{\rm Weyl}$ commute for each $i,j\in\{1,\ldots,n\}$. Thus, 
\begin{align*}
{\bf X}^{\ell_1} \cdots {\bf X}^{\ell_n} & = [{\textstyle \prod}_{v\in\mathcal{V}} {\bf X}_v^{{\rm a}_v(\ell_1)}]_{\rm Weyl} \cdots [{\textstyle \prod}_{v\in\mathcal{V}} {\bf X}_v^{{\rm a}_v(\ell_n)}]_{\rm Weyl} \\
& = [{\textstyle \prod}_{v\in\mathcal{V}} {\bf X}_v^{{\rm a}_v(\ell_1)+\cdots+{\rm a}_v(\ell_n)}]_{\rm Weyl} \quad (\because \mbox{Lem.\ref{lem:Weyl-ordering_basic}(A)}) \\
& = [{\textstyle \prod}_{v\in\mathcal{V}} {\bf X}_v^{{\rm a}_v(\ell_0)}]_{\rm Weyl} = {\bf X}^{\ell_0}.
\end{align*}
Likewise, if one writes $\ell'$ as a disjoint union $\ell' = \ell'_1 \cup \cdots \cup \ell'_r$ of single peripherals $\ell'_i$, the Laurent monomials ${\bf X}^{\ell'_1}, \ldots, {\bf X}^{\ell'_r}$ commute with each other, and ${\bf X}^{\ell'} = {\bf X}^{\ell'_1} \cdots {\bf X}^{\ell'_r}$. Note that Step 2 says that each of ${\bf X}^{\ell_i}$ commutes with each of ${\bf X}^{\ell'_j}$ also. Hence ${\bf X}^{\ell_0}$ and ${\bf X}^{\ell'}$ commute, as desired.

\vs

\ul{Step 4} : ${\bf X}^{\ell_0} {\bf X}^{\ell'} = {\bf X}^{\ell'} {\bf X}^{\ell_0}$ holds when $\ell_0$ is any peripheral and $\ell'$ is a weight-positive.

\vs

Write $\ell_0$ as a disjoint union $\ell_0 = \ell_1 \cup \ell_2 \cup \cdots \cup \ell_n$ of single peripherals, so that ${\bf X}^{\ell_0} = {\bf X}^{\ell_1} \cdots {\bf X}^{\ell_n}$ as seen in the proof of Step 3. Since each $\ell_i$ for $i=1,\ldots,n$ is a single peripheral, by Step 1 we have ${\bf X}^{\ell_i} {\bf X}^{\ell'} = {\bf X}^{\ell'} {\bf X}^{\ell_i}$. Thus
$$
{\bf X}^{\ell_0} {\bf X}^{\ell'} = {\bf X}^{\ell_1} {\bf X}^{\ell_2} \cdots {\bf X}^{\ell_n} {\bf X}^{\ell'}
= {\bf X}^{\ell_1} {\bf X}^{\ell_2} \cdots {\bf X}^{\ell'} {\bf X}^{\ell_n}
= \cdots
= {\bf X}^{\ell'} {\bf X}^{\ell_1} {\bf X}^{\ell_2} \cdots {\bf X}^{\ell_n}
= {\bf X}^{\ell'} {\bf X}^{\ell_0},
$$
where we moved ${\bf X}^{\ell'}$ to the left using the commutation relations ${\bf X}^{\ell_i} {\bf X}^{\ell'} = {\bf X}^{\ell'} {\bf X}^{\ell_i}$, $i=1,\ldots,n$.

\vs

\ul{Step 5} : ${\bf X}^{\ell_0} {\bf X}^{\ell'} = {\bf X}^{\ell'} {\bf X}^{\ell_0}$ holds when $\ell_0$ is any peripheral and $\ell'$ is any ${\rm SL}_3$-lamination.

\vs

First, let $\ell'_0$ be any peripheral such that the disjoint union $\ell'' := \ell' \cup \ell'_0$ is a weight-positive. Such $\ell'_0$ exists; for example, one can construct $\ell'_0$ as the ${\rm SL}_3$-lamination based on the ${\rm SL}_3$-web given by the union of all peripheral loops whose isotopy classes appear in $\ell'$, given big enough positive weights (this idea is used in \cite{Kim}). Following the notation in \cite{Kim}, denote by $-\ell'_0$ the ${\rm SL}_3$-lamination obtained from $\ell'_0$ by multiplying all weights by $-1$ (when $\ell'_0$ is represented by an ${\rm SL}_3$-web with weights on components). As in \cite[Lem.3.33]{Kim}, we then have ${\rm a}_v(-\ell'_0) = -{\rm a}_v(\ell'_0)$ for all $v\in \mathcal{V}$, and also $\ell' \cup \ell'_0 \cup (-\ell'_0) = \ell'$ holds as ${\rm SL}_3$-laminations (here the left-hand side $\ell' \cup \ell'_0 \cup (-\ell'_0)$ is a disjoint union of ${\rm SL}_3$-laminations). So we have $\ell'' \cup (-\ell'_0) = \ell'$, and hence ${\rm a}_v(\ell'')+{\rm a}_v(-\ell'_0)={\rm a}_v(\ell')$ for all $v\in \mathcal{V}$ (\cite[Lem.3.32]{Kim}).

\vs

Since $\ell_0$ is a peripheral and $\ell''$ is a weight-positive, Step 4 applies and yields
\begin{align}
\label{eq:step4_applied1}
{\bf X}^{\ell_0} {\bf X}^{\ell''} = {\bf X}^{\ell''} {\bf X}^{\ell_0}.
\end{align}
Since $-\ell'_0$ is a peripheral and $\ell''$ is a weight-positive, Step 4 applies and yields
\begin{align}
\nonumber
{\bf X}^{-\ell_0'} {\bf X}^{\ell''} = {\bf X}^{\ell''} {\bf X}^{-\ell_0'},
\end{align}
which, together with ${\rm a}_v(\ell'')+{\rm a}_v(-\ell'_0)={\rm a}_v(\ell')$ and Lem.\ref{lem:Weyl-ordering_basic}(A), implies that
\begin{align}
\label{eq:step5_eq1}
{\bf X}^{\ell'} = {\bf X}^{-\ell_0'}{\bf X}^{\ell''}.
\end{align}
Note now that
\begin{align*}
{\bf X}^{\ell_0} {\bf X}^{\ell'} & = {\bf X}^{\ell_0} {\bf X}^{-\ell'_0} {\bf X}^{\ell''} \quad (\because \mbox{eq.\eqref{eq:step5_eq1}}) \\
& = {\bf X}^{-\ell_0'} {\bf X}^{\ell_0} {\bf X}^{\ell''} \quad (\because \mbox{$\ell_0,-\ell_0'$ are peripherals; apply Step 3}) \\
& = {\bf X}^{-\ell_0'} {\bf X}^{\ell''} {\bf X}^{\ell_0} \quad (\because \mbox{eq.\eqref{eq:step4_applied1}}) \\
& = {\bf X}^{\ell'} {\bf X}^{\ell_0} \quad (\because \mbox{eq.\eqref{eq:step5_eq1}}),
\end{align*}
as desired. \qed \quad {\it [End of proof of Lem.\ref{lem:commutativity_of_Laurent_monomials_of_peripheral}]}

\begin{corollary}
\label{cor:peripheral_element_is_central}
Let $\ell_0$ be an ${\rm SL}_3$-lamination in $\frak{S}$ consisting only of peripheral loops (with arbitrary weights). Then ${\bf X}^{\ell_0} = [\prod_{v\in\mathcal{V}} {\bf X}_v^{{\rm a}_v(\ell_0)}]_{\rm Weyl}$ (eq.\eqref{eq:X_ell}) is in the center of the $\Delta$-balanced cube-root Fock-Goncharov algebra $\wh{\mathcal{Z}}^\omega_\Delta \subset \mathcal{Z}^\omega_\Delta$ (Def.\ref{def:balanced_subalgebras}). \end{corollary}

{\it Proof of Cor.\ref{cor:peripheral_element_is_central}.} 
In view of the definition of $\wh{\mathcal{Z}}^\omega_\Delta$ (Def.\ref{def:balanced_subalgebras}), it suffices to show that ${\bf X}^{\ell_0}$ commutes with each $\Delta$-balanced Laurent monomial $[\prod_{v\in\mathcal{V}} {\bf X}_v^{a_v}]_{\rm Weyl}$ for $\Delta$, with $(a_v)_{v\in\mathcal{V}}$ being an arbitrary element of $(\frac{1}{3}\mathbb{Z})^{\mathcal{V}}$ that is $\Delta$-balanced in the sense of Def.\ref{def:Delta-balanced_elements}. Pick any such $(a_v)_{v\in\mathcal{V}}$. Then by Prop.\ref{prop:coordinates_are_balanced}, there exists an ${\rm SL}_3$-lamination $\ell$ in $\frak{S}$ such that the tuple of its tropical coordinates $({\rm a}_v(\ell))_{v\in\mathcal{V}}$ coincides with the tuple $(a_v)_{v\in\mathcal{V}}$. Thus, $[\prod_{v\in\mathcal{V}} {\bf X}_v^{a_v}]_{\rm Weyl} = {\bf X}^\ell$, using the notation in eq.\eqref{eq:X_ell}. We showed in Lem.\ref{lem:commutativity_of_Laurent_monomials_of_peripheral} that ${\bf X}^{\ell_0}$ commutes with ${\bf X}^\ell$, as desired. \qed

\vs

We are now ready to prove Prop.\ref{prop:I_q_Delta_is_well-defined}, which is on the well-definedness of $\mathbb{I}^q_\Delta(\ell) \in \wh{\mathcal{Z}}^\omega_\Delta$.

\vs

{\it Proof of Prop.\ref{prop:I_q_Delta_is_well-defined}.} As mentioned, what needs to be checked is that each pair of factors $\mathbb{I}^q_\Delta(\ell_i)$ and $\mathbb{I}^q_\Delta(\ell_j)$ appearing in the product expression in the right-hand side of eq.\eqref{eq:I_q_Delta_ell_as_product} commute with each other. Recall that an ${\rm SL}_3$-lamination $\ell$ is expressed as a disjoint union $\ell_1\cup \ell_2 \cup \cdots \cup \ell_n$ of ${\rm SL}_3$-laminations, and that $\ell_i$ is represented by a single-component ${\rm SL}_3$-web $W_i$ in $\frak{S}$ with weight $k_i \in \mathbb{Z}\setminus \{0\}$; we can assume that $W_1,\ldots,W_n$ are disjoint from each other. We denoted by $\til{W}_i$ a constant-elevation lift in $\frak{S}\times {\bf I}$ of $W_i$ with upward vertical framing. Since we are dealing with a punctured surface, $\til{W}_i$ has no endpoints.

\vs

Since the projections of $\til{W}_i$ to $\frak{S}$ are mutually disjoint, one observes that the elements $[\til{W}_i,{\O}]$ of $\mathcal{S}^\omega_{\rm s}(\frak{S};\mathbb{Z})_{\rm red}$, for $i=1,\ldots,n$, mutually commute with respect to the superposition product (Def.\ref{def:stated_SL3-skein_algebra}). Since the ${\rm SL}_3$ quantum trace map ${\rm Tr}^\omega_\Delta : \mathcal{S}^\omega_{\rm s}(\frak{S};\mathbb{Z})_{\rm red} \to \mathcal{Z}^\omega_\Delta$ is an algebra homomorphism (Thm.\ref{thm:SL3_quantum_trace}), it follows that ${\rm Tr}^\omega_\Delta([\til{W}_i,{\O}])$, for $i=1,\ldots,n$, hence also their positive powers $({\rm Tr}^\omega_\Delta([\til{W}_i,{\O}]))^{k_i}$ (with $k_i\in\mathbb{Z}_{>0}$), mutually commute. It follows that $\mathbb{I}^q_\Delta(\ell_i)$ and $\mathbb{I}^q_\Delta(\ell_j)$ commute if neither $W_i$ nor $W_j$ is a peripheral loop; see eq.\eqref{eq:I_q_Delta_ell_i_non-peripheral}.

\vs

Suppose that $W_i$ is a peripheral loop. Then, since $\mathbb{I}^q_\Delta(\ell_j) \in \wh{\mathcal{Z}}^\omega_\Delta$ for all $j$, it follows from Cor.\ref{cor:peripheral_element_is_central} that $\mathbb{I}^q_\Delta(\ell_i) = {\bf X}^{\ell_i}$ commutes with $\mathbb{I}^q_\Delta(\ell_j)$ for all $j$.

\vs

So all $\mathbb{I}^q_\Delta(\ell_i)$, $i=1,\ldots,n$, mutually commute with each other, hence indeed $\mathbb{I}^q_\Delta(\ell)$ is well defined through eq.\eqref{eq:I_q_Delta_ell_as_product}. Moreover, since we saw that each $\mathbb{I}^q_\Delta(\ell_i)$ belongs to the algebra $\wh{\mathcal{Z}}^\omega_\Delta$, it follows that $\mathbb{I}^q_\Delta(\ell) \in \wh{\mathcal{Z}}^\omega_\Delta$. \qed

\vs

By applying Prop.\ref{prop:highest_term} to each $\mathbb{I}^q_\Delta(\ell_i)$ in eq.\eqref{eq:I_q_Delta_ell_as_product} and using arguments similar to above, together with the additivity of the tropical coordinates as in \cite[Lem.3.32]{Kim}, one can show that $\mathbb{I}^q_\Delta(\ell)$ has the unique highest term $[\prod_{v\in\mathcal{V}} {\bf X}_v^{{\rm a}_v(\ell)}]_{\rm Weyl}$ (\cite[Thm.5.83(3)]{Kim}). In case $\ell \in \mathscr{A}_{\rm L}(\frak{S};\mathbb{Z})$ belongs to $\mathscr{A}_{{\rm SL}_3,\frak{S}}(\mathbb{Z}^T)$, we have ${\rm a}_v(\ell) \in \mathbb{Z}$ for all $v\in \mathcal{V}$, so this highest term is an element of $\mathcal{X}^q_\Delta \subseteq \mathcal{Z}^\omega_\Delta$. Now, from Prop.\ref{prop:congruence_of_SL3_quantum_trace} it follows that the other Laurent monomial terms of $\mathbb{I}^q_\Delta(\ell) \in \mathcal{Z}^\omega_\Delta$ also belong to $\mathcal{X}^q_\Delta$, hence $\mathbb{I}^q_\Delta(\ell) \in \mathcal{X}^q_\Delta$ as desired, for $\ell \in \mathscr{A}_{{\rm SL}_3,\frak{S}}(\mathbb{Z}^T)$. So, indeed the image under $\mathbb{I}^q_\Delta:\mathscr{A}_{\rm L}(\frak{S};\mathbb{Z}) \to \wh{\mathcal{Z}}^\omega_\Delta$ of the subset $\mathscr{A}_{{\rm SL}_3,\frak{S}}(\mathbb{Z}^T)$ lies in $\mathcal{X}^q_\Delta$, so the restriction of $\mathbb{I}^q_\Delta$ on $\mathscr{A}_{{\rm SL}_3,\frak{S}}(\mathbb{Z}^T)$ can be written as $\mathbb{I}^q_\Delta : \mathscr{A}_{{\rm SL}_3,\frak{S}}(\mathbb{Z}^T) \to \mathcal{X}^q_\Delta$. Notice that the arguments so far, 
which are largely based on \cite{Kim}, prove that for $\ell \in \mathscr{A}_{{\rm SL}_3,\frak{S}}(\mathbb{Z}^T)$, $\mathbb{I}^q_\Delta(\ell)$ is a well-defined Laurent polynomial in the variables ${\bf X}_v$, $v\in \mathcal{V}=\mathcal{V}(Q_\Delta)$, with coefficients being in $\mathbb{Z}[\omega^{\pm 1/2}]=\mathbb{Z}[q^{\pm 1/18}]$. It is expected in \cite[\S5]{Kim} that the coefficients live in $\mathbb{Z}[q^{\pm 1/2}]$ (in fact in $\mathbb{Z}[q^{\pm 1}]$ when $\frak{S}$ is a punctured surface, without boundary, as we are assuming now), but this hasn't been proved and left as a conjecture.
\vs

We now turn into the proof of Thm.\ref{thm:main_consequence}, i.e. the naturality of $\mathbb{I}^q_\Delta$ with respect to the quantum coordinate change maps. In view of the above construction of $\mathbb{I}^q_\Delta$, Thm.\ref{thm:main} yields the naturality $\Theta^\omega_{\Delta\Delta'} \mathbb{I}^q_{\Delta'}(\ell_i) = \mathbb{I}^q_{\Delta}(\ell_i)$ for each $\ell_i$ that is not represented by a peripheral loop (eq.\eqref{eq:I_q_Delta_ell_i_non-peripheral}).  For peripheral loops, we need some more work; namely, the following proposition yields the naturality $\Theta^\omega_{\Delta\Delta'} \mathbb{I}^q_{\Delta'}(\ell_i) = \mathbb{I}^q_{\Delta}(\ell_i)$ for $\ell_i$ represented by a peripheral loop.

\begin{proposition}
\label{prop:compatibility_of_peripheral_term}
Let $\Delta$ be an ideal triangulation of a triangulable generalized marked surface $\frak{S}$ that has at least one puncture. Let $W$ be a constant-elevation (with upward vertical framing) lift in $\frak{S} \times {\bf I}$ of an oriented peripheral loop in $\frak{S}$ surrounding a puncture. 
\begin{enumerate}
\item[\rm (1)] ${\rm Tr}^\omega_\Delta([W,{\O}])$ is a sum of three Weyl-ordered $\Delta$-balanced Laurent monomials in $\mathcal{Z}^\omega_\Delta$ (Def.\ref{def:Weyl-ordered});

\item[\rm (2)] Among the three Laurent monomial terms of ${\rm Tr}^\omega_\Delta([W,{\O}]) \in \mathcal{Z}^\omega_\Delta$, the Laurent monomial of the highest partial ordering is $[\prod_{v\in \mathcal{V}(Q_\Delta)} {\bf X}_v^{{\rm a}_v(\pi(W))}]_{\rm Weyl}$.

\item[\rm (3)] For any other ideal triangulation $\Delta'$, one has
$$
[{\textstyle \prod}_{v\in \mathcal{V}(Q_\Delta)} {\bf X}_v^{{\rm a}_v(\pi(W))}]_{\rm Weyl} = \Theta^\omega_{\Delta\Delta'}([{\textstyle \prod}_{v\in \mathcal{V}(Q_{\Delta'})} ({\bf X}_v')^{{\rm a}_v'(\pi(W))}]_{\rm Weyl}),
$$
where ${\rm a}_v'(\pi(W))$ denote the tropical coordinates of $\pi(W)$ in terms of $\Delta'$.
\end{enumerate}
\end{proposition}

{\it Proof of Prop.\ref{prop:compatibility_of_peripheral_term}.} One may note that the two items (1) and (2) are proved for the classical setting $\omega^{1/2}=1$ in \cite[Prop.4.15]{Kim}. In the quantum setting, the arguments in the proof of this classical setting still go through, in case when the projection $\pi(W)$ of $W$ is a loop in $\frak{S}$ that meets each arc of $\Delta$ at most once. However, in a general case, one needs to come up with a more careful treatment, as we do now. First, one can isotope $W$, within the class of constant-elevation ${\rm SL}_3$-webs in $\frak{S}\times {\bf I}$, so that $W$ meets $\Delta \times {\bf I}$ in a minimal number. Then, $\Delta\times {\bf I}$ divides $W$ into left- or right-turn oriented edges living over triangles of $\Delta$, where these arcs are either all left turns or all right turns (see e.g. \cite[Lem.4.11]{Kim}). Assume that they are all left turns. The case of the all-right-turn can be taken care of with only a slight modification of the argument. Consider a split ideal triangulation $\wh{\Delta}$ for $\Delta$ (as explained immediately after Prop.\ref{prop:biangle_quantum_trace}), and assume that $W$ still meets $\wh{\Delta} \times {\bf I}$ in a minimal number, so that for each biangle $B$ of $\wh{\Delta}$, $W\cap (B\times {\bf I})$ consists of `parallel' arcs (i.e. non-intersecting simple arcs at the same elevation). Apply a vertical isotopy to $W$, so that for each triangle $\wh{t}$ of $\wh{\Delta}$, each of the components of $W\cap (\wh{t} \times {\bf I})$ is at a constant elevation at all times throughout the isotopy, and that in the end, for each triangle $\wh{t}$ of $\wh{\Delta}$, the components of $W\cap (\wh{t} \times {\bf I})$ are at mutually distinct elevations. So, $W$ would be in a `good position' with respect to $\wh{\Delta}$, and in fact in a `gool position', in the sense used in \cite[\S5.3]{Kim}. Still, for each biangle $B$ of $\wh{\Delta}$, the projection of $W\cap (B\times {\bf I})$ in $\frak{S}$ consists of parallel arcs in $B$ (i.e. non-intersecting simple arcs in $B$ connecting the two sides of $B$), but a component of $W\cap (B\times {\bf I})$ may not be at a constant elevation. In fact, what matters is the ordering of the elevations of the components of $W\cap (\wh{t} \times {\bf I})$ for each triangle $\wh{t}$. Let $t$ and $u$ be two triangles of $\Delta$ sharing a side, so that the corresponding triangles $\wh{t}$ and $\wh{u}$ of $\wh{\Delta}$ `share' a common biangle $B$. We say that the ordering of components of the part of $W$ living over $\wh{t}$ is \ul{\em compatible} with that for $\wh{u}$ \ul{\em at} this biangle $B$, if the ordering of elevations of the endpoints of $W\cap (\wh{t}\times {\bf I})$ lying over a side of $B$ and that of the endpoints of $W\cap (\wh{u}\times {\bf I})$ lying over the other side of $B$ correspond to each other by the connectedness relation by the arcs of $W\cap (B\times {\bf I})$. That is to say, the ${\rm SL}_3$-web $W\cap (B\times {\bf I})$ in $B\times {\bf I}$ can be isotoped by a vertical isotopy within the class of the ${\rm SL}_3$-webs in $B \times {\bf I}$ so that the components of $W\cap (B\times {\bf I})$ lie in mutually distinct elevation intervals. Another way to put it is that the element $[W\cap (B\times {\bf I})]$ of the ${\rm SL}_3$-skein algebra $\mathcal{S}^\omega(B;\mathbb{Z})$ (Def.\ref{def:stated_SL3-skein_algebra}) is given by the product of its constituent edges, each of which connects the two sides of $B$. It is proved in \cite[Thm.1.2]{CKKO} that the ordering of elevations of the components of $W\cap (\wh{t}\times {\bf I})$ for each triangle $\wh{t}$ can be chosen so that the above compatibility holds at all biangles of $\wh{\Delta}$. Let's use such an elevation ordering for each triangle $\wh{t}$.

\vs

Consider the {\em junctures} $W\cap (\wh{\Delta}\times {\bf I})$, and a juncture-state $J : W\cap(\wh{\Delta}\times {\bf I}) \to \{1,2,3\}$. First, the state-sum formula in eq.\eqref{eq:full_state-sum_formula} (\cite[\S5.3]{Kim}) yields
\begin{align}
\label{eq:state-sum_later}
{\rm Tr}^\omega_\Delta([W,{\O}]) = {\textstyle \sum}_J ({\textstyle \prod}_B {\rm Tr}^\omega_B([W\cap (B\times {\bf I}), J])) \, ({\textstyle \prod}_{\wh{t}} {\rm Tr}^\omega_{\wh{t}}([W\cap (\wh{t}\times {\bf I}),J]))
\end{align}
where the sum is over all juncture-states $J$, the product $\prod_B$ is over all biangles $B$ of $\wh{\Delta}$, and the product $\prod_{\wh{t}}$ is over all triangles $\wh{t}$ of $\wh{\Delta}$. Since $W$ was put into a good (or a gool) position, for each $\wh{t}$, ${\rm Tr}^\omega_{\wh{t}}([W\cap (\wh{t}\times {\bf I}),J])$ is a product of ${\rm Tr}^\omega_{\wh{t}}([W_{\wh{t};i}, J])$, where $W_{\wh{t};1}$, $W_{\wh{t};2}$, ... are components of $W\cap (\wh{t} \times {\bf I})$, each of which is a left-turn edge over $\wh{t}$. Because of the elevation compatibility at each biangle $B$, we see that ${\rm Tr}^\omega_B([W\cap (B\times {\bf I}),J])$ is a product of ${\rm Tr}^\omega_B([W_{B;j}, J])$, where $W_{B;1}$, $W_{B;2}$, ... are components of $W\cap (B\times {\bf I})$, each of which is a simple edge over $B$ connecting the two boundary walls of $B$. For each component $W_{B;j}$, by Prop.\ref{prop:biangle_quantum_trace}(BQT2),
\begin{align}
\label{eq:biangle_value_at_elementary_arc}
{\rm Tr}^\omega_B([W_{B;j},J]) = \left\{
\begin{array}{ll}
1 & \mbox{if $J$ assigns the same state values to the two endpoints of $W_{B;j}$,} \\
0 & \mbox{otherwise}.
\end{array}
\right.
\end{align}

\begin{lemma}
A juncture-state $J : W\cap (\wh{\Delta}\times {\bf I}) \to \{1,2,3\}$ has a nonzero contribution to the sum in eq.\eqref{eq:state-sum_later} if and only if $J$ is a constant juncture-state, i.e. assigns the same value to all junctures.
\end{lemma}
This follows from eq.\eqref{eq:state-sum_later}, eq.\eqref{eq:biangle_value_at_elementary_arc}, and Thm.\ref{thm:SL3_quantum_trace}(QT2) which says that the value of each ${\rm Tr}^\omega_{\wh{t}}([W_{\wh{t};i}, J])$ equals the $(\varepsilon_1,\varepsilon_2)$-th entry of the matrix in eq.\eqref{eq:QT2-1}, where this matrix is upper triangular. 

\vs

The $(1,1)$-th entry of ${\rm Tr}^\omega_{\wh{t}}([W_{\wh{t};i}, J])$ is of the highest partial ordering (among $(\varepsilon,\varepsilon)$-th entries), so the constant juncture-state with value $1$ yields the highest term of ${\rm Tr}^\omega_\Delta([W,{\O}])$. By the proof of \cite[Prop.4.15]{Kim}, which is the classical counterpart of the items (1) and (2) of the current proposition, we then obtain the items (1) and (2). In general, the quantum situation is more subtle and complicated than the classical situation, as the values of the biangle ${\rm SL}_3$ quantum trace ${\rm Tr}^\omega_B$ could be complicated, as they are essentially a Reshetikhin-Turaev invariant for the standard 3-dimensional representation of $\mathcal{U}_q(\frak{sl}_3)$ \cite{RT}, involving R-matrices. Here, the elevation compatibility of \cite{CKKO} allowed us to avoid such a complicated computation.

\vs

For the item (3), it suffices to show the statement in the case when $\Delta\leadsto \Delta'$ is a flip at an arc $k$. Let $e$ be an arc of $\Delta$ that is different from $k$ and that meets the peripheral loop $\pi(W)$. Suppose that $W$ is isotoped so that it satisfies the nice properties with respect to a split ideal triangulation $\wh{\Delta}$ of $\Delta$ as above, i.e. in a gool position and having the elevation ordering compatibility at biangles. Cut along $e$; let $\frak{S}_e$ be the resulting surface, $\Delta_e$ and $\Delta'_e$ the triangulations of $\frak{S}_e$ induced from $\Delta$ and $\Delta'$, and $W_e$ the ${\rm SL}_3$-web in $\frak{S}_e$, obtained by this cutting process (Def.\ref{def:cutting_process}). Pick one point $x$ in $W \cap (e\times {\bf I})$ (there can be at most two such points), and let $x_1$ and $x_2$ be the endpoints of $W_e$ corresponding to $x$. For a state $s_e : \partial W_e \to \{1,2,3\}$ of $W_e$ that is compatible with the original state $s : \partial W \to \{1,2,3\}$ of $W$, it must be that $s_e(x_1
)=s_e(x_2)$ (note $\partial W = {\O}$). For each $\varepsilon \in \{1,2,3\}$, denote by $(\varepsilon,\varepsilon)$ the state $s_e$ that assigns $\varepsilon$ to $x_1$ and $x_2$. By the cutting/gluing property (Thm.\ref{thm:SL3_quantum_trace}(QT1)) we have
$$
i_{\Delta,\Delta_e} {\rm Tr}^\omega_\Delta([W,{\O}]) = {\textstyle \sum}_{\varepsilon=1}^3 {\rm Tr}^\omega_{\Delta_e}([W_e,(\varepsilon,\varepsilon)]),
$$
and in view of the relationship between the above equation and the state-sum formula for  ${\rm Tr}^\omega_\Delta([W,{\O}])$ in eq.\eqref{eq:state-sum_later}, one can observe that the summand in eq.\eqref{eq:state-sum_later} corresponding to the constant juncture-state $J$ with value $\varepsilon$ is sent via the cutting map $i_{\Delta,\Delta_e}$ (Def.\ref{def:cutting_process}) to the term ${\rm Tr}^\omega_{\Delta_e}([W_e,(\varepsilon,\varepsilon)])$ on the right-hand side. Likewise for $\Delta'$. Hence we have
$$
i_{\Delta,\Delta_e}( [{\rm Tr}^\omega_\Delta([W,{\O}])]_{\rm high}) = {\rm Tr}^\omega_{\Delta_e}([W_e,(1,1)]), \quad
i_{\Delta',\Delta'_e}( [{\rm Tr}^\omega_{\Delta'}([W,{\O}])]_{\rm high}) = {\rm Tr}^\omega_{\Delta'_e}([W_e,(1,1)]),
$$
where $[\sim]_{\rm high}$ stands for the highest Laurent monomial term. Note now (in the balanced fraction algebras) that
\begin{align*}
i_{\Delta,\Delta_e}( \Theta^\omega_{\Delta\Delta'} [{\rm Tr}^\omega_{\Delta'}([W,{\O}])]_{\rm high}) 
& = \Theta^\omega_{\Delta_e \Delta'_e} i_{\Delta',\Delta'_e} ([{\rm Tr}^\omega_{\Delta'}([W,{\O}])]_{\rm high}) \quad (\because\mbox{Prop.\ref{prop:compatibility_of_Theta_under_cutting}}) \\
& = \Theta^\omega_{\Delta_e \Delta'_e} {\rm Tr}^\omega_{\Delta'_e}([W_e,(1,1)]) \\
& = {\rm Tr}^\omega_{\Delta_e}([W_e,(1,1)]) \quad (\because\mbox{Thm.\ref{thm:main}}) \\
& = i_{\Delta,\Delta_e}([{\rm Tr}^\omega_\Delta([W,{\O}])]_{\rm high}).
\end{align*}
Hence it follows from the injectivity of $i_{\Delta,\Delta_e}$ (Lem.\ref{lem:i_injective}) that $\Theta^\omega_{\Delta\Delta'} [{\rm Tr}^\omega_{\Delta'}([W,{\O}])]_{\rm high} = [{\rm Tr}^\omega_\Delta([W,{\O}])]_{\rm high}$. \qed

\begin{remark}
The proof of Prop.\ref{prop:compatibility_of_peripheral_term} also yields an analogous result for peripheral arcs (Def.\ref{def:SL3-lamination}), when the surface $\frak{S}$ has boundary.
\end{remark}

We finally arrive at a proof of Thm.\ref{thm:main_consequence}.

\vs

{\it Proof of Thm.\ref{thm:main_consequence}}. For a triangulable punctured surface $\frak{S}$, by the construction of the quantum duality map $\mathbb{I}^q_\Delta$, Thm.\ref{thm:main} and Prop.\ref{prop:compatibility_of_peripheral_term} together imply $\Theta^\omega_{\Delta\Delta'} \mathbb{I}^q_{\Delta'}(\ell_i) = \mathbb{I}^q_{\Delta}(\ell_i)$ for all $\ell_i$ in eq.\eqref{eq:I_q_Delta_ell_as_product}, hence $\Theta^\omega_{\Delta\Delta'} \mathbb{I}^q_{\Delta'}(\ell) = \mathbb{I}^q_{\Delta}(\ell)$ for all $\ell \in \mathscr{A}_{\rm L}(\frak{S};\mathbb{Z})$. When $\ell \in \mathscr{A}_{{\rm SL}_3,\frak{S}}(\mathbb{Z}^T)$, we saw that $\mathbb{I}^q_{\Delta'}(\ell) \in \mathcal{X}^q_{\Delta'}$ and $\mathbb{I}^q_\Delta(\ell) \in \mathcal{X}^q_\Delta$. Now, by Lem.\ref{lem:Theta_extends_Phi}, we can now write $\Phi^q_{\Delta\Delta'} \mathbb{I}^q_{\Delta'}(\ell) = \mathbb{I}^q_\Delta(\ell)$, for $\ell \in \mathscr{A}_{{\rm SL}_3,\frak{S}}(\mathbb{Z}^T)$, hence Thm.\ref{thm:main_consequence}.\qed

\vs

Thus, the quantum duality maps $\mathbb{I}^q_\Delta$ for all ideal triangulations $\Delta$ of a triangulable punctured surface $\frak{S}$ can be viewed as constituting a single quantum duality map
$$
\mathbb{I}^q : \mathscr{A}_{{\rm SL}_3,\frak{S}}(\mathbb{Z}^T) \to \mathcal{O}^q_{\rm tri}(\mathscr{X}_{{\rm PGL}_3,\frak{S}}),
$$
where $\mathcal{O}^q_{\rm tri}(\mathscr{X}_{{\rm PGL}_3,\frak{S}})$ stands for the ring of all elements that are `quantum $\mathcal{X}$-Laurent', i.e. are Laurent in the quantum $\mathscr{X}$-variables ${\bf X}_v$, $v\in \mathcal{V}(Q_\Delta)$, i.e. belong to $\mathcal{X}^q_\Delta$, for every ideal triangulation $\Delta$. That is, for each ideal triangulation $\Delta$ we define
$$
\mathcal{O}^q_\Delta(\mathscr{X}_{{\rm PGL}_3,\frak{S}}) := \bigcap_{\Delta'} \Phi^q_{\Delta\Delta'}(\mathcal{X}^q_{\Delta'}) \quad \subset \quad \mathcal{X}^q_\Delta \quad \subset \quad {\rm Frac}(\mathcal{X}^q_\Delta),
$$
where $\bigcap_{\Delta'}$ is over all ideal triangulations $\Delta'$; note that $\mathcal{O}^q_\Delta(\mathscr{X}_{{\rm PGL}_3,\frak{S}})$ for different triangulations are isomorphically identified through the (restrictions of) the maps $\Phi^q_{\Delta\Delta'}$. We denote these rings $\mathcal{O}^q_\Delta(\mathscr{X}_{{\rm PGL}_3,\frak{S}})$ for all possible ideal triangulations $\Delta$ collectively by $\mathcal{O}^q_{\rm tri}(\mathscr{X}_{{\rm PGL}_3,\frak{S}})$. Recall the classical ${\rm SL}_3$-${\rm PGL}_3$ duality map
$$
\mathbb{I} : \mathscr{A}_{{\rm SL}_3,\frak{S}}(\mathbb{Z}^T) \to \mathcal{O}(\mathscr{X}_{{\rm PGL}_3,\frak{S}})
$$
constructed in \cite{Kim}, whose image forms a basis of $\mathcal{O}(\mathscr{X}_{{\rm PGL}_3,\frak{S}})$, the ring of all regular functions on Fock and Goncharov's moduli stack $\mathscr{X}_{{\rm PGL}_3,\frak{S}}$ of `framed' ${\rm PGL}_3$-local systems on $\frak{S}$ \cite{FG06}. As mentioned in \S\ref{sec:introduction}, it is known from \cite[Thm.1.1]{Shen} that $\mathcal{O}(\mathscr{X}_{{\rm PGL}_3,\frak{S}})$ coincides with the ring $\mathcal{O}_{\rm cl}(\mathscr{X}_{{\rm PGL}_3,\frak{S}})$ of classical universally $\mathcal{X}$-Laurent elements, i.e. the elements of the Laurent polynomial ring $\mathcal{X}^1_\Delta \cong \mathbb{Z}[\{X_v^{\pm 1} : v\in \mathcal{V}(Q_\Delta)\}]$ that stay $\mathcal{X}$-Laurent after arbitrary sequence of mutations (that is, $\mathcal{X}$-Laurent not just for all triangulations, but for all possible cluster $\mathscr{X}$-seeds). Sending each basis element $\mathbb{I}(\ell) \in \mathcal{O}(\mathscr{X}_{{\rm PGL}_3,\frak{S}})$ to $\mathbb{I}^q(\ell) \in \mathcal{O}^q_{\rm tri}(\mathscr{X}_{{\rm PGL}_3,\frak{S}})$, one obtains a deformation quantization map
$$
\mathcal{O}(\mathscr{X}_{{\rm PGL}_3,\frak{S}}) \to \mathcal{O}^q_{\rm tri}(\mathscr{X}_{{\rm PGL}_3,\frak{S}})
$$
for the (Poisson) moduli space $\mathscr{X}_{{\rm PGL}_3,\frak{S}}$ for a triangulable punctured surface $\frak{S}$.

\subsection{Future perspectives}

We list some conjectures that are not mentioned in the introduction.
\begin{conjecture}
\label{conj:quantum_universally_Laurent}
Let $\frak{S}$ be a triangulable punctured surface. For each congruent ${\rm SL}_3$-lamination $\ell \in \mathscr{A}_{{\rm SL}_3,\frak{S}}(\mathbb{Z}^T)$, the corresponding element $\mathbb{I}^q(\ell) \in \mathcal{O}^q_{\rm tri}(\mathscr{X}_{{\rm PGL}_3,\frak{S}})$ belongs to $\mathcal{O}^q_{\rm cl}(\mathscr{X}_{{\rm PGL}_3,\frak{S}})$, i.e. is quantum $\mathcal{X}$-Laurent for {\em all} cluster $\mathscr{X}$-seeds for $\mathscr{X}_{{\rm PGL}_3,\frak{S}}$, not just for the seeds corresponding to ideal triangulations of $\frak{S}$.
\end{conjecture}
In the present paper, we had a glimpse of Conjecture \ref{conj:quantum_universally_Laurent}. Namely, our proofs and arguments can be used to show that this quantum $\mathcal{X}$-Laurent property holds for the cluster $\mathscr{X}$-seeds sitting `in between' the seeds for the ideal triangulations. One possible expectation is that a quantum version of \cite[Lem.2.2]{Shen} \cite[Thm.3.9]{GHK} would hold, which would say that if a quantum $\mathcal{X}$-Laurent element for any chosen seed stays quantum $\mathcal{X}$-Laurent after all possible single mutations from this seed, then this element is universally quantum $\mathcal{X}$-Laurent. The author was informed by Linhui Shen that one could prove this `expectation' for (the quantum version of) the cluster $\mathscr{X}$-variety $\mathscr{X}_{|Q|}$ associated to the mutation equivalence class $|Q|$ of {\it any} initial quiver $Q$. Here $\mathscr{X}_{|Q|}$ is the scheme made by gluing the split algebraic tori $\mathbb{G}_m^\mathcal{V} = {\rm Spec}(\mathbb{Q}[X_v^{\pm 1} : v\in \mathcal{V}])$ (where $\mathcal{V}$ is the set of nodes of $Q$) associated to cluster $\mathscr{X}$-seeds obtained from the initial seed by sequences of mutations and seed automorphisms, where the gluing maps are given by the corresponding rational coordinate change formulas. We present a rough argument of this proof as follows: first observe from \cite[Thm.5.1]{BZ} that this holds for (the quantum version of) the cluster $\mathscr{A}$-variety $\mathscr{A}_{|Q|}$ which is constructed with (quantum) cluster $\mathscr{A}$-mutations instead of (quantum) cluster $\mathscr{X}$-mutations, and then observe that one can embed the quantized $\mathscr{X}$-algebra $\mathcal{O}^q(\mathscr{X}_{|Q|})$ into a quantized $\mathscr{A}$-algebra $\mathcal{O}^q(\mathscr{A}_{|Q|})$ (see \cite{BZ}) as in \cite[Prop.18.5]{GS19} so that the situation reduces to \cite[Thm.5.1]{BZ}. Then, what would remain to show is whether $\mathbb{I}^q_\Delta(\ell)$ stays quantum $\mathcal{X}$-Laurent after mutating at a node lying in the interior of a triangle. 

\vs

Another direction of research related to Conjecture \ref{conj:quantum_universally_Laurent} is to extend the results of the present paper, as well as those of \cite{Kim}, so that they incorporate more general kinds of ideal triangulations without the non-self-foldedness (Def.\ref{def:regular_triangulation}) assumption \cite{JK}, and flips among them; we would allow self-folded triangles, and maybe we should also take into consideration the `tagged' ideal triangulations \cite{FST}.

\vs

Since $\mathcal{O}(\mathscr{X}_{{\rm PGL}_3,\frak{S}})$ coincides with the ring $\mathcal{O}_{\rm cl}(\mathscr{X}_{{\rm PGL}_3,\frak{S}}) \cong \mathcal{O}(\mathscr{X}_{|Q_\Delta|})$ of all rational functions on $\mathscr{X}_{{\rm PGL}_3,\frak{S}}$ that are regular (i.e. $\mathcal{X}$-Laurent) for all cluster $\mathscr{X}$-seeds (\cite[Thm.1.1]{Shen}), once one has Conjecture \ref{conj:quantum_universally_Laurent} then one can write the deformation quantization map as
$$
\mathcal{O}_{\rm cl}(\mathscr{X}_{{\rm PGL}_3,\frak{S}}) \to \mathcal{O}^q_{\rm cl}(\mathscr{X}_{{\rm PGL}_3,\frak{S}})
$$
or solely in terms of the cluster $\mathscr{X}$-variety $\mathscr{X}_{|Q_\Delta|}$ associated to the mutation-equivalence class $|Q_\Delta|$ of the quiver $Q_\Delta$
$$
\mathcal{O}(\mathscr{X}_{|Q_\Delta|}) \to \mathcal{O}^q(\mathscr{X}_{|Q_\Delta|}).
$$
After having Conjecture \ref{conj:quantum_universally_Laurent}, the natural next step is the following conjecture.
\begin{conjecture}
\label{conj:last}
Let $\frak{S}$ be a triangulable punctured surface. 
The elements $\mathbb{I}^q(\ell)$, $\ell \in \mathscr{A}_{{\rm SL}_3,\frak{S}}(\mathbb{Z}^T)$, form a basis of $\mathcal{O}^q_{\rm cl}(\mathscr{X}_{{\rm PGL}_3,\frak{S}})$.
\end{conjecture}
A natural approach to the above conjecture is to try to compare the elements $\mathbb{I}^q(\ell)$ with the quantum theta functions of Davison and Mandel \cite{DM}; this approach is already mentioned in \cite{Kim}, but an important tool for it that was missing in \cite{Kim} is precisely the main result of the present paper, Thm.\ref{thm:main} and Thm.\ref{thm:main_consequence}. Another approach to Conjecture \ref{conj:last} without going through \cite{DM} nor Conjecture \ref{conj:quantum_universally_Laurent} is studied in a work in progress, joint with Linhui Shen \cite{Kim-Shen}. In \cite{Kim-Shen} we will also extend the construction of the ${\rm SL}_3$-${\rm PGL}_3$ quantum duality map and proofs of its properties to the case of surfaces with boundary; we note that most of the arguments work similarly as for punctured surfaces, but there are some subtleties to deal with.

\vs

For readers' reference, as suggested by a referee, we note that the classical and quantum ${\rm SL}_3$-${\rm PGL}_3$ duality maps constructed in \cite{Kim} and studied in the present section are higher-rank analogs of the well-known ${\rm SL}_2$-${\rm PGL}_2$ duality maps, which we briefly review now. 

\vs

The classical ${\rm SL}_2$-${\rm PGL}_2$ duality map $\mathbb{I} : \mathscr{A}_{{\rm SL}_2,\frak{S}}(\mathbb{Z}^T) \to \mathcal{O}(\mathscr{X}_{{\rm PGL}_2,\frak{S}})$ is constructed by Fock and Goncharov \cite{FG06}, and a corresponding quantum ${\rm SL}_2$-${\rm PGL}_2$ duality map $\mathbb{I}^q : \mathscr{A}_{{\rm SL}_2,\frak{S}}(\mathbb{Z}^T) \to \mathcal{O}^q_{\rm tri}(\mathscr{X}_{{\rm PGL}_2,\frak{S}})$ is constructed by Allegretti and the author \cite{AK}, where the latter is based on Bonahon and Wong's ${\rm SL}_2$ quantum trace maps \cite{BW} going from the ${\rm SL}_2$ skein algebras to the Fock-Goncharov algebras for the $2$-triangulation quivers $Q^{[2]}_\Delta$ (Fig.\ref{fig:n-triangulation}) for ideal triangulations $\Delta$ of $\frak{S}$. The codomain ring $\mathcal{O}^q_{\rm tri}(\mathscr{X}_{{\rm PGL}_2,\frak{S}})$ of $\mathbb{I}^q$ is defined in a similar manner as $\mathcal{O}^q_{\rm tri}(\mathscr{X}_{{\rm PGL}_3,\frak{S}})$, using ideal triangulations of $\frak{S}$; we remark that, one can regard the quantum duality map $\mathbb{I}^q$ of \cite{AK} as $\mathbb{I}^q : \mathscr{A}_{{\rm SL}_2,\frak{S}}(\mathbb{Z}^T) \to \mathcal{O}^q_{\rm cl}(\mathscr{X}_{{\rm PGL}_2,\frak{S}})$ (note the codomain), thanks to the work of Mandel and Qin \cite{MQ} which relates the result of \cite{AK} and the work of Davison and Mandel \cite{DM}. In the meantime, the author was informed by Linhui Shen that it is possible to directly prove that $\mathcal{O}^q_{\rm tri}(\mathscr{X}_{{\rm PGL}_2,\frak{S}})$ coincides with $\mathcal{O}^q_{\rm cl}(\mathscr{X}_{{\rm PGL}_2,\frak{S}})$. Namely, we use the quantum version of \cite[Lem.2.2]{Shen} \cite[Thm.3.9]{GHK} of which we presented a brief proof above; that is, if a quantum $\mathcal{X}$-Laurent element for a seed stays quantum Laurent after all possible single mutations from that seed, then it is universally quantum $\mathcal{X}$-Laurent. In this ${\rm SL}_2$-${\rm PGL}_2$ setting, an ideal triangulation corresponds to a seed, and a mutation at a node means flip at an ideal arc of the ideal triangulation, in case this flipping also yields an ideal triangulation; we note that the coordinate change formula associated to a flip that involves a self-folded triangle is slightly different from the cluster mutation formula (see \cite[\S10.7]{FG06} \cite[\S9]{AB20} \cite{JK}). To prove $\mathcal{O}^q_{\rm tri}(\mathscr{X}_{{\rm PGL}_2,\frak{S}}) = \mathcal{O}^q_{\rm cl}(\mathscr{X}_{{\rm PGL}_2,\frak{S}})$, it suffices to show that there exists an ideal triangulation of $\frak{S}$ such that flipping at any of its constituent ideal arc yields an ideal triangulation; it is not hard to see that this is true for punctured surfaces $\frak{S}$, unless it is a thrice-punctured sphere. Indeed, for surfaces with exactly one puncture, any ideal triangulation satisfies this. For a fourth-punctured sphere it is easy to come up with an explicit example that works. Say we have such an ideal triangulation $\Delta$ for a surface with genus $g$ and $n$ punctures; introduce a new puncture in the interior of any of the ideal triangles of $\Delta$, and add to $\Delta$ three new ideal arcs in this triangle incident to the new puncture, to obtain an ideal triangulation $\Delta'$ for a surface with genus $g$ and $n+1$ punctures; one can see that $\Delta'$ satisfies the property, finishing the proof by induction.

\vs

Fock and Goncharov realized the abstract set $\mathscr{A}_{{\rm SL}_2,\frak{S}}(\mathbb{Z}^T)$, which is the domain of $\mathbb{I}$ and $\mathbb{I}^q$, geometrically as a set of `${\rm SL}_2$-laminations' in $\frak{S}$, where an ${\rm SL}_2$-lamination in a triangulable punctured surface $\frak{S}$ is a collection of isotopy classes of mutually disjoint unoriented simple essential loops in $\frak{S}$ with integer weights, subject to certain condition.

\vs

To a non-peripheral loop $\gamma$ in $\frak{S}$, Fock and Goncharov's duality map $\mathbb{I}$ associates the trace-of-monodromy function along $\gamma$ on $\mathscr{X}_{{\rm PGL}_2,\frak{S}}$, and to a peripheral loop $\gamma$ is associated a special monomial function (which is the `highest term' of the trace-of-monodromy along $\gamma$). To a non-peripheral loop $\gamma$ with weight $k \in \mathbb{Z}_{>0}$ is associated the trace of $[\gamma]^k=[\gamma^k] \in \pi_1(\frak{S})$ (i.e. $\gamma$ going $k$ times around), which can also be given as a result of applying the $k$-th Chebyshev polynomial $T_k$ of the first kind to the trace of $\gamma$, and to a peripheral loop $\gamma$ with weight $k\in \mathbb{Z}$ is associated the $k$-th power of the monomial function for $\gamma$. To an ${\rm SL}_2$-lamination, the duality map $\mathbb{I}$ associates the product of functions associated to its mutually non-homotopic components as just described. Here the polynomials $T_k(t) \in \mathbb{Z}[t]$ are defined recursively by $T_0(t)=1$, $T_1(t)=t$, $T_k(t) = t T_{k-1}(t) - T_{k-2}(t)$, $k\ge 2$; its characteristic property is ${\rm tr}(A^k) = T_k({\rm tr}(A))$ for $A\in {\rm SL}_2$. These polynomials are used in the construction of $\mathbb{I}$ by Fock and Goncharov \cite{FG06} and also in that of $\mathbb{I}^q$ by Allegretti and the author for non-peripheral loops with positive weights, and hence the corresponding bases (i.e. the images of the duality maps) are called {\it bracelets} bases. 

\vs

On the other hand, to a non-peripheral loop $\gamma$ with weight $k\in\mathbb{Z}_{>0}$ one can try associating the $k$-th power of the trace of $\gamma$, which is different from trace of $\gamma^k$, as $({\rm tr}(A))^k \neq {\rm tr}(A^k$) for $A\in {\rm SL}_2$ in general. Nevertheless, this also leads to duality maps and hence bases of the rings of regular functions, which are called {\it bangles} bases. The terms ``bracelets" and ``bangles" come from the literature on bases of the ${\rm SL}_2$ skein algebras (see \cite{T14}); a ``bracelet" represents a picture of a simple loop winding around multiple times (hence having some self-intersections), whereas a `bangle' represents a picture of multiple parallel copies of a simple loop. In terms of properties mostly related to the theory of cluster varieties, the bracelets bases have been considered as more canonical objects than the bangles bases. 

\vs

We note that the classical and quantum ${\rm SL}_3$-${\rm PGL}_3$ duality maps constructed in \cite{Kim} are higher-rank analogs of the bangles bases for ${\rm SL}_2$-${\rm PGL}_2$, as opposed to the bracelets bases for ${\rm SL}_2$-${\rm PGL}_2$. We expect that what can be called ${\rm SL}_3$-${\rm PGL}_3$ bracelets bases should be more canonical objects. For a non-peripheral loop component of an ${\rm SL}_3$-lamination, say represented by a loop $\gamma$ with weight $k \in \mathbb{Z}_{\ge 0}$, a bracelets duality map would associate the trace of $[\gamma]^k=[\gamma^k] \in \pi_1(\frak{S})$ (as already mentioned in \cite{Kim}), while the bangles duality map of \cite{Kim} associates the $k$-th power of $\gamma$ as seen in eq.\eqref{eq:I_q_Delta_ell_i_non-peripheral}. In an upcoming joint work with Thang L\^e \cite{KL}, one special property satisfied by these bracelets ${\rm SL}_3$-${\rm PGL}_3$ duality maps but not by the bangles duality maps will be studied; see also the work of Bonahon and Higgins \cite{BH23}. However, as mentioned in \cite[\S6]{Kim} it is not clear how the web components having 3-valent vertices should be dealt with, in what would be regarded as full `correct' higher analogs of the ${\rm SL}_2$-${\rm PGL}_2$ bracelets bases, and it is an interesting and important future problem to study this.


\begin{thebibliography}{CKKO20}

\bibitem[A20]{Akh} T. Akhmejanov, {\it Non-elliptic Webs and Convex Sets in the Affine Building}, Doc. Math. {\bf 25} (2020), 2413-2443. \quad arXiv:2004.13803

\bibitem[AB20]{AB20} D.G.L. Allegretti and T. Bridgeland, {\it The monodromy of meromorphic projective structures}, T. Am. Math. Soc. {\bf 373}(9) (2020), 6321--6367. \quad arXiv:1802.02505

\bibitem[AK17]{AK} D.G.L. Allegretti and H. Kim, {\it A duality map for quantum cluster varieties from surfaces}, Adv. Math. {\bf 306} (2017), 1164--1208. \quad arXiv:1509.01567

\bibitem[BL20]{BL20} W. Bloomquist and T.T.Q. L\^e, {\it The Chebyshev-Frobenius Homomorphism for stated skein modules of 3-manifolds}, Math. Z. {\bf 301} (2022), 1063--1105. \quad arXiv:2011.02130

\bibitem[BZ05]{BZ} A. Berenstein and A. Zelevinsky, {\it Quantum cluster algebras}, Adv. Math. {\bf 195}(2) (2005) 405-455. \quad arXiv:math.QA/0404446


\bibitem[BH23]{BH23} F. Bonahon and V. Higgins, {\it Central elements in the ${\rm SL}_d$-skein algebra of a surface}, arXiv:2308.13691

\bibitem[BW11]{BW} F. Bonahon and H. Wong, {\it Quantum traces for representations of surface groups in ${\rm SL}_2(\mathbb{C})$}, Geom. Topol. {\bf 15} (2011), 1569--1615. \quad arXiv:1003.5250

\bibitem[BW16]{BW16} F. Bonahon and H. Wong, {\it Representations of the Kauffman bracket skein algebra I: invariants and miraculous cancellations}, Invent. Math. {\bf 204} (2016) 195--243. \quad arXiv:1206.1638

\bibitem[CF99]{CF99} L.O. Chekhov and V.V. Fock, {\it A quantum Teichm\"uller space}, Theor. Math. Phys. {\bf 120} (1999), 1245--1259.

\bibitem[CKKO20]{CKKO} S. Cho, H. Kim, H. Kim and D. Oh, {\it Laurent positivity of quantized canonical bases for quantum cluster varieties from surfaces}, Commun. Math. Phys. {\bf 373} (2020), 655-705. \quad arXiv:1710.06217

\bibitem[C95]{Cohn} P.M. Cohn, {\it Skew Fields: Theory of General Division Rings}, Encyclopedia of Mathematics and its Applications {\bf 57}, Cambridge University Press, Cambridge, 1995. 

\bibitem[DM21]{DM} B. Davison and T. Mandel, {\it Strong positivity for quantum theta bases of quantum cluster algebras}, Invent. Math. {\bf 226} (2021), 725--843. \quad arXiv:1910.12915

\bibitem[D21]{Douglas1} D.C. Douglas, {\it Quantum traces for ${\rm SL}_n(\mathbb{C})$: the case $n=3$}, J. Pure Appl. Algebra {\bf 228}(7) (2024), 107652. \quad arXiv:2101.06817

\bibitem[DS20a]{DS1} D.C. Douglas and Z. Sun, {\it Tropical Fock-Goncharov coordinates for ${\rm SL}_3$-webs on surfaces I: construction}, Forum Math. Sigma {\bf 12} (2024), e5. \quad arXiv:2011.01768

\bibitem[DS20b]{DS2} D.C. Douglas and Z. Sun, {\it Tropical Fock-Goncharov coordinates for ${\rm SL}_3$-webs on surfaces II: naturality}, arXiv:2012.14202


\bibitem[FK94]{FaKa} L. D. Faddeev and R. M. Kashaev, {\it Quantum dilogarithm}, Modern Phys. Lett. A{\bf 9} (1994) 427--434

\bibitem[FG06a]{FG06} V.V. Fock and A.B. Goncharov, {\it Moduli spaces of local systems and higher Teichm\"uller theory}, Publ. Math. Inst. Hautes \'{E}tudes Sci. {\bf 103} (2006), 1--211. arXiv:math/0311149v4

\bibitem[FG06b]{FG06b} V.V. Fock and A.B. Goncharov, ``Cluster $\mathscr{X}$-varieties, amalgamation, and Poisson-Lie groups, in {\em Algebraic geometry and number theory}, Progr. Math. {\bf 253}, Birkh\"auser, 2006, 27--68.



\bibitem[FG09a]{FG09a} V.V. Fock and A.B. Goncharov, {\it Cluster ensembles, quantization and the dilogarithm}, Annales Scientifiques de l'\'Ecole Normale Sup\'erieure, Serie 4, vol. {\bf 42}(6) (2009), 865-930.

\bibitem[FG09b]{FG09b} V. V. Fock and A. B. Goncharov, {\it The quantum dilogarithm and representations of the quantum cluster varieties}, Invent. Math. {\bf 175}(2) (2009) 223--286.

\bibitem[FST08]{FST} S. Fomin, M. Shapiro and D. Thurston, {\it Cluster algebras and triangulated surfaces. Part I: cluster complexes}, Acta Math. {\bf 201}(1) (2008), 83--146. \quad arXiv:math/0608367

\bibitem[FS22]{FS} C. Frohman and A. Sikora, {\it ${\rm SU}(3)$-skein algebras and webs on surfaces}, Math. Z. {\bf 300} (2022), 33-56. \quad arXiv:2002.08151



\bibitem[GS15]{GS15} A.B. Goncharov and L. Shen, {\it Geometry of canonical bases and mirror symmetry}, Invent. Math. {\bf 202} (2015), 487--633. \quad arXiv:1309.5922

\bibitem[GS18]{GS18} A.B. Goncharov and L. Shen, {\it Donaldson-Thomas transformations of moduli spaces of ${\rm G}$-local systems}, Adv. Math. {\bf 327} (2018), 225--348. \quad arXiv:1602.06479

\bibitem[GS19]{GS19} A.B. Goncharov L. Shen, {\it Quantum geometry of moduli spaces of local systems and representation theory}, arXiv:1904.10491

\bibitem[GHK15]{GHK} M. Gross, P. Hacking, and S. Keel, {\it Birational geometry of cluster algebras}, Algebr. Geom. {\bf 2}(2) (2015), 137--175. \quad arXiv:1309.2573

\bibitem[GHKK18]{GHKK} M. Gross, P. Hacking, S. Keel, and M. Kontsevich, {\it Canonical bases for cluster algebras}, J. Am. Math. Soc. {\bf 31} no.2 (2018), 497--608. \quad arXiv:1411.1394

\bibitem[H10]{Hiatt} C. Hiatt, {\it Quantum traces in quantum Teichm\"uller theory}, Algebr. Geom. Topol. {\bf 10}(3) (2010), 1245--1283. \quad arXiv:0809.5118

\bibitem[H20]{Higgins} V. Higgins, {\it Triangular decomposition of ${\rm SL}_3$ skein algebras}, Quantum Topol. {\bf 14}(1) (2023), 1--63.

\bibitem[JK]{JK} Seung-Jo Jung and Hyun Kyu Kim, {\it On self-folded triangulations in ${\rm SL}_n$ cluster varieties for surfaces}, in preparation.

\bibitem[KN11]{KN} R.M. Kashaev and T. Nakanishi, {\it Classical and quantum dilogarithm identities}, SIGMA Symmetry Integrability Geom. Methods Appl. {\bf 7} (2011), Paper 102, 29 pages. \quad arXiv:1104.4630

\bibitem[K20]{Kim} H. Kim, {\it ${\rm SL}_3$-laminations as bases for ${\rm PGL}_3$ cluster varieties for surfaces}, to appear in Mem. Am. Math. Soc.\quad arXiv:2011.14765

\bibitem[K21]{Kim21} H. Kim, {\it Phase constants in the Fock-Goncharov quantum cluster varieties}, Anal. Math. Phys. {\bf 11}, 2 (2021). \quad arXiv:1602.00797

\bibitem[KL]{KL} H. Kim, T.T.Q. L\^e, and Z. Wang, {\it Frobenius homomorphisms for stated ${\rm SL}_n$ skein algebras}, in preparation.

\bibitem[KLS18]{KLS} H. Kim, T.T.Q. L\^e and M. Son, {\it ${\rm SL}_2$ quantum trace in quantum Teichm\"uller theory via writhe}, Algebr. Geom. Topol. {\bf 23}(1) (2023), 339-418. \quad arXiv:1812.11628

\bibitem[KS]{Kim-Shen} H. Kim and L. Shen, {\it On the canonicity of ${\rm SL}_3$-${\rm PGL}_3$ duality maps for cluster varieties for surfaces}, in preparation.


\bibitem[KY20]{KY} H. Kim and M. Yamazaki, {\it Comments on exchange graphs in cluster algebras}, Exp. Math. {\bf 29}(1) (2020), 79--100. \quad arXiv:1612.00145

\bibitem[KT99]{KT99} A. Knutson and T. Tao, {\it The honeycomb model of $GL_n(\mathbb{C})$ tensor products I: proof of the saturation conjecture}, J. Am. Math. Soc. {\bf 12}(4) (1999), 1055-1090.

\bibitem[K96]{Kuperberg} G. Kuperberg, {\it Spiders for Rank $2$ Lie Algebras}, Commun. Math. Phys. {\bf 180} (1996), 109--151. \quad arXiv:q-alg/9712003

\bibitem[L-F09]{LF} D. Labardini-Fragoso, {\it Quivers with potentials associated to triangulated surfaces, part II: arc representations}, arXiv:0909.4100

\bibitem[L16]{Le16} I. Le, {\it Higher laminations and affine buildings}, Geom. Topol. {\bf 20} (2016), 1673--1735. \quad arXiv:1209.0812

\bibitem[L18]{Le18} T.T.Q. L\^{e}, {\it Triangular decomposition of skein algebras}, Quantum Topol. {\bf 9} (2018), 591--632. \quad arXiv:1609.04987

\bibitem[L19]{Le17} T.T.Q. L\^{e}, {\it Quantum Teichm\"uller spaces and quantum trace map}, J. Inst. Math. Jussieu {\bf 18}(2) (2019), 249--291. \quad arXiv:1511.06054

\bibitem[LY23]{LY23} T.T.Q. L\^{e} and T. Yu, {\it Quantum traces for $SL_n$-skein algebras}, arXiv:2303.08082

\bibitem[MQ23]{MQ} T. Mandel and F. Qin, {\it Bracelets bases are theta bases}, arXiv:2301.11101

\bibitem[RT90]{RT} N.Y. Reshetikhin and V.G. Turaev. {\it 
Ribbon graphs and their invariants derived from quantum groups}, Commun. Math. Phys. {\bf 127}(1) (1990), 1--26. 

\bibitem[S20]{Shen} L. Shen, {\it Duals of semisimple Poisson-Lie groups and cluster theory of moduli spaces of $G$-local systems}, Int. Math. Res. Notices. {\bf 2022}(18) (2022), 14295--14318. \quad arXiv:2003.07901

\bibitem[S01]{S01} A.S. Sikora, {\it ${\rm SL}_n$-character varieties as spaces of graphs}, T. Am. Math. Soc. {\bf 353}(7) (2001), 2773--2804. \quad arXiv:math/9806016

\bibitem[S05]{S05} A.S. Sikora, {\it Skein theory for $SU(n)$-quantum invariants}, Algebr. Geom. Topol. {\bf 5} (2005), 865--897. \quad arXiv:math/0407299


\bibitem[SW07]{SW} A.S. Sikora and B.W. Westbury, {\it Confluence theory of graphs}, Algebr. Geom. Topol. {\bf 7} (2007), 439--478. \quad arXiv:math/0609832

\bibitem[S20]{Miri} M. Son, {\it Quantum coordinate change map for Chekhov-Fock square root algebras}, Master's thesis, Ewha Womans University, Seoul, 2020.

\bibitem[T14]{T14} D.P. Thurston, {\it Positive bases for surface skein algebras}, Proc. Natl. Acad. Sci. {\bf 111}(27) (2014), 9725--9732. \quad arXiv:1310.1959

\end{thebibliography}
\end{document}